%% file: main.tex
\documentclass[11pt, a4paper, oneside]{Thesis} % Paper size, default font size and one-sided paper

\graphicspath{{Pictures/}} % Specifies the directory where pictures are stored

\usepackage[square, numbers, comma, sort&compress]{natbib} % Use the natbib reference package - read up on this to edit the reference style; if you want text (e.g. Smith et al., 2012) for the in-text references (instead of numbers), remove 'numbers' 
\usepackage[spanish, english]{babel} % For translating the keyword Demostracion
\hypersetup{urlcolor=blue, colorlinks=true} % Colors hyperlinks in blue - change to black if annoying
\title{\ttitle} % Defines the thesis title - don't touch this
\usepackage[all]{xy}
\usepackage{nicefrac}
\usepackage[all]{xy}
\usepackage{stmaryrd}

\usepackage{enumitem}
\usepackage{comment}
\usepackage{empheq}
\usepackage{comment}
\usepackage[mathscr]{euscript}
\usepackage{textcomp}

\begin{document}

\input{mathCommands}

\frontmatter % Use roman page numbering style (i, ii, iii, iv...) for the pre-content pages

\setstretch{1.3} % Line spacing of 1.3

% Define the page headers using the FancyHdr package and set up for one-sided printing
\fancyhead{} % Clears all page headers and footers
\rhead{\thepage} % Sets the right side header to show the page number
\lhead{} % Clears the left side page header

\pagestyle{fancy} % Finally, use the "fancy" page style to implement the FancyHdr headers

\newcommand{\HRule}{\rule{\linewidth}{0.5mm}} % New command to make the lines in the title page

% PDF meta-data
\hypersetup{pdftitle={\ttitle}}
\hypersetup{pdfsubject=\subjectname}
\hypersetup{pdfauthor=\authornames}
\hypersetup{pdfkeywords=\keywordnames}

%----------------------------------------------------------------------------------------
%	TITLE PAGE
%----------------------------------------------------------------------------------------

\begin{titlepage}
\begin{center}

\textsc{\Huge \univname}\\[0.9cm] % University name
%\textsc{\Large Doctorate Thesis}\\[0.5cm] % Thesis type

\HRule \\[0.7cm] % Horizontal line

{\Huge \bfseries Topological Quantum Field Theories}

{\Huge \bfseries for Character Varieties}\\[0.5cm]

{\LARGE \bfseries Teor\'ias Topol\'ogicas de Campos Cu\'anticos}

{\LARGE \bfseries para Variedades de Caracteres}
\\[0.3cm] % Thesis title
\HRule \\[0.3cm] % Horizontal line

\huge{
\href{mailto:angel_gonzalez@ucm.es}{\authornames} % Author name - remove the \href bracket to remove the link
}\\[0.3cm]

\LARGE \bfseries{{Directores}}\\[0cm]

\Large
\href{mailto:marina.logares@plymouth.ac.uk}{\supnamefull}\\[-0.1cm]

\href{vicente.munoz@mat.ucm.es }{\supnamesecondfull}\\[0.8cm]

\includegraphics[scale=0.68]{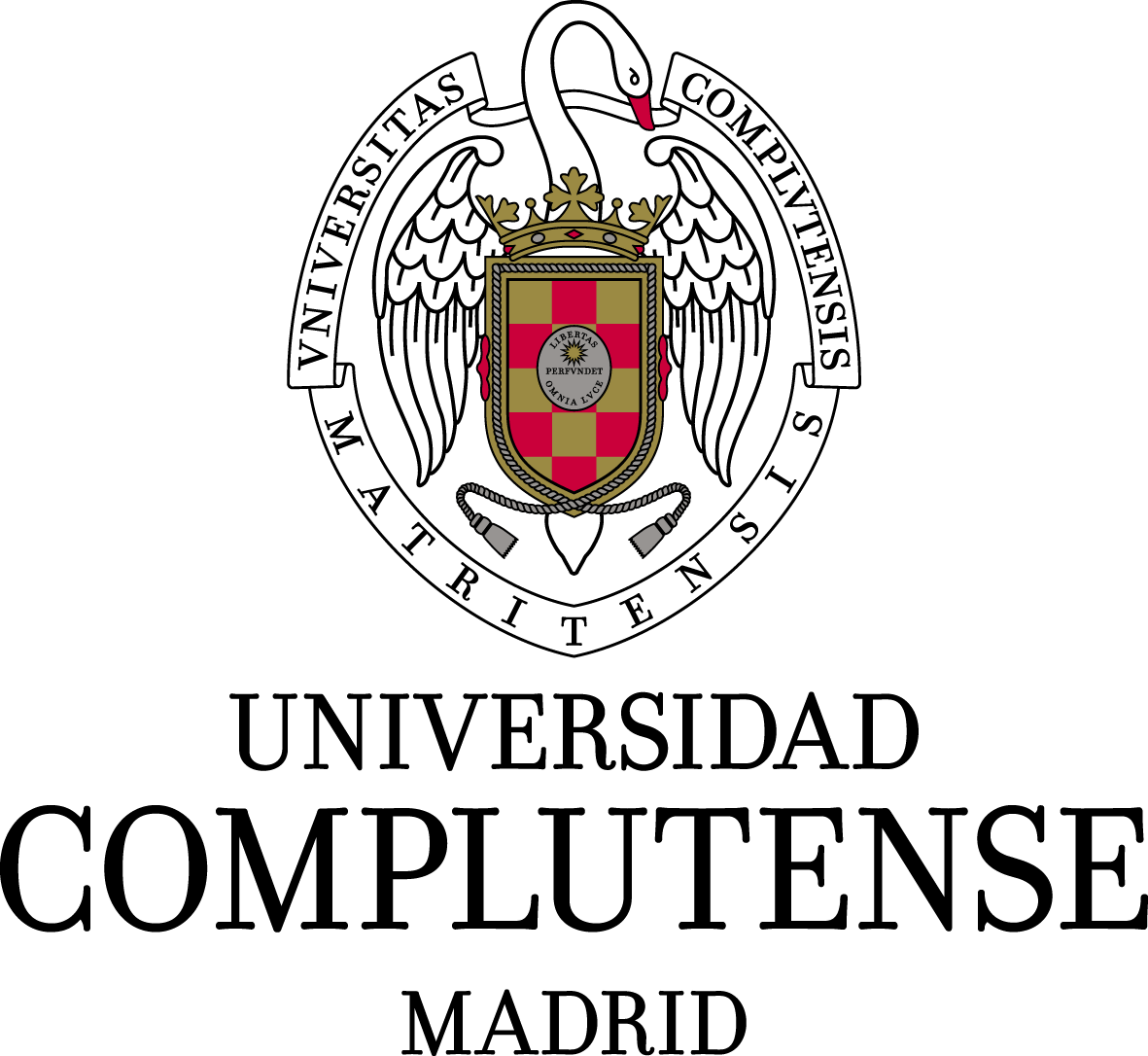}\\[0.6cm] % University/department logo - uncomment to place it

\large \textit{Tesis presentada para optar al grado de\\\Large{Doctor en Investigaci\'on Matem\'atica}}\\[0.2cm] % University requirement text

\Large{
\deptname\\[-0.1cm] % Research group name and department name
\facname\\[0.2cm] % Research group name and department name
}

{\Large Madrid, 2018}

%\vspace{-0.3cm}

%\flushleft{\small{\textcopyright\,\,Jos\'e \'Angel Gonz\'alez Prieto}}

\end{center}

\end{titlepage}

\newpage{}
\thispagestyle{empty}
\null

%----------------------------------------------------------------------------------------
%	TITLE PAGE 2
%----------------------------------------------------------------------------------------

\begin{titlepage}
\begin{center}

\textsc{\Huge \univname}\\[1.5cm] % University name
%\textsc{\Large Doctorate Thesis}\\[0.5cm] % Thesis type

\HRule \\[0.7cm] % Horizontal line

{\Huge \bfseries Topological Quantum Field Theories}

{\Huge \bfseries for Character Varieties}
\\[0.4cm] % Thesis title
\HRule \\[0.6cm] % Horizontal line

\huge{
\href{mailto:angel_gonzalez@ucm.es}{\authornames} % Author name - remove the \href bracket to remove the link
}\\[0.4cm]
 
\Large \bfseries{{Advisors: \href{mailto:marina.logares@plymouth.ac.uk}{\supname} and \href{vicente.munoz@mat.ucm.es }{\supnamesecond}}}\\[2cm]

\includegraphics[scale=0.65]{logo.png}\\[1.3cm] % University/department logo - uncomment to place it

\large \textit{Thesis submitted in fulfilment of the requirements for the degree of\\\Large{\degreename}}\\[0.3cm] % University requirement text

\Large{
\deptname\\[0cm] % Research group name and department name
\facname\\[0.8cm] % Research group name and department name
}

{\Large Madrid, 2018}

\vfill
\end{center}

\end{titlepage}

\pagestyle{empty} % No headers or footers for the following pages

\newpage{}
\null
\clearpage

\clearpage % Start a new page

\null\vfill % Add some space to move the quote down the page a bit

\begin{flushright}
\textit{``It's not who I am underneath,\\but what I do that defines me."}

Batman
%\textit{``If it's just turning the crank it's algebra;\\but if it's got an idea in it, it's topology"}

%Solomon Lefschetz
\end{flushright}

\vfill\vfill\vfill\vfill\vfill\vfill\null % Add some space at the bottom to position the quote just right

\newpage{}
\null
\clearpage

\clearpage % Start a new page

%----------------------------------------------------------------------------------------
%	ABSTRACT PAGE
%----------------------------------------------------------------------------------------

\addtotoc{Abstract} % Add the "Abstract" page entry to the Contents

\abstract{\addtocontents{toc}{\vspace{0cm}} % Add a gap in the Contents, for aesthetics

\vspace{0cm}
\Large{\textbf{Abstract en Espa\~nol}}
\normalsize

La presente tesis doctoral est\'a dedicada al estudio de las estructuras de Hodge de un tipo especial de variedades algebraicas complejas que reciben el nombre de variedades de caracteres. Para este fin, proponemos utilizar una poderosa herramienta de naturaleza \'algebro-geom\'etrica, proveniente de la f\'isica te\'orica, conocida como Teor\'ia Topol\'ogica de Campos Cu\'anticos (TQFT por sus siglas en ingl\'es). Con esta idea, en la presente tesis desarrollamos un formalismo que nos permite construir TQFTs a partir de dos piezas m\'as sencillas de informaci\'on: una teor\'ia de campos (informaci\'on geom\'etrica) y una cuantizaci\'on (informaci\'on algebraica). Como aplicaci\'on, construimos una TQFT que calcula las estructuras de Hodge de variedades de representaciones y la usamos para calcular expl\'icitamente los polinomios de Deligne-Hodge de $\mathrm{SL}_2(\mathbb{C})$-variedades de caracteres parab\'olicas.

\textbf{Palabras clave:} \emph{TQFT, variedades de caracteres, estructuras de Hodge, espacios de m\'oduli, GIT.}

\vspace{0.5cm}
\HRule

\vspace{0.7cm}

\Large{\textbf{Abstract in English}}
\normalsize

This PhD Thesis is devoted to the study of Hodge structures on a special type of complex algebraic varieties, the so-called character varieties. For this purpose, we propose to use a powerful algebro-geometric tool coming from theoretical physics, known as Topological Quantum Field Theory (TQFT). With this idea in mind, in the present Thesis we develop a formalism that allows us to construct TQFTs from two simpler pieces of data: a field theory (geometric data) and a quantisation (algebraic data). As an application, we construct a TQFT computing Hodge structures on representation varieties and we use it for computing explicity the Deligne-Hodge polynomials of parabolic $\mathrm{SL}_2(\mathbb{C})$-character varieties.

\vspace{0.0cm}


\vspace{-0.2cm}
 \vspace{-0.2cm}

\newpage{}
\null
\clearpage

\clearpage % Start a new page

%----------------------------------------------------------------------------------------
%	ACKNOWLEDGEMENTS
%----------------------------------------------------------------------------------------

\setstretch{1.3} % Reset the line-spacing to 1.3 for body text (if it has changed)

\acknowledgements{\addtocontents{toc}{\vspace{0em}} % Add a gap in the Contents, for aesthetics

El camino hacia la consecuci\'on del Doctorado en Investigaci\'on Matem\'atica que culmina con esta tesis ha sido largo, complejo y lleno de curvas. Sin embargo, de ninguna manera puedo decir que haya sido una senda dura o desagradable. Durante los tres a\~nos en los que he realizado este doctorado he sido, simple y llanamente, feliz. Me he podido dedicar en cuerpo y alma a lo que siempre hab\'ia deseado y, lejos de defraudarme, ha sido una experiencia enormemente enriquecedora y, por qu\'e no decirlo, divertida. Sin embargo, la causa \'ultima de toda esta felicidad no es sino la gente que me ha rodeado a lo largo de todo este proceso y, sin cuya ayuda y apoyo, esta tesis no habr\'ia sido posible. Normalmente en este tipo de escritos se suele argumentar el manido ``no mencionar\'e nombres para evitar omisiones". Yo creo que me arriesgar\'e.

En primer lugar, quiero dedicar esta tesis a mis padres. Durante toda mi vida, ellos han sido un gran apoyo, sustento y acicate para crecer y progresar. Si decid\'i hacerme cient\'ifico fue gracias a ellos. De mis padres aprend\'i que la esencia del ser humano es cuestionarse cosas y que la cultura nos hace libres. Ellos me ense\~naron que la oportunidad de aprender es la mayor fortuna que se puede tener.

No pod\'ia faltar aqu\'i una menci\'on a todos mis amigos que, de una u otra manera, me han acompa\~nado y empujado a lo largo de este viaje. Gracias a Adri, Alberto, Nico y Nuria por esas fant\'asticas cervezas en el Ferm\'in envueltas en risas y seguidas de una (mala) pel\'icula de terror. Tambi\'en quiero dar las gracias a mis amigos murcianos (de coraz\'on, solo dos de ellos de nacimiento) Germ\'an, Jaime, Luis, Miguel y Sara. Y a mi cu\~nada Cristina y a mis suegros, Juanjo y Mar\'ia, porque siempre me han apoyado y me han hecho sentir uno m\'as de la familia. Y, por supuesto, no me puedo olvidar de Elena, Kike y Tania por todas esas tardes en Principe P\'io, probablemente en el Burger King, en las que las horas se nos van volarnos sin darnos ni cuenta.

Pero no solo he contado con un gran apoyo fuera de la universidad, sino tambi\'en dentro de ella. Por eso, tambi\'en quiero dedicar esta tesis a \'Angel Ramos, Antonio Br\'u, Benjamin Ivorra, David G\'omez-Castro y Susana G\'omez por todas esas interesantes charlas, buenos consejos y risas en tantas y tantas comidas. Muchas gracias de todo coraz\'on. Gracias a vosotros, ir a la universidad cada d\'ia ha sido fant\'astico.

Tambi\'en quiero dedicar unas l\'ineas de estos agradecimientos a mis compa\~neros de departamento Giovanni, Juan \'Angel y Luc\'ia, por tantas largas conversaciones matem\'aticas y por siempre ayudarme a responder todas mis preguntas sobre spinores, \'algebras de Clifford y operadores de Dirac. Vosotros me hab\'eis ense\~nado que la matem\'atica no debe ser una actividad solitaria y que la verdadera comprensi\'on solo se alcanza cuando se discute sobre ella. Y, por supuesto, tambi\'en quer\'ia dedicar esta tesis a Antonio Vald\'es por todos estos a\~nos que hemos compartido mientras le ayudaba en la docencia. Quiero darte las gracias sinceramente porque siempre me has escuchado, valorado y tenido muy en cuenta mi opini\'on. A tu lado nunca me he sentido un mero ayudante, sino un profesor m\'as y eso me ha permitido encontrar en la docencia una verdadera vocaci\'on. Nuestros tes con porra y nuestras comidas de fin de cuatrimestre son siempre un aut\'entico placer y en ti he encontrado no solo un colega sino tambi\'en un amigo.

I would like to thank the reviews of this thesis, David Ben-Zvi and Peter Newstead, for their careful reading of the manuscript. Their very helpful comments and constructive suggestions have improved the exposition of this thesis.

I would like to specially thank Thomas Wasserman. I am greatly indebted for your invaluable help through the development of this thesis. It was you who taught me all the secrets of TQFTs and who gave me new tools and perspectives for addressing the problems I had found. I really enjoy all our talks about maths, physics and life in general. This thesis would have been impossible without your help and support.

Por descontado, tambi\'en me gustar\'ia agradecer a mis directores, Marina Logares y Vicente Mu\~noz, todo el apoyo, cari\~no y confianza que han depositado en m\'i durante todos estos a\~nos. Sinceramente, creo que no podr\'ia haber tenido mejores gu\'ias para esta ruta. Gracias por haber escuchado con gran atenci\'on todas y cada una de mis ideas, por alocadas que fuesen, y por siempre haber sido capaces de extraer algo bueno de ellas. Nunca me hab\'eis desalentado y siempre me hab\'eis animado a no desistir de mis intentos por atacar un problema. Creo sinceramente que eso es lo que distingue a un mero supervisor de un buen director. Ojal\'a pueda estar a la altura cuando me llegue el momento de ser yo el director. Cuando os conoc\'i, podr\'ia saber m\'as o menos matem\'aticas (m\'as bien menos), pero vosotros me hab\'eis ense\~nado a ser Matem\'atico, con may\'usculas. Sois mucho m\'as que mis directores, sois mis maestros.

Finalmente, mis \'ultimas l\'ineas de esta secci\'on son, como no podr\'ia ser de otro modo, para Laura. Gracias. Gracias por todo tu amor, tu apoyo, tu cari\~no y tu comprensi\'on durante toda la realizaci\'on de esta tesis (y antes). T\'u has sido mi energ\'ia, mi motor y mi fuerza para seguir estudiando, trabajando y avanzando durante todo el doctorado. T\'u siempre confiaste en m\'i, incluso cuando ni yo mismo lo hac\'ia y los teoremas no sal\'ian. Siempre cre\'iste en m\'i. Simplemente, gracias por existir y por estar a mi lado.

%\textit{Durante el desarrollo de esta tesis, he sido financiado por una beca de la Fundaci\'on "La Caixa" destinada a cursar estudios de doctorado en universidades españolas.}
 
}

\clearpage % Start a new page

%----------------------------------------------------------------------------------------
%	LIST OF CONTENTS/FIGURES/TABLES PAGES
%----------------------------------------------------------------------------------------

\pagestyle{fancy} % The page style headers have been "empty" all this time, now use the "fancy" headers as defined before to bring them back

\lhead{\emph{Table of Contents}} % Set the left side page header to "Contents"
\tableofcontents % Write out the Table of Contents

\setstretch{1.3} % Return the line spacing back to 1.3

\pagestyle{empty} % Page style needs to be empty for this page

\null
\vspace{3cm}

\begin{flushright}
\textit{Para Laura,\\ mi quintaesencia.}
\end{flushright}

\nonumber
\thispagestyle{empty}

%\dedicatory{} % Dedication text

%\addtocontents{toc}{\vspace{1em}} % Add a gap in the Contents, for aesthetics

\newpage{}
\null
\clearpage

\clearpage % Start a new page

%----------------------------------------------------------------------------------------
%	INTRODUCTION
%----------------------------------------------------------------------------------------

\setstretch{1.3} % Reset the line-spacing to 1.3 for body text (if it has changed)

\introduction{\addtocontents{toc}{\vspace{0em}} % Add a gap in the Contents, for aesthetics
}

\input{Chapters/Introduction}

\clearpage % Start a new page

\newpage{}
\null
\clearpage

\clearpage % Start a new page

%----------------------------------------------------------------------------------------
%	INTRODUCTION IN SPANISH
%----------------------------------------------------------------------------------------

\setstretch{1.3} % Reset the line-spacing to 1.3 for body text (if it has changed)

\introduccion{\addtocontents{toc}{\vspace{0em}} % Add a gap in the Contents, for aesthetics
}

\input{Chapters/Introduccion}

\clearpage

%----------------------------------------------------------------------------------------
%	THESIS CONTENT - CHAPTERS
%----------------------------------------------------------------------------------------

\mainmatter % Begin numeric (1,2,3...) page numbering

\pagestyle{fancy} % Return the page headers back to the "fancy" style

% Include the chapters of the thesis as separate files from the Chapters folder
% Uncomment the lines as you write the chapters

\input{Chapters/Tqft}

\newpage{}
\thispagestyle{empty}
\null

\input{Chapters/Hodge}

\newpage{}
\thispagestyle{empty}
\null

\input{Chapters/RepresentationVarieties}

\input{Chapters/TopGIT}

%----------------------------------------------------------------------------------------
%	THESIS CONTENT - FUTURE WORK
%----------------------------------------------------------------------------------------

\clearpage

%\lhead{\emph{Future work}} % Change the page header to say "Future work"
\futurework{\addtocontents{toc}{\vspace{1em}} % Add a gap in the Contents, for aesthetics
}

\input{Chapters/FutureWork}

%----------------------------------------------------------------------------------------
%	THESIS CONTENT - APPENDICES
%----------------------------------------------------------------------------------------

%\addtocontents{toc}{\vspace{2em}} % Add a gap in the Contents, for aesthetics

%\appendix % Cue to tell LaTeX that the following 'chapters' are Appendices

% Include the appendices of the thesis as separate files from the Appendices folder
% Uncomment the lines as you write the Appendices

%\input{Appendices/ComplexGeometry}

%\addtocontents{toc}{\vspace{2em}} % Add a gap in the Contents, for aesthetics

%----------------------------------------------------------------------------------------
%	BIBLIOGRAPHY
%----------------------------------------------------------------------------------------
\backmatter
\label{Bibliography}

\lhead{\emph{Bibliography}} % Change the page header to say "Bibliography"

\bibliographystyle{thesisstyle} % Use the "unsrtnat" BibTeX style for formatting the Bibliography

%\nocite{*}
%

%
\bibliography{bibliography}

%\begin{thebibliography}{20}
%\end{thebibliography}

\end{document}

%% file: mathCommands.tex
\theoremstyle{teorema}
\newtheorem{thm}{Theorem}[section]
\newtheorem*{thm*}{Theorem}
\newtheorem{prop}[thm]{Proposition}
\newtheorem{lem}[thm]{Lemma}
\newtheorem{cor}[thm]{Corollary}
\newtheorem{conjecture}{Conjecture}

\theoremstyle{definition}
\newtheorem{defn}[thm]{Definition}
\newtheorem{ex}[thm]{Example}
\newtheorem{as}{Assumption}

\theoremstyle{remark}
\newtheorem{rmk}[thm]{Remark}

\theoremstyle{remark}
\newtheorem*{prf}{Remark}
\newtheorem*{prfskt}{Sketch of Proof\textbf{}}
\renewcommand{\qedsymbol}{$\blacksquare$}

\newcommand{\iacute}{\'{\i}} %i con acento%
\newcommand{\norm}[1]{\lVert#1\rVert} %norma%

\newcommand{\lto}{\longrightarrow}
\newcommand{\hra}{\hookrightarrow}

\newcommand{\suchthat}{\;\;|\;\;}
\newcommand{\dbar}{\overline{\partial}}

\newcommand{\cA}{\mathcal{A}}
\newcommand{\cC}{\mathcal{C}}
\newcommand{\cD}{\mathcal{D}}
\newcommand{\cE}{\mathcal{E}}
\newcommand{\cF}{\mathcal{F}}
\newcommand{\cG}{\mathcal{G}} %Gauge group%
\newcommand{\cI}{\mathcal{I}} %Gauge group%
\newcommand{\cO}{\mathcal{O}} %Holomorphic functions sheaf%
\newcommand{\cM}{\mathcal{M}} %Moduli space% %moduli of parabolic bundles% %moduli of U(p,q) bundles%
\newcommand{\cN}{\mathcal{N}} %Space of minimal points of the Morse function%
\newcommand{\cP}{\mathcal{P}} %Moduli of K(D) pairs%
\newcommand{\cQ}{\mathcal{Q}} %Moduli of K(D) pairs%
\newcommand{\cS}{\mathcal{S}} %Moduli of solutions of Hitchin's equations, contructed by Konno%
\newcommand{\cU}{\mathcal{U}} %Moduli of stable U(p,q) parabolic Higgs bundles%
\newcommand{\cJ}{\mathcal{J}}
\newcommand{\cX}{\mathcal{X}}
\newcommand{\cT}{\mathcal{T}}
\newcommand{\cV}{\mathcal{V}}
\newcommand{\cW}{\mathcal{W}}
\newcommand{\cB}{\mathcal{B}}
\newcommand{\cR}{\mathcal{R}}
\newcommand{\cH}{\mathcal{H}}
\newcommand{\cZ}{\mathcal{Z}}
\newcommand{\D}{\bar{B}}

\newcommand{\Ker}{\textrm{Ker}\,}
\newcommand{\coker}{\textrm{coker}\,}

\newcommand{\ext}{\mathrm{ext}} % an extension%
\newcommand{\x}{\times}

\newcommand{\mM}{\mathscr{M}} %Meromorphic function sheaf%

\newcommand{\CC}{\mathbb{C}} %Complex numbers%
\newcommand{\QQ}{\mathbb{Q}} %Rational numbers%
\newcommand{\FF}{\mathbb{F}} %Rational numbers%
\newcommand{\PP}{\mathbb{P}} %projective space%
\newcommand{\HH}{\mathbb{H}} %Hypercohomology, quaternions..%
\newcommand{\RR}{\mathbb{R}} %Real numbers%
\newcommand{\ZZ}{\mathbb{Z}} %Integer numbers%
\newcommand{\NN}{\mathbb{N}} %Natural numbers%
\newcommand{\DD}{\mathbb{D}} %Natural numbers%

\renewcommand{\lg}{\mathfrak{g}} %Lie algebra of G%
\newcommand{\lh}{\mathfrak{h}} %Lie algebra of H%
\newcommand{\lu}{\mathfrak{u}} %Lie algebra of U%
\newcommand{\la}{\mathfrak{a}} %Lie algebra of A%
\newcommand{\lb}{\mathfrak{b}} %Lie algebra of B%
\newcommand{\lm}{\mathfrak{m}} %Lie algebra of M%
\newcommand{\lgl}{\mathfrak{gl}} %Lie algebra of GL%
\newcommand{\lZ}{\mathfrak{Z}} %Almost-TQFT%

\newcommand{\too}{\longrightarrow}
\newcommand{\imat}{\sqrt{-1}} %i%
\newcommand{\tinyclk}{{\scriptscriptstyle \Taschenuhr}} %Tiny Clock Symbol from ifsym package%
\newcommand\restr[2]{\left.#1\right|_{#2}}
\newcommand\rtorus[1]{{\mathbb{T}^{#1}}}
\newcommand\actPartial{\overline{\partial}}
\newcommand\handle[2]{\mathcal{A}^{#1}_{#2}}

\newcommand\pastingArea[2]{S^{#1-1} \times \bar{B}^{#2-#1}}
\newcommand\pastingAreaPlus[2]{S^{#1} \times \bar{B}^{#2-#1-1}}
\newcommand\coordvector[2]{\left.\frac{\partial}{\partial {#1}}\right|_{#2}}

% Categories
\newcommand\Sets{\textbf{Set}}
\newcommand\Cat{\textbf{Cat}}
\newcommand\Top{\textbf{Top}}
\newcommand\TopHLC{\textbf{Top}_{hlc}}
\newcommand\TopS{\textbf{Top}_\star}
\newcommand\Diff{\textbf{Diff}}
\newcommand\Diffc{\textbf{Diff}_c}

\newcommand\CBord[3]{\mathbf{Bd}_{{#1 #3}}^{#2}}
\newcommand\Bord[1]{\CBord{#1}{}{}}
\newcommand\EBord[2]{\CBord{#1}{#2}{}}
\newcommand\Bordo[1]{\CBord{#1}{or}{}}
\newcommand\Bordp[1]{\mathbf{Bdp}_{{#1}}}
\newcommand\Bordpar[2]{\mathbf{Bd}_{{#1}}(#2)}
\newcommand\Bordppar[2]{\mathbf{Bdp}_{{#1}}(#2)}

\newcommand\Tubo[1]{\mathbf{Tb}_{#1}^0}
\newcommand\Tub[1]{\mathbf{Tb}_{#1}}
\newcommand\ETub[2]{\mathbf{Tb}_{#1}^{#2}}
\newcommand\Tubp[1]{\mathbf{Tbp}_{#1}}
\newcommand\Tubpo[1]{\mathbf{Tbp}_{#1}^0}
\newcommand\Tubppar[2]{\mathbf{Tbp}_{#1}(#2)}

\newcommand\Embc{\textbf{Emb}_c}
\newcommand\EEmbc[1]{\textbf{Emb}_c^{#1}}
\newcommand\Embpc{\textbf{Embp}_c}
\newcommand\Embparc[1]{\textbf{Emb}_c(#1)}
\newcommand\Embpparc[1]{\textbf{Embp}_c(#1)}

\newcommand\cSp{\cS_p}
\newcommand\cSpar[1]{\cS_{#1}}

\newcommand\Diffpc{\textbf{Diff}_c}

\newcommand\PVar[1]{\mathbf{PVar}_{#1}}
\newcommand\PVarC{\mathbf{PVar}_{\CC}}
\newcommand\Sr[1]{\mathrm{S}{#1}}

\newcommand\Bordpo[1]{\mathbf{Bdp}_{{#1}}^{or}}
\newcommand\ClBordp[1]{\mathbf{l}\CBord{#1}{}{}}
\newcommand\CTub[3]{\mathbf{Tb}_{{#1 #3}}^{#2}}
\newcommand\CTubp[1]{\CTub{#1}{}{}}
\newcommand\CTubpp[1]{\mathbf{Tbp}_{#1}{}{}}
\newcommand\CTubppo[1]{\mathbf{Tbp}_{#1}^0}
\newcommand\CClose[3]{\mathbf{Cl}_{{#1 #3}}^{#2}}
\newcommand\CClosep[1]{\CClose{#1}{}{}}
\newcommand\Obj[1]{\mathrm{Obj}(#1)}
\newcommand\Mor[1]{\mathrm{Mor}(#1)}
\newcommand\Vect[1]{{#1}\textrm{-}\mathbf{Vect}}
\newcommand\Vecto[1]{{#1}\textrm{-}\mathbf{Vect}_0}
\newcommand\Mod[1]{{#1}\textrm{-}\mathbf{Mod}}
\newcommand\Modt[1]{{#1}\textrm{-}\mathbf{Mod}_t}
\newcommand\Rng{\mathbf{Ring}}
\newcommand\Grp{\mathbf{Grp}}
\newcommand\Grpd{\mathbf{Grpd}}
\newcommand\Grpdo{\mathbf{Grpd}_0}
\newcommand\HS[1]{\mathbf{HS}^{#1}}
\newcommand\MHS[1]{\mathbf{MHS}}
\newcommand\PHS[2]{\mathbf{HS}^{#1}_{#2}}
\newcommand\MHSq{\mathbf{HS}}
\newcommand\Sch{\mathbf{Sch}}
\newcommand\Sh[1]{\mathbf{Sh}\left(#1\right)}
\newcommand\QSh[1]{\mathbf{QSh}\left(#1\right)}
\newcommand\Var[1]{\mathbf{Var}_{#1}}
\newcommand\CVar{\Var{\CC}}
\newcommand\PHM[2]{\cM_{#1}^p(#2)}
\newcommand\MHM[1]{\cM_{#1}}
\newcommand\HM[2]{\textrm{HM}^{#1}(#2)}
\newcommand\HMW[1]{\textrm{HMW}(#1)}
\newcommand\VMHS[1]{VMHS({#1})}
\newcommand\geoVMHS[1]{VMHS_g({#1})}
\newcommand\goodVMHS[1]{\mathrm{VMHS}_0({#1})}
\newcommand\Par[1]{\mathrm{Par}({#1})}
\newcommand\K[1]{\mathrm{K}#1}
\newcommand\Ko[1]{\mathrm{K}{#1}_0}
\newcommand\KM[1]{\mathrm{K}\MHM{#1}}
\newcommand\KMo[1]{\mathrm{K}{\MHM{#1}}_0}
\newcommand\Ab{\mathbf{Ab}}
\newcommand\CPP{\cP\cP}
\newcommand\Bim[1]{{#1}\textrm{-}\mathbf{Bim}}
\newcommand\Span[1]{\mathrm{Span}({#1})}
\newcommand\Spano[1]{\mathrm{Span}^{op}({#1})}
\newcommand\GL[1]{\mathrm{GL}_{#1}}
\newcommand\SL[1]{\mathrm{SL}_{#1}}
\newcommand\PGL[1]{\mathrm{PGL}_{#1}}
\newcommand\Rep[1]{\mathfrak{X}_{#1}}

\newcommand\Qtm[1]{\mathcal{Q}_{#1}}
\newcommand\sQtm[1]{\mathcal{Q}_{#1}^0}
\newcommand\Fld[1]{\mathcal{F}_{#1}}

% Hodge theory
\newcommand\DelHod[1]{e\left(#1\right)}
\newcommand\RDelHod{e}
\newcommand\eVect{\mathcal{E}}
\newcommand\e[1]{\eVect\left(#1\right)}
\newcommand\intMor[2]{\int_{#1}\,#2}

% Representation Varieties
\newcommand\Dom[1]{\mathcal{D}_{#1}}

\newcommand\Xf[1]{{X}_{#1}}					% Free non-parabolic
\newcommand\Xs[1]{\mathfrak{X}_{#1}}							% Surface group non-parabolic

\newcommand\Xft[2]{\overline{{X}}_{#1, #2}} % Free tr=2
\newcommand\Xst[2]{\overline{\mathfrak{X}}_{#1, #2}}			% Surface group tr=2

\newcommand\Xfp[2]{X_{#1, #2}}			% Free J+
\newcommand\Xsp[2]{\mathfrak{X}_{#1, #2}}						% Surface group J+

\newcommand\Xfd[2]{X_{#1; #2}}			% Free diagonal
\newcommand\Xsd[2]{\mathfrak{X}_{#1; #2}}						% Surface group diagonal

\newcommand\Xfm[3]{\mathcal{X}_{#1, #2; #3}}		% Free mixed
\newcommand\Xsm[3]{X_{#1, #2; #3}}					% Surface group mixed

\newcommand\XD[1]{#1^{D}}
\newcommand\XDh[1]{#1^{\delta}}
\newcommand\XU[1]{#1^{UT}}
\newcommand\XP[1]{#1^{U}}
\newcommand\XPh[1]{#1^{\upsilon}}
\newcommand\XI[1]{#1^{\iota}}
\newcommand\XTilde[1]{#1^{\varrho}}
\newcommand\Xred[1]{#1^{r}}
\newcommand\Xirred[1]{#1^{ir}}

\newcommand\Char[1]{\cR_{#1}}
\newcommand\Chars[1]{\cR_{#1}}
\newcommand\CharW[1]{\mathscr{R}_{#1}}

% Derived Category
\newcommand\Ch[1]{\textrm{Ch}\,{#1}}
\newcommand\Chp[1]{\textrm{Ch}^+{#1}}
\newcommand\Chm[1]{\textrm{Ch}^-{#1}}
\newcommand\Chb[1]{\textrm{Ch}^b{#1}}
\newcommand\Der[1]{\textrm{D}{#1}}
\newcommand\Dp[1]{\textrm{D}^+{#1}}
\newcommand\Dm[1]{\textrm{D}^-{#1}}
\newcommand\Db[1]{\textrm{D}^b{#1}}

% TQFT
\newcommand\Gs{\cG}
\newcommand\Gq{\cG_q}
\newcommand\Gg{\cG_c}
\newcommand\Zs[1]{Z_{#1}}
\newcommand\Zg[1]{Z^{gm}_{#1}}
\newcommand\cZg[1]{\cZ^{gm}_{#1}}

% Miscelany
\newcommand\RM[2]{R\left(\left.#1\right|#2\right)}
\newcommand\RMc[3]{R_{#1}\left(\left.#2\right|#3\right)}
\newcommand\set[1]{\left\{#1\right\}}
\newcommand{\Stab}{\textrm{Stab}\,} %\Stabilizer%
\newcommand{\tr}{\textrm{tr}\,}             %Trace Tr%
\newcommand\EuChS{E}             %Trace Tr%
\newcommand\EuCh[1]{E\left(#1\right)}             %Trace Tr%

\newcommand{\Hom}{\textrm{Hom}\,}           %Hom%
\newcommand{\End}{\textrm{End}\,} 
\newcommand{\Aut}{\textrm{Aut}} 
\newcommand{\Sym}{\textrm{Sym}\,} 
\newcommand{\Ann}{\textrm{Ann}\,}
\newcommand{\Rad}{\textrm{Rad}\,}  
\newcommand{\img}{\textrm{im}\,} 
\newcommand\supp[1]{\mathrm{supp}{(#1)}}
\newcommand\coh[1]{\left[H_c^\bullet\hspace{-0.05cm}\left(#1\right)\right]}
\newcommand\Bt{B_t}
\newcommand\Be{B_e}
\newcommand\re{\textrm{Re}\,}
\newcommand\imag{\textrm{Im}\,}
\newcommand\Kahc{\textbf{K\"ah}_c}
\newcommand\Id{\textrm{Id}}
\newcommand\colExt[2]{\textrm{Coll}(#1, #2)}

% Mixed Hodge modules
\newcommand\Ccs[1]{C_{cs}(#1)}
\newcommand\can{\textrm{can}}
\newcommand\var{\textrm{var}}
\newcommand\Perv[1]{\textrm{Perv}(#1)}
\newcommand\DR[1]{\textrm{DR}{#1}}
\newcommand\Gr[2]{\textrm{Gr}_{#1}^{#2}\,}
\newcommand\ChV{\textrm{Ch}\,}
\newcommand\VerD{^{\textrm{Ve}}\DD}
\newcommand\DHol[1]{\textrm{D}^b_{\textrm{hol}}(\cD_{#1})}
\newcommand\RegHol[1]{\Mod{\cD_{#1}}_{\textrm{rh}}}
\newcommand\Drh[1]{\textrm{D}^b_{\textrm{rh}}(\cD_{#1})}
\newcommand\Dcs[2]{\textrm{D}^b_{\textrm{cs}}({#1}; {#2})}
\newcommand{\rat}{\mathrm{rat}}
\newcommand{\dmod}{\mathrm{Dmod}}
\newcommand{\bimapcat}[1]{{#1}^{bi}}

%% file: Chapters/Introduction.tex
\label{introduction} % For referencing the chapter elsewhere, use \ref{Chapter1} 

\lhead{Introduction} % This is for the header on each page - perhaps a shortened title

%----------------------------------------------------------------------------------------

Let $M$ be a compact differentiable manifold, possibly with boundary, and let $G$ be a complex algebraic group. The set of representations $\rho: \pi_1(M) \to G$ can be endowed with a structure of complex algebraic variety, the so-called representation variety of $M$ into $G$, denoted $\Rep{G}(M)$. Moreover, the group $G$ itself acts on $\Rep{G}(M)$ by conjugation so, taking the Geometric Invariant Theory (GIT) quotient of $\Rep{G}(W)$ by this action, we obtain the so-called character variety
$$
	\cR_G(M) = \Rep{G}(M) \sslash G.
$$
It is the moduli space of representations of $\pi_1(M)$ into $G$, as treated in \cite{Nakamoto}, also known as the Betti moduli space. Even in the simplest cases, the topology and algebraic structure of these character varieties is extremely rich and has been objective of studies for the last twenty years.

One of the main reasons for studying character varieties is their prominent role in the non-abelian Hodge correspondence. This beautiful theory states that character varieties are just one of the three faces of the same object. The first incarnation of this principle is the so-called Riemann-Hilbert correspondence (\cite{SimpsonI} \cite{SimpsonII}). For $M=\Sigma$ a compact Riemann surface and $G = \GL{n}(\CC)$ (resp.\ $G = \SL{n}(\CC)$), an element of $\Char{G}(\Sigma)$ defines a $G$-local system and, thus, a rank $n$ algebraic vector bundle $E \to \Sigma$ of degree zero (resp.\ and fixed determinant) with a flat connection on it. In this way, the Riemann-Hilbert correspondence gives a real analytic correspondence between $\Char{G}(\Sigma)$ and the moduli space of flat bundles on $\Sigma$ of rank $n$ and degree zero (resp.\ and fixed determinant), usually called the de Rham moduli space. Furthermore, via the Hitchin-Kobayashi correspondence (\cite{Simpson:1992} \cite{Corlette:1988}), we also have that, for $G=\GL{n}( \CC)$ (resp.\ $G=\SL{n}(\CC)$), the Betti moduli space $\Char{G}(\Sigma)$ is real analytic equivalent to the Dolbeault moduli space, that is, the moduli space of rank $n$ and degree zero (resp.\ and fixed determinant) $G$-Higgs bundles i.e.\ vector bundles $E \to \Sigma$ together with a field $\Phi: E \to E \otimes K_\Sigma$ called the Higgs field.

Using these correspondences, it is possible to compute the Poincar\'e polynomial of character varieties by means of Morse theory. Following these ideas, Hitchin, in the seminal paper \cite{Hitchin}, gave the Poincar\'e polynomial for $G=\SL{2}(\CC)$. Gothen also computed it for $G=\SL{3}(\CC)$ in \cite{Gothen} and Garc\'ia-Prada, Heinloth and Schmitt for $G=\GL{4}(\CC)$ in \cite{GP-Heinloth-Schmitt}.

However, these correspondences from non-abelian Hodge theory are far from being algebraic. For this reason, the study of the algebraic structure on the character varieties turns important. In particular, it is interesting to analyze an extra linear structure that appears on the cohomology of these varieties, the so called mixed Hodge structure. The origin of these invariants goes back to classical Hodge theory, where they are related with something as subtle as the harmonic forms on a variety. However, it was not until the work of Deligne in \cite{DeligneI:1971}, \cite{DeligneII:1971} and \cite{DeligneI:1971}, in the context of the Weil conjectures, when Hodge structures gained lot of visibility as important algebraic invariants. The point was that, in those papers, Deligne proved that the cohomology of every complex algebraic variety is endowed with a natural Hodge structure.

Roughly speaking, a Hodge structure on a complex variety $X$ is a decomposition of $H_c^k(X;\CC)$ as a lattice $H_c^{k;p,q}(X)$ with $p,q \in \ZZ$. The dimensions of these Hodge pieces, $h_c^{k;p,q}(X) = \dim_\CC H_c^{k;p,q}(X)$, are the so-called Hodge numbers and can be thought as an algebraic refinement of the Betti numbers. In particular, from them we can construct the Deligne-Hodge polynomial, or $E$-polynomial, of $X$ as a kind of Euler characteristic with the Hodge numbers
$$
	\DelHod{X} = \sum_{k, p, q} (-1)^k\,h_c^{k;p,q}(X)\,u^pv^q \in \ZZ[u^{\pm 1}, v^{\pm 1}].
$$
However, in general, the Deligne-Hodge polynomial of character varieties is unknown. One of the most important advances in this direction was accomplished by Hausel and Rodr\'iguez-Villegas in \cite{Hausel-Rodriguez-Villegas:2008} by means of a theorem of Katz of arithmetic flavour, based on the Weil conjectures and the Lefschetz principle. In that paper, they computed an expression of the $E$-polynomial for $\GL{n}(\CC)$-twisted character varieties on surfaces in terms of generating functions. This technique was later expanded in \cite{Mereb} for $G=\SL{n}(\CC)$ and, recently, in \cite{Baraglia-Hekmati:2016} for the untwisted case and orientable surfaces ($G=\GL{3}(\CC)$ and $\SL{3}(\CC)$) and for non-orientable surfaces ($G=\GL{2}(\CC)$ and $\SL{2}(\CC)$). However, this approach suffers of a major problem since, in order to use Katz' theorem, it is necessary to count the number of points of the character varieties over finite fields and that requieres to use character tables of the finite groups $\GL{n}(\mathbb{F}_q)$. However, for high rank, these tables cannot be explicitly computed so only an expression of the $E$-polynomial in terms of generating functions can be obtained.

Trying to overcome this problem, a new approach was initiated by Logares, Mu\~noz and Newstead in \cite{LMN}. The strategy there was to focus on the computation of the Deligne-Hodge polynomial of the representation variety by chopping it into simpler strata for which this polynomial can be computed and to sum all the contributions. After that, they studied the identification that took place in the GIT quotient and, from these data, they computed the Deligne-Hodge polynomial of the corresponding character variety. Following this idea, explicit expressions of the $E$-polynomials for genus $g=1,2$ and $G=\SL{2}(\CC)$ were computed in \cite{LMN} and, later, extended for arbitrary genus in \cite{MM} and, for $G=\textrm{PGL}_2(\CC)$, in \cite{Martinez:2017}.

A step forward in representation varieties is to consider a parabolic structure on a surface $\Sigma$, which is a finite collection of pairs $Q=\left\{(p_i, \lambda_i)\right\}_{i=1}^s$, where $p_i \in \Sigma$ and $\lambda_i \subseteq G$ is a subvariety closed under internal automorphisms of $G$, called the holonomy of $p_i$. In that case, the parabolic representation variety $\Rep{G}(\Sigma, Q)$ is the set of representations $\rho: \pi_1(\Sigma - \left\{p_1, \ldots, p_s\right\}) \to G$ such that the image of the positive loop around $p_i$ lies in $\lambda_i$.

The importance of these parabolic counterparts is that the non-abelian Hodge correspondence extends to the parabolic context. For example, for $G=\GL{n}(\CC)$ and $Q$ a single marked point with holonomy a primitive root of unit $e^{2 \pi i d/n}$ (the so-called twisted case), we have that $\cR_{\GL{n}(\CC)}(\Sigma, Q)$ is diffeomorphic to the moduli space of rank $n$ and degree $d$ Higgs bundles and to the moduli space of rank $n$ logarithmic flat bundles with a simple pole with residue $-\frac{d}{n}\Id$. Analogous correspondences also hold for generic semi-simple conjugacy classes \cite{Simpson:parabolic}.
Using these correspondences, Boden and Yokogawa calculated the Poincar\'e polynomial of character varieties in \cite{Boden-Yokogawa} for $G=\SL{2}(\CC)$ and generic semi-simple conjugacy classes, and Garc\'ia-Prada, Gothen and Mu\~noz for $G=\GL{3}(\CC)$ and $G=\SL{3}(\CC)$ in \cite{GP-Gothen-Munoz}. In the general twisted case, in \cite{Schiffmann:2016} and \cite{Mozgovoy-Schiffmann} (see also \cite{Mellit}) a combinatorial formula is given for arbitrary $G=\GL{n}(\CC)$ provided that $n$ and $d$ are coprime.
Incidentally, other correspondences can also appear, as for $G=\SL{2}(\CC)$, $\Sigma$ an elliptic curve and $Q$ two marked points with different semi-simple conjugacy classes not containing a multiple of identity, in which $\cR_G(\Sigma, Q)$ is diffeomorphic to the moduli space of doubly periodic instantons through the Nahm transform \cite{Biquard-Jardim} \cite{Jardim}.

However, the calculation of $E$-polynomials for parabolic character varieties is much harder than for the non-parabolic case due to the extra perturbation generated by the punctures. Some advances were given in \cite{Hausel-Letellier-Villegas}, following the arithmetic method, for $G=\GL{n}(\CC)$ and generic semi-simple marked points (but the result is not explicit) and in \cite{LM} with the geometric method but only at most two marked points. In the general case there is not a satisfactory solution.

In this PhD Thesis, we propose a new approach to the problem of computing Hodge structures on character varieties based on Topological Quantum Field Theories (TQFTs). The seed of this idea is the paper \cite{MM}, in which the authors computed the Deligne-Hodge polynomials of $\SL{2}(\CC)$-representation varieties by means of the geometric method. However, beyond the calculation itself, the important point is that they proved that, if $\Sigma_g$ is the compact genus $g$ surface, there exists a recursive algorithm that computes the Deligne-Hodge polynomial of $\Rep{\SL{2}(\CC)}(\Sigma_g)$ in terms of the $E$-polynomial of $\Rep{\SL{2}(\CC)}(\Sigma_{g-1})$ and an extra piece of data called the Hodge monodromy representation. This recursive nature can also be observed in other works as in and \cite{Mozgovoy:2012}, \cite{Hausel-Letellier-Villegas:2013}, \cite{Diaconescu:2017} and \cite{Carlsson-Rodriguez-Villegas}. This pattern brings us back to TQFTs. 

The origin of TQFTs dates back to the seminar paper of Witten \cite{Witten:1988} in which he showed that the Jones polynomial (a knot invariant) can be obtained from a quantum field theory known as Chern-Simons theory. Aware of the importance of this discovery, Atiyah in \cite{Atiyah:1988} gave a full definition of this kind of theories and coined the term Topological Quantum Field Theory. There, he emphasized in the categorical nature of this kind of constructions and formulated them as monoidal functors $Z: \Bord{n} \to \Mod{R}$ from the category of $n$-dimensional bordisms to the category of $R$-modules.

\begin{figure}[h]
\vspace{0cm}
\begin{center}
\includegraphics[scale=0.55]{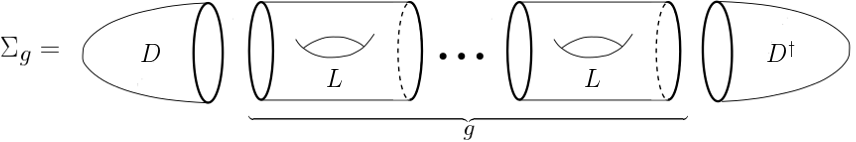}
\caption{}
\label{img:decomp-tubes-intro}
\vspace{-0.6cm}
\end{center}
\end{figure}

However, for our purposes, the most important application of TQFTs is that they give effective methods of computation of algebraic invariants. Suppose that we are interesed in some algebraic invariant that takes values on a ring $R$, and suppose that have a TQFT, $Z: \Bord{n} \to \Mod{R}$, such that, if $M$ is a closed $n$-dimensional manifold, seen as a bordism $M: \emptyset \to \emptyset$, we have that
$Z(M)(1) \in R$ is the desired invariant of $M$. Observe that, as $Z$ is monoidal, $Z(\emptyset)=R$ so $Z(M)$ is a $R$-linear automorphism of $R$ i.e.\ it is multiplication by some fixed element of $R$. In that case, we can split $M$ into simpler pieces in order to obtain such a invariant in a recursive way. For example, suppose that $M= \Sigma_g$ and decompose $M = D^\dag \circ L^g \circ D$, where $D: \emptyset \to S^1$ is a disc, $L: S^1 \to S^1$ is a torus with two discs removed and $D^\dag: S^1 \to \emptyset$ is a disc in the other way around, as depicted in Figure \ref{img:decomp-tubes-intro}.

In that case, the TQFT gives us a decomposition $Z(\Sigma_g) = Z(D^\dag) \circ Z(L)^g \circ Z(D)$, where $Z(D),Z(D^\dag)$ and $Z(L)$ are fixed $R$-module homomorphisms, independent of $g$. Hence, as $Z(\Sigma_g)(1)$ is the desired invariant, it can be computed automatically for all the surfaces $\Sigma_g$ from the knowledge of just three homomorphisms. In that sense, a TQFT can be thought as the dynamic version of something static as an algebraic invariant.

This PhD Thesis is structured in four chapters. Chapter \ref{chap:tqft} is devoted to Topological Quantum Field Theories. In Section \ref{sec:categories-tqft}, especially \ref{sec:low-dim-tqft}, we will review the definitions concerning TQFTs and some of their fundamental properties, including the fact that, for all $M \in \Bord{n}$, the module $Z(M)$ is forced to be finitely generated. This restriction imposes a major obstruction to the construction of TQFTs since, in most cases, the natural assigment for $Z(M)$ is not finitely generated. In order to overcome this problem, in Section \ref{sec:lax-monoidality} we propose to relax the monoidality condition and to consider lax monoidal Topological Quantum Field Theories, that is, functors $Z: \Bord{n} \to \Mod{R}$ such that, for all $M, N \in \Bord{n}$, there exists a homomorphism $\Delta: Z(M) \otimes Z(N) \to Z(M \sqcup N)$ which is no longer an isomorphism.

Under these relaxed conditions, life becomes easier and, actually, there are plenty of TQFTs that can be constructed. In this direction, in Section \ref{sec:lax-monoidal-tqft} we describe a new procedure, of physical inspiration, that creates a lax monoidal TQFT from two simpler pieces: a field theory (a geometric data) and a quantisation (an algebraic data). The idea is to consider a category $\cC$ with final object and pullbacks that is going to play the role of a `category of fields' (in the physical sense), and to split our TQFT as a composition
$$
	\Bord{n} \stackrel{\cF}{\longrightarrow} \Span{\cC} \stackrel{\cQ}{\longrightarrow} \Mod{R},
$$
where $\Span{\cC}$ is the category of spans in $\cC$ (see Section \ref{sec:lax-monoidal-tqft}). The functor $\cF$ is the `field theory' and can be constructed from a functor $\Embc \to \cC$ that preserves gluings (Section \ref{sec:field-theory}). Here, $\Embc$ stands for the wide subcategory of compact differentiable manifolds but only embeddings between them. In the other hand, the functor $\cQ$ is the `quantisation' and can be made from a piece of algebraic data called a $\cC$-algebra (Sections \ref{sec:C-algebras} and \ref{sec:quantisation}). On order to finish Chapter \ref{chap:tqft}, in Section \ref{sec:other-versions-tqft} we describe how this procedure can be expanded to construct other types of Topological Quantum Field Theories as TQFTs preserving some extra structure (Section \ref{sec:TQFT-over-sheaf}), TQFTs only on tubes instead of general bordisms (Section \ref{sec:almost-TQFT-tubes}), simplified versions of TQFTs from an heuristic datum (Section \ref{sec:reduction-TQFT}) and a kind of weaker version of TQFTs, called soft TQFTs, when only half of the $\cC$-algebra is available (Section \ref{sec:soft-TQFT}). These alternative versions of TQFTs will be very useful in Chapter \ref{chap:representation}, where we will construct a TQFT computing Hodge structures on representation varieties.

In Chapter \ref{chap:hodge}, we will focus on Hodge theory, as the framework in which the algebraic invariants that we want to study fit. After some preliminaries about $K$-theory and Grothendieck rings (Section \ref{sec:panoramic-Ktheory}), we will review Hodge structures in Section \ref{sec:hodge-structures}, both pure Hodge structures (as in the compact K\"ahler case) and mixed Hodge structures (in the general case).

However, the core of Chapter \ref{chap:hodge} is the Section \ref{sec:mixed-hodge-modules}, in which we review Saito's theory of mixed Hodge modules. These mixed Hodge modules are cornerstones in our construction of a TQFT computing Hodge structures of representation varieties, since the Grothendieck rings of mixed Hodge modules will be the $\CVar$-algebra choosen for the quantisation. As preparation, in Section \ref{sec:riemann-hilbert} we will review $\cD$-modules and perverse sheaves as well as their interplay by means of the Riemann-Hilbert correspondence. These objects are needed because they actually are the building blocks of the category of mixed Hodge modules, as explained in Section \ref{sec:saito-mixed-hodge-modules}. Roughly speaking, the category of mixed Hodge modules is just a suitable subcategory of the category of filtered pairs of $\cD$-modules and perverse sheaves. However, the selection rule to decide whether one of such a pairs is a mixed Hodge module is extremely involved and requieres lots of subtle functors controlling the behaviour of several filtrations. No original work is intended to be shown in this section, but only a reinterpretation of Saito's theory as a black box suitable for our purposes.

In fact, the definition of mixed Hodge modules is so hard that almost no example of a mixed Hodge module can be constructed by hand. However, in Section \ref{sec:monodromy-as-mhm}, we will review some of the fundamental results of Saito that state that the category of variations of Hodge structures (i.e.\ sheaves whose stalks are Hodge structures) can be seen as a subcategory of the category of mixed Hodge modules. With this result at hand, given a fibration $f: X \to B$ of smooth complex varieties locally trivial in the analytic topology, we will construct a mixed Hodge module on $B$ that controls the monodromy of $f$. Such a Hodge module is called the Hodge monodromy representation of $f$ and, as we will show, generalizes the Hodge monodromy representation introduced in \cite{LMN}.

Chapter \ref{chap:representation} is devoted to representation varieties. First of all, in Section \ref{sec:intro-representation} we will review some of the fundamentals of representation and character varieties and we will put them in context by describing briefly the non-abelian Hodge theory (Section \ref{sec:non-abelian-hodge}). However, the heart of this Thesis is Section \ref{sec:TQFT-for-repr}. In that section, we will prove the following fundamental result (Theorem \ref{thm:existence-s-LTQFT}).

\begin{thm*}
Let $G$ be a complex algebraic group. For every $n \geq 1$, there exists a lax monoidal TQFT of pairs, $\Zs{G}: \Bordp{n} \to \Modt{\K{\MHSq}}$, such that, for every $n$-dimensional connected closed manifold $W$ and every non-empty finite subset $A \subseteq W$, we have
$$
	\Zs{G}(W, A)\left(1\right) = \coh{\Rep{G}(W); \QQ} \otimes \coh{G; \QQ}^{|A|-1} \in \K{\MHSq}.
$$ 
Here, $\MHSq$ stands for the category of (mixed) Hodge structures with $1 \in \K{\MHSq}$ the unit of its Grothendieck ring and $\Modt{\K{\MHSq}}$ is an enharced version of the usual category of modules with an extra $2$-category structure (see Example \ref{ex:categories})
\end{thm*}

This result, as well as the construction methods of TQFTs of Section \ref{sec:lax-monoidal-tqft}, appeared for the first time in the paper \cite{GPLM-2017} (but with a remarkable lower degree of generality). For the construction of $\Zs{G}$, we will take as category of fields $\cC = \CVar$, the category of complex algebraic varieties. The field theory will be constructed from the functor $\Rep{G}: \Embpc \to \CVar$ that assigns, to every pair $(M,A)$ with $M$ a compact manifold and $A \subseteq M$ a finite set, the generalized representation variety $\Rep{G}(M,A)=\Hom(\Pi(M, A), G)$. Recall that $\Pi(M, A)$ is the fundamental groupoid of $(M,A)$, that is, the groupoid of homotopy classes of paths on $M$ between points of $A$. However, the core of the construction is the quantisation part. For it, we will use the Grothendieck rings of mixed Hodge modules that, as we will have proved in Section \ref{sec:saito-mixed-hodge-modules}, can be made a $\CVar$-algebra. Moreover, this construction can be extended to the parabolic case without further modifications, as shown in Section \ref{sec:parabolic-case-TQFT}.

Despite that the TQFT introduced above captures the Hodge structures of representation varieties, with a view towards practical calculations it is not the best option. With this idea, in Section \ref{sec:other-tqft-repr} we will show several variants of the previous TQFT. The most useful on for our purposes will be the geometric version introduced in Section \ref{sec:geometric-reduced-TQFT}. For that, we will need to introduce the category of piecewise algebraic varieties, $\PVar{k}$, as explained in Section \ref{sec:piecewise-regular-varieties}. Roughly speaking, this category has, as objects, the elements of the Grothendieck semi-ring of algebraic varieties and, as morphisms, maps that can be decomposed as disjoint union of regular morphisms.

Finally, in Section \ref{sec:sl2-repr-var}, we will show how this geometric version of TQFT, $\Zg{G}$, can be effectively used for computing Hodge structures of representation varieties. This result can also be found in the paper \cite{GP-2018}. We will focus on the case $G = \SL{2}(\CC)$ and parabolic structures with holonomies in the conjugacy classes of the Jordan type elements, $J_\pm \in \SL{2}(\CC)$ and $\lambda = [J_\pm]$ (see Section \ref{sec:sl2-repr-var} for an introduction to this group). We will show that, in fact, all the computations of the TQFT can be done in a certain finitely generated submodule $\cW$ and we will compute explicitly the homomorphisms $\Zg{\SL{2}(\CC)}(D): \K{\MHSq} \to \cW$, $\Zg{\SL{2}(\CC)}(D^\dag): \cW \to \K{\MHSq}$ and $\Zg{\SL{2}(\CC)}(L_\lambda): \cW \to \cW$ needed for the computation. Here, $L_\lambda: (S^1, \star) \to (S^1, \star)$ is the trivial tube with a puncture in the class $\lambda$. Actually, as shown in Section \ref{sec:discs-tubes}, the first two homomorphisms are rather trivial so the unique challenge is the third one. Section \ref{sec:tube-J+} is devoted to the computation of this morphism $\Zg{\SL{2}(\CC)}(L_\lambda)$. Finally, in Section \ref{sec:genus-tube} we will show that the algorithmic procedure of \cite{MM} is nothing but the calculation of (a slightly modified version of) $\Zg{\SL{2}(\CC)}(L)$. Therefore, putting together the calculations of \cite{MM} and the ones of this Thesis, we conclude the following result (Theorem \ref{thm:Hodge-repr-sl2}).

\begin{thm*}
Let $Q$ be a parabolic structure with $r_+$ punctures with holonomies $[J_+]$, $r_-$ punctures with holonomies $[J_-]$ and $t$ punctures with holonomies in $-\Id$. Denote $r = r_+ + r_-$ and $\sigma = (-1)^{r_- + t}$. Then, the image of the mixed Hodge structure on the cohomology of $\Rep{\SL{2}(\CC)}(\Sigma_g, Q)$ in $\K{\MHSq}$ is:
\begin{itemize}
	\item If $\sigma = 1$, then
\begin{align*}
	\coh{\Rep{\SL{2}(\CC)}(\Sigma_g, Q)} =& \, {\left(q^2 - 1\right)}^{2g + r - 1} q^{2g - 1} +
\frac{1}{2} \, {\left(q -
1\right)}^{2g + r - 1}q^{2g -
1}(q+1){\left({2^{2g} + q - 3}\right)} 
\\ &+ \frac{\left(-1\right)^{r}}{2} \,
{\left(q + 1\right)}^{2g + r - 1} q^{2g - 1} (q-1){\left({2^{2g} +q -1}\right)}.
\end{align*}
	\item If $\sigma = -1$, then
\begin{align*}
	\coh{\Rep{\SL{2}(\CC)}(\Sigma_g, Q)} = &\, {\left(q - 1\right)}^{2g + r - 1} (q+1)q^{2g - 1}{{\left( {\left(q + 1\right)}^{2 \,
g + r-2}+2^{2g-1}-1\right)} } \\
&+ \left(-1\right)^{r + 1}2^{2g - 1}  {\left(q + 1\right)}^{2g + r
- 1} {\left(q - 1\right)} q^{2g - 1}.
\end{align*}
Here, $q = \QQ(-1)=\coh{\CC}$ is the one-dimensional pure Hodge structure concentrated in the piece $(1,1)$.
\end{itemize}
\end{thm*}

However, the final object of desire is not the Hodge structures of representation varieties but of their GIT quotients, the character varieties. In order to fill that gap, Chapter \ref{chap:git} is devoted to GIT-like techniques that allow us to compute Deligne-Hodge polynomials of character varieties from the ones of representation varieties. In order to do so, we will need to consider weak quotients that go beyond the scope of classical GIT (as reviewed in Section \ref{sec:review-git}). The results of this chapter can also be found in \cite{GP-2018-2}.

The problem is the following. Suppose that $X$ is an algebraic variety with an action of an algebraic group $G$ such that the GIT quotient, $\pi: X \to X \sslash G$, is good. Suppose that we can decompose $X =  Y \sqcup U$ with $Y \subseteq X$ closed and $U \subseteq X$ open, both invariant for the action. In general, $X \sslash G \neq (Y \sslash G) \sqcup (U \sslash G)$ so the GIT quotient is not well-behaved under decompositions. The underlying problem is that, while $U \to \pi(U)$ is still a good quotient (so $\pi(U) = U \sslash G$), the closed part $\pi|_{Y}: Y \to \pi(Y)$ may no longer be a good quotient since they could exist $G$-invariant functions on $Y$ that do not factorize through $\pi|_{Y}$. 

Nonetheless, the closed part $Y \to \pi(Y)$ still holds all the topological properties of good quotients. For this reason, in Section \ref{section:pseudo-quotients}, we introduce the notion of pseudo-quotient as a kind of weak quotient regarding only the topological properties of good quotients. Using this approach, we will prove in Section \ref{section:stratification} (Theorem \ref{prop:decomposition-quotient}) that, if $X = Y \sqcup U$ as above and $\pi: X \to \overline{X}$ is a pseudo-quotient for the action of $G$, then $\pi: Y \to \pi(Y)$ and $\pi: U \to \pi(U)$ are pseudo-quotients. Moreover, pseudo-quotients allow us to obtain a sort of quotient even in the case that $G$ is not reductive. In this sense, pseudo-quotients are a simpler tool than the sophisticated techniques used for constructing non-reductive GIT quotients, as in \cite{Berczi-Dolan-Hawes-Kirwan:2016}, \cite{Dolan-Kirwan:2007} or \cite{Kirwan:2009}.

On the other hand, pseudo-quotients are no longer unique but, as we will show in Section \ref{sec:uniqueness-pseudo-quot}, their classes in the Grothendieck ring of complex algebraic varieties are so. In particular, this is enough to ensure that their $E$-polynomials agree. Putting together these two facts, we obtain the following result.

\begin{thm*}
Let $X$ be a complex algebraic variety with a linear action of a reductive group $G$. If we decompose $X = Y \sqcup U$ with $Y \subseteq X$ closed and $U$ orbitwise-closed (see Definition \ref{defn:orbitwise-closed}), then
$$
	\DelHod{X \sslash G} = \DelHod{Y \sslash G} + \DelHod{U \sslash G}.
$$
\end{thm*}

As an application to character varieties, there is a natural decomposition of the representation variety as $\Rep{G}(\Gamma) = \Xred{\Rep{G}}(\Gamma) \sqcup \Xirred{\Rep{G}}(\Gamma)$, where $\Xred{\Rep{G}}(\Gamma)$ denotes the set of reducible representations and $\Xirred{\Rep{G}}(\Gamma)$ the set of irreducible ones. In that case, the results of Section \ref{sec:rep-var} imply that
$$
	\DelHod{\Char{G}(\Gamma)} = \DelHod{\Xred{\Rep{G}}(\Gamma) \sslash G} + \DelHod{\Xirred{\Rep{G}}(\Gamma) \sslash G}.
$$

In this way, each stratum can be analyzed separately. For the stratum $\Xirred{\Rep{G}}(\Gamma)$, the situation is quite simple since the action on it is closed and (essentially) free, so $\DelHod{\Xirred{\Rep{G}}(\Gamma) \sslash G}$ is just the quotient of the $E$-polynomial of $\Xirred{\Rep{G}}(\Gamma)$ over the $E$-polynomial of $G/G^0$, where $G^0 \subseteq G$ is the center of $G$.

For the stratum $\Xred{\Rep{G}}(\Gamma)$, the situation is a bit harder. The idea here is that the closures of the orbits of elements of $\Xred{\Rep{G}}(\Gamma)$ always intersect the subvariety of diagonal representations. This is a geometric situation is treated in Proposition \ref{prop:core} and implies that the GIT quotient is isomorphic to the quotient of the diagonal representations under permutation of eigenvalues. Therefore, for this stratum, the calculation reduces to the analysis of a quotient by a finite group.

Using these ideas, in Section \ref{subsec:reducible-rep} we will recompute, for $G = \SL{2}(\CC)$, the Deligne-Hodge polynomial of character varieties of free groups and surface groups (i.e.\ fundamental groups of closed orientable surfaces) from the ones of the corresponding representation variety, reproving the results of \cite{Cavazos-Lawton:2014} and \cite{MM:2016} respectively. Moreover, in Section \ref{sec:parabolic-rep}, we will explore the parabolic case and we will compute the Deligne-Hodge polynomial of $\SL{2}(\CC)$-parabolic character varieties of free and surface groups with any number of punctures with holonomies of Jordan type. In the case of surface groups, we obtain the following result.

\begin{thm} Let $Q$ be a parabolic structure on $\Sigma_g$ with $r_+$ punctures with holonomies $[J_+]$, $r_-$ punctures with holonomies $[J_-]$ and $t$ punctures with holonomies in $-\Id$. Denote $r = r_+ + r_-$ and $\sigma = (-1)^{r_- + t}$. Then, the Deligne-Hodge polynomials of parabolic $\SL{2}(\CC)$-character varieties are:
\begin{itemize}
	\item If $\sigma = 1$, then
	\begin{align*}
	\DelHod{\Char{\SL{2}(\CC)}(\Sigma_g, Q)} =& \,{\left(q^2 - 1\right)}^{2g + r -
2} q^{2g - 2} +\left(-1\right)^{r} 2^{2g}  {\left(q - 1\right)} q^{2g - 2}
{\left({1-\left(1-q\right)}^{r - 1}\right)}\\
&+ \frac{1}{2}{\left(q - 1\right)}^{2g +r - 2} q^{2g - 2} \, {\left(2^{2g} + q - 3\right)}  \\
&+ \frac{1}{2} {\left(q + 1\right)}^{2g +
r - 2} q^{2g - 2}\,
\left(2^{2g} + q - 1\right).
	\end{align*}
	\item If $\sigma = -1$, then
	\begin{align*}
	\hspace{-1.75cm}\DelHod{\Char{\SL{2}(\CC)}(\Sigma_g, Q)} =& \left(-1\right)^{r-1}2^{2g - 1}  {\left(q + 1\right)}^{2g + r - 2} q^{2g - 2} \\ &+ {\left(q - 1\right)}^{2g + r - 2} q^{2g - 2}\left( {\left(q + 1\right)}^{2g + r - 2} + 2^{2g - 1} - 1\right).
\end{align*} 
\end{itemize}
\end{thm}

From the results of this PhD Thesis, several lines of reseach are opened, as sketched in the final section of this document. Let us review some of them here. First of all, in this Thesis we have focused on the case of Jordan type punctures. The reason is that, when we want to introduce conjugacy classes of semi-simple type, a new interaction phenomenon appears that complicates the calculations. However, we expect that, using the techniques introduced in this Thesis, a detailed study of these interferences will allow us to extend the calculations also to this case.

Moreover, so far we have only centered our attention to the case $G = \SL{2}(\CC)$. However, the TQFTs constructed in this Thesis are valid for any complex algebraic group so a next step would be to consider other groups. As a starting point, it would be interesting to study the case $G = \SL{r}(\CC)$ with $r \geq 2$ and to understand how the TQFT and the GIT quotient behave when we increase the rank.

Furthermore, along this Thesis, we have always focused on character varieties, forgetting about the other faces of the non-abelian Hodge theory. However, if we still use the mixed Hodge modules as quantisation, it may be expectable that some suitable fields theories will allow us to obtain similar results for the moduli spaces of Higgs bundles and of flat connections. That would be very interesting since, in that case, we could capture the whole hyperk\"ahler structure of these moduli spaces.

Another exciting exploration point would be to study the monoidality restriction of the described TQFT. As we mentioned above, the constructed TQFT is no longer a strict monoidal functor but only a lax monoidal functor. However, as the known cases show, all the computations can be performed in a finitely generated submodule. This shows that, morally, the TQFT wants to be monoidal but there are some obstructions that prevent it from being so. For this reason, we expect that the previous construction could be modified to give rise to a strict monoidal TQFT, maybe using a different quantisation or varying the monoidal structures.

Finally, another framework in which character varieties are central is the geometric Langlands program (see \cite{Beilinson-Drinfeld}). In this setting, the Hitchin fibration for the Dolbeault moduli space satisfies the Strominger-Yau-Zaslow conditions of mirror symmetry for Calabi-Yau manifolds (see \cite{Strominger-Yau-Zaslow}), from which it arises several questions about relations between $E$-polynomials of character varieties for Langlands dual groups, as conjectured in \cite{Hausel:2005} and \cite{Hausel-Rodriguez-Villegas:2008}. The validity of these conjectures has been discussed in some cases as in \cite{LMN} and \cite{Martinez:2017}. Despite of that, the general case remains unsolved. We hope that the ideas introduced in this PhD Thesis could be useful to shed light upon these questions.

%% file: Chapters/Introduccion.tex
\label{introduccion} % For referencing the chapter elsewhere, use \ref{Cap\'itulo1} 

\lhead{Introducci\'on} % This is for the header on each page - perhaps a shortened title

%----------------------------------------------------------------------------------------

Sea $M$ una variedad diferenciable compacta, posiblemente con frontera, y sea $G$ un grupo algebraico complejo. El conjunto de representaciones $\rho: \pi_1(M) \to G$ puede ser dotado de la estructura de una variedad algebraica, la llamada variedad de representaciones de $M$ en $G$ y denotada por $\Rep{G}(M)$. Adem\'as, el grupo $G$ act\'ua en $\Rep{G}(M)$ por conjugaci\'on por lo que, tomando el cociente de Teor\'ia de Invariantes Geom\'etricos (GIT, por sus siglas en ingl\'es) de $\Rep{G}(W)$ bajo esta acci\'on, obtenemos la variedad de caracteres
$$
	\cR_G(M) = \Rep{G}(M) \sslash G.
$$
La variedad de caracteres es el espacio de m\'oduli de representaciones de $\pi_1(M)$ en $G$, como se trata en \cite{Nakamoto} y, en este contexto, este espacio es conocido como el espacio de m\'oduli de Betti. Incluso en los casos m\'as sencillos, la topolog\'ia y estructura algebraica de estas variedades de caracteres es extremada rica y ha sido muy estudiada durante los \'ultimos veinte a\~nos.

Uno de los principales motivos para estudiar variedades de caracteres es su importante papel en la correspondencia de Hodge no abeliana. Esta bella teor\'ia muestra que las variedades de caracteres son solo una de las tres caras de una misma realidad. La primera encarnaci\'on de esta teor\'ia es la correspondencia de Riemann-Hilbert (\cite{SimpsonI} \cite{SimpsonII}). Para $M=\Sigma$ una superficie de Riemann compacta y $G = \GL{n}(\CC)$ (resp.\ $G = \SL{n}(\CC)$), un elemento de $\Char{G}(\Sigma)$ define un $G$-sistema local y, por tanto, un fibrado vectorial algebraico de rango $n$, $E \to \Sigma$, de grado cero (resp.\ y determinante fijado) con una conexi\'on plana en \'el. En esta filosof\'ia, la correspondencia de Riemann-Hilbert nos da una una correspondencia anal\'itica-real entre $\Char{G}(\Sigma)$ y el espacio de m\'oduli de fibrados planos en $\Sigma$ de rango $n$ y grado cero (resp.\ y determinante fijado), el conocido como espacio de m\'oduli de de Rham. M\'as a\'un, v\'ia la correspondencia de Hitchin-Kobayashi (\cite{Simpson:1992} \cite{Corlette:1988}), tambi\'en tenemos que, para $G=\GL{n}(\CC)$ (resp.\ $G=\SL{n}(\CC)$), el espacio de m\'oduli de Betti $\Char{G}(\Sigma)$ es equivalente anal\'itico-real al espacio de m\'oduli de Dolbeault, esto es, el espacio de m\'oduli de $G$-fibrados de Higgs de rango $n$ y grado cero (resp.\ y determinante fijado) i.e.\ fibrados vectoriales $E \to \Sigma$ junto con un campo $\Phi: E \to E \otimes K_\Sigma$, llamado el campo de Higgs.

Por medio de estas correspondencias, es posible calcular el polinomio de Poincar\'e de las variedades de caracteres a trav\'es de teor\'ia de Morse. Siguiendo estas ideas, Hitchin, en el paper seminal \cite{Hitchin}, dio el polinomio de Poincar\'e para $G=\SL{2}(\CC)$. Posteriormente, Gothen tambi\'en lo calcul\'o para $G=\SL{3}(\CC)$ en \cite{Gothen} y Garc\'ia-Prada, Heinloth y Schmitt para $G=\GL{4}(\CC)$ en \cite{GP-Heinloth-Schmitt}. 

Sin embargo, estas correspondencias de la teor\'ia de Hodge no abeliana est\'an muy lejos de ser algebraicas. Por esta raz\'on, es importante estudiar la estructura algebraica de las variedades de caracteres. En particular, es interesante analizar una estructura lineal extra que aparece en la cohomolog\'ia de estas variedades, la llamada estructura de Hodge mixta. El origen de estos invariantes se retrotrae a teor\'ia de Hodge cl\'asica, donde se relacionan con cuestiones tan sutilies como las formas arm\'onicas en la variedad. Sin embargo, no fue hasta los trabajos de Deligne en \cite{DeligneI:1971}, \cite{DeligneII:1971} y \cite{DeligneI:1971}, en el contexto de las conjeturas de Weil, cuando las estructuras de Hodge cobraron importancia como invariantes algebraicos. La raz\'on de ello es que, en esos art\'iculos, Deligne prob\'o que la cohomolog\'ia de una variedad algebraica compleja est\'a equipada, de forma natural, con una estructura de Hodge.

Sin entrar en detalles, una estructura de Hodge en una variedad compleja $X$ es una descomposici\'on de la cohomolog\'ia con soporte compacto, $H_c^k(X;\CC)$, como un ret\'iculo de piezas $H_c^{k;p,q}(X)$ con $p,q \in \ZZ$. La dimensi\'on de estas piezas, $h_c^{k;p,q}(X) = \dim_\CC H_c^{k;p,q}(X)$, son los llamados n\'umeros de Hodge y pueden pensarse como un refinamiento de los n\'umeros de Betti. En particular, a partir de estos n\'umeros de Hodge, podemos construir el polinomio de Deligne-Hodge, o $E$-polinomio, de $X$ como una suma alternada al estilo de una caracter\'istica de Euler
$$
	\DelHod{X} = \sum_{k, p, q} (-1)^k\,h_c^{k;p,q}(X)\,u^pv^q \in \ZZ[u^{\pm 1}, v^{\pm 1}].
$$
Sin embargo, en general, el polinomio de Deligne-Hodge de las variedades de caracteres es desconocido. Uno de los mayores avances en este sentido fue llevado a cabo por Hausel y Rodr\'iguez-Villegas en \cite{Hausel-Rodriguez-Villegas:2008} mediante un teorema de Katz de sabor aritm\'etico basado en las conjeturas de Weil y el principio de Lefschetz. En ese art\'iculo, los autores calcularon una expresi\'on del $E$-polinomio para $\GL{n}(\CC)$-variedades de caracteres torcidas (twisted) en superficies compactas en t\'erminos de funciones generatrices. Esta t\'ecnica fue posteriormente expandida en \cite{Mereb} para $G=\SL{n}(\CC)$ y, recientemente, en \cite{Baraglia-Hekmati:2016} para el caso no-torcido y superficies orientables ($G=\GL{3}(\CC)$ y $\SL{3}(\CC)$) y para variedades no orientables ($G=\GL{2}(\CC)$ y $\SL{2}(\CC)$). No obstante, esta aproximaci\'on adolece de un grave problema porque, para poder usar el teorema de Katz, es necesario contar el n\'umero de puntos de variedades de caracteres sobre cuerpos finitos y esto requiere usar la tabla de caracteres de los grupos finitos $\GL{n}(\mathbb{F}_q)$. Sin embargo, para rango alto, estas tablas no pueden calcularse expl\'icitamente y, por tanto, solo se puede dar una expresi\'on del $E$-polinomio en t\'erminos de funciones generatrices.

Tratando de superar este escollo, Logares, Mu\~noz y Newstead en \cite{LMN} iniciaron una nueva aproximaci\'on a este problema. En ese art\'iculo, la estrategia consiste en centrarse en el c\'alculo del polinomio de Deligne-Hodge de las variedades de representaciones. Para ello, la variedad se trocea en estratos m\'as sencillos de los que s\'i se conoce su $E$-polinomio y luego se suman todas esas contribuciones. Tras ello, estudiaron qu\'e identificaciones deben realizarse para pasar al cociente GIT y, a partir de esta informaci\'on, calcularon el polinomio de Deligne-Hodge de las variedades de caracteres correspondientes. Siguiendo este principio, se han calculado expresiones expl\'icitas de los $E$-polinomios para g\'enero $g=1,2$ y $G=\SL{2}(\CC)$ en \cite{LMN} y, posteriormente, se han extendido a g\'enero arbitrario en \cite{MM} y, para $G=\textrm{PGL}_2(\CC)$, en \cite{Martinez:2017}.

Un paso m\'as all\'a en las variedades de representaciones es considerar una estructura parab\'olica sobre una superficie $\Sigma$. \'Estas no son m\'as que una colecci\'on finita de pares $Q=\left\{(p_i, \lambda_i)\right\}_{i=1}^s$, donde $p_i \in \Sigma$ son los llamados puntos marcados y $\lambda_i \subseteq G$ es una subvariedad cerrada bajo conjugaci\'on por $G$, llamada la holonom\'ia de $p_i$. En ese caso, la variedad de representaciones parab\'olica $\Rep{G}(\Sigma, Q)$ es el conjunto de representaciones $\rho: \pi_1(\Sigma - \left\{p_1, \ldots, p_s\right\}) \to G$ tales que la imagen de los lazos positivos alrededor de los $p_i$ viven en $\lambda_i$.

Estas contrapartidas parab\'olicas son importantes porque la correspondencia de Hodge no abeliana tambi\'en se extiende a este contexto. Por ejemplo, para $G=\GL{n}(\CC)$ y $Q$ un solo punto marcado con holonom\'ia una ra\'iz primitiva de la unidad $e^{2 \pi i d/n}$ (el llamado caso torcido), se tiene que $\cR_{\GL{n}(\CC)}(\Sigma, Q)$ es difeomorfo al espacio de m\'oduli de fibrados de Higgs de rango $n$ y grado $d$, y espacio de m\'oduli de fibrados logar\'itmicos planos de rango $n$ con un polo simple de residuo $-\frac{d}{n}\Id$. An\'alogas correspondencias existen tambi\'en para clases de conjugaci\'on semi-simples gen\'ericas, \cite{Simpson:parabolic}. A trav\'es de estas correspondencias, Boden y Yokogawa calcularon el polinomio de Poincar\'e de variedades de caracteres parab\'olicas en \cite{Boden-Yokogawa} para $G=\SL{2}(\CC)$ y clases de conjugaci\'on semisimples, y Garc\'ia-Prada, Gothen y Mu\~noz para $G=\GL{3}(\CC)$ y $G=\SL{3}(\CC)$ en \cite{GP-Gothen-Munoz}. En el caso general torcido, en \cite{Schiffmann:2016} y \cite{Mozgovoy-Schiffmann} (v\'ease tambi\'en \cite{Mellit}) se calcula una f\'ormula combinatoria para rango arbitrario $G=\GL{n}(\CC)$ siempre que $n$ y $d$ sean coprimos. Incidentalmente, otras correspondencias excepcionales pueden aparecer, como para $G=\SL{2}(\CC)$, $\Sigma$ una curva el\'iptica y $Q$ dos puntos marcados con diferente holonom\'ia semi-simple, en cuyo caso $\cR_G(\Sigma, Q)$ es difeomorfo al espacio de m\'oduli de instantones doblemente peri\'odicos a trav\'es de la transformada de Nahm \cite{Biquard-Jardim} \cite{Jardim}.

No obstante, el calculo de $E$-polinomios para variedades de caracteres parab\'olicas es mucho m\'as dif\'icil que para el caso no parab\'olico debido a las nuevas perturbaciones generadas por los puntos eliminados. Algunos avances en este sentido se han producido en \cite{Hausel-Letellier-Villegas}, mediante el m\'etodo aritm\'etico, para $G=\GL{n}(\CC)$ y puntos marcados gen\'ericos semi-simples (pero el resultado no es expl\'icito) y en \cite{LM} con el m\'etodo geom\'etrico pero, a lo sumo, dos puntos marcados. El caso general sigue abierto.

En esta tesis doctoral, proponemos una nueva aproximaci\'on al problema de calcular las estructuras de Hodge de variedades de representaciones basada en Teor\'ias Topol\'ogicas de Campos Cu\'anticos (Topological Quantum Field Theories o TQFTs). El germen de esta idea es el art\'iculo \cite{MM}, en el que los autores calculan el polinomio de Deligne-Hodge de $\SL{2}(\CC)$-variedades de representaciones por medio del m\'etodo geom\'etrico. Sin embargo, m\'as all\'a del c\'alculo en s\'i mismo, el punto fundamental es que probaron que, si $\Sigma_g$ es la superficie compacta de g\'enero $g$, existe un algoritmo recursivo que calcula el polinomio de Deligne-Hodge de $\Rep{\SL{2}(\CC)}(\Sigma_g)$ en t\'erminos del de $\Rep{\SL{2}(\CC)}(\Sigma_{g-1})$ y una pieza de informaci\'on extra conocida como representaci\'on de monodrom\'ia de Hodge. Esta naturaleza recursiva puede tambi\'en ser observada en otros trabajos, como en \cite{Mozgovoy:2012}, \cite{Hausel-Letellier-Villegas:2013}, \cite{Diaconescu:2017} y \cite{Carlsson-Rodriguez-Villegas}. Este patr\'on evoca fuertemente a las TQFTs.

El origen de las TQFTs se remonta al art\'iculo seminal de Witten \cite{Witten:1988}, en el que mostr\'o que el polinomio de Jones (un invariante de nudos) puede obtenerse a trav\'es de una teor\'ia de campos cu\'anticos conocida como teor\'ia de Chern-Simons. Consciente de la importancia de este descubrimiento, Atiyah en \cite{Atiyah:1988} dio una fundamentaci\'on completa de este tipo de teor\'ias y acu\~n\'o el t\'ermino Topological Quantum Field Theory. En ese art\'iculo, se enfatiza la naturaleza categ\'orica de este tipo de construcciones y las formula como functores monoidales $Z: \Bord{n} \to \Mod{R}$ de la categor\'ia de bordismos $n$-dimensionales a la categor\'ia de $R$-m\'odulos.

\begin{figure}[h]
\vspace{0cm}
\begin{center}
\includegraphics[scale=0.55]{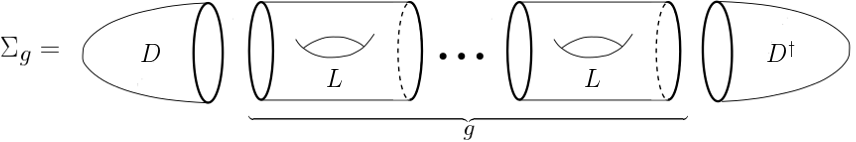}
\caption{}
\label{img:decomp-tubes-intro-spanish}
\vspace{-0.6cm}
\end{center}
\end{figure}

Sin embargo, para nuestros prop\'ositos, la aplicaci\'on m\'as importante de las TQFTs es que proporcionan un m\'etodo efectivo de c\'alculo de invarientes algebraicos. Para verlo, supongamos que estamos interesados en un invariante algebraico que toma valores en un cierto anillo $R$, y supongamos que disponemos de una TQFT, $Z: \Bord{n} \to \Mod{R}$, tal que, si $M$ es una variedad cerrada $n$-dimensional, vista como un bordismo $M: \emptyset \to \emptyset$, tenemos que $Z(M)(1) \in R$ es el invariante de $M$ buscado. Obs\'ervese que, como $Z$ es monoidal, $Z(\emptyset)=R$ por lo que $Z(M)$ es un automorfismo $R$-lineal de $R$ i.e.\ es multiplicar por un cierto elemento fijo de $R$. En ese caso, podemos descomponer $M$ en piezas m\'as sencillas para obtener el invariante de forma recursiva. Por ejemplo, supongamos que $M= \Sigma_g$ y descompongamos $M = D^\dag \circ L^g \circ D$, donde $D: \emptyset \to S^1$ es un disco, $L: S^1 \to S^1$ es un toro con dos discos recortados y $D^\dag: S^1 \to \emptyset$ es un disco en el sentido opuesto a $D$, tal y como se muestra en la Figura \ref{img:decomp-tubes-intro-spanish}.

\renewcommand{\figurename}{Figura}

\renewcommand{\figurename}{Figure}

En ese caso, la TQFT nos aporta una descomposici\'on $Z(\Sigma_g) = Z(D^\dag) \circ Z(L)^g \circ Z(D)$, donde $Z(D),Z(D^\dag)$ y $Z(L)$ son homomorfismos de $R$-m\'odulos fijados e independientes de $g$. De este modo, el invariante deseado $Z(\Sigma_g)(1)$ puede ser calculado autom\'aticamente para todas las superficies $\Sigma_g$ sin m\'as que el entender estos tres homomorfismos. En este sentido, una TQFT puede pensarse como la versi\'on din\'amica de algo est\'atico como un invariante algebraico.

Esta tesis doctoral se estructura en cuatro cap\'itulos. El Cap\'itulo \ref{chap:tqft} est\'a dedicado a las Teor\'ias Topol\'ogicas de Campos Cu\'anticos. En la Secci\'on \ref{sec:categories-tqft}, y especialmente en \ref{sec:low-dim-tqft}, revisaremos algunas de la definiciones y propiedades fundamentales de las TQFTs, incluido el hecho de que, para todo $M \in \Bord{n}$, el m\'odulo $Z(M)$ est\'a obligado a ser finitamente generado. Esta restricci\'on es crucial en tanto que impone serveras obstrucciones a la construcci\'on de TQFTs pues, en muchos casos, la asignaci\'on natural para $Z(M)$ no es finitamente generada. Para salvar este problema, en la Secci\'on \ref{sec:lax-monoidality} proponemos relajar la condici\'on de monoidalidad y considerar TQFTs monoidales laxas (lax monoidal), esto es, functores $Z: \Bord{n} \to \Mod{R}$ tales que, para todo $M, N \in \Bord{n}$, existe un homomorfismo $\Delta: Z(M) \otimes Z(N) \to Z(M \sqcup N)$ que podr\'ia no ser un isomorfismo.

Bajo estas condiciones relajadas, la vida se vuelve mucho m\'as sencilla y podemos construir multitud de TQFTs. En esta direcci\'on, en la Secci\'on \ref{sec:lax-monoidal-tqft} describimos un nuevo procedimiento, de inspiraci\'on f\'isica, que crea una TQFT monoidal laxa a partir de dos piezas m\'as sencillas: una teor\'ia de campos (informaci\'on geom\'etrica) y una cuantizaci\'on (informaci\'on algebraica). Para ello, la idea es considerar una categor\'ia $\cC$ con objeto final y pullbacks que jugar\'a el papel de una cierta ``categor\'ia de campos'' (en el sentido f\'isico), y descomponer nuestra TQFT como una composici\'on
$$
	\Bord{n} \stackrel{\cF}{\longrightarrow} \Span{\cC} \stackrel{\cQ}{\longrightarrow} \Mod{R},
$$
donde $\Span{\cC}$ es la categor\'ia de los spans de $\cC$ (v\'ease Secci\'on \ref{sec:lax-monoidal-tqft}). El functor $\cF$ es la ``teor\'ia de campos'' y puede ser construido a partir de otro functor $\Embc \to \cC$ que preserve pegados (Secci\'on \ref{sec:field-theory}). Obs\'ervese que, aqu\'i, $\Embc$ es la categor\'ia amplia (wide subcategory) de variedades diferenciables compactas pero \'unicamente embebimientos entre ellas. Por otra parte, el functor $\cQ$ es la ``cuantizaci\'on" y puede construirse a partir de una pieza de informaci\'on algebraica conocida como $\cC$-\'algebra (Secciones \ref{sec:C-algebras} y \ref{sec:quantisation}). Para terminar el Cap\'itulo \ref{chap:tqft}, en la Secci\'on \ref{sec:other-versions-tqft} describimos c\'omo este procedimiento puede extenderse para construir otros tipos de Teor\'ias Topol\'ogicas de Campos Cu\'anticos, como TQFTs que preservan alg\'un tipo de estructura extra (Secci\'on \ref{sec:TQFT-over-sheaf}), TQFTs definidas \'unicamente en tubos y no en bordismos generales (Secci\'on \ref{sec:almost-TQFT-tubes}), versiones simplificadas de TQFTs a partir de informaci\'on heur\'istica (Secci\'on \ref{sec:reduction-TQFT}) y un tipo de versi\'on d\'ebil de TQFTs conocidas como TQFTs blandas (soft TQFTs), disponibles cuando solo la mitad de la $\cC$-\'algebra est\'a a nuestro alcance (Secci\'on \ref{sec:soft-TQFT}). Estas versiones alternativas de TQFTs ser\'an de gran utilidad en el Cap\'itulo \ref{chap:representation}, donde construiremos una TQFT que calcule estructuras de Hodge en variedades de representaciones.

En el Cap\'itulo \ref{chap:hodge}, nos centraremos en la teor\'ia de Hodge como el marco natural en el que encajan los invariantes algebraicos que queremos estudiar. Tras algunos preliminares sobre $K$-teor\'ia y anillos de Grothendieck (Secci\'on \ref{sec:panoramic-Ktheory}), revisaremos las estructuras de Hodge en la Secci\'on \ref{sec:hodge-structures}, tanto de tipo puro (como en el caso de variedades K\"ahler compactas) como de tipo mixto (en el caso general).

Sin embargo, el n\'ucleo del Cap\'itulo \ref{chap:hodge} es la Secci\'on \ref{sec:mixed-hodge-modules}, en la cual revisamos la teor\'ia de Saito de m\'odulos de Hodge mixtos. Estos m\'odulos de Hodge mixtos son la piedra angular de nuestra construcci\'on de una TQFT que calcule las estructuras de Hodge de variedades de representaciones, en tanto que sus anillos de Grothendieck ser\'an la $\CVar$-\'algebra escogida para la cuantizaci\'on. Como preparaci\'on, en la Secci\'on \ref{sec:riemann-hilbert} revisaremos la teor\'ia de $\cD$-m\'odulos y haces perversos (perverse sheaves) as\'i como sus interrelaciones a trav\'es de la correspondencia de Riemann-Hilbert. Estos objetos son enormemente necesarios porque, como explicamos en la Secci\'on \ref{sec:saito-mixed-hodge-modules}, son precisamente los bloques b\'asicos de la construcci\'on de m\'odulos de Hodge mixtos. Grosso modo, la categor\'ia de m\'odulos de Hodge mixtos no es m\'as que una subcategor\'ia apropiada de la categor\'ia de pares filtrados de $\cD$-m\'odulos y haces perversos. No obstante, la regla de selecci\'on que decide si uno de esos pares merece la denominaci\'on de m\'odulo de Hodge mixto es extremadamente compleja y requiere el uso de m\'ultiples functores que controlan, de forma muy sutil, el comportamiento de diversas filtraciones. En esta secci\'on, no pretendemos mostrar ning\'un tipo de trabajo original, sino \'unicamente una reinterpretaci\'on de la teor\'ia de Saito como una caja negra con las funciones apropiadas para nuestros prop\'ositos.

De hecho, la definici\'on de los m\'odulos de Hodge mixtos es tan dif\'icil que pr\'acticamente no se pueden construir m\'odulos de Hodge mixtos a mano. Sin embargo, en la Secci\'on \ref{sec:monodromy-as-mhm} revisaremos algunos de los resultados b\'asicos de Saito que afirman que la categor\'ia de variaciones de estructuras de Hodge (i.e.\ haces cuyos tallos (stalks) son estructuras de Hodge) puede verse como una subcategor\'ia de la categor\'ia de los m\'odulos de Hodge mixtos. A partir de este resultado, dada una fibraci\'on $f: X \to B$ de variedades complejas lisas que es localmente trivial en la topolog\'ia anal\'itica, construiremos un m\'odulo de Hodge mixto en $B$ que controle la monodrom\'ia de $f$. Llamaremos a ese m\'odulo de Hodge la representaci\'on de monodrom\'ia de Hodge de $f$ y, como veremos, generaliza las representaciones de monodrom\'ia de Hodge introducidas en \cite{LMN}.

El Cap\'itulo \ref{chap:representation} est\'a dedicado a las variedades de representaciones. En primer lugar, en la Secci\'on \ref{sec:intro-representation} revisaremos los fundamentos de la variedades de representaciones y los situaremos en su contexto, esbozando los hitos de la teor\'ia de Hodge no abeliana (Secci\'on \ref{sec:non-abelian-hodge}). Sin embargo, el coraz\'on de esta disertaci\'on es la Secci\'on \ref{sec:TQFT-for-repr}, en la cual probamos el siguiente resultado fundamental (Teorema \ref{thm:existence-s-LTQFT}).

\begin{thm*}
Sea $G$ un grupo algebraico complejo. Para todo $n \geq 1$, existe una TQFT monoidal laxa de pares, $\Zs{G}: \Bordp{n} \to \Modt{\K{\MHSq}}$, tal que, para toda variedad $n$-dimensional cerrada y conexa $W$ y todo subconjunto finito no vac\'io $A \subseteq W$, tenemos
$$
	\Zs{G}(W, A)\left(1\right) = \coh{\Rep{G}(W); \QQ} \otimes \coh{G; \QQ}^{|A|-1} \in \K{\MHSq}.
$$ 
Aqu\'i, $\MHSq$ denota la categor\'ia de las estructuras de Hodge (mixtas) y cuya unidad de su anillo de Grothendieck es $1 \in \K{\MHSq}$ y $\Modt{\K{\MHSq}}$ es una versi\'on ampliada de la categor\'ia de m\'odulos usuales con una estructura de $2$-categor\'ia extra (v\'ease Ejemplo \ref{ex:categories})
\end{thm*}

Este resultado, as\'i como los m\'etodos de construcci\'on de TQFTs de la Secci\'on \ref{sec:lax-monoidal-tqft}, aparecieron por primera vez en el art\'iculo \cite{GPLM-2017} (aunque con un notablemente menor grado de generalidad). Para la construcci\'on de $\Zs{G}$, tomaremos como categor\'ia de campos $\cC = \CVar$, la categor\'ia de las variedades algebraicas complejas. La teor\'ia de campos se construye a partir del functor $\Rep{G}: \Embpc \to \CVar$ que asigna, a cada par $(M,A)$ con $M$ una variedad compacta y $A \subseteq M$ un subconjunto finito, la variedad de representaciones generalizada $\Rep{G}(M,A)=\Hom(\Pi(M, A), G)$. Recu\'erdese que $\Pi(M, A)$ denota el grupoide fundamental del par $(M,A)$, esto es, el grupoide de clases de homotop\'ia de caminos de $M$ entre puntos de $A$. Sin embargo, la clave de esta construcci\'on es la cuantizaci\'on. Para ella, usaremos el anillo de Grothendieck de m\'odulos de Hodge mixtos que, como habremos probado en la Secci\'on \ref{sec:saito-mixed-hodge-modules}, forma una $\CVar$-\'algebra. M\'as a\'un, esta construcci\'on puede extenderse al caso parab\'olico sin grandes modificaciones, como mostramos en la Secci\'on \ref{sec:parabolic-case-TQFT}.

A pesar de que esta TQFT captura las estructuras de Hodge de las variedades de representaciones, con vista a los c\'alculos expl\'icitos no es la mejor opci\'on. Con esta idea, en la Secci\'on \ref{sec:other-tqft-repr} mostraremos diversas variantes de la TQFT anterior. La m\'as \'util para nuestros prop\'ositos ser\'a la versi\'on geom\'etrica introducida en la  Secci\'on \ref{sec:geometric-reduced-TQFT}. Para construirla, necesitaremos introducir la categor\'ia de las variedades algebraicas a trozos, $\PVar{k}$, como explicamos en la Secci\'on \ref{sec:piecewise-regular-varieties}. En l\'ineas generales, esta categor\'ia tiene, como objetos, los elementos del semi-anillo de Grothendieck de variedades algebraicas y, como morfismos, aplicaciones que pueden descomponerse como uni\'on disjunta de morfismos regulares.

Finalmente, en la Secci\'on \ref{sec:sl2-repr-var} mostreremos que la versi\'on geom\'etrica de la TQFT, $\Zg{G}$, puede ser utilizada para realizar c\'alculos efectivos de estructuras de Hodge en variedades de representaciones. Este resultado tambi\'en puede encontrarse en el art\'iculo \cite{GP-2018}. Nos centraremos en el caso $G = \SL{2}(\CC)$ y estructuras parab\'olicas con holonom\'ias en clases de conjugaci\'on de elementos de tipo Jordan, $J_\pm \in \SL{2}(\CC)$ y $\lambda = [J_\pm]$ (v\'ease Secci\'on \ref{sec:sl2-repr-var} para una revisi\'on de este grupo). Mostraremos que, de hecho, todos los c\'alculos de la TQFT pueden realizarse en un cierto subm\'odulo $\cW$ finitamente generado y calcularemos expl\'icitamente los homomorfismos $\Zg{\SL{2}(\CC)}(D): \K{\MHSq} \to \cW$, $\Zg{\SL{2}(\CC)}(D^\dag): \cW \to \K{\MHSq}$ y $\Zg{\SL{2}(\CC)}(L_\lambda): \cW \to \cW$ necesarios para el c\'omputo. Aqu\'i, $L_\lambda: (S^1, \star) \to (S^1, \star)$ denota el tubo trivial con un punto marcado con holonom\'ia en la clase $\lambda$. De hecho, como se muestra en la Secci\'on \ref{sec:discs-tubes}, los dos primeros homomorfismos son bastante triviales, por lo que el verdadero desaf\'io se encuentra en el tercero de ellos. En este sentido, la Secci\'on \ref{sec:tube-J+} est\'a dedicada al c\'alculo de este morfismo $\Zg{\SL{2}(\CC)}(L_\lambda)$. Finalmente, en la Secci\'on \ref{sec:genus-tube} explicaremos que el m\'etodo algor\'itmico de \cite{MM} no es sino el c\'alculo de (una versi\'on ligeramente modificada de) $\Zg{\SL{2}(\CC)}(L)$. De este modo, juntado los c\'alculos de \cite{MM} con los de esta tesis doctoral, concluimos el siguiente resultado (Teorema \ref{thm:Hodge-repr-sl2}).

\begin{thm*}
Sea $Q$ una estructura parab\'olica con $r_+$ puntos marcados con holonom\'ias $[J_+]$, $r_-$ puntos marcados con holonom\'ias $[J_-]$ y $t$ puntos marcados con holonom\'ias $-\Id$. Denotemos $r = r_+ + r_-$ y $\sigma = (-1)^{r_- + t}$. Entonces, la imagen en $\K{\MHSq}$ de la estructura de Hodge mixta en la cohomolog\'ia de $\Rep{\SL{2}(\CC)}(\Sigma_g, Q)$ es:
\begin{itemize}
	\item Si $\sigma = 1$, entonces
\begin{align*}
	\coh{\Rep{\SL{2}(\CC)}(\Sigma_g, Q)} =& \, {\left(q^2 - 1\right)}^{2g + r - 1} q^{2g - 1} +
\frac{1}{2} \, {\left(q -
1\right)}^{2g + r - 1}q^{2g -
1}(q+1){\left({2^{2g} + q - 3}\right)} 
\\ &+ \frac{\left(-1\right)^{r}}{2} \,
{\left(q + 1\right)}^{2g + r - 1} q^{2g - 1} (q-1){\left({2^{2g} +q -1}\right)}.
\end{align*}
	\item Si $\sigma = -1$, entonces
\begin{align*}
	\coh{\Rep{\SL{2}(\CC)}(\Sigma_g, Q)} = &\, {\left(q - 1\right)}^{2g + r - 1} (q+1)q^{2g - 1}{{\left( {\left(q + 1\right)}^{2 \,
g + r-2}+2^{2g-1}-1\right)} } \\
&+ \left(-1\right)^{r + 1}2^{2g - 1}  {\left(q + 1\right)}^{2g + r
- 1} {\left(q - 1\right)} q^{2g - 1}.
\end{align*}
Aqu\'i, $q = \QQ(-1)=\coh{\CC}$ es la estructura de Hodge pura unidimensional concentrada en la pieza $(1,1)$.
\end{itemize}
\end{thm*}

Sin embargo, recordemos que el objeto del deseo no son las estructuras de Hodge de variedades de representaciones, sino de su cociente GIT, las variedades de caracteres. Con el fin de rellenar ese hueco, el Cap\'itulo \ref{chap:git} est\'a dedicado a t\'ecnicas de tipo GIT que nos permitan calcular el polinomio de Deligne-Hodge de variedades de caracteres a partir del de variedades de representaciones. Para ello, necesitaremos considerar versiones d\'ebiles de cocientes que van m\'as all\'a del \'ambito de la Teor\'ia de Invariantes Geom\'etricos cl\'asica (que revisamos en la Secci\'on \ref{sec:review-git}). Los resultados de este cap\'itulo tambi\'en se encuentran recogidos en el art\'iculo \cite{GP-2018-2}.

El problema fundamental es el siguiente. Supongamos que $X$ es una variedad algebraica con una acci\'on de un grupo algebraico $G$ tal que su cociente GIT, $\pi: X \to X \sslash G$, es bueno. Supongamos que descomponemos $X =  Y \sqcup U$ con $Y \subseteq X$ cerrado y $U \subseteq X$ abierto, ambos invariantes bajo la acci\'on. En general, $X \sslash G \neq (Y \sslash G) \sqcup (U \sslash G)$ por lo que el cociente GIT no se comporta bien respecto a descomposiciones. El problema subyacente es que, mientras que $U \to \pi(U)$ sigue siendo un buen cociente, (y, por tanto, $\pi(U) = U \sslash G$), la parte cerrada $\pi|_{Y}: Y \to \pi(Y)$ podr\'ia dejar de ser un buen cociente porque podr\'ian existir funciones $G$-invariantes en $Y$ que no factorizasen a trav\'es de $\pi|_{Y}$. 

A pesar de ello, la parte cerrada $Y \to \pi(Y)$ a\'un conserva todas las propiedades topol\'ogicas de los buenos cocientes. Por esta raz\'on, en la Secci\'on \ref{section:pseudo-quotients}, introducimos la noci\'on de pseudo-cociente como un tipo de cociente d\'ebil que \'unicamente tiene en cuenta las propiedades top\'ologicas de los buenos cocientes. Usando esta aproximaci\'on, en la Secci\'on \ref{section:stratification} (Teorema \ref{prop:decomposition-quotient}) probaremos que, si $X = Y \sqcup U$ como antes y $\pi: X \to \overline{X}$ es un pseudo-cociente para la acci\'on de $G$, entonces $\pi: Y \to \pi(Y)$ y $\pi: U \to \pi(U)$ son pseudo-cocientes. M\'as a\'un, los pseudo-cocientes nos permiten obtener una suerte de cociente incluso cuando el grupo $G$ no es reductivo. En este sentido, los pseudo-cocientes son herramientas m\'as sencillas que las sofisticadas t\'ecnicas empleadas para construir cocientes GIT no reductivos, como en \cite{Berczi-Dolan-Hawes-Kirwan:2016}, \cite{Dolan-Kirwan:2007} o \cite{Kirwan:2009}.

En contrapartida, los pseudo-cocientes no son \'unicos pero, como vermos en la Secci\'on \ref{sec:uniqueness-pseudo-quot}, sus clases en el anillo de Grothendieck de variedades complejas s\'i lo son, lo cual es suficiente para garantizar que sus $E$-polinomios coinciden. Por ello, uniendo estos dos hechos, obtenemos el siguiente resultado.

\begin{thm*}
Sea $X$ una variedad compleja algebraica con una acci\'on lineal de un grupo reductivo $G$. Si descomponemos $X = Y \sqcup U$ con $Y \subseteq X$ cerrado y $U$ cerrado por \'orbitas (orbitwise-closed, v\'ease Definici\'on \ref{defn:orbitwise-closed}), entonces
$$
	\DelHod{X \sslash G} = \DelHod{Y \sslash G} + \DelHod{U \sslash G}.
$$
\end{thm*}

Como aplicaci\'on a variedades de caracteres, tenemos una descomposici\'on natural de la variedad de representaciones como $\Rep{G}(\Gamma) = \Xred{\Rep{G}}(\Gamma) \sqcup \Xirred{\Rep{G}}(\Gamma)$, donde $\Xred{\Rep{G}}(\Gamma)$ denota el conjunto de representaciones reducibles y $\Xirred{\Rep{G}}(\Gamma)$ el conjunto de las irreducibles. En este caso, los resultados de la Secci\'on \ref{sec:rep-var} implican que
$$
	\DelHod{\Char{G}(\Gamma)} = \DelHod{\Xred{\Rep{G}}(\Gamma) \sslash G} + \DelHod{\Xirred{\Rep{G}}(\Gamma) \sslash G}.
$$
De este modo, cada estrato puede ser analizado separadamente. Para el estrato $\Xirred{\Rep{G}}(\Gamma)$, la situaci\'on es bastante simple porque la acci\'on en \'el es (esencialmente) libre, luego $\DelHod{\Xirred{\Rep{G}}(\Gamma) \sslash G}$ no es m\'as que el cociente del $E$-polinomio de $\Xirred{\Rep{G}}(\Gamma)$ sobre el $E$-polinomio de $G/G^0$, donde $G^0 \subseteq G$ es el centro de $G$.

Para el estrato $\Xred{\Rep{G}}(\Gamma)$, la situaci\'on es un poco m\'as complicada. La idea aqu\'i es que la clausura de las \'orbitas de los elementos de $\Xred{\Rep{G}}(\Gamma)$ siempre interseca a la subvariedad de las representaciones diagonales. \'Esta es una situaci\'on geom\'etrica que es tratada en la Proposici\'on \ref{prop:core} e implica que el cociente GIT sea isomorfo al cociente de las representaciones diagonales bajo la permutaci\'on de sus autovalores. De este modo, para este estrato, el c\'alculo se reduce al an\'alisis del cociente por un grupo finito.

Usando estas ideas, en la Secci\'on \ref{subsec:reducible-rep} recalcularemos, para $G = \SL{2}(\CC)$, el polinomio de Deligne-Hodge de variedades de caracteres de grupos libres y grupos fundamentales de superficies cerradas orientables (grupos de superficie) a partir de los correspondientes para variedades de representaciones, reprobando los resultados de \cite{Cavazos-Lawton:2014} y \cite{MM:2016} respectivamente. M\'as a\'un, en la Secci\'on \ref{sec:parabolic-rep}, exploraremos el caso parab\'olico y calcularemos el polinomio de Deligne-Hodge de $\SL{2}(\CC)$-variedades de caracteres parab\'olicas de grupos libres y grupos de superficie con cualquier n\'umero de puntos marcados con holonom\'ias de tipo Jordan. En el caso de grupos de superficie, obtenemos el siguiente resultado.

\begin{thm} Sea $Q$ una estructura parab\'olica sobre $\Sigma_g$ con $r_+$ puntos marcados con holonom\'ias $[J_+]$, $r_-$ puntos marcados con holonom\'ias $[J_-]$ y $t$ puntos marcados con holonom\'ias $-\Id$. Denotemos $r = r_+ + r_-$ y $\sigma = (-1)^{r_- + t}$. Entonces, el polinomio de Deligne-Hodge de $\SL{2}(\CC)$-variedades de caracteres parab\'olicas son:
\begin{itemize}
	\item Si $\sigma = 1$, entonces
	\begin{align*}
	\DelHod{\Char{\SL{2}(\CC)}(\Sigma_g, Q)} =& \,{\left(q^2 - 1\right)}^{2g + r -
2} q^{2g - 2} +\left(-1\right)^{r} 2^{2g}  {\left(q - 1\right)} q^{2g - 2}
{\left({1-\left(1-q\right)}^{r - 1}\right)}\\
&+ \frac{1}{2}{\left(q - 1\right)}^{2g +r - 2} q^{2g - 2} \, {\left(2^{2g} + q - 3\right)}  \\
&+ \frac{1}{2} {\left(q + 1\right)}^{2g +
r - 2} q^{2g - 2}\,
\left(2^{2g} + q - 1\right).
	\end{align*}
	\item Si $\sigma = -1$, entonces
	\begin{align*}
	\hspace{-1.75cm}\DelHod{\Char{\SL{2}(\CC)}(\Sigma_g, Q)} =& \left(-1\right)^{r-1}2^{2g - 1}  {\left(q + 1\right)}^{2g + r - 2} q^{2g - 2} \\ &+ {\left(q - 1\right)}^{2g + r - 2} q^{2g - 2}\left( {\left(q + 1\right)}^{2g + r - 2} + 2^{2g - 1} - 1\right).
\end{align*} 
\end{itemize}
\end{thm}

A partir de los resultados de esta tesis doctoral, varias lineas de investigaci\'on se abren, tal y como son esbozadas en la \'ultima secci\'on de este documento. Revisemos brevemente algunas de ellas aqu\'i. En primer lugar, en esta tesis doctoral nos hemos centrado en el caso de puntos marcados de tipo Jordan. La raz\'on para ello es que, en el caso de clases de conjugaci\'on semi-simples, surge un nuevo fen\'omeno de interacci\'on que complica los c\'alculos. Sin embargo, es de esperar que, con las t\'ecnicas introducidas en esta tesis doctoral, sea posible extender los c\'alculos tambi\'en a este caso a trav\'es de un estudio detalado de estas interferencias.

Adem\'as, hasta ahora hemos centrado nuestra atenci\'on en el caso $G = \SL{2}(\CC)$. Sin embargo, la TQFT construida en esta disertaci\'on es v\'alida para cualquier grupo algebraico complejo, por lo que el siguiente paso natural ser\'ia considerar otros grupos. Como punto de partida, ser\'ia interesante estudiar el caso $G = \SL{r}(\CC)$ con $r \geq 2$ y entender c\'omo la TQFT y el cociente GIT se comportan cuando incrementamos en rango.

M\'a a\'un, a lo largo de esta tesis doctoral, nos hemos centrado en variedades de caracteres, olvid\'andonos completamente de las otras caras de la teor\'ia de Hodge no abeliana. Sin embargo, si seguimos usando m\'odulos de Hodge mixtos como cuantizaci\'on, es de esperar que, con teor\'ias de campos apropiadas, seamos capaces de obtener resultados similares para los espacios de m\'oduli de fibrados de Higgs y para el de conexiones planas. Este estudio ser\'ia muy interesante porque, en ese caso, podr\'iamos capturar la estructura hyperk\"ahler de estos espacio de m\'oduli.

Otra v\'ia de exploraci\'on interesante ser\'ia estudiar las restricciones de monoidalidad de que presenta la TQFT descrita aqu\'i. Como mencionamos anteriormente, la TQFT construida no es un functor monoidal estricto, sino \'unicamente monoidal laxo. Sin embargo, en los casos conocidos, siempre los c\'alculos han podido ser realizados en un subm\'odulo finitamente generado. Moralmente, eso significa que la TQFT quiere ser monoidal, pero existen algunas obstrucciones que se lo inpiden. Por este motivo, esperamos que las construcciones previas puedan ser modificadas para dar lugar a una TQFT monoidal estricta, para lo cual quiz\'a sea necesario cambiar la cuantizaci\'on o variar las estructuras monoidales.

Finalmente, otro contexto en el que las variedades de caracteres son centrales es el programa de Langlands geom\'etrico (v\'ease \cite{Beilinson-Drinfeld}). En este escenario, la fibraci\'on de Hitchin para el espacio de m\'oduli de Dolbeault satisface las condiciones de Strominger-Yau-Zaslow de simetr\'ia especular (mirror symmetry) para variedades Calabi-Yau (v\'ease \cite{Strominger-Yau-Zaslow}). A partir de esta observaci\'on, surgen varias preguntas que relacionan los $E$-polinomios de variedades de caracteres para grupos duales Langlands, tal y como se conjetur\'o en \cite{Hausel:2005} y \cite{Hausel-Rodriguez-Villegas:2008}. La validez de estas conjeturas ha sido estudiada en algunos casos, como en \cite{LMN} y \cite{Martinez:2017}. A pesar de ello, el caso general permanece sin resolver. Por ello, esperamos que las ideas introducidas en esta tesis doctoral puedan ser \'utiles para arrojar algo de luz sobre estas cuestiones.

%% file: Chapters/Tqft.tex
% Chapter 2

\chapter{Topological Quantum Field Theories} % Main chapter title

\label{chap:tqft} % For referencing the chapter elsewhere, use \ref{Chapter1} 

\lhead{Chapter 1. \emph{Topological Quantum Field Theories}} % This is for the header on each page - perhaps a shortened title

%----------------------------------------------------------------------------------------

\section{Review of category theory and TQFTs}
\label{sec:categories-tqft}

In order to fix notation, let us review some fundamental concepts of category theory that is going to appear along this thesis. The basic reference for all the categorical machinery that we will use is \cite{MacLane}.

Recall that a category $\cC$ is comprised by a class $\Obj{\cC}$ (in the sense of Von Neumann-Bernays-Gödel set theory, see \cite{Bernays}), called the \emph{class of objects} of $\cC$, and, for any $a, b \in \Obj{\cC}$, a class $\Hom_{\cC}(a, b)$, called the \emph{class of morphisms} between $a$ and $b$.
For any $a \in \Obj{\cC}$, the class $\Hom_{\cC}(a, a)$ has a distinguished element $1_a$, called the \emph{identity morphism}, and, for any $a, b, c \in \Obj{\cC}$, there is a map $\circ: \Hom_{\cC}(b, c) \times \Hom_{\cC}(a, b) \to \Hom_{\cC}(a,c)$ called \emph{composition}. This composition is associative in the sense that $h \circ (g \circ f) = (h \circ g) \circ f$, for any $f \in \Hom_{\cC}(a,b)$, $g \in \Hom_{\cC}(b,c)$ and $h \in \Hom_{\cC}(c,d)$; and $1_a \in \Hom_{\cC}(a,a)$ is an identity in the sense that $f \circ 1_a = f$ and $1_a \circ g = g$ for any $f \in \Hom_{\cC}(a,b)$ and $g \in \Hom_{\cC}(c, a)$. A morphism $f \in \Hom_{\cC}(a,b)$ is called an \emph{isomorphism} if it has a right and left inverse, i.e.\ if there exists $g \in \Hom_{\cC}(b,a)$ such that $g \circ f = 1_a$ and $f \circ g = 1_b$. In order to lighten notation, we will shorten $a \in \Obj{\cC}$ by $a \in \cC$ and a morphism $f \in \Hom_{\cC}(a,b)$ will be denoted $f: a \to b$. When the underlying category is clear from the context, $\Hom_{\cC}(a,b)$ will be denoted just by $\Hom(a,b)$.

The role of morphisms between categories will be played by functors. A (covariant) \emph{functor} $F: \cC \to \cD$ between categories $\cC$ and $\cD$ is an assignment rule that, for any $a \in \cC$, gives an object $F(a) \in \cD$ and, to any morphism $f: a \to b$ in $\cC$, assigns a morphism $F(f): F(a) \to F(b)$ in $\cD$. Such an assignment must preserve composition in the sense that $F(g \circ f) = F(g) \circ F(f)$ for any $f: a \to b$ and $g: b \to c$ (in Section \ref{sec:lax-monoidality} we will weaken this notion to obtain lax functors).

Given two categories $\cC$ and $\cD$, it is said that $\cD$ is a \emph{subcategory} of $\cC$ if $\Obj{\cD} \subseteq \Obj{\cC}$ and $\Hom_{\cD}(a,b) \subseteq \Hom_{\cC}(a,b)$ for any $a, b \in \cD$. In that case, we have a \emph{forgetful functor}, $\cD \to \cC$ that 'forgets' the extra data given by $\cD$. In particular, the forgetful functor from a category $\cC$ to itself as a subcategory will be denoted by $1_{\cC}: \cC \to \cC$.
Given $\cD$ a subcategory of $\cC$, if, for any $a,b \in \cD$, we have $\Hom_{\cD}(a,b) = \Hom_{\cC}(a,b)$ it is customary to call $\cD$ a \emph{full subcategory} of $\cC$. On the other hand, if $\Obj{\cD} = \Obj{\cC}$, $\cD$ is usually referred to as a \emph{lluf} or a \emph{wide subcategory}.
	
\begin{rmk}
	Given a category $\cC$, we can construct its \emph{opposite category}, $\cC^{op}$ having the same objects as $\cC$ and the same morphisms but with the arrows reversed (i.e.\ $\Hom_{\cC^{op}}(a,b) = \Hom_{\cC}(b,a)$). A functor $F: \cC^{op} \to \cD$ is also called a \emph{contravariant functor}, $F: \cC \to \cD$. If we have a functor $F: \cC \to \cD$, we can construct the opposite functor $F^{op}: \cC^{op} \to \cD^{op}$.
\end{rmk}

Analogously, the role of morphisms between functors will be played by natural transformations. Given functors $F, G: \cC \to \cD$, a \emph{natural transformation} $\tau: F \Rightarrow G$ is an assignment, to every $a \in \cC$, of a morphism $\tau_a: F(a) \to G(a)$ in $\cD$, called the components. Such an assignment must satisfy that, for any morphism $f: a \to b$ in $\cC$, the following diagram commutes
\[
\begin{displaystyle}
   \xymatrix
   {
  	 F(a) \ar[r]^{\tau_a} \ar[d]_{F(f)}& G(a) \ar[d]^{G(f)}\\
  	 F(b) \ar[r]_{\tau_b} & G(b) 
   }
\end{displaystyle}
\]
A natural transformation is called a \emph{natural isomorphism} if all the morphisms $\tau_a$, for $a \in \cC$, are isomorphisms.

\begin{ex}
\begin{itemize}
	\item Some examples of categories are the empty category $\emptyset$, with no objects and no morphisms, the category $\Sets$, whose objects are sets and whose morphisms are maps between sets, and the category $\Cat$, for which the objects are small categories (i.e. categories whose class of objects is actually a set) and whose morphisms are functors between them.

	\item In a more geometric context, we find the category of topological spaces and continuous maps between them, $\Top$, or the category of compact differentiable manifolds (maybe with boundary) with differentiable maps, $\Diffc$. Another useful category will be $\Embc$, which is the wide subcategory of $\Diffc$ whose morphisms are tame differentiable embeddings between manifolds of the same dimension (see Section \ref{sec:field-theory} for a more precise description).
	\item From algebraic geometry, we also have the category of varieties over an (algebraically closed) field $k$ with regular morphisms between them, $\Var{k}$. Here, and all along this document, by an algebraic variety over $k$ we mean a reduced separated scheme of finite type over $k$. In particular, a variety may not be irreducible. By a projective variety we mean a variety admitting a closed embedding into the projective space $\PP^N$ for $N$ large enough. A quasi-projective variety is an open subset of a projective variety. More generally, we can also consider the category of schemes, $\Sch$.

	\item In an algebraic setting, an archetypal category is, fixed a field $k$, the category of $k$-vector spaces with linear maps, $\Vect{k}$ or, more generally, the category of modules with module homomorphisms over a ground ring $R$, $\Mod{R}$. Another useful categories would be the category of groups and group homomorphisms, $\Grp$, the category of abelian groups, $\Ab$, and the category of commutative rings with unit, $\Rng$.
	
	\item Given two categories $\cC$ and $\cD$, we can build its product category $\cC \times \cD$. It is the category whose objects class is the cartesian product $\Obj{\cC} \times \Obj{\cD}$ and whose morphisms are $\Hom_{\cC \times \cD}((c, d), (c', d')) = \Hom_{\cC}(c, c') \times \Hom_{\cD}(d, d')$ for $c, c' \in \cC$ and $d, d' \in \cD$.
\end{itemize}
\end{ex}

Recall that, given a category $\cJ$ (called the \emph{diagram}) and a functor $F: \cJ \to \cC$, a \emph{cone} of $F$ is a pair $(c, \psi)$, where $c \in \cC$ and $\psi$ is a collection of morphisms $\psi_{j}: c \to F(j)$ for $j \in \cJ$ such that $F(f) \circ \psi_j = \psi_{j'}$ for any morphism $f: j \to j'$ in $\cJ$. A \emph{limit} is a `universal cone' $(\ell, \phi)$, in the sense that, for any other cone $(c, \psi)$ of $F$, there exists an unique morphism $\varphi: c \to \ell$ such that $\phi_j \circ \varphi = \psi_j$ for all $j \in \cJ$.
Analogously, given a functor $F: \cC \to \cJ$, a \emph{cocone} (resp.\ \emph{colimit}) is a cone (resp.\ limit) for $F^{op}: \cC^{op} \to \cJ^{op}$.

\begin{rmk}
\begin{itemize}
	\item By the universal property defining them, if a limit or colimit exists, it is unique up to isomorphism.
	\item Taking $\cJ = \emptyset$, a limit of the unique functor $\emptyset \to \cC$ is called a \emph{final object} of $\cC$ and a colimit is called an \emph{initial object} of $\cC$. If there exists an object of $\cC$ which is both initial and final, it is called the \emph{zero object} and it is denoted by $0 \in \cC$. In particular, for any $a \in \cC$, there are unique morphisms $a \to 0$ and $0 \to a$ which are both denoted by $0$.
	\item Let $\cJ$ be a subcategory of $\cC$ with only identity morphisms. A limit of the forgetful map $\cJ \to \cC$ is called a \emph{product}, and it is denoted by $\prod_{a \in \cJ} a$, and a colimit is called a \emph{coproduct} and it is denoted $\coprod_{a \in \cJ} a$. If $\cC$ has products (resp.\ coproducts) for any finite full subcategories (including the empty one), then $\cC$ is said to have \emph{finite products} (resp.\ \emph{finite coproducts}).
	\item Let $\cJ$ be the subcategory with three objects $a_1, a_2, b \in \cC$ and two morphisms $f_1: b \to a_1$, $f_2: b \to a_2$ (and the identities between them). A \emph{pushout} (resp.\ \emph{pullback}) is a limit (resp.\ colimit) of the forgetful functor $F: \cJ \to \cC$. In other words, a pushout is an object, usually denoted $a_1 \sqcup_b a_2 \in \cC$ which is a limit of the diagram
\[
\begin{displaystyle}
   \xymatrix
   {
	b \ar[r]^{f_1}\ar[d]_{f_2} & a_1 \ar@{--{>}}[d]^{g_1} \\
	a_2 \ar@{--{>}}[r]_{g_2\;\;\;} & a_1 \sqcup_b a_2
   }
\end{displaystyle}
\]
where $g_1: a_1 \to a_1 \sqcup_b a_2$ and $g_2: a_2 \to a_1 \sqcup_b a_2$ are the morphisms of the cone defining the limit. Analogously, it is customary to denote the pullback by $a_1 \times_b a_2$.
	\item If $\cC$ has an final object $\star \in \cC$, then the product of $a_1$ and $a_2$ is the same as the pullback of $a_1$ and $a_2$ over $\star$. Analogously, a coproduct is a pushout over the initial object. For short, we will denote them $a_1 \star a_2$ and $a_1 \sqcup a_2$.
\end{itemize}
\end{rmk}

\subsection{Monoidal categories}

A main concept along this chapter are the so-called monoidal categories. These categories can be thought as a category with an extra binary operation in the spirit of tensor product on modules. For further information about monoidal categories, see \cite{MacLane1963} or \cite{MacLane}.

\begin{defn}\label{defn:monoidal-cat}
A \emph{monoidal category} is a category $\cC$, with a functor
$$
	\otimes: \cC \times \cC \to \cC
$$
and a distinguished object $I \in \cC$ such that:
\begin{itemize}
	\item (Identity) There are natural isomorphisms
	$$
		\lambda: I \otimes - \Rightarrow 1_{\cC}, \hspace{1cm} \rho: - \otimes I \Rightarrow 1_{\cC}.
	$$
	\item (Associtivity) There is a natural isomorphism $\alpha: (- \otimes -) \otimes - \Rightarrow - \otimes (- \otimes -)$.
	\item (Triangle identity) For any $a, b \in \cC$, the following diagram commutes
	\[
\begin{displaystyle}
   \xymatrix
   {
  	 (a \otimes I) \otimes b \ar[rr]^{\alpha_{a, I, b}} \ar[rd]_{\rho_a \otimes 1_{b}} & & a \otimes (I \otimes b) \ar[ld]^{1_a \otimes \lambda_b}\\
  	 & a \otimes b &
   }
\end{displaystyle}
\]
	\item (Pentagon identity) For any $a, b, c, d \in \cC$, the following diagram commutes
	\[
\begin{displaystyle}
   \xymatrix
   {
  	  & (a \otimes b) \otimes (c \otimes d) \ar[rd]^{\alpha_{a, b, c \otimes d}} & \\
  	 ((a \otimes b) \otimes c) \otimes d \ar[d]_{\alpha_{a,b,c} \otimes 1_d} \ar[ru]^{\alpha_{a \otimes b, c, d}} &  & a \otimes (b \otimes (c \otimes d)) \\
  	  (a \otimes (b \otimes c)) \otimes d \ar[rr]_{\alpha_{a,b \otimes c, d}}& & a \otimes ((b \otimes c) \otimes d) \ar[u]_{1_a \otimes \alpha_{b,c,d}}
   }
\end{displaystyle}
\]
\end{itemize}
A monoidal category is said to be \emph{braided} if there exists a natural isomorphism $B_{a,b}: a \otimes b \to b \otimes a$ for $a,b \in \cC$, called the braiding. If $B_{b,a} \circ B_{a,b} = 1_{a,b}$, the category is said to be \emph{symmetric}.
\end{defn}

\begin{rmk}
\begin{itemize}
	\item Despite that, a priori, the coherence conditions for monoidal categories could comprise much more relations, Kelly coherence theorem (see \cite{Kelly1964} shows that the triangle and the pentagon identity are enough generate them. Hence, on a monoidal category, all the possible placements of brackets in a steam of objects of $\cC$ give the same object.
	\item To be precise, some extra coherence conditions are needed for braided categories, in the spirit of the triangle of pentagon identities. For simplicity, we will not quote them here, but they can be checked in \cite{Joyal-Street}.
\end{itemize}
\end{rmk}

\begin{ex}\label{ex:cartesian-monoidal}
The model example of monoidal category is $\Vect{k}$ with tensor product $\otimes$ and unit $k \in \Vect{k}$ or, more generally, $\Mod{R}$ with tensor product $\otimes_R$ and unit $R$. Another example is $\Sets$ with the cartesian product and the singleton set $\star \in \Sets$ as the unit, or $\Cat$ with the product of categories and unit the category with a single object and only the identity as morphism.

In general, if a category has a final object and pullbacks, it can be endowed with a monoidal structure where the unit is the final object and the tensor product is the pullback of two objects over the final object. In this fashion, such a category is called a \emph{cartesian monoidal category}. Analogously, if the category has initial object and pushout, it becomes a monoidal category usually referred to as the \emph{cocartesian monoidal category}.
\end{ex}

\begin{figure}[h]
	\begin{center}
	\includegraphics[scale=0.2]{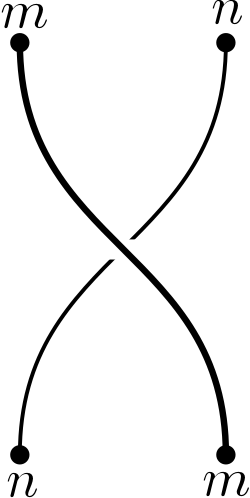}
	\caption{}
	\label{img:crossing-stands}
	\end{center}
\vspace{-1cm}
\end{figure}

The reason of the name `braiding' comes from the category of braids. Roughly speaking, this category, $\cB$, has $\Obj{\cB} = \NN$ and, given a natural number $n$, $\Hom_{\cB}(n,n)$ is the braid group in $n$ elements (see \cite{Wilhelm} for a survey about braid groups). It has no morphisms between different natural numbers. The category $\cB$ is monoidal with the disjoint union of braids (i.e.\ given natural numbers $n,m$, $n \otimes m = n + m$ and the tensor product of two braids is just to juxtapose them). Furthermore, it is a braided category with the braiding $B_{n,m}: n \otimes m \to m \otimes n$ given by 'crossing strands', as shown in Figure \ref{img:crossing-stands}.

Obviously, $B_{n,m}$ is an isomorphism (all the morphism are so) but $B_{m,n} \circ B_{n,m} \neq 1_{n + m}$ so $\cB$ is not a symmetric monoidal category.

The most important example of monoidal category for our purposes is the so-called \emph{category of $n$-dimensional oriented bordisms}, $\Bordo{n}$. Given $n \geq 1$, this category is given by the following data:
\begin{itemize}
	\item Objects: The objects of $\Bordo{n}$ are $(n-1)$-dimensional closed (i.e.\ compact without boundary) oriented differentiable manifolds, maybe empty.
	\item Morphisms: Given objects $M_1, M_2 \in \Bordo{n}$, a morphism $M_1 \to M_2$ is a class of oriented bordisms between $M_1$ and $M_2$, i.e.\ of compact oriented differentiable manifolds $W$ with boundary $\partial W = M_1 \sqcup M_2$ in such a way that the inclusion $i_1: M_1 \hookrightarrow \partial W$ (resp.\ $i_2: M_2 \hookrightarrow \partial W$) preserves (resp.\ inverts) orientation over its image. It is customary to denote this situation by $\partial W = M_1 \sqcup \overline{M_2}$. Two bordisms are in the same class if there exists an orientation-preserving diffeomorphism of bordisms (i.e.\ fixing the boundaries) between them. For a more precise definition of bordisms, see \cite{Milnor:1965}.
	\item Composition: Given $W: M_1 \to M_2$ and $W': M_2 \to M_3$, we define $W' \circ W$ as the morphism $W \cup_{M_2} W': M_1 \to M_3$ where $W \cup_{M_2} W'$ is the usual gluing of bordisms along $M_2$.
	\item Monoidal structure: It is given by the bifunctor $\sqcup: \Bordo{n} \times \Bordo{n} \to \Bordo{n}$ that takes disjoint union of both objects and bordisms. The unit of the monoidal structure is the empty set manifold $\emptyset \in \Bordo{n}$.
\end{itemize}
For further information about this category, see \cite{Kock:2004}.

\begin{defn}
Let $(\cC, \otimes_\cC, I_\cC)$ and $(\cD, \otimes_\cD, I_\cD)$ be monoidal categories. A functor $F: \cC \to \cD$ is called \emph{monoidal} if:
\begin{itemize}
	\item There is an isomorphism $\alpha: I_{\cD} \to F(I_\cC)$.
	\item There is a natural isomorphism $\Delta: F(-) \otimes_\cD F(-) \Rightarrow F(- \otimes_\cC -)$.
\end{itemize}
If $\cC$ and $\cD$ are symmetric monoidal categories, then $F$ is also called \emph{symmetric} if $\Delta_{b, a} \circ B_{F(a), F(b)} = F(B_{a,b}) \circ \Delta_{a,b}$.
\end{defn}

\subsection{Duality properties and low dimensional classification of TQFTs}
\label{sec:low-dim-tqft}

With these definitions at hand, we can define the central concept of this chapter (and of this thesis), the so-called Topological Quantum Field Theories. Introducted by Segal in \cite{Segal:1988} from the physical intuition behind Feynman's path integral, they were formalized by Atiyah in \cite{Atiyah:1988} and Witten \cite{Witten:1988}. Here, we will focus on a description in terms of categories (for further information, see \cite{Kock:2004}).

\begin{defn}
Let $R$ be a commutative ring with unit. An (oriented) \emph{Topological Quantum Field Theory}, shorten as TQFT, is a symmetric monoidal functor $Z: \Bordo{n} \to \Mod{R}$.
\end{defn}

Given two TQFTs, $Z, Z': \Bordo{n} \to \Mod{R}$, we will say that they are isomorphic if there exists a natural isomorphism $\tau: Z \Rightarrow Z'$ such that the following two diagrams commute
\[
\begin{displaystyle}
   \xymatrix
   {
	Z(M_1) \otimes_R Z(M_2) \ar[rr]^{\tau_{M_1} \otimes_R \tau_{M_2}} \ar[d]_{\Delta_{M_1, M_2}} && Z'(M_1) \otimes_R Z'(M_2) \ar[d]^{\Delta_{M_1, M_2}'} \\
	Z(M_1 \sqcup M_2)\ar[rr]_{\tau_{M_1 \sqcup M_2}} && Z'(M_1 \sqcup M_2) \\
   }
   \hspace{1cm}
   \xymatrix
   {
	& R \ar[ld]_\alpha \ar[rd]^\alpha' & \\
	Z(\emptyset) \ar[rr]_{\tau_\emptyset} & & Z'(\emptyset)
   }
\end{displaystyle}
\]
In general, such a natural transformation between monoidal functors is called a \emph{natural monoidal transformation}.

Fix $M \in \Bordo{n}$ and denote by $\overline{M}$ the same manifold but with the opposite orientation. Among the cobordisms of $M$, the simplest one is the oriented bordism $M \times [0,1]$ whose boundary is two copies of $M$, where $[0,1]$ has the natural orientation. Depending on the orientation of the boundaries, we obtain three different morphisms of $\Bordo{n}$. We have the unit morphism $1_M: M \to M$ if both boundaries have the same orientation, $\mu_M: M \sqcup \overline{M} \to \emptyset$ if both boundaries have the induced orientation and $\epsilon_M: \emptyset \to \overline{M} \sqcup M$ if both boundaries have the opposite orientation to the induced one. Another important bordism is the swap tube $S_M: M \sqcup \overline{M} \to \overline{M} \sqcup M$ given by the twist of two identity tubes.
\begin{figure}[h]
	\begin{center}
	\includegraphics[scale=0.5]{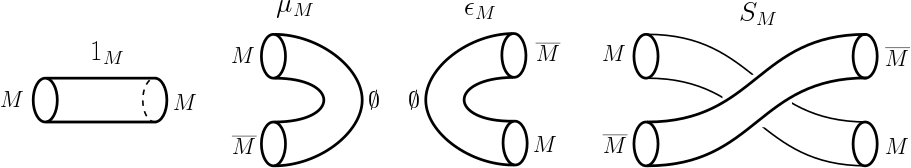}
	\end{center}
\vspace{-0.3cm}
\end{figure}

\begin{lem}[Zorro's property]\label{lem:zorro}
Let $Z: \Bordo{n} \to \Mod{R}$ be a Topological Quantum Field Theory and $M \in \Bordo{n}$.
\begin{itemize}
	\item[i)] It holds $\left(Z(\mu_M) \otimes 1_{Z(M)}\right)  \circ \left(1_{Z(M)} \otimes Z(\epsilon_M)\right) = 1_{Z(M)}$.
	\item[ii)] The $R$-module $Z(M)$ is finitely generated.
	\item[iii)] The morphism $Z(\mu_M): Z(M) \otimes Z(\overline{M}) \to R$ is non-degenerate.
	\item[iv)] The morphism $Z(S_M): Z(M) \otimes Z(\overline{M}) \to Z(\overline{M}) \otimes Z(M)$ is an isomorphism.
\end{itemize}
\begin{proof}
In the category $\Bordo{n}$, we have that $(\mu_M \sqcup 1_M) \circ (1_M \sqcup \epsilon_M) = 1_M$ as depicted in Figure \ref{img:zorro-property}. For that reason, $i)$ follows inmediately after applying $Z$.

\begin{figure}[h]
	\begin{center}
	\includegraphics[scale=0.3]{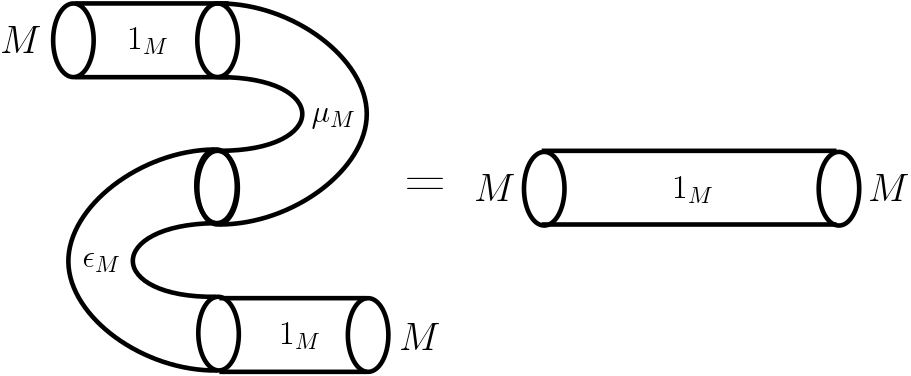}
	\caption{}
	\label{img:zorro-property}
	\end{center}
\vspace{-0.4cm}
\end{figure}

For $ii)$ observe that, since $Z(\epsilon_M): R \to Z(\overline{M}) \otimes Z(M)$, we can write $Z(\epsilon_M)(1) = \displaystyle{\sum_{i=1}^N}\, y_i \otimes x_i$ for some $x_i \in Z(M)$ and $y_i \in Z(\overline{M})$. Given $x \in Z(M)$, from $i)$ we have that
$$
	x = 1_{Z(M)}(x) = \left(Z(\mu_M) \otimes 1_{Z(M)}\right)\left(\sum_{i=1}^N x \otimes y_i \otimes x_i\right) = \sum_{i=1}^N \mu(x \otimes y_i) x_i.
$$
Hence, the elements $x_1, \ldots, x_N$ generate $Z(M)$, as claimed. The same calculation also proves that $\mu$ is non-degenerate in the first variable since, if $x \neq 0$, some of the products $\mu(x \otimes y_i) \neq 0$.

For $iv)$, observe that $S_{\overline{M}} \circ S_M = 1_{M \sqcup \overline{M}}$ and $S_M \circ S_{\overline{M}} = 1_{\overline{M} \sqcup M}$. Hence, $S_M$ is an isomorphism in $\Bordo{n}$ so $Z(S_M)$ is so in $\Mod{R}$. In order to finish, observe that $\mu_M \circ S_M = \mu_{\overline{M}}$ so the non-degeneracy in the second component of $\mu_M$ follows from the one in the first component of $\mu_{\overline{M}}$.
\end{proof}
\end{lem}

\begin{rmk}\label{rmk:finite-dim-TQFT}
\begin{itemize}
	\item The cornerstone in the previous proof of the finite generation of $Z(M)$ was the elbow $\epsilon_M: \emptyset \to \overline{M} \sqcup M$. Actually, this bordism itself determines a set of generators of $Z(M)$ since, as shown in the proof above, if $Z(\epsilon_M)(1) = \sum y_i \otimes x_i$ then the elements $x_i$ generate $Z(M)$.
	\item Suppose that the ground ring $R$ is a field $k$. In that case, by the previous lemma, $Z(M)$ is a finite dimensional vector space and $Z(\mu_M)$ is a perfect pairing, so we have a natural identification $Z(\overline{M}) = Z(M)^*$. Under this identification, $Z(\mu_M) = \textrm{ev}: Z(M) \otimes Z(M)^* \to k$ is the evaluation map and $Z(\epsilon_M) = \textrm{coev}: k \to Z(M)^* \otimes Z(M)$ is the coevaluation map, that is, $Z(\epsilon_M)(1) = \sum v_i^* \otimes v_i$ where the vectors $v_i$ form a basis of $Z(M)$ and $v_i^*$ is the dual basis of $Z(M)^*$.
\end{itemize}
\end{rmk}

\begin{ex}\label{ex:low-dim-TQFT}
Along these examples, suppose that $R=k$ is a field.
\begin{itemize}
	\item Let $n=1$. In that case, the objects of $\Bordo{1}$ are disjoint union of three basic objects: $\emptyset$ and the point with its two possible orientations, $+$ and $-$. Moreover, morphisms are just paths, so they are disjoint union of the four previous ones and the circle $S^1: \emptyset \to \emptyset$ (and their opposite orientations). Hence, by the previous Proposition, a $1$-TQFT, $Z: \Bordo{1} \to \Vect{k}$, is completely determined by the finite dimensional vector space $Z(+)$. Even more, since $S^1 = \mu_+ \circ \epsilon_+$, we have that
	$$
	Z(S^1)(1) = Z(\mu_+) \circ Z(\epsilon_+)(1) = Z(\mu_+)\left(\sum_{i=1}^d v_i^* \otimes v_i\right) = d, 
	$$
where $v_1, \ldots, v_d$ is a basis of $Z(+)$. Hence, $Z(S^1)$ is just multiplication by the dimension of $Z(+)$. It is a general fact that $Z(W)$ for a closed $n$-dimensional manifold $W: \emptyset \to \emptyset$ usually is multiplication by some important datum.
	\item For $n=2$, a similar classification result holds, since a $2$-TQFT is equivalent to a Frobenius algebra. Recall that a Frobenius algebra $A$ is a commutative $k$-algebra of finite type with a non-degenerate bilinear form $B: A \otimes A \to k$ such that $B(ab,c)=B(a,bc)$ for all $a,b,c \in A$. We are not goint to present here a full proof of this result, that can be found in \cite{Kock:2004}. However, at least we will show how the Frobenius algebra structure appears. The category $\Bordo{2}$ has two basis objects, $\emptyset$ and $S^1$, beeing the other objects disjoint unions of them or their opposite orientation. Let $Z: \Bordo{n} \to \Vect{k}$ be a TQFT and set $A = Z(S^1)$. We have two special types of bordisms, the pair of pants $\Delta: S^1 \sqcup \overline{S^1} \to S^1$, as depicted in Figure \ref{img:pair-of-pants}, and $\mu_{S^1}: S^1 \sqcup \overline{S^1} \to \emptyset$.
\begin{figure}[h]
	\begin{center}
	\includegraphics[scale=0.18]{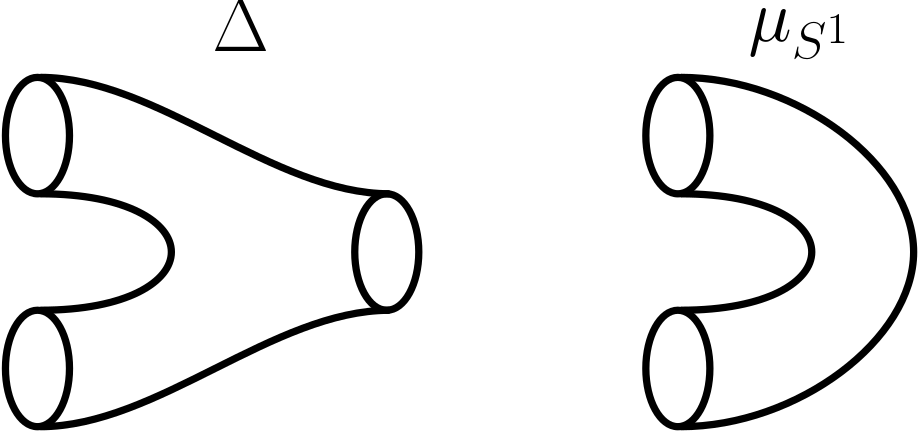}
	\caption{}
	\label{img:pair-of-pants}
	\end{center}
\vspace{-0.8cm}
\end{figure}

After identifying $Z(S^1)$ and $Z(\overline{S^1})$, the image of $\Delta$ is an homomorphism $Z(\Delta): A \otimes A \to A$, which gives the algebra structure. The bilinear form is nothing but $B = Z(\mu_{S^1}): A \otimes A \to k$. Using the relations between these bordisms in $\Bordo{2}$, it can be proven that these data give to $A$ a Frobenius algebra structure. Conversely, it can be given a full description of the morphism of $\Bordo{2}$ in terms of generators and relations (which are nothing but compact orientable surfaces minus discs, which are completely classified). From them, it can be proven that, just giving the assignment of $\Delta$ and $\mu_{S^1}$ under the TQFT, the rest of bordisms are completely determined.
\end{itemize}
\end{ex}

As a final comment, observe that we may remove the orientability condition on $\Bordo{n}$. In that case, we obtain the \emph{category of $n$-dimensional bordisms}, $\Bord{n}$ that has, as objects, $n$-dimensional closed manifolds and, as morphisms, bordisms between them up to boundary preserving isomorphism.

\begin{defn}
Let $R$ be a commutative ring with unit. A \emph{Topological Quantum Field Theory} is a symmetric monoidal functor $Z: \Bord{n} \to \Mod{R}$.
\end{defn}

Since there is a forgetful functor $\Bordo{n} \to \Bord{n}$ just by ignoring the orientation, any non-oriented TQFT induces an oriented one. Moreover, observe that the proof of the Zorro's property \ref{lem:zorro} does not requiere indeed the orientation on $M$ so it can be adjusted to prove the same result for non-oriented TQFT. For this reason, since all the morphisms in $\Bord{1}$ are orientable, it can be shown that an analogous classification for $1$-TQFTs holds on the non-oriented case. However, the classification of $2$-TQFT is not preserved in the non-oriented setting since there are non-oriented bordisms that must be analyzed.

\begin{rmk}
Historically, Topological Quantum Field Theories were introduced as a phenomenological description of path integrals. Roughly speaking, the idea is the following. Suppose that we have a field (in the physical sense) on a boundary of a bordism. We can propagate it to the opposite boundary by means of the integration of the action of all the fields on the whole bordism that restrict to the initial field in the incoming boundary. If such an action does not depend on the metric of the bordism, the obtained propagation method is purely topological and satisfies the same conditions as a TQFT. For a more detailed introduction, see \cite{Carqueville-Runkel:2017}.

However, this origin is extremely informal and, after the papers of Witten (\cite{Witten:1988}) and, specially, Atiyah (\cite{Atiyah:1988}), this physical background was pulled apart and the mathematical framework presented here was introduced. Nonetheless, most of the mathematical constructions of TQFT follow similar lines to the original construction for path integrals. Actually, in the following section we give a construction of TQFT that deeply reminds the historical one.
\end{rmk}

\subsection{Lax monoidality and $2$-categories}
\label{sec:lax-monoidality}

As we showed in the previous section, monoidality of a TQFT implies very strong algebraic restrictions on the TQFT. As shown in Zorro's lemma \ref{lem:zorro}, if $Z$ is a strict monoidal TQFT, then $Z(M)$ is a finitely generated module. As explained in Remark \ref{rmk:finite-dim-TQFT}, the culprit of this phenomenon is the elbow bordism $\epsilon_M: \emptyset \to M \sqcup M$. For this reason, if we want to bypass Zorro's lemma and allow infinitely generated modules, we need to focus on $\epsilon_M$. There are (at least) two obvious approaches to this problem.
\begin{itemize}
	\item Restrict the morphisms. Consider an appropiate subcategory $\cB$ of $\Bord{n}$ that no longer contains $\epsilon_M$. In that case, Zorro's lemma cannot be accomplished in $\cB$ so a monoidal functor $Z: \cB \to \Mod{R}$ is not forced to land into finitely generated modules. A classical choice for $\cB$ is to restrict to bordisms with non-empty ingoing boundary. In the extended framework, this gives rise to the so-called non-compact TQFTs, as studied in \cite{Lurie-eTQFT:2009}, Section 4.2. In this thesis, we also propose to take $\cB=\Tub{n}$ as the category of tubes i.e.\ bordisms with connected (or empty) ingoing and outgoing boundaries. A symmetric monoidal functor $\lZ: \Tub{n} \to \Mod{R}$ is called an \emph{almost-TQFT} and they will be studied in Section \ref{sec:almost-TQFT-tubes}. This choice of $\cB$ has the advantage that closed manifolds remain as morphisms, so almost-TQFTs are good enough for computational purposes.
%In the non-compact case, there is an analogous classification of extended TQFTs by means of Calabi-Yau objects (Theorem 4.2.11)
	\item Relax monoidality. Suppose that the functor $Z: \Bord{n} \to \Mod{R}$ is no longer monoidal. In that case, it may happen $Z(M \sqcup M) \neq Z(M) \otimes Z(M)$ so, in the proof of Zorro's lemma, we cannot write $Z(\epsilon_M)(1)=\sum y_i \otimes x_i$ with $y_i,x_i \in Z(M)$. For this reason, the proof fails and $Z(M)$ is not forced to by finitely generated (actually, the lax monoidal TQFT constructed in Chapter \ref{chap:representation} will be an example of this phenomenon). This approach has the advantage that we impose no restrictions on the allowed bordisms, so the source category remains unchanged. 
\end{itemize} 

Along this section, we will focus on the second method and we will consider the so-called lax monoidal TQFTs. In some sense, the TQFTs considered in this thesis arise naturaly in this context.

\begin{defn}\label{defn:lax-monoidal}
Let $(\cC, \otimes_\cC, I_\cC)$ and $(\cD, \otimes_\cD, I_\cD)$ be monoidal categories. A functor $F: \cC \to \cD$ is said to be \emph{lax monoidal} if there exists:
	\begin{itemize}
		\item A morphism $\alpha: I_{\cD} \to F(I_{\cC})$.
		\item A natural transformation $\Delta: F(-) \otimes_{\cD} F(-) \Rightarrow F(- \otimes_\cC -)$.		
	\end{itemize}
If $\alpha$ and $\Delta$ are isomorphisms, $F$ is said to be \emph{pseudo-monoidal} (or just monoidal) and, if they are identity morphisms, $F$ is called \emph{strict monoidal}.
\end{defn}

\begin{rmk}
\begin{itemize}
	\item To be precise, the transformations $\alpha$ and $\Delta$ have to satisfy a set of coherence conditions similar to the ones of Definition \ref{defn:monoidal-cat}. We will not state them here but, for a precise definition, see Definition 1.2.10 of \cite{Leinster}.
	\item If a functor $F: \cC \to \cD$ is a lax monoidal functor when seen in the opposite categories $F: \cC^{op} \to \cD^{op}$ it is customary to call $F$ oplax monoidal.
\end{itemize}
\end{rmk}

\begin{defn}
Let $R$ be a commutative ring with unit. A \emph{lax monoidal Topological Quantum Field Theory} is a symmetric lax monoidal functor $Z: \Bord{n} \to \Mod{R}$.
\end{defn}

%\begin{rmk}
%Since a lax monoidal TQFT is not strict monoidal, the proof of Zorro's lemma \ref{lem:zorro} does not hold since $Z(M \sqcup N) \neq Z(M) \otimes_R Z(N)$. For this reason, $Z(M)$ may be infinitely generated. Actually, the lax monoidal TQFT constructed in Chapter \ref{chap:representation} will be an example of this phenomenon.
%\end{rmk}

Another important concept in category theory is the notion of a $2$-category. The idea is very simple: a $2$-category is a category with an extra layer of arrows between morphisms, called $2$-morphisms.
However, at the moment of defining them, we find several technical difficulties concerning the interplay between $1$-morphisms and $2$-morphisms. For these reasons, there exists several (non-equivalent) formulations of $2$-categories depending on the level of strictness required.

In this section, we will first introduce the notion of a strict $2$-category as a particular case of an enriched category. It is somehow the cleanest way of introducing $2$-categories but, sadly, they are too restictive for our purposes. In order to solve it, we will introduce the weak $2$-categories just by requiring that the coherence conditions are satisfied up to invertible $2$-morphism.

\begin{defn}\label{defn:enriched}
Given a monoidal category $(\cV, \otimes, I)$, a \emph{$\cV$-enriched category} $\cC$ is a category such that:
\begin{itemize}
	\item[$i)$] For any $a,b \in \cC$, the morphisms between $a$ and $b$, $\Hom_{\cC}(a,b)$, is an object of $\cV$.
	\item[$ii)$] The composition is given by morphisms of $\cV$, $\circ: \Hom_{\cC}(b,c) \otimes \Hom_{\cC}(a,b) \to \Hom_{\cC}(a,c)$
	\item[$iii)$] For every $a \in \cC$, there exists a morphism $j_a: I \to \Hom_{\cC}(a,a)$.
	\item[$iv)$] For any $a,b,c,d \in \cC$, the following diagram commutes
	\[
\begin{displaystyle}
   \xymatrix
   {
	(\Hom_{\cC}(c,d)\otimes \Hom_{\cC}(b,c)) \otimes \Hom_{\cC}(a,b) \ar[d]_{\circ \otimes 1}
      \ar[r]^{\alpha} 
     &\Hom_{\cC}(c,d) \otimes (\Hom_{\cC}(b,c) \otimes \Hom_{\cC}(a,b)) \ar[d]^{1 \otimes \circ}
     \\ 
     \Hom_{\cC}(b,d)\otimes \Hom_{\cC}(a,b) \ar[d]_{\circ}&
     \Hom_{\cC}(c,d) \otimes \Hom_{\cC}(a,c) \ar[dl]^{\circ} \\
      \Hom_{\cC}(a,d) &
   }
\end{displaystyle}
\]
	\item[$v)$] For any $a,b \in \cC$, the following diagram commutes
	\[
\begin{displaystyle}
   \xymatrix
   {
    \Hom_{\cC}(b,b)\otimes \Hom_{\cC}(a,b) \ar[r]^{\hspace{1.2cm}\circ}
    &
    \Hom_{\cC}(a,b)
    &
    \Hom_{\cC}(a,b) \otimes \Hom_{\cC}(a,a) \ar[l]_{\circ\hspace{1cm}}
    \\
    I \otimes \Hom_{\cC}(a,b) \ar[u]^{j_b \otimes 1} \ar[ur]_{\lambda}
    &&
    \Hom_{\cC}(a,b) \otimes I \ar[u]_{1 \otimes j_a} \ar[ul]^{\rho}
   }
\end{displaystyle}
\]
\end{itemize}
\end{defn}

\begin{ex}
Most of the algebraic categories are enriched categories. For example, the category $\Mod{R}$ is a $\Mod{R}$-enriched category (in general, categories enriched over $\Mod{R}$ are called $R$-linear categories). A category enriched over $\Sets$ is the same as a locally small category (i.e.\ a category whose morphisms between two objects are genuine sets).
\end{ex}

A category $\cC$ enriched over $\Cat$ is called a \emph{strict $2$-category}. Given $a,b \in \cC$, the objects of $\Hom_{\cC}(a,b)$ are usually referred to as $1$-morphisms. However, as $\Cat$-enriched category, $\Hom_{\cC}(a,b)$ itself has a category structure. Hence, given $f, g \in \Hom_{\cC}(a,b)$, we have also morphism between $f$ and $g$, denoted as $f \Rightarrow g$ and called $2$-morphisms or $2$-cells.

However, for most applications, strict $2$-categories are too restrictive. In order to solve this problem, we can consider the notion of \emph{weak $2$-categories}. These categories satisfy the same conditions as items $i)$-$iii)$ of Definition \ref{defn:enriched} for $\cV = \Cat$ (in particular, there is a notion of $2$-morphisms). However, items $iv)$-$v)$ are satisfied only up to invertible $2$-morphism. For a complete introduction to $2$-categories, see \cite{Benabou}. Except when explicitly said, from now on $2$-category will mean weak $2$-category.

\begin{rmk}
The incarnation of weak $2$-category given here actually corresponds to the one of bicategories, as decribed in \cite{Benabou}. Other formulations of weak $2$-categories can be given in the context of higher category theory (see \cite{Lurie:2009}).
\end{rmk}

\begin{ex}\label{ex:categories}
\begin{itemize}
	\item $\Cat$ is a (strict) $2$-category. The $1$-morphisms are functors between categories and the $2$-morphisms are natural transformations between functors.
	
	\item The category $\Top$ can be endowed with a $2$-category structure by taking as $2$-cell $f \Rightarrow g$ between continuous maps an homotopy between $f$ and $g$ (actually, homotopy classes of homotopies in order to have a well-behaved composition).
	
	\item Let $\cC$ be a monoidal category. We build the category $\textbf{B}\cC$ with a single object $\star$ and $\Hom_{\textbf{B}\cC}(\star, \star)=\cC$. Composition of two $1$-morphisms $a,b \in \cC$ is given by tensor product $a \otimes b \in \cC$ and composition of two $2$-morphisms $f: a \to b$ and $g: b \to c$ is given by the usual composition $g \circ f: a \to c$. Hence, $\textbf{B}\cC$ is a $2$-category which is usually referred to as the delooping of $\cC$. Actually, the construction is invertible so a monoidal category is the same as a $2$-category with a single object.
	
	\item Let $R$ be a ring. We define the $2$-category $\Bim{R}$ of $R$-algebras and bimodules whose objects are commutative $R$-algebras with unit and, given algebras $A$ and $B$, a $1$-morphism $A \to B$ is a $(A,B)$-bimodule. By convention, an $(A,B)$-bimodule is a set $M$ with a left $A$-module and a right $B$-module compatible structures, usually denoted ${_A}M_{B}$. Composition of $M: A \to B$ and $N: B \to C$ is given by ${_A}(M \otimes_B N)_C$.
	
With this definition, the set $\Hom_{\Bim{R}}(A,B)$ is naturally endowed with a category structure, namely, the category of $(A,B)$-bimodules. Hence, a $2$-morphism $M \Rightarrow N$ between $(A,B)$-bimodules is a bimodule homomorphism $f: M \to N$. Therefore, $\Bim{R}$ is a monoidal $2$-category with tensor product over $R$.
	\item Let $R$ be a fixed ring (commutative and with unit). Given two homomorphisms of $R$-modules $f,g: M \to N$, we say that $g$ is an \emph{immediate twist} of $f$ if there exists an $R$-module $D$, homomorphisms $f_1: M \to D$, $f_2: D \to N$ and $\psi: D \to D$ such that $f = f_2 \circ f_1$ and $g = f_2 \circ \psi \circ f_1$.
\[
\begin{displaystyle}
   \xymatrix
   {	M \ar[r]^{f_1} \ar@/_1pc/[rr]_g & D \ar[r]^{f_2} \ar@(ul,ur)^{\psi} & N
   }
\end{displaystyle}   
\]
In general, given $f,g: M \to N$ two $R$-module homomorphisms, we say that $g$ is a \emph{twist} of $f$ if there exists a finite sequence $f=f_0, f_1, \ldots, f_r = g: M \to M$ of homomorphisms such that $f_{i+1}$ is an immediate twist of $f_i$.

In that case, we define the \emph{category of $R$-modules with twists}, $\Modt{R}$ as the category whose objects are $R$-modules, its $1$-morphisms are $R$-modules homomorphisms and, given homomorphisms $f,g: A \to B$, a $2$-morphism $f \Rightarrow g$ is a twist from $f$ to $g$ (i.e.\ a sequence of immediate twists). Composition of $2$-cells is juxtaposition of twists. With this definition, $\Modt{R}$ has a $2$-category structure. Moreover it is a monoidal category with the usual tensor product.
\end{itemize}
\end{ex}

Exactly as for usual categories we have a notion of functor, for $2$-categories there is an analogous version, called a $2$-functor.

\begin{defn}\label{defn:lax-2-funct}
A \emph{lax 2-functor} between $2$-categories, $F: \cC \to \cD$, is an assignment that:
\begin{itemize}
	\item For each object $x \in \cC$, it gives an object $F(x) \in \cD$.
	\item For each pair of objects $x, y$ of $\cC$, we have a functor
	$$
		F_{x,y}: \Hom_\cC(x, y) \to \Hom_{\cD}(F(x), F(y)).
	$$
	\item For each object $x \in \cC$, we have a $2$-morphism $F_{1_x}: 1_{F(x)} \Rightarrow F_{x,x}(1_x)$.
	\item For each triple $x,y,z \in \cC$ and every $f: x \to y$ and $g: y \to z$, we have a $2$-morphism $F_{x,y,z}(g,f): F_{y,z}(g) \circ F_{x,y}(f) \Rightarrow F_{x,z}(g \circ f)$, natural in $f$ and $g$.
\end{itemize}
Also, some technical conditions, namely the coherence conditions, have to be satisfied (see \cite{Benabou} or \cite{MacLane} for a complete definition). If the $2$-morphisms $F_{1_x}$ and $F_{x,y,z}$ are isomorphisms, it is said that $F$ is a \emph{pseudo}-functor or a weak functor (or even simply a $2$-functor) and, if they are the identity $2$-morphism, $F$ is called a \emph{strict} $2$-functor.
\end{defn}

Finally, we also have a notion of natural equivalence between $2$-functors.
\begin{defn}
Let $F, G: \cC \to \cD$ be (maybe lax) $2$-functors between $2$-categories. A lax natural transformation is a collection of assignments:
\begin{itemize}
	\item For any $a \in \cC$, a morphism $\tau_a: F(a) \to G(a)$.
	\item For any morphism $f: a \to b$ of $\cC$, a $2$-morphism $\tau_f: G(f) \circ \tau_a \Rightarrow \tau_b \circ F(f)$.
\end{itemize}
such that $\tau_f: (\tau_{a})^* \circ G_{a, b} \Rightarrow (\tau_{b})_* \circ F_{a, b}$ is a natural transformation (of usual categories). Here $(\tau_{a})^*: \Hom_{\cD}(G(a), G(b)) \to \Hom_{\cD}(F(a), G(b))$ is the functor given by $f \mapsto f \circ \tau_a$ and $(\tau_{b})^*: \Hom_{\cD}(F(a), F(b)) \to \Hom_{\cD}(F(a), G(b))$ is given by $f \mapsto \tau_b \circ f$. 

As in the case of $2$-functors, if the $2$-cells $\tau_f$ are invertible, $\tau$ is called pseudo-natural (or just natural) and, if they are identities, $\tau$ is called strict natural.
\end{defn}

\begin{rmk}
\begin{itemize}
	\item If the $2$-cells in the definitions of lax natural transformation and lax functor go in the other way around, it is customary to call them oplax.
	\item In higher category theory, the functors that appear between $2$-categories are usually pseudo-functors. That is the reason why, by default, $2$-functor stands for pseudo-functor.
\end{itemize}
\end{rmk}

\section{Lax monoidal Topological Quantum Field Theories}
\label{sec:lax-monoidal-tqft}

In this section, we will describe a general recipe for constructing lax monoidal TQFTs from two simpler pieces of data: one of geometric nature (a field theory) and one of algebraic nature (a quantisation). This is a very natural construction that, in fact, has been widely used in the literature in some related form (see, for example, \cite{Freed:1994}, \cite{Freed-Hopkins-Lurie-Teleman:2010}, \cite{Freed-Hopkins-Teleman:2010}, \cite{Haugseng}, \cite{Ben-Zvi-Nadler:2016} or \cite{Ben-Zvi-Gunningham-Nadler} among others), sometimes referred to as the `push-pull construction'. In this thesis, we recast this construction to identify the requiered input data in a simple way that will be useful for applications. For example, it fits perfectly with Saito's theory of mixed Hodge modules (see Chapter \ref{chap:hodge}).

The idea of the construction is to consider an auxiliar category $\cC$ with pullbacks and final object, that is going to play the role of a category of fields (in the physical sense). Then, we are going to split our functor $Z$ as a composition
$$
	\Bord{n} \stackrel{\Fld{}}{\longrightarrow} \Span{\cC} \stackrel{\Qtm{}}{\longrightarrow} \Mod{R},
$$
where $\Span{\cC}$ is the category of spans of $\cC$ (see Section \ref{sec:quantisation}). The first arrow, $\Fld{}$, is the field theory and we will describe how to build it in Section \ref{sec:field-theory}. The second arrow, $\Qtm{}$, is the quantisation part. It is, in some sense, the most subtle piece of data. It will be constructed by means of something called a $\cC$-algebra (see Section \ref{sec:quantisation}). A $\cC$-algebra can be thought as collection of rings parametrized by $\cC$ with a pair of induced homomorphisms from every morphism of $\cC$.

\subsection{$\cC$-algebras}
\label{sec:C-algebras}
Let $\cC$ be a cartesian monoidal category (see \ref{ex:cartesian-monoidal}) and let $A: \cC \to \Rng$ be a contravariant functor. If $\star \in \cC$ is the final object, the ring $A(\star)$ plays a special role since, for every $a \in \cC$, we have an unique map $c_a: a \to \star$ which gives rise to a ring homomorphism $A(c_a): A(\star) \to A(a)$. Hence, we can see $A(a)$ as a $A(\star)$-module in a natural way. Such a module structure is the one considered in the first condition of Definition \ref{def:C-algebra} below.

Given $a, b,d \in \cC$ with morphisms $a \to d$ and $b \to d$, let $p_1: a \times_d b \to a$ and $p_2: a \times_d b \to b$ be the corresponding projections. We define the \emph{external product} over $d$
$$
	\boxtimes_d: A(a) \otimes_{A(d)} A(b) \to A(a \times_d b),
$$ by $z \boxtimes_d w = A(p_1)(z) \cdot A(p_2)(w)$ for $z \in A(a)$ and $w \in A(b)$. The external product over the final object will be denoted just by $\boxtimes = \boxtimes_\star: A(a) \otimes_{A(\star)} A(b) \to A(a \times b)$.

\begin{rmk}\label{rmk:module-exterior-prod}
For $b = \star$, we have that $p_1 = \lambda: a \times \star \to a$ is the unital isomorphism of the monoidal structure (see Definition \ref{defn:monoidal-cat}) so it gives raise to an isomorphism $A(\lambda): A(a) \to A(a \times \star)$ of rings and $A(\star)$-modules. Under this isomorphism, the external product $\boxtimes: A(a) \otimes A(\star) \to \cA(a \times \star) \cong A(a)$ is the given $A(\star)$-module structure on $A(a)$.
\end{rmk}

\begin{prop}\label{prop:A-preserves-external-prod}
The external product $\boxtimes: A(-) \otimes_{A(\star)} A(-) \Rightarrow A(- \times -)$ is a natural tranformation.
\begin{proof}
Let $f: a \to b$ and $g: a' \to b'$ be morphisms in $\cC$. In order to prove that $\boxtimes$ is natural, we have to check that $A(f \times g)(z \boxtimes w) = A(f)z \boxtimes A(g)w$ for $z \in A(b)$ and $w \in A(b')$. Let $p_1: a \times a' \to a$, $p_2: a \times a' \to a'$, $\rho_1: b \times b' \to b$ and $\rho_2: b \times b' \to b'$ be the corresponding projections. Observe that $\rho_1 \circ (f \times g) = f \circ p_1$ and $\rho_2 \circ (f \times g) = g \circ p_2$ so
\begin{align*}
	A(f)z \boxtimes A(g)w &= [A(p_1) \circ A(f)z] \cdot [A(p_2) \circ A(g)w] = A(f \circ p_1)z \cdot A(g \circ p_2)w \\
	&= A(\rho_1 \circ (f \times g))z \cdot A(\rho_2 \circ (f \times g))w\\
	& = [A(f \times g) \circ A(\rho_1) z] \cdot [A(f \times g) \circ A(\rho_2)w] = A(f \times g)(z \boxtimes w),
\end{align*}
as we wanted to prove.
\end{proof}
\end{prop}

\begin{defn}\label{def:C-algebra}
Let $(\cC, \times, \star)$ be a cartesian monoidal category. A \emph{$\cC$-algebra}, $\cA$, is a pair of functors
$$
	A: \cC^{op} \to \Rng, \hspace{1cm} B: \cC \to \Mod{A(\star)}
$$
such that:
\begin{itemize}
	\item They agree on objects, that is, $A(a) = B(a)$ for all $a \in \cC$, as $A(\star)$-modules.
	\item They satisfy the Beck-Chevaley condition (also known as the base change formula), that is, given $a_1, a_2, b \in \cC$ and a pullback diagram
\[
\begin{displaystyle}
   \xymatrix
   {
	d \ar[r]^{\;g'}\ar[d]_{f'} & a_1\ar[d]^{f} \\
	a_2 \ar[r]_{g} & b
   }
\end{displaystyle}
\]
we have that $A(g) \circ B(f) = B(f') \circ A(g')$.
	\item The external product $\boxtimes: B(-) \otimes_{A(\star)} B(-) \Rightarrow B(- \times -)$ is a natural transformation with respect to $B$.
\end{itemize}
\end{defn}

\begin{rmk}\label{rmk:notation-grothendieck-six}
\begin{itemize}
	\item A $\cC$-algebra can be thought as a collection of rings parametrized by $\cC$. For this reason, we will denote $\cA_a = A(a)$ for $a \in \cC$.
	\item The Beck-Chevalley condition appears naturally in the context of Grothendieck's yoga of six functors $f_*, f^*, f_!, f^!, \otimes$ and $\mathbb{D}$ in which $(f_*, f^*)$ and $(f_!, f^!)$ are adjoints, and $f^*$ and $f_!$ satisfy the Beck-Chevalley condition. 
In this context, we can take $A$ to be the functor $f \mapsto f^*$ and $B$ the functor $f \mapsto f_!$. Moreover, in order to get in touch with this framework, we will denote $A(f) = f^*$ and $B(f) = f_!$. For further infomation about Grothendieck' six functors, see for example \cite{Fausk-May:2003} or \cite{Ayoub}.
\end{itemize}
\end{rmk}

Using the covariant functor $B$, for every object $a \in \cC$, we obtain a $\cA_\star$-module homomorphism $(c_a)_!: \cA_a \to \cA_\star$. The special element $\mu(a) = (c_a)_!(1) \in \cA_\star$, where $1 \in \cA_a$ is the unit of the ring, will be called the \emph{measure} of $a$.

\begin{ex}
A toy example of these type of algebras can be given from commutative algebra. Let us denote by $\Rng_{fl}$ the subcategory of $\Rng$ of rings and flat ring homomorphisms i.e.\ a ring homomorphism $R \to R'$ is a morphism in $\Rng_{fl}$ if $R'$ is flat as $R$-module. We define the $\Rng_{fl}^{op}$-algebra $\K{(\Mod{(-)})}=(A,B)$ as follows:
	\begin{itemize}
		\item The contravariant functor $A: \Rng_{fl}^{op} \to \Rng$ takes a ring $R$ and assigns the Grothendieck ring $\K{(\Mod{R})}$ (see Section \ref{sec:grothendieck-groups}). If we have a morphism $f: R' \to R$ (which corresponds to a genuine flat ring homomorphism $f: R \to R'$), consider the functor $\Mod{R} \to \Mod{R'}$ given by $M \mapsto M \otimes_R R'$ for $M \in \Mod{R}$. Observe that, as $R'$ is a flat $R$-module, this functor descends to Grothendieck ring to give a ring homomorphism $A(f): \K{(\Mod{R})} \to \K{(\Mod{R'})}$. Observe that $\ZZ$ is the initial object of $\Rng_{fl}$ so it is the final object of $\Rng_{fl}^{op}$ and, thus, the ground ring of $\cA$ is $\K{(\Mod{\ZZ})} = \K{\Ab}$.
		\item The functor $B: \Rng_{fl}^{op} \to \Mod{\K{\Ab}}$ is given by restriction of scalars. That is, $B(R) = \K{(\Mod{R})}$, seen as a $\K{\Ab}$-module, and, given a morphism $f: R' \to R$ in $\Rng_{fl}^{op}$, it assigns the $K$-theory morphism associated to the functor $\Mod{R'} \to \Mod{R}$ given by restriction of scalars from $R'$ to $R$ via $f$ i.e.\ $M \mapsto M_R$ for $M \in \Mod{R'}$.
		\end{itemize}
By construction, these functors agree on objects. Moreover, pushouts in the category $\Rng$ are given by tensor product. Therefore, a cartesian square in $\Rng_{fl}^{op}$ (i.e.\ a pushout in $\Rng_{fl}$) is given by
\[
\begin{displaystyle}
   \xymatrix
   {
	S \otimes_R T \ar[r] \ar[d] & S \ar[d] \\
	T \ar[r] & R
   }
\end{displaystyle}
\]
Then, $A$ and $B$ satisfy the Beck-Chevalley condition since, for all $M \in \Mod{S}$, $M_R \otimes_R T \cong \left(M \otimes_S \left[S \otimes_R T\right]\right)_T$. In the same spirit, $B$ preserves external product, so $\K{(\Mod{(-)})}$ is a $\Rng_{fl}^{op}$-algebra.
\end{ex}

\begin{ex}\label{ex:algebra-sheaves}
Given a locally compact Hausdorff topological space $X$, let us consider the category $\Sh{X}$ of sheaves (of rational vector spaces) on $X$ with sheaf transformations between them (i.e.\ natural transformations). It is an abelian monoidal category with monoidal structure given by tensor product of sheaves. Given a continuous map $f: X \to Y$, we can induce two maps at the level of sheaves. The first one, called the \emph{inverse image} or pullback, $f^*: \Sh{Y} \to \Sh{X}$ takes a sheaf $\cF$ on $Y$ and assigns the sheaf on $X$ associated to the presheaf $U \mapsto \lim\limits_{\longrightarrow} \cF(V)$, where the direct limit is taken over all the open sets $V \subseteq Y$ such that $f(U) \subseteq V$. It is an exact functor (see \cite{Iversen:1986}, Section II.4).
	
	On the other way around, we also have a functor $f_!: \Sh{X} \to \Sh{Y}$, called the \emph{direct image with compact support} or the exceptional pushout. It is the sheaf associated to the presheaf that, given an open set $U \subseteq Y$, assigns the set of sections $s \in f^{-1}(U)$ such that $f: \supp{s} \to U$ is proper. In this case, $f_!$ is only left exact, so we can consider its derived functor $Rf_!: \Sh{X} \to \Dp{\Sh{Y}}$ (see Section \ref{sec:derived-categories}). In order to have some flavour about this map, observe that, for any sheaf $\cF$ on $X$ the stalk of $R^kf_!(\cF)$ at $y \in Y$ is $(R^kf_!(\cF))_y = H_c^k(f^{-1}(y), \cF)$. In this context, the base change theorem with compact support (see \cite{Dimca:2004}, Theorem 2.3.27) implies that, for any pullback diagram of locally Hausdorff spaces
\[
\begin{displaystyle}
   \xymatrix
   {
	X \times_Z Y \ar[r]^{\;\;\;\;\;\;g'} \ar[d]_{f'} & X \ar[d]^{f} \\
	Y \ar[r]_{g} & Z
   }
\end{displaystyle}
\]
there is a natural isomorphism $g^* \circ Rf_! \cong Rf'_! \circ g'^*$. Even more, $Rf_!$ preserves the external product.

This situation can be exploited to obtain a $\cC$-algebra. However, we must surround the annoying difficulty that, for a general topological space $X$, $R^kf_!\cF$ may not vanish for arbitrary large $k$. In order to solve this problem, let us restrict to the full subcategory $\Top_0$ of $\Top$ of locally compact Hausdorff topological spaces with finite cohomological dimension (or, even simpler, to $\Diff$), see \cite{Bredon:1997}. In this subcategory, we do have $Rf_!: \Sh{X} \to \Db{\Sh{Y}}$ for $f: X \to Y$ continuous. Hence, $Rf_!$ induces a map in $K$-theory $f_!: \K{\Sh{X}} \to \K{\Sh{Y}}$ (see Proposition \ref{prop:k-theory-algebra}) and analogously for $f^*: \K{\Sh{Y}} \to \K{\Sh{X}}$. Even more, $f^*$ is a ring homomorphism and $f_!$ is a module homomorphism over $\K{\Sh{\star}} = \K{(\Vect{\QQ})} = \ZZ$ (see Example \ref{ex:grothendieck-groups}). Thus, by the previous properties, $\K{\Sh{-}}$ is a $\Top_0$-algebra with $A: f \mapsto f^*$ and $B: f \mapsto f_!$. Observe that the unit object in $\K{\Sh{X}}$ is the image of constant sheaf $\underline{\QQ}_X$ on $X$ with stalk $\QQ$.

Moreover, given $X \in \Top_0$, let $c_X: X \to \star$ the projection onto the singleton set. Then, the measure of $X$ is the object $(c_X)_!(\underline{\QQ}_X)$ which is a sheaf whose unique stalk is $(c_X)_!(\underline{\QQ}_X)_\star = [H_c^\bullet(X; \QQ)] = \chi(X) \in \K{(\Vect{\QQ})} = \ZZ$. Hence, the measure of $X$ is nothing but the Euler characteristic of $X$ (with compact support), which is a kind of data compression of $X$. This justifies the fancy name `measure' of $X$.
\end{ex}

\begin{rmk}
\begin{itemize}
	\item The construction of the $\Top_0$-algebra in Example \ref{ex:algebra-sheaves} is paradigmatic in the sense that, in Theorem \ref{thm:mixed-hodge-modules-cvar-alg}, we will construct the $\Var{\CC}$-algebra of mixed Hodge modules following the same recipe.
	\item Analogous construction can be built by considering the category of algebraic varieties (or, more generally, quasi-compact and quasi-separated schemes) and quasi-coherent sheaves on them. 
	\item We can also consider the usual direct image functor $f_*\cF(U) = \cF(f^{-1}(U))$ for $U \subseteq Y$. It is related with the one with compact support via a sheaf morphism $f_! \to f_*$ that is an isomorphism for $f$ proper. For usual direct image, the Beck-Chevalley condition $g^* \circ Rf_* \cong Rf'_* \circ g'^*$ holds if $X$ is Hausdorff, $Y$ is locally compact Hausdorff and $f$ is universally closed (\cite{Milne:1980}, Theorem 17.3). In the context of algebraic geometry, base change also holds if $f$ is proper between noetherian schemes and flat sheaves \cite{Hartshorne}. The base change formula has also been extended to the context of stacks in \cite{Ben-Zvi-Francis-Nadler:2010}.
\end{itemize}
\end{rmk}

\begin{lem}\label{lem:iso-C-algebra}
Let $\cA$ be a $\cC$-algebra and let $f: a \to b$ be an isomorphism in $\cC$. Then, $f^* = (f^{-1})_!$. In particular, $f_!$ is also a ring homomorphism, $f^* \circ f_! = 1_a$ and $f_! \circ f^* = 1_b$.
\begin{proof}
Since $f$ is an isomorphism, the following square is a pullback
\[
\begin{displaystyle}
   \xymatrix
   {
	a \ar[r]^{1_a}\ar[d]_{1_a} & a \ar[d]^{f} \\
	a \ar[r]_{f} & b
   }
\end{displaystyle}
\]
Hence, $f^* \circ f_! = (1_a)_! \circ (1_a)^* = 1_{\cA_a}$. Since $f_!$ has a two-sided inverse $(f^{-1})_!$, we have that $f^* = (f^{-1})_!$, as claimed.
\end{proof}
\end{lem}

\begin{rmk}\label{rmk:C-algebra-as-correspondence}
There is an equivalent way of defining a $\cC$-algebra. Recall that a bivariant functor $\Phi: \cC \to \cD$ is an assigment that associates, to every $a \in \cC$ an object $\Phi(a) \in \cD$, and to every morphism $f: a \to b$ of $\cC$ a pair of morphisms $\Phi(f)_!: \Phi(a) \to \Phi(b)$ and $\Phi(f)^*: \Phi(b) \to \Phi(a)$. In this way, a $\cC$-algebra $\cA$ can be seen as a bivariant functor $\cA: \cC \to \Mod{R}$ such that $\Phi(a)$ has a prescribed ring structure for all $a \in \cC$, $\Phi(\star)=R$, $\Phi(f)^*$ is a ring homomorphism for all morphism $f$ of $\cC$ and $\Phi$ satisfies the six functors properties. In this fashion, the notion of $\cC$-algebra is the non-categorified analog of the sheaf theories of \cite{Gaitsgory-Rozenblyum}, Part III.

However, with a view towards its application in Chapters \ref{chap:hodge} and \ref{chap:representation}, it is very useful to consider separely the $A$ and the $B$ part of the $\cC$-algebra. In Chapter \ref{chap:hodge}, we will show that a $\cC$-algebra can be constructed via the $K$-theory of a derived functor. In that formulation, the functors $f^*$ and $f_!$ play, philosophically, different roles (actually, $f^*$ is exact and $f_!$ is not). Moreover, in some contexts, as for soft TQFTs (see Section \ref{sec:soft-TQFT}), only half of the $\cC$-algebra structure is available.
\end{rmk}

\subsection{Quantisation}
\label{sec:quantisation}

Given a category $\cC$ with pullbacks, we can construct the $2$-category of \emph{spans of $\cC$}, $\Span{\cC}$.  As described in \cite{Benabou}, the objects of $\Span{\cC}$ are the same as the ones of $\cC$. A morphism $a \to b$ in $\Span{\cC}$ is a span, that is, a triple $(d, f,g)$ of morphisms
$$
\xymatrix{
&\ar[dl]_{f} d \ar[dr]^{g}&\\
a&& b
}
$$
where $d \in \cC$. Given two spans $(d_1, f_1, g_1): a \to b$ and $(d_2, f_2, g_2): b \to c$, we define the composition $(d_2, f_2, g_2) \circ (d_1, f_1,g_1) = (d_1 \times_b d_2, f_1 \circ f_2', g_2 \circ g_1')$, where $f_2', g_1'$ are the morphisms in the pullback diagram
$$
\xymatrix{
&&\ar[dl]_{f'_2}d_1 \times_b d_2\ar[dr]^{g'_1}&&\\
&\ar[dl]_{f_1}d_1\ar[dr]^{g_1}&&\ar[dl]_{f_2}d_2\ar[dr]^{g_2}&\\
a&&b&&c
}
$$
Finally, a $2$-morphism $(d, f,g) \Rightarrow (d', f',g')$ between $a \stackrel{\hspace{3pt}f}{\leftarrow} d \stackrel{g}{\rightarrow} b$ and $a \stackrel{\hspace{5pt}f'}{\leftarrow} d' \stackrel{g'}{\rightarrow} b$ is a morphism $\alpha: d' \to d$ (notice the inverted arrow!) such that the following diagram commutes
\[
\begin{displaystyle}
   \xymatrix
   {	& d' \ar[rd]^{g'} \ar[ld]_{f'} \ar[dd]^\alpha & \\
   		a &  & b\\
	&d  \ar[ru]_{g} \ar[lu]^{f}& 
   }
\end{displaystyle}   
\]
Moreover, if $\cC$ is a monoidal category, $\Span{\cC}$ inherits a monoidal structure by tensor product on objects and morphisms.

\begin{prop}\label{prop-quantisation}
Let $\cA$ be a $\cC$-algebra, where $\cC$ is a category with final object $\star$ and pullbacks. Then, there exists a lax monoidal (strict) $2$-functor $\Qtm{\cA}: \Span{\cC} \to \Modt{\cA_\star}$ such that
$$
	\Qtm{\cA}(a) = \cA_a \hspace{1.5cm} \Qtm{\cA}(d, f, g) = g_! \circ f^*: \cA_a \to \cA_b,
$$
for $a, b \in \cC$ and a span $a \stackrel{f}{\leftarrow} d \stackrel{g}{\to} b$. The functor $\Qtm{\cA}$ is called the \emph{quantisation} of $\cA$.
\begin{proof}
More detailed, the functor $\Qtm{\cA}: \Span{\cC} \to \Modt{\cA_\star}$ is given as follows:
\begin{itemize}
	\item For any $a \in \cC$ we define $\Qtm{\cA}(a) = \cA_a$.
	\item Fixed $a, b \in \cC$, we define the functor
	$$
		(\Qtm{\cA})_{a,b}: \Hom_{\Span{\cC}}(a, b) \to \Hom_{\Modt{\cA_\star}}(\cA_a, \cA_b)$$ by:
	\begin{itemize}
		\item For a $1$-morphism $a \stackrel{f}{\leftarrow} d \stackrel{g}{\rightarrow} b$ we define $\Qtm{\cA}(d,f,g) = g_! \circ f^*: \cA_a \to \cA_d \to \cA_b$.
		\item For a $2$-morphism $\alpha: (d, f,g) \Rightarrow (d', f',g')$ given by $\alpha: d' \to d$, we define the endomorphism $\psi = \alpha_! \alpha^*: \cA_d \to \cA_{d}$. Since $\alpha$ is a $2$-cell in $\Span{\cC}$ we have that $g_! \circ \psi \circ f^* = (g')_! \circ (f')^*$ so $\psi$ is a twist $g_! \circ f^* \Rightarrow (g')_! \circ (f')^*$. Observe that, if $\alpha$ is an isomorphism, then $\psi = 1_{\cA_d}$ by Lemma \ref{lem:iso-C-algebra}.
	\end{itemize}
	\item For $(\Qtm{\cA})_{1_a}$, $a \in \cC$, we take the identity $2$-cell.
	\item Given $1$-morphisms $(d_1, f_1, g_1): a \to b$ and $(d_2, f_2, g_2): b \to c$ we have $(\Qtm{\cA})_{b,c}(d_2, f_2,g_2) \circ (\Qtm{\cA})_{a,b}(d_1, f_1,g_1) = (g_2)_!(f_2)^*(g_1)_!(f_1)^*$. On the other hand, $(\Qtm{\cA})_{a,c}((d_2, f_2,g_2) \circ (d_1, f_1, g_1))  =  (g_2)_! (g_1')_!  (f_2')^*  (f_1)^*$, where $g_1'$ and $f_2'$ are the maps in the pullback
	$$
\xymatrix{
&&\ar[dl]_{f'_2}d_1\times_b d_2\ar[dr]^{g'_1}&&\\
&\ar[dl]_{f_1}d_1\ar[dr]^{g_1}&&\ar[dl]_{f_2}d_2\ar[dr]^{g_2}&\\
a&&b&&c.
}
$$
By the Beck-Chevalley condition we have $(g_1')_!(f_2')^* = (f_2)^*(g_1)_!$ and the two morphisms agree. Thus, we can take the $2$-cell $(\Qtm{\cA})_{a,b,c}$ as the identity.
\end{itemize}

This proves that $\Qtm{\cA}$ is a strict $2$-functor. Furthermore, $\Qtm{\cA}$ is also lax monoidal taking $\Delta_{a, b}$ to be the external product $\boxtimes: \cA_a \otimes \cA_b \to \cA_{a \times b}$ as in Section \ref{sec:C-algebras}. In order to check that, suppose that we have spans $(d, f, g): a \to b$ and $(d', f', g'): a' \to b'$. Then, since both the pullback (Proposition \ref{prop:A-preserves-external-prod}) and the pushout maps (Definition \ref{def:C-algebra}) preserve the tensor product we have, for all $z \in \cA_{a}, w \in \cA_{a'}$, $g_!f^*z \boxtimes g'_!f'^*w = (g \times g')_!(f^*z \boxtimes f'^*w) = (g \times g')_!(f \times f')^*(z \boxtimes w)$. Therefore, the following diagram commutes
\[
\begin{displaystyle}
   \xymatrix
   {
   	\cA_a \otimes \cA_{a'} \ar[r]^{\boxtimes} \ar[d]_{g_!f^* \otimes g'_!f'^*} & \cA_{a \times a'} \ar[d]^{(g \times g')_!(f \times f')^*}\\
   	\cA_{b} \otimes \cA_{b'} \ar[r]_{\boxtimes} & \cA_{b \times b'}
   }
\end{displaystyle}   
\]
Hence, $\Delta$ is natural, as we wanted.
\end{proof}
\end{prop}

\begin{rmk}
\begin{itemize}
	\item If the functor $A$ defining a $\cC$-algebra $\cA = (A,B)$ is monoidal, then $\Qtm{\cA}$ is strict monoidal.
	\item If the Beck-Chevalley condition was satisfied up to natural isomorphism, then we would have equality after a pair of automorphisms in $\cA_{d_1}$ and $\cA_{d_2}$ which correspond to an invertible $2$-cell in $\Modt{\cA_\star}$. In this case, $\Qtm{\cA}$ would be a pseudo-functor. 
\end{itemize}
\end{rmk}

\begin{rmk}\label{rmk:non-invertible-2-cells}
If $\alpha: (d,f,g) \Rightarrow (d',f',g')$ is an invertible $2$-cell between spans, then $\Qtm{\cA}(d,f,g)= g_! \alpha_!\alpha^*f^* = (g')_!(f')^* = \Qtm{\cA}(d',f',g')$ since $\alpha_!\alpha^* = 1_{\cA_d}$ by Lemma \ref{lem:iso-C-algebra}. Thus, it seems reasonable to restrict the $2$-cells in $\Span{\cC}$ to isomorphisms and to drop out the cumbersome twisting in $\Mod{\cA_\star}$. However, there is a deep reason for keeping track the twisting. In Section \ref{sec:TQFT-over-sheaf}, we will consider the category of bordisms with an extra structure, $\EBord{n}{\cS}$, and a $2$-cell between two bordisms will be given by a morphism of the extra data. In this way, $\EBord{n}{\cS}$ has a natural $2$-category structure and, in general, its $2$-cells are non-invertible. Hence, under a TQFT, they become non-trivial twists.
\end{rmk}

The construction described in this section is strongly related to the sheaf theoretic formalism of \cite{Gaitsgory-Rozenblyum}, Part III (specially, Chapters 7 and 8). In that paper, it is proven that a bivalent functor $\Phi: \cC \to \cS$ satisfying the six functors formalism (which is the analogous of a $\cC$-algebra, see Remark \ref{rmk:C-algebra-as-correspondence}) is equivalent to a functor out of the $2$-category of correspondences $\textbf{Corr}_{vert, horiz}^{adm}(\cC)$ (Theorem 2.13, Chapter 7).

The category $\textbf{Corr}_{vert, horiz}^{adm}(\cC)$ is very similar to our category $\Span{\cC}$. The idea of this category is that we choose three subsemi-groups of morphisms $vert, horiz, adm$ of $\cC$, with $adm \subseteq vert \cap horiz$. Then, the objects of $\textbf{Corr}_{vert, horiz}^{adm}(\cC)$ are the ones of $\cC$, a $1$-morphism is a span $a \stackrel{f}{\leftarrow} d \stackrel{g}{\rightarrow} b$ with $f \in vert$ and $g \in horiz$, and a $2$-morphism is a morphism of spans given by a morphism in $adm$. In this way, $\textbf{Corr}_{vert, horiz}^{adm}(\cC)$ is a subcategory of $\Span{\cC}$ where the maps that take place in the cells are restricted. Under this interpretation, our Proposition \ref{prop-quantisation} is parallel to half of Theorem 2.13 of Chapter 7 of \cite{Gaitsgory-Rozenblyum}.

However, in the present thesis, we want to emphasize in the $\cC$-algebra part rather than its formulation as correspondences. One of the reasons for this is that $\cC$-algebras may be thought precisely as a quantized analog of the invariant under study. In this way, it is more natural to reformulate Saito's theory of mixed Hodge modules (see Chapter \ref{chap:hodge}) as a $\CVar$-algebra than as a functor out of correspondences. Actually, this interpretation underlies many of the references in Hodge theory, like \cite{Schurmann:2011} or \cite{Saito:2017}.

In any case, it would be very interesting to study in further detail the relation between the sheaf theoretic framework of \cite{Gaitsgory-Rozenblyum} and the $\cC$-algebras developed in this thesis. We hope that the precise statement of this interplay will allow us to formulate some of the $\cC$-algebras considered in this thesis, as mixed Hodge modules, in a sheaf theoretic way. That would produce a deeper insight in the (very tangled) category of mixed Hodge modules.

\subsection{Field theory}
\label{sec:field-theory}

\begin{figure}[h]
	\begin{center}
	\includegraphics[scale=0.23]{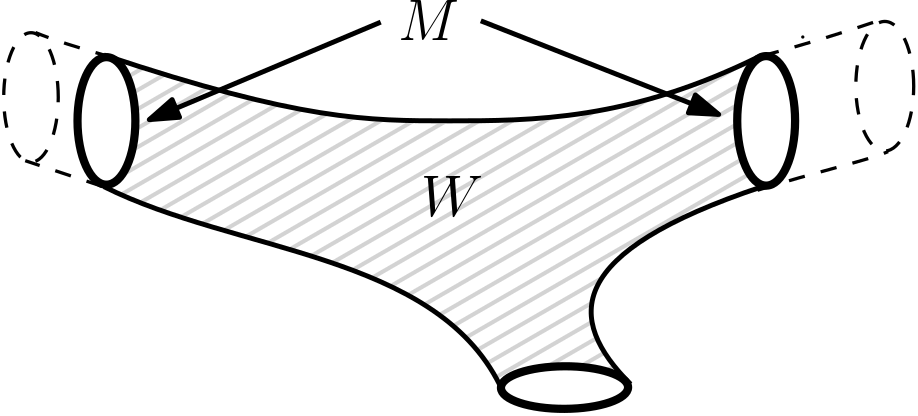}
	\caption{Collar extension $W \cup_M \left(M \times [0,1]\right)$.}
	\label{img:collar-extension}
	\end{center}
\vspace{-0.8cm}
\end{figure}

Let $W$ be a compact manifold with boundary and let $M \subseteq \partial W$ be a union of connected components of $\partial W$. The \emph{collar extension} of $W$ along $M$ is the open manifold (with boundary) $W \cup_M \left(M \times [0,1]\right)$, where the pasting is performed along $M \times \left\{0\right\}$ (see Figure \ref{img:collar-extension}).
Now, let $W_1$ and $W_2$ be smooth manifolds (maybe with boundary) of the same dimension. A \emph{tame smooth embedding} $f: W_1 \to W_2$ is a smooth map that, for some union of connected components $M \subseteq \partial W_1$, extends to a smooth embedding $\tilde{f}: W_1 \cup_M \left(M \times [0,1]\right) \to W_2$ with $\tilde{f}^{-1}(\partial W_2) = \partial W_1 - M$.

\begin{rmk}
Observe that $\partial \left[W_1 \cup_M \left( M \times [0,1]\right)\right] = \partial W_1 - M$ so the last condition means that $\tilde{f}$ sends the boundary of the collar extension $W_1 \cup_M \left( M \times [0,1]\right)$ into the boundary of $W_2$. We actually have that $M = \partial W_1 - f^{-1}(\partial W_2)$, that is, $M$ is the set of connected components of $\partial W_1$ that do not meet $\partial W_2$ under $f$. In this way, the tameness condition says precisely that $f$ can be slightly extended around these boundary components to give an open embedding.
\end{rmk}

Let us consider the category $\Embc$ whose objects are compact differentiable manifolds, maybe with boundary, and, given compact manifolds $M_1$ and $M_2$ of the same dimension, a morphism $M_1 \to M_2$ in $\Embc$ is a class of tame smooth embeddings $f: M_1 \to M_2$. Two such an embeddings $f, f': M_1 \to M_2$ are equivalent if there exists an ambient diffeotopy $h$ between them i.e.\ a smooth map $h: M_2 \times [0, 1] \to M_2$ such that $h_t = h(-,t)$ is a diffeomorphism of $M_2$ for all $0 \leq t \leq 1$, $h_0 = \Id_{M_2}$, $h_t|_{\partial M_2} = \Id_{\partial M_2}$ for all $t$, and $h_1 \circ f = f'$. There are no morphisms in $\Embc$ between manifolds of different dimensions.

\begin{ex}\label{ex:embc-map-boundary}
Let $W$ be a compact manifold and let $M \subseteq \partial W$ be a union of connected components of $\partial W$. Let $U \subseteq W$ be an open collaring around $M$, that is, an open subset of $W$ such that there exists a diffeomorphism $\varphi: M \times [0,1) \to U$. In that case, $\varphi: M \times [0, 1/2] \to U \subseteq W$ defines a morphism in $\Embc$. Moreover, such a morphism does not depend on the chosen collaring as any two collars are ambient diffeotopic. For this reason, we will denote this map $M \times [0, 1/2] \to W$ in $\Embc$ just by $M \to W$.
\end{ex}

Let $W_1$ and $W_2$ be $n$-dimensional compact manifolds and let $M$ be a $(n-1)$-dimensional closed manifold. Suppose that there exist tame embeddings $f_1: M \to W_1$ and $f_2: M \to W_2$ as in Example \ref{ex:embc-map-boundary}. This means that, for $k = 1,2$, $M$ is diffeomorphic to a union of connected components of $\partial W_k$ and $f_k$ is the collaring $f_k: M \times [0,1/2] \to W_k$ along that boundary.
In that case, the pushout of $f_1$ and $f_2$ in $\Embc$, $W_1 \cup_M W_2$, exists and it is given by the gluing of $W_1$ and $W_2$ along $M$.
\begin{figure}[h]
	\begin{center}
	\includegraphics[scale=0.35]{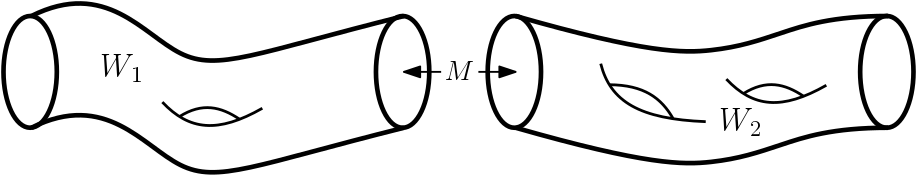}
	\end{center}
\vspace{-0.8cm}
\end{figure}

For this reason, we will call this special situation a \emph{gluing pushout}.

\begin{rmk}\label{rmk:gluing-embc}
On the category $\Diffc$, pushouts may or may not exist. For example, take $f_1,f_2: \left\{1\right\} \to S^1$ as the inclusions onto $1 \in S^1$. In the category $\Top$, the pushout of these diagram is the wedge sum of two $S^1$ which is no longer a differentiable manifold. However, the forgetful functor $\Diffc \to \Top$ does not preserves colimits, so we cannot assure that such a wedge sum has to be the pushout. In order to solve it, the category of Fr\"olicher spaces, $\textbf{Fr\"ol}$, can be used as an intermediate category between $\Diffc$ (or even $\Diffc$) and $\Top$ (see \cite{Kriegl-Michor:1997}, Theorem 23.2). This category always has limits and colimits and the forgetful functor $\Embc \to \textbf{Fr\"ol}$ preserves them. It can be proven directly that, in $\textbf{Fr\"ol}$, the pushout of $f_1, f_2$ is the wedge sum of two copies of $S^1$, completing the proof that such a pushout cannot exist in $\Diffc$. In the same spirit, pushouts in $\Embc$ may not exist.
\end{rmk}

\begin{defn}
Given a category $\cC$ with final object, a contravariant functor $F: \Embc \to \cC$ is said to have the \emph{Seifert-van Kampen property} if $F$ sends gluing pushouts into pullbacks and sends the initial object of $\Embc$ (i.e.\ $\emptyset$) into the final object of $\cC$.
\end{defn}

\begin{prop}\label{prop:field-theory}
Let $\cC$ be a category with final object and pullbacks and let $F: \Embc \to \cC$ be a contravariant functor satisfying the Seifert-van Kampen property. Then, there exists a monoidal functor $\cF_F: \Bord{n} \to \Span{\cC}$ such that
$$
	\cF_F(M) = F(M), \hspace{1cm} F(M_1) \stackrel{F(i_1)}{\leftarrow} F(W) \stackrel{F(i_2)}{\rightarrow} F(M_2),
$$
for all objects $M, M_1, M_2 \in \Bord{n}$ and bordisms $W: M_1 \to M_2$ where $i_k: M_k \to W$ are the inclusions. In this situation, the functor $F$ is called the \emph{geometrisation} and $\cF_F$ is called the \emph{field theory} of $F$.
\begin{proof}
The complete definition of $\cF$ is given as follows. Given $M \in \Bord{n}$, $M$ is a compact manifold so we can define $\cF_F(M) = F(M)$. With respect to morphisms, given a bordism $W: M_1 \to M_2$ in $\Bord{n}$, let $i_1: M_1 \to W$ and $i_2: M_2 \to W$ be the inclusions of $M_1$ and $M_2$ as boundaries of $W$ (as in Example \ref{ex:embc-map-boundary}). To $W$, $\cF_F$ assigns the span
\[
\begin{displaystyle}
   \xymatrix
   {
   	F(M_1) & F(W) \ar[r]^{F(i_2)}\ar[l]_{F(i_1)} & F(M_2) 
   }
\end{displaystyle}   
\]
This assignment is a functor. In order to check it, let $W: M_1 \to M_2$ and $W': M_2 \to M_3$ be two bordisms with inclusions $i_k: M_k \to W$ and $i_k': M_k \to W'$. As we mentioned in Remark \ref{rmk:gluing-embc}, as a manifold, $W' \circ W$ is the gluing pushout of $\Embc$
\[
\begin{displaystyle}
   \xymatrix
   {
  	M_2 \ar[d]_{i_2'}\ar[r]^{i_2} & W \ar[d]^{j} \\
  	W' \ar[r]_{j'\;\;\;\;}  & W' \circ W
   }
\end{displaystyle}   
\]
and, since $F$ sends gluing pushouts into pullbacks, the following diagram is a pullback in $\cC$
\[
\begin{displaystyle}
   \xymatrix
   {
  	F(W' \circ W) \ar[r]^{F(j)} \ar[d]_{F(j')} & F(W) \ar[d]^{F(i_2)} \\
  	F(W') \ar[r]_{F(i_2')} & F(M_2)
   }
\end{displaystyle}   
\]
Therefore, $\cF_F(W') \circ \cF_F(W)$ is given by the span
	$$
\xymatrix{
&&\ar[dl]_{F(j)} F(W' \circ W) \ar[dr]^{F(j')}&&\\
&\ar[dl]_{F(i_1)}F(W)\ar[dr]^{F(i_2)}&&\ar[dl]_{F(i_2')}F(W')\ar[dr]^{F(i_3')}&\\
F(M_1)&&F(M_2)&&F(M_3).
}
$$
Since $i_1 \circ j$ and $i_3' \circ j'$ are the inclusions onto $W' \circ W$ of $M_1$ and $M_3$ respectively, the previous span is also $\cF_F(W' \circ W)$, as we wanted.

For the monoidality, let $M_1, M_2 \in \Bord{n}$. As the coproduct $M_1 \sqcup M_2$ can be seen as a gluing pushout along $\emptyset$, $F$ sends coproducts in $\Bord{n}$ into products on $\cC$. Hence, $F(M_1 \sqcup M_2) = F(M_1) \times F(M_2)$. Since the monoidal structure on $\Span{\cC}$ is given by products on $\cC$, monoidality holds for objects. For morphisms, the argument is analogous.
\end{proof}
\end{prop}

Putting together Proposition \ref{prop-quantisation} and Proposition \ref{prop:field-theory}, we obtain the main theorem of this section.

\begin{thm}[Lagrangian formulation of TQFTs]
\label{thm:physical-constr}
Let $\cC$ be a category with final object $\star$ and pullbacks. Given a functor $F: \Embc \to \cC$ with the Seifert-van Kampen property and a $\cC$-algebra $\cA$, there exists a lax monoidal Topological Quantum Field Theory
$$
	Z_{F, \cA}: \Bord{n} \to \Mod{\cA_\star}.
$$
\end{thm}

\begin{rmk}
\label{rmk:expression-TQFT}
From the explicit construction given in the proofs of Propositions \ref{prop-quantisation} and \ref{prop:field-theory}, the functor $Z_{F, \cA}$ satisfies:
\begin{itemize}
	\item For an object $M \in \Bord{n}$, it assigns $Z_{F, \cA}(M) = \cA_{F(M)}$.
	\item For a bordism $W: M_1 \to M_2$, it assings $Z_{F, \cA}(W) = F(i_2)_! \circ F(i_1)^*: \cA_{F(M_1)} \to \cA_{F(M_2)}$.
	\item For a closed $n$-dimensional manifold $W$, seen as a bordism $W: \emptyset \to \emptyset$, the homomorphism $Z_{F, \cA}(W): \cA_\star \to \cA_\star$ is given by multiplication by the measure $\mu({F(W)}) \in \cA_\star$. This follows from the fact that, since $F(\emptyset)=\star$, the inclusion $i: \emptyset \to W$ gives the projection $c_W = F(i): F(W) \to \star$. Hence, we have, $Z_{F, \cA}(W)(1_\star) = (c_W)_! \circ (c_W)^*(1_\star) = (c_W)_!(1_{F(W)}) = \mu(F(W))$. Here, $1_\star$ and $1_{F(W)}$ denote the unit in $\cA_\star$ and $\cA_{F(W)}$ respectively.
\end{itemize}
\end{rmk}

\begin{rmk}\label{rmk:field-theory-very-lax}
We can also consider the case in which the geometrisation functor $F: \Embc \to \cC$ no longer has the Seifert-van Kampen property, but it still maps the initial object into the final object. In that case, the image of a gluing pushout under $F$ is not a pullback but, as for any functor, it is a cone. Suppose that we have two bordisms $W: M_1 \to M_2$ and $W': M_2 \to M_3$ that fit in the gluing pushout in $\Embc$
\[
\begin{displaystyle}
   \xymatrix
   {
   		M_2 \ar[r]\ar[d] & W\ar[d] \\
   		W' \ar[r] & W' \circ W
   }
\end{displaystyle}   
\]
By definition of pullback, there exists an unique morphism $\phi: F(W' \circ W) \to F(W') \times_{F(M_2)} F(W)$ in $\cC$ such that the induced diagram commutes
\[
\begin{displaystyle}
   \xymatrix
   {
   		F(W' \circ W) \ar@{--{>}}[rd]^{\phi} \ar@/^1pc/[rrd] \ar@/_1pc/[ddr] & & \\
   		&F(W') \times_{F(M_2)} F(W)  \ar[r]\ar[d]  & F(W') \ar[d] \\
   		&F(W)  \ar[r] & F(M_2)
   }
\end{displaystyle}  
\]
This morphism $\phi$ induces a $2$-morphism $\cF_F(W') \circ \cF_F(W) \Rightarrow \cF_F(W' \circ W)$. Therefore, in this case, $\cF: \Bord{n} \to \Span{\cC}$ is no longer a functor but a lax $2$-functor (see Definition \ref{defn:lax-2-funct}). Thus, the induced functor
$$
	Z_{F, \cA}: \Bord{n} \to \Modt{\cA_{\star}}
$$
is a lax monoidal symmetric lax $2$-functor (recall that $\Qtm{\cA}: \Span{\cC} \to \Modt{\cA_{\star}}$ is a $2$-functor). We will call such functors \emph{very lax Topological Quantum Field Theories}.
\end{rmk}

\begin{rmk}
The first version of this theorem was described in \cite{GPLM-2017}. However, there the method is built \textit{ad hoc} for the lax monoidal TQFT that computes mixed Hodge structures (Section 4.1). In that paper, the field theory procedure described here corresponds to the geometrisation part of the construction. We thank D. Ben-Zvi for suggesting us the name `field theory' and `formulation of lagrangian field theory' based on the physical interpretation of TQFTs.
\end{rmk}

As a final remark, it is customary in the literature to focus on the field theory as a functor $\Fld{}: \Bord{n} \to \Span{\cC}$ and on the quantisation as a functor $\Qtm{}: \Span{\cC} \to \Mod{R}$, and to forget about the geometrisation and the $\cC$-algebra, despite that they underlie the whole construction (see, for example, \cite{Freed-Hopkins-Lurie-Teleman:2010} or \cite{Freed-Hopkins-Teleman:2010}). In this form, the Seifert-van Kampen property and the base change property stay encoded in the functoriality of $\Fld{}$ of $\Qtm{}$, respectively, and many properties follow in a simpler way so this description gives a deeper insight in the properties of the TQFT.

However, in this thesis we want to emphasize the role of the geometrisation $F$ and the $\cC$-algebra $\cA$ in the construction. The reason is that we are constructing TQFTs with a view towards the creation of new effective computational methods of algebraic invariants. In this way, $\cA$ determines the algebraic invariant under study and $F$ determines the object for which we are going to compute the invariant. As an example of this principle, in Chapter \ref{chap:representation}, we will build a TQFT that computes the Hodge structure of representation varieties. As ingredients for that TQFT, we will take $F$ as the representation variety functor and $\cA$ will be the $\CVar$-algebra of mixed Hodge modules, that might be thought as a variational version of Hodge structures. This principle opens the door to the development of new computational methods based on TQFTs, just by modifying the input data: the geometrisation and the $\cC$-algebra.

\section{Other versions of Topological Quantum Field Theories}
\label{sec:other-versions-tqft}

The previous formulation of TQFTs can be refined in order to obtain several variants. Such modified versions will be useful in Chapter \ref{chap:representation} since, there, the TQFT that suits the problem is an extended version of the one constructed in the previus section.

In this section, we will discuss four types of variations. The first one concerns TQFTs that preserve a sheaf that can be thought as an extra structure on the category of bordisms. 
The second variant are the so-called almost-TQFT. Recall that the TQFTs constructed by the Lagrangian formulation are, in general, far from being monoidal and the problem is that the elbow bordism breaks this monoidality (see Remark \ref{rmk:finite-dim-TQFT}). However, if we restrict ourselves to the subcategory $\Tub{n}$ of $\Bord{n}$ of tubes (i.e.\ of bordisms with only at most one connected component for each boundary) then we remove all the sources of problems and the resultant functor happens to be monoidal.

The third construction will be very relevant in Chapter \ref{chap:representation}. The idea is that a TQFT might be too involved for computational purposes because the field $\Fld{F}(M)$ is too complex. However, if we have a collection of morphisms that send $\Fld{F}(M)$ to a simplest field, the complexity of the calculation can be reduced drastically. As payback, we lose the strict functoriality property, but it can be fixed by means of some extra data.

To finish the section, we will introduce soft Topological Quantum Field Theories. The point is that the construction of a $\cC$-algebra may be too hard for some applications, but having only `half of a $\cC$-algebra' may be feasible. In that case, instead of a quantisation, we can build a pre-quantisation and, by means of it, construct a sort of TQFT with vales in the category of bimodules.

\subsection{TQFTs over a sheaf}
\label{sec:TQFT-over-sheaf}

In this section, we will consider a larger class of Topological Quantum Field Theories. For that purpose, we will endow the bordisms with an extra algebraic structure, described in terms of a sheaf. This idea is recurrent in the literature and, actually, the classical definition of a TQFT requires not pure bordisms, but oriented bordisms \cite{Atiyah:1988}. In fact, the celebrated Lurie's cobordism hypothesis is stated for framed bordisms, i.e.\ bordisms with a chosen trivialization of the trivial bundle or, in general, with a chosen $(X, \xi)$-structure, where $X$ is a CW-complex and $\xi$ is a real vector bundle of rank $r$ (see Theorem 2.4.18 of \cite{Lurie-eTQFT:2009}). The formulation of these extra structures as sheaves goes back to \cite{Madsen-Weiss:2007} (see also \cite{Galatius-Tillmann-Madsen-Weiss:2009} and \cite{Ayala:2009}).

The consideration of extra structures on bordisms has two advantages. The first one is that these TQFTs have an additional structure of $2$-functor that gives extra data about the topological problem. The second one is that this enriched version of TQFTs may be built in contexts in which a classical TQFTs cannot be constructed (at least, in a natural way).

Let us consider a sheaf $\cS: \Embc \to \Cat$ with $\cS(\emptyset)=\emptyset$. By a sheaf we mean the following (c.f.\ \cite{Madsen-Weiss:2007}). Let $W$ be a compact manifold and suppose that we have a collection of tame morphisms $i_\alpha: M_\alpha \to W$, for $\alpha \in \Lambda$, such that $W = \cup_\alpha i_\alpha(M_\alpha)$. Then, a sheaf $\cS$ is a contravariant functor such that, for any collection $s_\alpha \in \cS(M_\alpha)$ with $\cS(i_{\alpha,\beta})(s_\alpha) = \cS(i_{\beta,\alpha})(s_\beta)$ for all $\alpha,\beta \in \Lambda$ (where $i_{\alpha,\beta}: M_\alpha \cap i_\alpha^{-1} \left(i_\beta (M_\beta)\right) \to M_\alpha$ are the restrictions) there exists an unique $s \in \cS(M)$ such that $\cS(i_\alpha)(s) = s_\alpha$ for all $\alpha \in \Lambda$.

\begin{rmk}
We insist that, in contrast with usual sheaves, the coverings considered here are closed. However, they also give an open covering since the embeddings $i_\alpha$ have to extend to a small open set around $M_\alpha$.
\end{rmk}

Given such a sheaf $\cS$, we define the \emph{category of embeddings over $\cS$}, or just $\cS$-embeddings, $\EEmbc{\cS}$. The objects of this category are pairs $(M, s)$ with $M \in \Embc$ and $s \in \cS(M)$ (and the initial object $\emptyset$). Given objects $(M_1, s_1), (M_2, s_2) \in \EEmbc{\cS}$, a morphism between them is a pair $(f, \alpha)$ where $f: M_1 \to M_2$ is a morphism in $\Embc$ and $\alpha: s_1 \to \cS(f)(s_2)$ is a morphism of $\cS(M_1)$. The composition of two morphisms $(f, \alpha)$ and $(f', \alpha')$ is $(f' \circ f, \cS(f)(\alpha') \circ \alpha)$.

In analogy with Section \ref{sec:field-theory}, given two morphisms $(f_1, \alpha_1): (M_1, s_1) \to (W, s)$ and $(f_2, \alpha_2): (M_2, s_2) \to (W, s)$ in $\EEmbc{\cS}$, we will call them a \emph{gluing pushout} if $f_1$ and $f_2$ are gluing pushouts of $\Embc$, $\alpha_1 = 1_{s_1}$ and $\alpha_2 = 1_{s_2}$. Hence, if $\cC$ is a category with final object and pullbacks and $F: \EEmbc{\cS} \to \cC$ is a contravariant functor, we will say that $F$ satisfies the \emph{Seifert-van Kampen property} if it sends the initial object of $\EEmbc{\cS}$ into the final object of $\cC$ and assigns gluing pushouts to pullbacks.

Analogously, given $n \geq 1$, we define the \emph{category of $n$-bordisms over $\cS$}, $\EBord{n}{\cS}$. It is a $2$-category given by the following data:
\begin{itemize}
	\item Objects: The objects of $\EBord{n}{\cS}$ are pairs $(M, s)$ where $M$ is a $(n-1)$-dimensional closed manifold and $s \in \cS(M)$. Observe that, if $M \neq \emptyset$ has $\cS(M)=\emptyset$, then $M$ does not appear as object of $\EBord{n}{\cS}$.
	\item $1$-morphisms: Given objects $(M_1, s_1)$, $(M_2, s_2)$ of $\EBord{n}{\cS}$, a morphism $(M_1, s_1) \to (M_2, s_2)$ is a class of pairs $(W, s)$, where $W: M_1 \to M_2$ is a morphism in $\Bord{n}$ such that $\cS(i_1)(s)=s_1$ and $\cS(i_2)(s)=s_2$, $i_k: M_k \to W$ being the inclusions. Two bordisms $(W, s)$ and $(W', s')$ are in the same class if there exists a boundary preserving diffeomorphism $f: W \to W'$ such that $\cS(f)(s')=s$.
	\\
With respect to the composition, given $(W, s): (M_1, s_1) \to (M_2, s_2)$ and $(W', s'): (M_2, s_2) \to (M_3, s_3)$, we define $(W',s') \circ (W,s)$ as the morphism $(W \cup_{M_2} W', s \cup s'): (M_1, s_1) \to (M_3, s_3)$ where $W \cup_{M_2} W'$ is the usual gluing of bordisms along $M_2$ and $s \cup s' \in \cS(W \cup_{M_2} W')$ is the object given by gluing $s$ and $s'$ with the sheaf property of $\cS$.
	\item $2$-morphisms: Given two $1$-morphisms $(W, s), (W',s'): (M_1, s_1) \to (M_2, s_2)$, a $2$-cell $(W, s) \Rightarrow (W',s')$ is a pair $(f, \alpha)$, where $f$ is a diffeomorphism of bordisms $f: W \to W'$ and $\alpha$ is a morphism $\alpha: s \to \cS(f)(s')$. Composition of $2$-cells $(f, \alpha)$ and $(f', \alpha')$ is just $(f' \circ f, \cS(f)(\alpha') \circ \alpha)$.
\end{itemize}

\begin{rmk}\label{rmk:cobord-hypoth}
The equipment of a $2$-category structure on the category of bordisms is a recurrent idea for improving the understanding of the topological nature of the TQFT. A similar idea appears in Extended Topological Quantum Field Theories (see \cite{Lurie-eTQFT:2009} or \cite{Baez-Dolan:1995}), where the extra structure allows their full classification via cobordism's hypothesis. 
\end{rmk}

In this form, $\EBord{n}{\cS}$ is not exactly a category since, for $(M,s) \in \EBord{n}{\cS}$, it may be no unit morphism in the category $\Hom_{\EBord{n}{\cS}} ((M,s),(M,s))$. In that case, it can be solved by weakening slightly the notion of bordism, allowing that $M$ itself could be seen as a bordism $M: M \to M$. With this modification, $(M,s): (M,s) \to (M,s)$ is the desired unit and it is a straightforward check to see that $\EBord{n}{\cS}$ is a (strict) $2$-category. Furthermore, as for $\Bord{n}$, it has a natural monoidal structure by means of disjoint union.

\begin{defn}
Let $\cS: \Embc \to \Cat$ be a sheaf and let $R$ be a ring. A \emph{(lax monoidal) Topological Quantum Field Theory over $\cS$}, shorten as (lax monoidal) $\cS$-TQFT, is a (lax) monoidal symmetric (strict) $2$-functor $Z: \EBord{n}{\cS} \to \Modt{R}$.
\end{defn}

In analogy with the procedure of Theorem \ref{thm:physical-constr}, TQFTs over a sheaf also has a physical-inspired construction method. For that, we only need to adjust the previous construction to take into account the sheaf.

\begin{prop}
Let $\cC$ be a category with final object and pullbacks and let $F: \EEmbc{\cS} \to \cC$ be a contravariant functor satisfying the Seifert-van Kampen property. Then, there exists a monoidal $2$-functor $\cF_F: \EBord{n}{\cS} \to \Span{\cC}$ such that
$$
	\cF_F(M,s) = F(M,s), \hspace{1cm} F(M_1,s_1) \stackrel{F(i_1,1_{s_1})}{\longleftarrow} F(W,s) \stackrel{F(i_2,1_{s_2})}{\longrightarrow} F(M_2,s_2),
$$
for all objects $(M,s), (M_1, s_1), (M_2, s_2) \in \EBord{n}{\cS}$ and bordisms $(W, s): (M_1, s_1) \to (M_2, s_2)$ where $i_k: M_k \to W$ are the inclusions.
\begin{proof}
The exhaustive definition of $\cF_F$ is analogous to the one of Proposition \ref{prop:field-theory} and the proof that it is a monoidal functor works verbatim in this context. For the $2$-functor structure, suppose that $(f, \alpha)$ is a $2$-cell between $1$-morphisms $(W, s), (W', s'): (M_1, s_1) \to (M_2, s_2)$. In that case, $(f, \alpha)$ is also a morphism in $\EEmbc{\cS}$ so we obtain a morphism $F(f, \alpha): (W', s') \to (W, s)$ fitting in the commutative diagram
\[
\begin{displaystyle}
   \xymatrix
   {
   	& F(W', s') \ar[rd]^{F(i_2', 1)} \ar[dd]_{F(f,\alpha)} \ar[ld]_{F(i_1', 1)} & \\
	F(M_1, s_1) & & F(M_2, s_2) \\
	& F(W, s) \ar[ru]_{F(i_2, 1)} \ar[lu]^{F(i_1, 1)} &
   }
\end{displaystyle}
\]
This produces the desired $2$-morphism in $\Span{\cC}$.
\end{proof}
\end{prop}

Combining this result with Proposition \ref{prop-quantisation} we obtain the following result, analogous to Theorem \ref{thm:physical-constr}.

\begin{thm}\label{thm:physical-constr-sheaf}
Let $\cC$ be a category with final object $\star$ and pullbacks. Given a functor $F: \EEmbc{\cS} \to \cC$ with Seifert-van Kampen property and a $\cC$-algebra $\cA$, there exists a lax monoidal Topological Quantum Field Theory over $\cS$
$$
	Z_{F, \cA}: \EBord{n}{\cS} \to \Modt{\cA_\star}.
$$
\end{thm}

The description of the obtained TQFT given in Remark \ref{rmk:expression-TQFT} also holds, with the corresponding modifications. In particular, if $W$ is a closed $n$-manifold and $s \in \cS(W)$, then $Z_{F, \cA}(W,s): \cA_\star \to \cA_\star$ is given by multiplication by the measure $\mu({F(W, s)}) \in \cA_\star$. With respect to the $2$-functor structure, given a $2$-morphism $(f, \alpha): (W, s) \Rightarrow (W', s')$ in $\EBord{n}{\cS}$, it assigns the twist given by the endomorphism $F(f, \alpha)_! \circ F(f, \alpha)^*: \cA_{F(W,s)} \to \cA_{F(W,s)}$.

\begin{rmk}\label{rmk:field-theory-very-lax-sheaf}
As in Remark \ref{rmk:field-theory-very-lax}, if $F: \EEmbc{\cS} \to \cC$ no longer has the Seifert-van Kampen property but it still maps the initial object into the final object, then $\Fld{F}: \EBord{n}{\cS} \to \Span{\cC}$ is a lax $2$-functor (see Definition \ref{defn:lax-2-funct}). Thus, the induced functor $Z_{F, \cA}: \EBord{n}{\cS} \to \Modt{\cA_{\star}}$ is a lax monoidal symmetric lax $2$-functor, i.e.\ a very lax Topological Quantum Field Theory.
\end{rmk}

In the context of bordisms with sheaves, the requirement of the whimsical twisting structure on $\Modt{\cA_\star}$ becomes evident. The existence of a non-invertible $2$-cell $(f,\alpha): (W, s) \Rightarrow (W, s')$ reflects a non-invertible morphism $\alpha: s \to \cS(f)(s')$ that can be interpreted as a restriction of the extra structure sheaf (e.g.\ in the case of unordered configurations of points, as in the following section). As explained in Remark \ref{rmk:non-invertible-2-cells}, these non-invertible cells become non-trivial twists in $\Modt{\cA_\star}$ that compare how the homomorphism changes under morphisms of sheaves.
	
Moreover, as a future work, we might extend the category of bordisms to also comprise manifolds with singular points. In this way, degenerations of bordisms can be encoded into non-invertible $2$-cells. This idea opens the door to formulate some observed degeneration phenomena in the language of TQFTs, as the work on monodromy manifolds and character varieties in \cite{Chekhov-Mazzocco-Rubtsov:2017} and \cite{Chekhov-Mazzocco-Rubtsov:2017b}. In those papers, it is proven that the well-known confluence scheme of Painlev\'e equations is a consequence of the degeneration of decorated character varieties over holed spheres with nodal singularities under the so-called `chewing-gum' operations on surfaces. We hope that the kind of ideas sketched here will be useful to formulate this confluence scheme as the image, under a TQFT, of a diagram of non-invertible $2$-cells in the category of $2$-dimensional bordisms.

\subsubsection*{Some useful extra structures}
\label{sec:some-extra-struct}

In this section, we will describe in detail some examples of sheaves, some of which will be used in Chapter \ref{chap:representation}. As we will see, the extra structure given by these sheaves is necessary for defining a field theory that captures the geometric nature of representation varieties.

The first example is the sheaf $\cS^{or}: \Embc \to \Cat$ of orientations. Given a compact orientable manifold $W$ with $m$ connected components, the category $\cS^{or}(W)$ has $2^m$ objects, one for each possible orientation of $W$, and no morphisms different than the units. If $W$ is non-orientable, then $\cS^{or}(W)$ is the empty category. For a tame embedding $f: W_1 \to W_2$, the functor $\cS^{or}(f)(W_2) \to \cS^{or}(W_1)$ sends an orientation in $W_2$ into the induced orientation in $W_1$ (recall that they have the same dimension). Observe that, in this case, $\EEmbc{\cS^{or}}$ has, as objects, oriented compact manifolds and, as morphisms, orientation-preserving differentiable maps. In this way, $\EBord{n}{\cS^{or}}$ is nothing but the usual category of oriented bordisms, $\Bordo{n}$ and an $\cS^{or}$-TQFT is the same as an oriented TQFT.

Another important example is the sheaf $\cSp: \Embc \to \Cat$ of unordered configurations of points (or the sheaf of pairs). For a compact manifold $W$, the category $\cSp(W)$ has, as objects, finite subsets $A \subseteq W$ meeting every connected component and every boundary component of $W$. Given two finite subsets $A_1, A_2 \subseteq W$, a morphism $A_1 \to A_2$ in $\cSp(W)$ is an inclusion $A_1 \subseteq A_2$. With respect to morphisms, if $f: W_1 \to W_2$ is an embedding, the functor $\cSp(f): \cSp(W_2) \to \cSp(W_1)$ is given as follows. For $A \subseteq W_2$, it assigns $\cSp(f)(A) = f^{-1}(A) \in \cSp(W_1)$ and, if we have an inclusion $A_1 \subseteq A_2$ then it gives the inclusion $f^{-1}(A_1) \subseteq f^{-1}(A_2)$ as a morphism in $\cSp(W_1)$. It is straighforward to check that $\cSp$ is a sheaf.

The associated categories $\EEmbc{\cSp}$ and $\EBord{n}{\cSp}$ will be called the \emph{category of embeddings of pairs} and the \emph{category of $n$-bordisms of pairs} and will be denoted by $\Embpc$ and $\Bordp{n}$, respectively. Explicitly, the later is given by the following data:
\begin{itemize}
	\item Objects: The objects of $\Bordp{n}$ are pairs $(M, A)$ where $M$ is a $(n-1)$-dimensional closed manifold together with a finite subset of points $A \subseteq M$ that intersects every connected component of $M$.
	\item $1$-morphisms: Given objects $(M_1, A_1)$, $(M_2, A_2)$ of $\Bordp{n}$, a morphism $(M_1, A_1) \to (M_2, A_2)$ is a class of pairs $(W, A)$ where $W: M_1 \to M_2$ is a bordism and $A \subseteq W$ is a finite set meetting every connected component of $W$ such that $M_1 \cap A = A_1$ and $M_2 \cap A = A_2$. Two pairs $(W, A), (W',A')$ are in the same class if there exists a diffeomorphism of bordisms  $f: W \to W'$ such that $f(A)=A'$. The composition of $(W, A): (M_1, A_1) \to (M_2, A_2)$ and $(W', A'): (M_2, A_2) \to (M_3, A_3)$ is given by the morphism $(W \cup_{M_2} W', A \cup A'): (M_1, A_1) \to (M_3, A_3)$.
	\item $2$-morphisms: Given two $1$-morphisms $(W, A), (W',A'): (M_1, A_1) \to (M_2, A_2)$, a $2$-cell $(W, A) \Rightarrow (W',A')$ is a diffeomorphism of bordisms $f: W \to W'$ such that $A \subseteq f^{-1}(A') $.
\end{itemize}

A (lax monoidal) $\cSp$-TQFT will be referred to as a (lax monoidal) \emph{Topological Quantum Field Theory of pairs}.

Another important sheaf for applications is the so-called sheaf of parabolic structures. The starting point is a fixed set $\Lambda$ that we will call the \emph{parabolic data}. We define $\cSpar{\Lambda}: \Embc \to \Cat$ as the following functor. For a compact manifold $W$, $\cSpar{\Lambda}(W)$ is the category whose objects are (maybe empty) finite sets $Q = \left\{(S_1, \lambda_1), \ldots, (S_r, \lambda_r)\right\}$, with $\lambda_i \in \Lambda$, called \emph{parabolic structures} on $W$. The $S_i$ are pairwise disjoint compact submanifolds of $W$ of codimension $2$ with a co-orientation (i.e.\ an orientation of their normal bundle) such that $S_i \cap \partial M = \partial S_i$ transversally. A morphism $Q \to Q'$ between two parabolic structures in $\cSpar{\Lambda}(W)$ is just an inclusion $Q \subseteq Q'$.

Moreover, suppose that we have a tame embedding $f: W_1 \to W_2$ in $\Embc$ and let $Q \in \cSpar{\Lambda}(W_2)$. Given $(S, \lambda) \in Q$, if $S \cap f(\partial W_1)$ transversally then the intersection has the expected dimension and, thus, $f^{-1}(S)$ is a codimension $2$ submanifold of $W_1$. Furthermore, the co-orientation of $S$ induces a co-orientation on $f^{-1}(S)$ by pullback. Hence, we can define $\cSpar{\Lambda}(f)(Q)$ as the set of pairs $(f^{-1}(S), \lambda)$ for $(S, \lambda) \in Q$ with $S \cap f(\partial W_1)$ transversal. For short, we will denote $Q|_{W_1} = \cSpar{\Lambda}(f)(Q)$. Obviously, if $Q \subseteq Q'$ then $Q|_{W_1} \subseteq Q'|_{W_1}$ so $\cSpar{\Lambda}(f): \cSpar{\Lambda}(W_2) \to \cSpar{\Lambda}(W_1)$ is a functor. With this definition, $\cSpar{\Lambda}$ is a sheaf, called the \emph{sheaf of parabolic structures over $\Lambda$}.

\begin{ex}\label{ex:parab-struct}
\begin{itemize}
	\item Since a $1$-dimensional manifold cannot have codimension $2$ subvarieties, there do not exist non-empty parabolic structures over a curve.
	\item In the case of surfaces, a parabolic structure over a surface $W$ is a set $Q = \left\{(p_1, \lambda_1), \ldots, (p_r, \lambda_r)\right\}$ with $\lambda_i \in \Lambda$ and $p_i \in W$ points with a preferred orientation of a small disc around them.
	\item For a $3$-fold $W$, a parabolic structure is given by $Q = \left\{(S_1, \lambda_1), \ldots, (S_r, \lambda_r)\right\}$, where the $S_i$ are pairwise disjoint paths in $W$ with endpoints in $\partial M$ and reaching the boundary transversally. If we have an embedding $f: W' \to W$ with $W'$ another $3$-fold, the restriction $Q|_{W'}$ is given by the restriction to $W'$ of those paths $S_i$ that are not tangent to the surface $\partial W'$. In Figure \ref{img:ex-parabolic-structure-3-fold} it is depicted and example of this phenomenon where $W_1$ is the shadowed open set with $\partial W_1 = \Sigma_2$, $W_2 = S^3$ and $S \subseteq S^3$ is a knot. The induced parabolic structure is the segment $f^{-1}(S) \subseteq W_1$.
	\vspace{-0.1cm}
\begin{figure}[h]
	\begin{center}
	\includegraphics[scale=0.23]{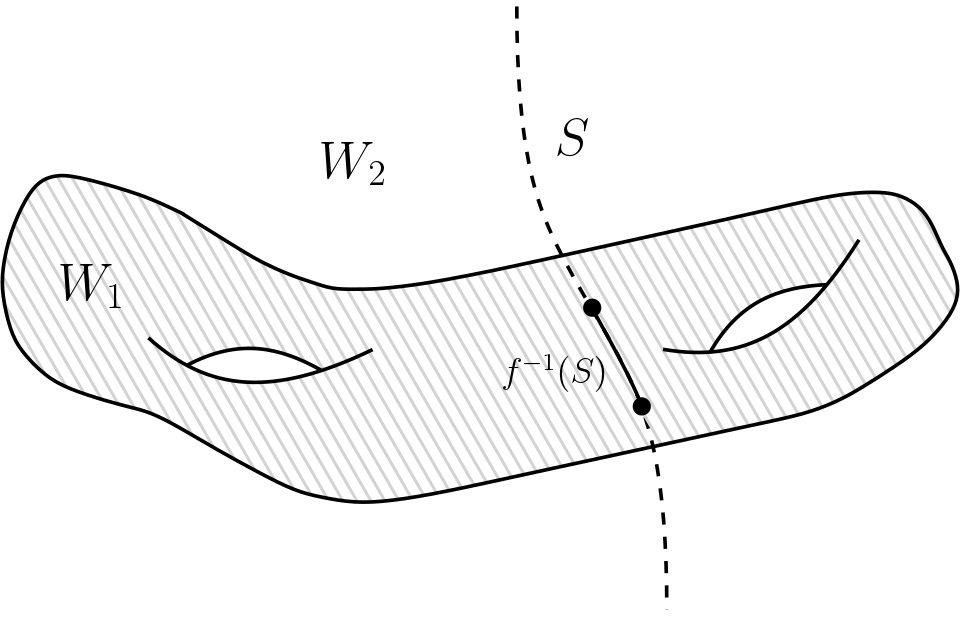}
	\vspace{-0.5cm}
	\caption{Example of induced parabolic structure.}
	\label{img:ex-parabolic-structure-3-fold}
	\end{center}
\vspace{-1cm}
\end{figure}
\end{itemize}
\end{ex}

As for the sheaf of unordered points, the sheaf $\cSpar{\Lambda}$ gives us categories $\EEmbc{\cSpar{\Lambda}}$ and $\EBord{n}{\cSpar{\Lambda}}$, that we will shorten $\Embparc{\Lambda}$ and $\Bordpar{n}{\Lambda}$. Even more, we can combine the two previous sheaves and to consider the sheaf $\cS_{p, \Lambda} = (\cSpar{\Lambda} \times \cSp) \circ \Delta: \Embc \to \Cat$, where $\Delta: \Embc \to \Embc \times \Embc$ is the diagonal functor. In that case, we will denote by $\Embpparc{\Lambda}$ and by $\Bordppar{n}{\Lambda}$ the categories of $\cS_{p, \Lambda}$-embeddings and $\cS_{p, \Lambda}$-bordisms, respectively.

For the convenience of Chapter \ref{chap:representation}, let us write down the definition of $\Bordppar{n}{\Lambda}$. The objects of $\Bordppar{n}{\Lambda}$ are triples $(M, A, Q)$ with $(M, A) \in \Bordp{n}$ and $Q$ a parabolic structure on $M$. A morphism $(M_1, A_1, Q_1) \to (M_2, A_2, Q_2)$ is a class of triples $(W, A, Q)$ where $W: M_1 \to M_2$ is a bordism, $A \cap M_1 = A_1$, $A \cap M_2 = A_2$ and $Q$ is a parabolic structure on $W$ such that $Q|_{M_1} = Q_1$ and $Q|_{M_2} = Q_2$. Two triples $(W, A, Q)$ and $(W', A', Q')$ are in the same class if there exists a boundary preserving diffeomorphism $F: W \to W'$ such that $F(A)=A'$ and $(S, \lambda) \in Q$ if and only if $(F(S), \lambda) \in Q'$. A composition of morphisms in $\Bordp{n}(\Lambda)$ is given by $(W', A', Q') \circ (W, A, Q) = (W \cup W', A \cup A', Q \cup Q')$. In the same spirit, we have a $2$-morphism $(W, A, Q) \Rightarrow (W', A', Q')$ if there exists a boundary preserving diffeomorphism $F: W \to W'$ such that $F(A) \subseteq A'$ and $Q \subseteq \cSpar{\Lambda}(f)(Q')$.

In this case, a (lax monoidal) $\cS_{p, \Lambda}$-TQFT will be called a (lax monoidal) \emph{parabolic Topological Quantum Field Theory of pairs}.

\subsection{Almost-TQFTs and the category of tubes}
\label{sec:almost-TQFT-tubes}

Recall that Zorro's lemma \ref{lem:zorro} forces a strict monoidal TQFT to land into finitely generated modules. This fact imposes very strong algebraic restrictions that avoid the construction of several TQFTs in a natural way. In order to overcome this problem, a solution is to relax the monoidality condition and to consider instead lax monoidal functors, as done in Section \ref{sec:lax-monoidal-tqft}.
However, as mentioned in the introduction of Section \ref{sec:lax-monoidality}, another approach to this problem is to restrict the source category of the TQFT, dropping out the elbow bordism. This is the method explored in this section.

Consider the subcategory $\Tub{n}$ of $\Bord{n}$ of tubes. Roughly speaking, it is a subcategory with the same objects but with all the elbow-like bordisms removed. To be precise, let us consider the subcategory $\Tubo{n} \subseteq \Bord{n}$ of \emph{strict tubes}. An object $M \in \Bord{n}$ is an object of $\Tubo{n}$ if $M$ is connected or empty and there exists a compact $n$-dimensional manifold $W$ such that $\partial W = M$. Given objects $M_1$ and $M_2$ of $\Tubo{n}$, a morphism $W: M_1 \to M_2$ of $\Bord{n}$ is in $\Tubo{n}$ if $W$ is connected. We will call such morphisms strict tubes. From this category, we consider the \emph{category of tubes}, $\Tub{n}$, as the subcategory of $\Bord{n}$ whose objects and morphisms are disjoint unions of the ones of $\Tubo{n}$. That is, $\Tub{n}$ is the monoidal closure of $\Tubo{n}$ in $\Bord{n}$. Observe that, in particular, $\Tub{n}$ is a wide subcategory of $\Bord{n}$. As in the case of $\Bord{n}$, $\Tub{n}$ is a monoidal category with the disjoint union.

\begin{defn}
Let $R$ be a ring. An almost-TQFT is a symmetric monoidal functor $\lZ: \Tub{n} \to \Mod{R}$.
\end{defn}

\begin{rmk}
Since $\Tub{n}$ does not contain the elbow bordisms, dualizing arguments no longer hold for almost-TQFTs so $\lZ(M)$ is not forced to be finitely generated for $M \in \Tub{n}$. Moreover, for $n=2$ the pair of pants is not a tube, so, in contrast with strict monoidal TQFTs, we cannot assure that $\lZ(S^1)$ is a Frobenius algebra, as it happened in Section \ref{ex:low-dim-TQFT}.
\end{rmk}

An almost-TQFT gives an effective way of computing invariants as follows. Fix $n \geq 1$ and suppose that we have a set of generators $\Delta$ for the morphisms of $\Tubo{n}$, i.e.\ up to boundary preserving diffeomorphisms, every morphism of $\Tubo{n}$ is a compositions of elements of $\Delta$. For example, for $n=2$, these generators can be obtained by means of Morse theory (see Section 1.4 of \cite{Kock:2004}).

Suppose that we want to compute an invariant that, for a closed connected orientable $n$-dimensional manifold $W$, is given by $\lZ(W)(1)$. In that case, seeing $W$ as a morphism $W: \emptyset \to \emptyset$, we can decompose $W = W_{s} \circ \ldots \circ W_{1}$ with $W_{i} \in \Delta$. Thus, for $\lZ(W): R \to R$ we have
$$
	\lZ(W)(1) = \lZ(W_{s}) \circ \ldots \circ \lZ(W_{1})(1)
$$
Hence, the knowledge of $\lZ(W_i)$ for $W_i \in \Delta$ is enough to compute that invariant for all the closed manifolds.

Given a lax monoidal TQFT, $Z: \Bord{n} \to \Mod{R}$, we can build an almost-TQFT from it, $\lZ = \lZ_{Z}: \Tub{n} \to \Mod{R}$. For that, notice that the restriction of $Z$ to $\Tubo{n}$ gives us a functor $Z: \Tubo{n} \to \Mod{R}$. Now, we define $\lZ$ by:
\begin{itemize}
	\item Given objects $M^i$ of $\Tubo{n}$, we take
	$$
		\lZ\left(\bigsqcup_i \,M^i\right) = {\bigotimes_i}\, Z(M^i)
	$$
where the tensor product is taken over $R$.
	\item Given strict tubes $W^i: M_{1}^i \to M_2^i$, we define
	$$
		\lZ\left(\bigsqcup_i \,W^i\right) = {\bigotimes_i}\, Z(W^j): {\bigotimes_i}\, Z(M_1^i) \to {\bigotimes_i}\, Z(M_2^i)
	$$
\end{itemize}

\begin{rmk}
The apparently artificial definition of $\lZ$ can be better understood in terms of the corresponding map of $Z$. To see it, let $W: M_1 \to M_2$ and $W': M_1' \to M_2'$ be two tubes. Recall that, since $Z$ is lax monoidal, we have a natural trasformation $\Delta: Z(-) \otimes_R Z(-) \Rightarrow Z(- \sqcup -)$ (see Definition \ref{defn:lax-monoidal}). Then, $\lZ\left(W \sqcup W'\right)$ is just a lift
	\[
\begin{displaystyle}
   \xymatrix
   {
      Z(M_1) \otimes_R Z(M_1') \ar[d]_{\Delta_{M_1, M_1'}} \ar@{--{>}}[rrr]^{\lZ\left(W \sqcup W'\right)} &&& Z(M_2) \otimes_R Z(M_2') \ar[d]^{{\Delta_{M_2, M_2'}}} \\
   Z(M_1 \sqcup M_1') \ar[rrr]_{Z\left(W \sqcup W'\right)} &&& Z(M_2 \sqcup M_2')
   }
\end{displaystyle}   
\]
Since composition of tubes is performed by componentwise gluing, these lifts behave well with composition, which implies that $\lZ$ is well defined. This property no longer holds for $Z$ since there exists bordisms mixing several components.
\end{rmk}

\begin{figure}[h]
	\begin{center}
	\includegraphics[scale=0.35]{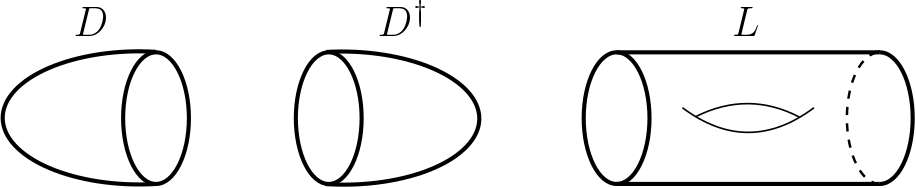}
	\caption{}
	\label{img:generator-tubes}
	\end{center}
\vspace{-0.7cm}
\end{figure}

In the case $n=2$, we can go a step further and describe explicitly the functor $\lZ$. Regarding the objects of $\Tubo{2}$, recall that every non-empty $1$-dimensional closed connected manifold is diffeomorphic to $S^1$. Regarding morphisms, observe that, adapting the proof of Proposition 1.4.13 of \cite{Kock:2004} for the generators of $\Bord{2}$, we have that a set of generators of $\Tub{2}$ is $\Delta = \left\{D, D^\dag, L\right\}$, where $D: \emptyset \to S^1$ is the disc, $D^\dag: S^1 \to \emptyset$ is the opposite disc and $L: S^1 \to S^1$ is the torus with two holes. They are depicted in Figure \ref{img:generator-tubes}.

Suppose that $\lZ$ comes from a lax monoidal TQFT constructed by means of a geometrisation $F: \Bord{n} \to \cC$ and a $\cC$-algebra $\cA$. If $i: S^1 \hookrightarrow D$ denotes the inclusion of $S^1$ as a boundary in $D$ and $j_1, j_2: S^1 \hookrightarrow L$ denote the inclusions of the two boundaries of $T$, the image under the field theory are the spans
$$
   \cF_F(D): \star \longleftarrow F(D) \stackrel{F(i)}{\longrightarrow} F(S^1), \hspace{1cm} \cF_F(D^\dag): F(S^1) \stackrel{F(i)}{\longleftarrow} F(D) \longrightarrow \star,
$$
$$
	\cF_F(L): F(S^1) \stackrel{F(j_1)}{\longleftarrow} F(L) \stackrel{F(j_2)}{\longrightarrow} F(S^1).
$$
Therefore, the image under the almost-TQFT are the ring homomorphisms
$$
	\lZ(D) = i_! \circ c^*: \cA_\star \to \cA_{F(S^1)}, \hspace{1cm} \lZ(D^\dag) =  c_! \circ i^*: \cA_{F(S^1)} \to \cA_\star,
$$
$$
	\lZ(L) = (j_2)_! \circ (j_1)^*: \cA_{F(S^1)} \to \cA_{F(S^1)},
$$
where we have shorten $i^* = F(i)^*$, $i_! = F(i)_!$ and so on, and $c = c_{F(D)}: F(D) \to \star$ is the projection of $F(D)$ onto the final object.

Analogous considerations can be done if we consider a sheaf $\cS$ given some extra data. In that case, the the preimage of $\Tub{n} \subseteq \Bord{n}$ under the forgetful functor $\EBord{n}{\cS} \to \Bord{n}$ will be called the category of $\cS$-tubes and will be denoted by $\ETub{n}{\cS} \subseteq \EBord{n}{\cS}$. As before, we will short $\Tubp{n} = \ETub{n}{\cSp}$ and $\Tubppar{n}{\Lambda} = \ETub{n}{\cS_{p,\Lambda}}$.

\subsection{Reduction of a TQFT}
\label{sec:reduction-TQFT}

Let $\cC$ be a category with final object and pullbacks, let $F: \Embc \to \cC$ be a contravariant functor with the Seifert-van Kampen property and let $\cA$ be a $\cC$-algebra. By Theorem \ref{thm:physical-constr}, these data give rise to a lax monoidal TQFT $Z = Z_{F, \cA}: \Bord{n} \to \Mod{\cA_\star}$ and, as explained in Section \ref{sec:almost-TQFT-tubes}, theoretically $Z$ can be used for computing algebraic invariants. However, for $M \in \Bord{n}$, the module $Z(M) = \cA_{F(M)}$ may be very complicated (e.g.\ it might be infinitely generated). Hence, from a computational point of view, it is not feasible to explicitly compute the homomorphisms associated to bordisms.

This problem can be mitigated if the objects $F(M)$ have some kind of symmetry that may be exploited. For example, suppose that there is a natural action of a group $G$ on $F(M)$ so that the categorical quotient object $F'(M) = F(M)/G \in \cC$ can be defined (e.g.\ if $\cC = \Var{k}$ and $G$ is an algebraic reductive group, see Chapter \ref{chap:git}). In the general case, such a symmetry can be modelled by an assignment $\tau$ that, for any $M \in \Bord{n}$, gives a pair $(F'(M), \tau_{M})$, where $F'(M) \in \cC$ and $\tau_{M}: F(M) \to F'(M)$ is a morphism in $\cC$. Suppose also that $F'(\emptyset) = F(\emptyset)$.

In that case, we can `reduce' the field theory and to consider $\Fld{F, \tau}: \Bord{n} \to \Span{\cC}$ such that $\Fld{F, \tau}(M)=F'(M)$, for $M \in \Bord{n}$, and, for a bordism $W: M_1 \to M_2$, $\Fld{F, \tau}(W)$ is the span
\[
\begin{displaystyle}
   \xymatrix
   {	
	F'(M_1) && F(W) \ar[rr]^{\tau_{M_2} \circ F(i_2)} \ar[ll]_{\tau_{M_1} \circ F(i_1)} && F'(M_2)
   }
\end{displaystyle}.
\]
With this field theory, we form $Z_\tau = \Qtm{\cA} \circ \Fld{F}'$, called the \emph{prereduction} of $Z$ by $\tau$.

However, even if $F$ had the Seifert-van Kampen property, $\Fld{F, \tau}$ may not be a functor so $Z_\tau$ will be, in general, a very lax TQFT. Indeed, the prereduction $Z_\tau$ is simpler than $Z$ but, as it is not a functor, we can no longer use it for a computational method. In this section we will show that, under some mild conditions, we can slightly modify $Z_\tau$ in order to obtain a (strict) almost-TQFT, $\cZ_\tau$, called the \emph{reduction}, with essentially the same complexity as $Z_\tau$. In this way, $\cZ_\tau$ can be used to give an effective computational method.

In order to do so, for any $M \in \Tubo{n}$, let us denote by $\cV_{M}$ be the submodule of $\cA_{F'(M)}$ generated by the images $Z_\tau(W_r) \circ \ldots \circ Z_\tau(W_1)(1) \in \cA_{F(M)}$ of all the sequences of strict tubes $W_k: M_{k-1} \to M_{k}$ with $M_0 = \emptyset$. Notice that, by definition, the submodules $\cV_M$ are invariant for strict tubes, that is, $Z_\tau(W)(\cV_M) \subseteq \cV_N$ for any strict tube $W: M \to N$.

\begin{prop}\label{prop:existence-descend-reduction}
Suppose that $(\tau_M)_! \circ (\tau_M)^*(\cV_M) \subseteq \cV_M$ for all $M$ and that the extended morphisms $\eta_M = (\tau_M)_! \circ (\tau_M)^*: \cV_M \to \cV_M$ are invertible. Then, for every strict tube $W: M \to N$, there exists an unique morphism $\cZ_\tau(W): \cV_M \to \cV_N$ such that the diagram commute
\[
\begin{displaystyle}
   \xymatrix
   {	
		\tau_M^*\cV_M \ar[r]^{Z(W)\;\;\;\;\;\;\;\;\;}\ar[d]_{(\tau_M)_!}  & Z(W)\left(\tau_M^*{\cV_M}\right) \ar[d]^{(\tau_N)_!}\\
		\cV_M \ar@{--{>}}[r]_{\cZ_\tau(W)} & \cV_N
   }
\end{displaystyle}   
\]
\begin{proof}
First of all, observe that the diagram above makes sense since $(\tau_N)_! \circ Z(W) \circ (\tau_M)^*(\cV_M) = Z_\tau(W)(\cV_M) \subseteq \cV_N$. Suppose that we have a morphism $f: \cV_M \to \cV_N$ such that $f \circ (\tau_M)_! = (\tau_N)_! \circ Z(W)$. Then, pre-composing with $\tau_M^* \circ \eta_M^{-1}: \cV_M \to \tau_M^*\cV_M$ we obtain that
$$
	(\tau_N)_! \circ Z(W) \circ \tau_M^* \circ \eta_M^{-1} = f \circ (\tau_M)_! \circ \tau_M^* \circ \eta_M^{-1} = f \circ \eta_M \circ \eta_M^{-1} = f.
$$
Hence, if such a map exists, it is unique. Actually, this calculation shows that we must take $\cZ_\tau(W) = (\tau_N)_! \circ Z(W) \circ \tau_M^* \circ \eta_M^{-1} = Z_\tau(W) \circ \eta_M^{-1}$ and it is a straighforward check that $\cZ_\tau(W)$ has the desired property.
\end{proof}
\end{prop}

\begin{cor}\label{cor:existence-reduction-TQFT}
In the situation of Proposition \ref{prop:existence-descend-reduction}, there exists an (strict) almost-TQFT
$$
	\cZ_\tau: \Tub{n} \to \Mod{\cA_\star}
$$
such that $\cZ_\tau(M) = \cV_M$ for all $M \in \Tub{n}$ connected and $\cZ_\tau(W) (\tau_M)_! = (\tau_N)_! \circ Z(W)$ for all strict tubes $W: M \to N$. This TQFT is called the \emph{reduction} of $Z$ via $\tau$.
\begin{proof}
Recall that, in order to define an almost TQFT, it is enough to define it on strict tubes. Hence, for $M \in \Tubo{n}$, we assign $\cZ_\tau(M) = \cV_M$ and, for a strict tube $W: M \to N$, we assign the morphism $\cZ_\tau(W): \cV_M \to \cV_N$ previously constructed.
In order to prove that $\cZ_\tau$ is a functor, suppose that $W: M \to N$ and $W': N \to Z$ are strict tubes. Then, by the previous proposition, we have a commutative diagram
\[
\begin{displaystyle}
   \xymatrix
   {	
		\tau_M^*\cV_M \ar[r]^{Z(W)\;\;\;\;\;\;\;\;\;}\ar[d]_{(\tau_M)_!}  & Z(W)\left(\tau_M^*\cV_M\right) \ar[d]^{(\tau_N)_!} \ar[r]^{Z(W')\;\;\;\;\;}& Z(W' \circ W)\left(\tau_M^*\cV_M\right) \ar[d]^{(\tau_Z)_!}\\
		\cV_M \ar[r]_{\cZ_\tau(W)} & \cV_N \ar[r]_{\cZ_\tau(W')} & \cV_Z
   }
\end{displaystyle}   
\]
Therefore, we have $\cZ_\tau(W') \circ \cZ_\tau(W) \circ (\tau_M)_! = (\tau_Z)_! \circ Z(W') \circ Z(W) = (\tau_Z)_! \circ Z(W' \circ W)$ and, by uniqueness, this implies that $\cZ_\tau(W') \circ \cZ_\tau(W) = \cZ_\tau(W' \circ W)$.
\end{proof}
\end{cor}

\begin{rmk}\label{rmk:reduction-TQFT}
\begin{itemize}
	\item The almost-TQFT, $\cZ_\tau: \Tub{n} \to \Mod{\cA_\star}$, satisfies that, for all closed $n$-dimensional manifolds $W$
$$
	\cZ_\tau(W)(1) = \cZ_\tau(W) \circ (\tau_\emptyset)_!(1) = (\tau_\emptyset)_! \circ Z(W)(1) = Z(W)(1),
$$
where we have used that $\tau_\emptyset = 1_\star$. Hence, $\cZ_\tau$ and $Z$ compute the same invariant.
	\item It may happen, and it will be the case for representation varieties, that $\eta_M: \cV_M \to \cV_M$ are not invertible as $\cA_\star$-module. However, it could happen that, for some fixed multiplicative system $S \subseteq \cA_\star$, all the extensions to the localizations $\eta_M: S^{-1}\cV_M \to S^{-1}\cV_M$ are invertible. In that case, we can fix the problem by localizing all the modules and morphisms of the original TQFT to obtain another TQFT, $Z: \Bord{n} \to \Mod{S^{-1}\cA_\star}$, to which we can apply the previous construction. 
	\item If $Z(W)(\tau_M^*\cV_M) \subseteq \cV_N$ for all the strict tubes $W: M \to N$, then the original almost-TQFT, $\lZ: \Tub{n} \to \Mod{\cA_\star}$, can be restricted to assign $\lZ(M)= \tau_M^*\cV_M$. In that case, the reduction $\tau$ defines a natural transformation between the restricted almost-TQFT and the reduced one, $\lZ \Rightarrow \cZ_\tau$.
	\item Instead of the submodule $\cV_M$ generated by the elements $Z_\tau(W_r) \circ \ldots \circ Z_\tau(W_1)(1) \in \cA_{F(M)}$ with $W_r \circ \ldots \circ W_1: \emptyset \to M$ strict tubes, we can consider the submodule $\cV_M'$ generated by the elements $Z_\tau(W_r) \circ \ldots \circ Z_\tau(W_1)(1) \in \cA_{F(M)}$ with $W_r \circ \ldots \circ W_1: \emptyset \to M$ any bordism. In that case, if $\cV_M'$ is invariant under $\eta_M = (\tau_M)_! \circ (\tau_M)^*$ and $\eta_M$ is invertible there, the same procedure as above will allow us to constuct a genuine TQFT, $\cZ_\tau: \Bord{n} \to \Mod{\cA_\star}$, with $\cZ_\tau(M) = \cV_M'$. However, in general this condition is much harder to check, so we will focus on the version presented above, which is enough for computations.
	\item With the appropiate modifications, the same results of this section hold for TQFTs over a sheaf.
\end{itemize}
\end{rmk}

\begin{ex}\label{ex:reduction-TQFT-surfaces}
For $n=2$ (surfaces), the unique non-empty connected object of $\Bord{2}$ is $S^1$. Hence, it is enough to consider $\cV = \cV_{S^1}$ and $\tau = \tau_{S^1}: F(S^1) \to F'(S^1)$. Actually, $\cV$ is the submodule generated by the elements $Z_\tau(L)^g \circ Z_\tau(D)(1)$ for $g \geq 0$, where $L: S^1 \to S^1$ is the holed torus and $D: \emptyset \to S^1$ is the disk (i.e.\ by the image of all the compact surfaces with a disc removed). In that case, the only condition we need to check is that $\eta = \tau_! \circ \tau^*: \cV \to \cV$ is invertible.
\end{ex}

\begin{rmk}\label{rmk:left-reduction}
Instead of the diagram of Proposition \ref{prop:existence-descend-reduction}, we could look for a morphism $\cZ_\tau'(W): \cV_M \to \cV_N$ such the diagram
	\[
\begin{displaystyle}
   \xymatrix
   {	
		\tau_M^*\cV_M \ar[r]^{Z(W)\;\;\;\;\;\;\;\;\;}  & Z(W)\left(\tau_M^*{\cV_M}\right) \\
		\cV_M \ar[u]^{\tau_M^*} \ar@{--{>}}[r]_{\cZ_\tau'(W)} & \cV_N \ar[u]_{\tau_N^*}
   }
\end{displaystyle}   
\]
commutes. In that case, in the same conditions as in Proposition \ref{prop:existence-descend-reduction}, $\cZ_\tau'$ exists, it is an almost TQFT and $\cZ_\tau'(W)=\eta^{-1}_N \circ Z_\tau(W)$. We will call $\cZ_\tau'$ the \emph{left $\tau$-reduction}.
\end{rmk}

\subsection{Soft TQFTs}
\label{sec:soft-TQFT}

In this section, we will study a lighter version of Topological Quantum Field Theories called soft TQFTs. As we will see, these TQFTs can be built with less data than genuine TQFTs so they are available in more general circumstances. The construction is similar to the $(1,2)$-parts of extended TQFTs, as studied in \cite{Lauda-Pfeiffer}, \cite{Schommer-Pries:2009} or \cite{Ben-Zvi-Francis-Nadler:2010}. It is a prospective work to establish the link between the two approaches and to try to improve these soft TQFTs to extended TQFTs. 

Recall from Example \ref{ex:categories} that, given a ground commutative ring $R$ with identity, we can define the $2$-category $\Bim{R}$ of $R$-algebras and bimodules whose objects are commutative $R$-algebras, $1$-morphisms are bimodules and $2$-morphisms are bimodule homomorphisms. The category $\Bim{R}$ is a monoidal $2$-category with tensor product over $R$.
Observe that, since $R$ is commutative, a $(A,B)$-bimodule is the same as a $(B,A)$-bimodule. Hence, we have an involutive functor $(-)^\dag: \Bim{R} \to \Bim{R}$ which is the identity on objects and, for $M: A \to B$ a $(A,B)$-bimodule, it gives $M^\dag: B \to A$, seen as a $(B,A)$-bimodule (in general, a category with a functor $(-)^\dagger$ as above is called a dagger category).

\begin{defn}\label{defn:soft-TQFT}
Let $R$ be a commutative ring with unit. A \emph{soft Topological Quantum Field Theory} is a lax functor $\mathscr{Z}: \Bord{n} \to \Bim{R}$ such that:
\begin{itemize}
	\item There is an isomorphism $\alpha: R \to F(\emptyset)$.
	\item For any $M, N \in \Bord{n}$, there exists a $1$-morphism $\Delta_{M,N}: Z(M) \otimes_R Z(N) \to Z(M \sqcup N)$. These morphisms satisfy that, for any morphisms $W: M \to N$ and $W': M' \to N'$, there is a $2$-morphism $\tau_{W,W'}: \Delta_{N, N'} \circ (Z(W) \otimes_R Z(W')) \circ \Delta_{M,M'}^\dag \Rightarrow Z(W \sqcup W')$.
		\[
\begin{displaystyle}
   \xymatrix
   {
     Z(M) \otimes_R Z(M') \ar[dd]_{Z(W) \otimes_R Z(W')} & Z(M \sqcup M') \ar[dd]^{Z(W \sqcup W')} \ar[l]_{\;\;\;\;\;\;\;\;\Delta_{M,M'}^\dag}\\
     \;\;\;\;\;\; \ar@{={>}}[r]^{\tau_{W,W'}} &  \;\;\;\;\;\;\\
     Z(N) \otimes_R Z(N') \ar[r]_{\;\;\;\;\;\;\;\Delta_{N,N'}} & Z(N \sqcup N')
   }
\end{displaystyle}   
\]
\end{itemize}
\end{defn}

A general recipe for building a soft TQFT from simpler data can be given analogously as for genuine TQFT decomposing the functor as $
	\mathscr{Z}: \Bord{n} \to \Span{\cC} \to \Bim{R}$.
As above, the field theory $\cF: \Bord{n} \to \Span{\cC}$ can be determined from a contravariant functor $F: \Embc \to \cC$ with the Seifert-van Kampen property. However, for the `quantum' part $\Span{\cC} \to \Bim{R}$, we do not need a full $\cC$-algebra but only half of it.
To be precise, let $A: \cC \to \Rng$ be a contravariant functor where $\cC$ is a category with pullbacks and final object $\star$. In analogy with the notation of Section \ref{sec:C-algebras} (in particular, Remark \ref{rmk:notation-grothendieck-six}), we denote $A(a) = A_a$ and $A(f)=f^*$ for $a \in \cC$ and a morphism $f$ in $\cC$. We define the pre-quantisation, $\sQtm{A}: \cC \to \Bim{A_\star}$, as follows:
\begin{itemize}
	\item For any $a \in \cC$, we set $\sQtm{A}(a) = A_a$. The $A_\star$-algebra structure on $A_a$ is given by the morphism $A_\star \to A_a$ image of the final map $a \to \star$ (c.f.\ Section \ref{sec:C-algebras}).
	\item Let us fix $a,b \in \cC$. We define 
	$$
	\sQtm{A}: \Hom_{\Span{\cC}}(a,b) \to \Hom_{\Bim{A_\star}}(A_a, A_b)
	$$
as the functor that:
		\begin{itemize}
			\item For any $1$-morphism $a \stackrel{f}{\leftarrow} d \stackrel{g}{\rightarrow} b$ it assigns $\sQtm{A}(d,f,g) = A_d$ with the $(A_a,  A_b)$-bimodule structure given by $\lambda z = f^*(\lambda) \cdot z$ and $z\mu = z \cdot g^*(\mu)$, for $\lambda \in A_a$, $\mu \in A_b$ and $z \in A_d$.
			
			\item For a $2$-morphism $\alpha: (d, f, g) \Rightarrow (d',f',g')$ given by a morphism $\alpha: d' \to d$, we define $\sQtm{A}(\alpha) = \alpha^*: A_d \to A_{d'}$. Since $A$ is a functor to rings, $\alpha^*$ is a bimodule homomorphism.
		\end{itemize}
	\item The $2$-morphism $1_{A_a} \Rightarrow \sQtm{A}(1_a)$ is the identity for all $a \in \cC$.
	\item Let $a \stackrel{f_1}{\leftarrow} d_1 \stackrel{g_1}{\rightarrow} b$ and $b \stackrel{f_2}{\leftarrow} d_2 \stackrel{g_2}{\rightarrow} c$. By definition of $\sQtm{A}$,
\begin{align*}
\sQtm{A}(d_2,f_2,g_2) \circ \sQtm{A}(d_1,f_1, g_1) &= {_{A_b}}[A_{d_2}]_{A_c} \circ {_{A_a}}[A_{d_1}]_{A_b} \\
	&= {_{A_a}}\left[A_{d_1} \otimes_{A_b} A_{d_2}\right]_{A_c}, \\
	& \\[-0.4cm]
\sQtm{A}((d_2,f_2,g_2) \circ (d_1,f_1, g_1)) &= {_{A_a}}\left[A_{d_1 \times_b d_2}\right]_{A_c}.
\end{align*}
Hence, the $2$-morphism $\sQtm{A}(d_2, f_2, g_2) \circ \sQtm{A}(d_1, f_1, g_1)) \Rightarrow \sQtm{A}((d_2, f_2, g_2), (d_1, f_1, g_1))$ is the $(A_a, A_b)$-bimodule homomorphism given by external product
$\boxtimes_b: A_{d_1} \otimes_{A_b} A_{d_2} \to A_{d_1 \times_b d_2}
$.
\end{itemize}

Furthermore, with the construction above, the lax functor $\sQtm{A}: \cC \to \Bim{A_\star}$ is also monoidal in the sense of Definition \ref{defn:soft-TQFT}. Since $\sQtm{A}$ preserves the units of the monoidal structures it is enough to define the morphisms $\Delta_{a,b}: A_a \otimes_{A_\star} A_b \to A_{a \times b}$. For that, just take as bimodule the own ring $A_{a \times b}$, where the left $(A_a \otimes_{A_\star} A_b)$-module structure comes from the external product $\boxtimes: A_a \otimes_{A_\star} A_b \to A_{a \times b}$. In that case, if $(d, f, g): a \to b$ and $(d', f', g'): a' \to b'$ are two spans, we have a commutative diagram
\[
\begin{displaystyle}
   \xymatrix
   {
     A_a \otimes A_{a'} \ar[dd]_{A_{d} \otimes A_{d'}} & A_{a \times a'} \ar[dd]^{A_{d \times d'}} \ar[l]_{\;\;\;\;\;\;A_{a \times a'}}\\
     \;\;\;\;\;\; \ar@{={>}}[r]^{\tau} & \;\;\;\;\;\;\\
     A_{b} \otimes A_{b'} \ar[r]_{\;\;\;\;\;A_{b \times b'}} & A_{b \times b'}
   }
\end{displaystyle}   
\]
where $\tau$ is the bimodule homomorphism $\tau: A_{a' \times b'} \otimes_{(A_{a'} \otimes A_{b'})} [A_{d} \otimes A_{d'}] \otimes_{(A_a \otimes A_b)} A_{a \times b} \to A_{d \times d'}$ given by $\tau(x,z,w,y) = A(g \times g')(y) \cdot [z \boxtimes w] \cdot A(f \times f')(x)$ for $x \in A_{a \times a'}, y \in A_{b \times b'}$, $z \in A_d$ and $w \in A_{d'}$.

\begin{rmk}
If the functor $A: \cC \to \Rng$ is monoidal ($\cC$ with the cartesian monoidal structure), then $\sQtm{A}$ is strict monoidal. Moreover, if it satisfies that $A_{d_1} \otimes_{A_b} A_{d_2} = A_{d_1 \times_b d_2}$ (resp.\ isomorphic) then it is also a strict (resp.\ pseudo) functor.
\end{rmk}

Therefore, with this construction at hand, given a contravariant functor $F: \Embc \to \cC$ with the Seifert-van Kampen property and a contravariant functor $A: \cC \to \Rng$, we can build a soft TQFT, $\mathscr{Z}_{F,A}$, by
$$
	\mathscr{Z}_{F,A} = \sQtm{A} \circ \cF_F: \Bord{n} \to \Bim{R}.
$$
Analogously to Section \ref{sec:TQFT-over-sheaf}, given a sheaf $\cS$, a lax monoidal (in the sense of Definition \ref{defn:soft-TQFT}) lax functor $\mathscr{Z}: \EBord{n}{\cS} \to \Bim{R}$ will be called a \emph{$\cS$-soft TQFT}. Since the construction above does not depend on the field theory, combined with the methods of Section \ref{sec:field-theory}, we can obtain $\cS$-soft TQFT from the geometrisation $F: \EEmbc{\cS}$ and $A$. We will use this approach in Remark \ref{rmk:Deligne-Hodge-soft-TQFT} in order to define a soft TQFT computing Hodge structures of representation varieties.

%% file: Chapters/Hodge.tex
% Chapter 1

\chapter{Hodge Theory} % Main chapter title

\label{chap:hodge} % For referencing the chapter elsewhere, use \ref{Chapter1} 

\lhead{Chapter 2. \emph{Hodge Theory}} % This is for the header on each page - perhaps a shortened title

%----------------------------------------------------------------------------------------

\section{A panoramic view of $K$-theory}
\label{sec:panoramic-Ktheory}
Let us consider a category $\cA$ with a zero object (i.e.\ an object that is both initial and final). Given a morphism $f: A \to B$ in $\cA$, the \emph{kernel} of $f$, denoted $i: \ker{f} \to A$, is the pullback of the diagram
\[
\begin{displaystyle}
   \xymatrix
   {
	\ker{f} \ar@{--{>}}[r] \ar@{--{>}}[d]_{i}& 0 \ar[d] \\
	A \ar[r] & B
   }
\end{displaystyle}
\]
	Analogously, a \emph{cokernel}, denoted $p: B \to \coker{f}$, is the pushout of the diagram
\[
\begin{displaystyle}
   \xymatrix
   {
	A \ar[r]^f \ar[d] & B \ar@{--{>}}[d]^{p} \\
	0 \ar@{--{>}}[r] & \coker{f}
   }
\end{displaystyle}
\]
A morphism $f$ with $\ker{f} = 0$ is called a \emph{monomorphism} and a morphism with $\coker{f} = 0$ is called an \emph{epimorphism}.

\begin{rmk}
Even if $\cA$ is a subcategory of $\Sets$, not every monomorphism has to be injective. However, every injective map is a monomorphism.
\end{rmk}

With these concepts at hand, we can define a very special kind of categories that behave like the usual `algebraic categories' of modules, called abelian categories. For a complete introduction to abelian categories, see Chapter VIII of \cite{MacLane} or Appendix A.4 of \cite{Weibel}.

\begin{defn}
A category $\cA$ with zero object is an \emph{abelian category} if:
\begin{itemize}
	\item[$i)$] For any $A,B \in \cA$, $\Hom_\cC(A,B)$ has an abelian group structure (i.e.\ $\cA$ is an $\Ab$-enriched category, see Definition \ref{defn:enriched}) and the composition is bilinear.
	\item[$ii)$] $\cA$ has finite coproducts.
	\item[$iii)$] Every morphism of $\cA$ has a kernel and a cokernel.
	\item[$iv)$] For every morphism $f: A \to B$ of $\cA$, there exists a decomposition $f = i \circ \bar{f} \circ p$, where
	$$
		A \stackrel{p}{\longrightarrow} \coker{(\ker{f})} \stackrel{\bar{f}}{\longrightarrow} \ker{(\coker{f})} \stackrel{i}{\longrightarrow} B,
	$$
with $\bar{f}$ an isomorphism. 
\end{itemize}
\end{defn}

\begin{rmk}
\begin{itemize}
	\item The prototypical example of an abelian category is the category of abelian groups, $\Ab$, or, more generally, of $R$-modules, $\Mod{R}$. In this spirit, the category of sheaves over a topological space $X$, $\Sh{X}$, is abelian. However, the full subcategory of vector bundles on $X$, $\textbf{Vb}(X)$, (i.e.\ locally free sheaves over the sheaf of continuous functions on $X$) is no longer abelian, since the kernel of a vector bundle morphism may be no longer a vector bundle. 
	\item As proven in \cite{MacLane}, condition $iv)$ is equivalent to the fact that every monomorphism is a kernel and every epimorphism is a cokernel.
	\item In an abelian category, finite products and coproducts agree. It is customary to denote the product of $\left\{A_i\right\}_{i \in \Lambda}$ just by $\bigoplus_{i \in \Lambda} A_i$ and this operation is called the \emph{direct sum}.
	\item A functor $F: \cA \to \cB$ between abelian categories is said to be \emph{additive} if it preserves finite products. Hence, if $F$ is additive, $F(0_\cA) \cong 0_\cB$ and $F(A \oplus B) \cong F(A) \oplus F(B)$ for any $A,B \in \cA$.
\end{itemize}
\end{rmk}

Let $\cA$ be an abelian category and let $A \in \cA$. An increasing \emph{filtration} $F^\bullet$ of $A$ is a sequence of subobjects of $A$
$$
	0 \subseteq \ldots \subseteq F^{k-1}A \subseteq F^kA \subseteq F^{k+1}A \subseteq \ldots \subseteq A.
$$
The filtration $F^\bullet$ is called finite if $F^kA = 0$ for $k$ small enough and $F^k A = A$ for $k$ large enough.
Given two filtered objects, $(A, F^\bullet)$ and $(B, G^\bullet)$, a morphism $f: A \to B$ in $\cA$ is said to be filtered (or to preserve the filtrations) if $f(F^kA) \subseteq f(G^kB)$ for all $k \in \ZZ$.

\begin{rmk}\label{rmk:filtration-subobject}
\begin{itemize}
	\item Analogously, a decreasing filtration $F_\bullet$ of $A$ is a sequence of subobjects with $F_k A \supseteq F_{k+1}A$.
	\item Strictly speaking, in a general abelian category, there is no notion of inclusion and it has to be replaced by the notion of subobject. Given $A \in \cA$ and monomorphisms $i: B \to A$ and $i': B' \to A$, we say that $i \subseteq i'$ if there exists $\psi: B \to B'$ such that $i = i' \circ \psi$. Two such monomorphisms are declared to be equivalent if $i \subseteq i'$ and $i' \subseteq i$ (observe that being equivalent does not imply that $B$ and $B'$ are isomorphic). In this way, the subobjects of $A$ are nothing but these equivalence classes, which are a partially ordered class with $\subseteq$, and a filtration of $A$ is a chain (i.e.\ a totally ordered subset) of this partial order.
\end{itemize}
\end{rmk}

Recall that, in an abelian category, a pair of morphisms
$$
	B \stackrel{f}{\longrightarrow} A \stackrel{g}{\longrightarrow} C
$$
is said to be exact in $A$ if $\img{f} = \ker{g}$ (in the sense that they are the same subobject of $A$, see Remark \ref{rmk:filtration-subobject}). If a sequence of morphisms is exact in all its objects, such a sequence is called \emph{exact}. An important class of exact sequences are the so-called \emph{short exact sequences} that have the form
$$
	0 \longrightarrow B \longrightarrow A \longrightarrow C \longrightarrow 0.
$$

A functor $F: \cA \to \cB$ between abelian categories is said to be \emph{left exact} if, given a short exact sequence $0 \to B \stackrel{f}{\to} A \stackrel{g}{\to} C \to 0$, the induced sequence 
$$
	0 \longrightarrow F(B)  \stackrel{F(f)}{\longrightarrow} F(A)  \stackrel{F(g)}{\longrightarrow} F(C)
$$
is exact. By duality, a \emph{right exact} functor is a functor that sends short exact sequences into sequences that are exact everywhere except in their leftmost term. A functor which is both left and right exact is called \emph{exact}.

\begin{rmk}
There is a weaker notion than abelian categories, the so-called \emph{triangulated categories}. Consider a category $\cA$ with an invertible functor $T: \cA \to \cA$, called the translation, which is tipically denoted $TA = A[1]$, for $A \in \cA$. In this setting, we can consider triangles, which are triples of morphisms of the form
$$
	B \longrightarrow A \longrightarrow C \longrightarrow B[1].
$$
In that case, a triangulated category is a pair $(\cA, T)$ with a choose of a special collection of triangles, called the \emph{distinguished triangles}, satisfying a set of coherence axioms (see \cite{Weibel}, Section 10.2). The importance of these triangulated categories is that the derived category (see the following section) is no longer an abelian category, but a triangulated one.
\end{rmk}

\subsection{Review of derived categories}
\label{sec:derived-categories}

Given an abelian category, we can extend it in order to consider the categories of cochain complexes of objects on it. These are, roughly speaking, the derived categories. Such a constructions are extremely important in homological algebra and most of the homological machinery can be restated in the context of derived categories. In order to fix notation and for the sake of completeness, we will review some fundamental notions of derived categories that will be useful for future constructions. For a complete exposition of derived categories, see \cite{Huybrechts:2006} and \cite{Weibel}.

Let $\cA$ be an abelian category. The category $\Ch{\cA}$ is the category whose objects are cochain complexes
$$
	\ldots \stackrel{d}{\longrightarrow} A^{k-1} \stackrel{d}{\longrightarrow} A^k \stackrel{d}{\longrightarrow} A^{k+1} \stackrel{d}{\longrightarrow} \ldots,
$$
with $d^2 = 0$. Such a chain complex is usually shortened $(A^\bullet,d) \in \Ch{\cA}$ or just by $A^\bullet$ when the morphisms $d$ are understood. A morphism of chain complexes $f: A^\bullet \to B^\bullet$ is a set of morphisms of $\cA$, $f^k: A^k \to B^k$, commuting with $d$, which is determined up to homotopy of chain complexes (see \cite{Weibel} for the definition). The most important operation that can be performed in a chain complex $A^\bullet$ is to consider its cohomology
$$
	H^k(A^\bullet) = \frac{\ker{d: A^k \to A^{k+1}}}{\img{d: A^{k-1} \to A^k}}.
$$
Observe that, by construction, a morphism of chain complexex $f: A^\bullet \to B^\bullet$ induces a map on cohomology $H^k(f): H^k(A^\bullet) \to H^k(B^\bullet)$.

\begin{rmk}
Sometimes in the literature (notably \cite{Huybrechts:2006}), $\Ch{\cA}$ stands for the category of cochain complexes before identifying homotopic chain morphisms. In that case, the category obtained after the identification of homotopic morphisms is denoted by $\mathbf{K}{\cA}$. However, we will reserve such notation for the Grothendieck group of $\cA$ (see Section \ref{sec:grothendieck-groups}).
\end{rmk}

From this category, one can build the \emph{derived category} of $\cA$, denoted $\Der{\cA}$. It has the same objects as $\Ch{\cA}$ but we localize the morphism of $\Ch{\cA}$ by quasi-isomorphisms. Recall that a quasi-isomorphism is a chain complex map $\phi: A^\bullet \to B^\bullet$ such that the induced maps in cohomology $H^k(\phi): H^k(A^\bullet) \to H^k(B^\bullet)$ are isomorphism. Hence, a morphism $A^\bullet \to B^\bullet$ in $\Der{\cA}$ is a span
$$
	A^\bullet \stackrel{\,\,\phi}{\longleftarrow} D^\bullet {\longrightarrow}\, B^\bullet,
$$
where $\phi$ is a quasi-isomorphism. Obviously, there is a projection $\Ch{\cA} \to \Der{\cA}$.

The full subcategory of $\Der{\cA}$ of chain complexes $A^\bullet$ such that $A^k = 0$ for $k$ small enough (resp.\ $k$ big enough) will be denoted by $\Dp{\cA}$ (resp.\ $\Dm{\cA}$) and will be called the \emph{bounded below derived category} (resp.\ \emph{bounded above derived category}). The subcategory of both bounded above and below cochain complexes is called the \emph{bounded derived category} and is denoted by $\Db{\cA}$. Analogous considerations hold for $\Chp{\cA}$, $\Chm{\cA}$ and $\Chb{\cA}$.

\begin{rmk}
The derived category $\Der{\cA}$ of an abelian category is no longer an abelian category. However, it has a natural structure of a triangulated category (see \cite{Weibel}, Section 10.2). The translation operator is the shifting operator of the chain complex to the left i.e.\ for $A^\bullet \in \Der{\cA}$ with differential $d$, $A^\bullet[1]$ is the chain complex with $A^\bullet[1]^k = A^{k+1}$ and differential $-d$. In order to define the distinguished triangules, suppose that we have an exact sequence $0 \to B^\bullet \to A^\bullet \to C^\bullet \to 0$ of chain complexes and suppose that, for all $k$, the sequence $0 \to B^k \to A^k \to C^k \to 0$ in $\cA$ splits (i.e. $A^k \cong B^k \oplus C^k$). In that case, a chain morphism $h: C^\bullet \to B^\bullet[1]$ can be built from the sections of the splittings, called the homotopy invariant. Then, we declare that $B^\bullet \to A^\bullet \to C^\bullet \stackrel{h}{\to} B^\bullet[1]$ is a distinguished triangle.
\end{rmk}

Derived categories are the natural framework in which derived functors can be formulated. Suppose that our abelian category $\cA$ has enough injectives, that is, for any $A \in \cA$ there is an injective object $I \in \cA$ and a monomorphism $A \to I$ (see \cite{Weibel}, Section 2.3, for the definition of an injective object). In that case, given a cochain complex $A^\bullet \in \Dp{\cA}$, there exists a quasi-isomorphism of chain complexes $A^\bullet \to I^\bullet$, where $I^\bullet \in \Dp{\cA}$ is a complex of injective objects (see \cite{Huybrechts:2006}, Proposition 2.35). Such a quasi-isomorphism is called an \emph{injective resolution}. Injective resolutions of a chain complex are unique up to homotopy equivalence. Even more, if we have a chain morphism $f: A^\bullet \to B^\bullet$ and $A^\bullet \to I^\bullet$, $B^\bullet \to J^\bullet$ are injective resolutions, there exists a chain morphism $\tilde{f}: I^\bullet \to J^\bullet$, unique up to homotopy, such that the following diagram commute
\[
\begin{displaystyle}
   \xymatrix
   {
	A^\bullet \ar[r]^f \ar[d] & B^\bullet \ar[d]\\
	I^\bullet \ar[r]_{\tilde{f}} & J^\bullet
   }
\end{displaystyle}
\]

In that context, suppose that we have a left exact functor $F:\cA \to \cB$ between abelian categories. In that case, we define the \emph{derived functor} of $F$, $RF: \Dp{\cA} \to \Dp{\cB}$, as follows.
\begin{itemize}
	\item For $A^\bullet \in \Dp{\cA}$, it assigns $RF(A^\bullet)= F(I^\bullet)$, where $A^\bullet \to I^\bullet$ is any injective resolution.
	\item Given a morphism of chain complexes $f: A^\bullet \to B^\bullet$, it assigns $RF(f) = F(\tilde{f}): F(I^\bullet) \to F(J^\bullet)$ where $A^\bullet \to I^\bullet$ and $B^\bullet \to J^\bullet$ are injective resolutions.
\end{itemize}
Finally, for $k \in \ZZ$, the \emph{k-th derived functor} is the $k$-th cohomology of $RF$, that is, the functor $R^kF = H^k(RF): \cA \to \cB$. A functor such that $RF(\Db{\cA})\subseteq \Db{\cB}$ is said to have a \emph{bounded derived functor}.

The importance of the derived functor is that, given a short exact sequence $0 \to B \to A \to C \to 0$ in $\cA$, the $k$-th derived functors fit in a long exact sequence
$$
\xymatrix{
    0 \ar[r] & F(A) \ar[r] & F(B) \ar[r] & F(C)
                \ar@{->} `r/8pt[d] `/10pt[l] `^dl[ll]|{} `^r/3pt[dll] [dll] \\
             & R^1F(A) \ar[r] & R^1F(B) \ar[r] & R^1F(C)
                \ar@{->} `r/8pt[d] `/10pt[l] `^dl[ll]|{} `^r/3pt[dll] [dll] \\
             & R^2F(A) \ar[r] & R^2F(B) \ar[r] & R^2F(C) \ar[r] & \ldots
                \\
}
$$

\begin{rmk}
\begin{itemize}
	\item If we consider just one object $A \in \cA$ as a cochain complex concentrated in degree $0$, we can also consider $RF(A)$. Observe that, in that case, an injective resolution $A \to I^\bullet$ is given by a sequence $A \to I^0 \to I^1 \to \ldots$ with the first arrow a monomorphism. Hence, $R^0F(A) = A$ for all $A \in \cA$, so $R^0F=F$.
	\item If $F: \cA \to \cB$ is exact, it sends quasi-isomorphisms into quasi-isomorphisms. Hence, for an injective resolution $A^\bullet \to I^\bullet$, $RF(A^\bullet) = F(I^\bullet)$ is quasi-isomorphic to $F(A^\bullet)$. Thus, $RF$ is isomorphic to $F$ itself (when extended to complexes).
	\item If $F$ is right exact, a dual construction can be given by means of projective objects. In that case, the derived functor is called the left derived functor, $LF: \Dm{\cA} \to \Dm{\cB}$. Analogous considerations can be done for contravariant functors.
	\item Strictly speaking, the derived functor is defined only up to isomorphism, since we have to choose an injective resolution. It becomes clearer if we consider the induced morphism $i: \Chp{(\textbf{Inj}(\cA))} \to \Dp{\cA}$, where $\textbf{Inj}(\cA)$ is the full subcategory of injective objects of $\cA$. It can be proven that all the quasi-isomorphisms in $\Chp{(\textbf{Inj}(\cA))}$ are actually isomorphisms so $i$ is an equivalence of categories (i.e.\ $i$ has an inverse up to natural isomorphism). In that case, if we choose an inverse $i^{-1}: \Dp{\cA} \to \Chp{(\textbf{Inj}(\cA))}$, the derived functor is given by the composition
$$
	RF: \Dp{\cA} \stackrel{i^{-1}}{\to} \Chp{(\textbf{Inj}(\cA))} \stackrel{F}{\to} \Chp{\cB} \to \Dp{\cB}.
$$
Hence, the definition of $RF$ depends on the chosen inverse $i^{-1}$.
\end{itemize}
\end{rmk}

\begin{ex}
\begin{itemize}
	\item Let $X$ be a topological space and let $\Gamma: \Sh{X} \to \Vect{\QQ}$ be the functor of global sections of sheaves on $X$, $\Gamma(\cF)=\cF(X)$ for $\cF \in \Sh{X}$. The $k$-the derived functor of $\Gamma$ is the sheaf cohomology of $\cF$, $R^k\Gamma(\cF) = H^k(X; \cF)$. Analogously, if $X$ is locally compact, we can consider the functor of global sections with compact support $\Gamma_c: \Sh{X} \to \Vect{\QQ}$. In that case, $R^k\Gamma_c(\cF) = H^k_c(X; \cF)$ is the compactly supported $k$-th sheaf cohomology.
	\item Fix a ring $R$ and, for a $R$-module $M$, consider the functor $- \otimes_R M: \Mod{R} \to \Mod{R}$. It is a right exact functor and its $k$-th derived functor $L^k(- \otimes_R M)(N) = \textrm{Tor}_k^R(N,M)$ is called the \emph{Tor functor}. Observe that $M$ is flat exactly when $ \textrm{Tor}_1^R(N,M) = 0$ for all $N \in \Mod{R}$. On the other hand, the functor $\Hom(M, -): \Mod{R} \to \Mod{R}$ is right exact and its derived functors $R^k\Hom(M, -)(N)=\textrm{Ext}^k_R(M,N)$ are called the \emph{Ext functors}.
\end{itemize}
\end{ex}

\subsection{Grothendieck group}
\label{sec:grothendieck-groups}

From an abelian category, we can construct an abelian group, the so-called Grothendieck group, that will be extremelly useful for our purposes. Here, we will sketch the fundamental properties of this construction. For a complete description of this important theory, see \cite{Weibel:2013}.

\begin{defn}
Let $\cA$ be an abelian category which is skeletally small (i.e.\ the class of isomorphism classes of objects of $\cA$ is a genuine set). We define the \emph{Grothendieck group} or the \emph{$K$-theory} of $\cA$, $\K{\cA}$, as the abelian group generated by the isomorphism classes $[A]$ of objects $A \in \cA$ subject to the restriction that $[A] = [B] + [C]$ whenever there exists a short exact sequence in $\cA$
$$
	0 \longrightarrow B \longrightarrow A \longrightarrow C \longrightarrow 0.
$$
The group structure is given by $[A] + [B] = [A \oplus B]$ for $A, B \in \cA$.
\end{defn}

Obviously, in that situation we have a function $\Obj{\cA} \to \K{\cA}$ that sends $A \in \cA$ to the image of its isomorphism class $[A] \in \K{\cA}$.

\begin{rmk}\label{rmk:grothendieck-semiring}
\begin{itemize}
	\item If the category $\cA$ is not abelian but only triangulated, we can still define the Grothendieck group of $\cA$ by declaring that $[A]=[B]+[C]$ for every distinguished triangle $B \to A \to C \to B[1]$. In particular, we can consider the Grothendieck group of the derived category of an abelian category.
	\item As an intermediate step, we can also consider the \emph{Grothendieck semiring} of a skeletally small abelian category $\cA$, $\Sr{\cA}$, as the abelian semi-group generated by the isomorphism classes of objects of $\cA$ after decreeing that all the short exact sequences split. From $\Sr{\cA}$, the Grothendick group $\K{\cA}$ can be constructed by introducing formal inverses to all the elements of $\Sr{\cA}$, in the same spirit as $\ZZ$ is constructed from $\NN$.
	\item It may be considered higher versions of $K$-theory, as described in \cite{Weibel:2013}. For this reason, it is customary in the literature to denote the Grothendieck group of $\cA$ by $K_0\cA$ in order to distinguish it from the higher ones. However, we will not use these higher versions so $\K{\cA}$ always stands for Grothendieck groups.
\end{itemize}
\end{rmk}

It can be proven (\cite{Weibel:2013}, 6.1.2) that $\K{\cA}$ satisfies the following universal property. Suppose that $G$ is an abelian group and that $f: \Obj{\cA} \to G$ is a function with $f(A) = f(B) + f(C)$ whenever we have a short exact sequence
$
	0 \to B \to A \to C \to 0.
$
Then, there exists an unique group homomorphism $\tilde{f}: \K{\cA} \to G$ with $\tilde{f} \circ i = f$.
	\[
\begin{displaystyle}
   \xymatrix
   {
  	 \Obj{\cA} \ar[rd]^{f} \ar[d]_{i} & \\
  	 \K{\cA} \ar@{--{>}}[r]_{\tilde{f}} & G
   }
\end{displaystyle}
\]

Moreover, suppose that we have a bounded chain complex $A^\bullet \in \Chb{\cA}$. We define its \emph{Euler characteristic}, $\chi(A^\bullet)$, as the element
$$
	\chi(A^\bullet) = \sum_k (-1)^k \,[H^k(A^\bullet)] = \sum_k (-1)^k \,[A^k] \in \K{\cA},
$$
where the last equality can be shown by splitting the long exact sequence in cohomology into sort exact sequences.
In particular, quasi-isomorphic chain complexes map to the same element of $\K{\cA}$ and $\chi(A^\bullet)=0$ if $A^\bullet$ is exact. Hence, the Euler characteristic descends to the bounded derived category to give a map
$$
	\chi: \Db{\cA} \to \K{\cA}.
$$
Using the universal property of $K$-theory, this map descends to the $K$-theory of $\Db{\cA}$ to give a map $\tilde{\chi}: \K{(\Db{\cA})} \to \K{\cA}$. This $\tilde{\chi}$ is actually an isomorphism so we have a natural identification $\K{(\Db{\cA})} = \K{\cA}$.

\begin{prop}\label{prop:induced-functor-K}
Let $\cA$ and $\cB$ be abelian categories and let $F: \cA \to \cB$ be a left-exact functor with bounded derived functor $RF: \cA \to \Db{\cB}$. Then, $F$ descends to a map $\K{F}: \K{\cA} \to \K{\cB}$ given by
$$
	\K{F}([A]) = \chi(RF(A)) = \sum_k (-1)^k\,[R^kF(A)].
$$
\begin{proof}
Consider the map $h: \Obj{\cA} \to \K{\cB}$ that sends $A \mapsto \chi(RF(A))$. By the long exact sequence in cohomology of the derived functor, for any short exact sequence in $\cA$,
$
	0 \to B \to A \to C \to 0
$
, we have that $\chi(RF(A)) = \chi(RF(B)) + \chi(RF(C))$. Hence, by the universal property of $\K{\cA}$, $h$ descends to $K$-theory to give the desired map $\K{F} = \tilde{h}: \K{\cA} \to \K{\cB}$. 
\end{proof}
\end{prop}

Furthermore, suppose that our skeletally small abelian category $\cA$ has also a symmetric monodal structure $(\cA, \otimes, I)$. If $\otimes: \cA \times \cA \to \cA$ is additive on each argument separately, $\cA$ is said to be an \emph{abelian monoidal category}. If, in addition, for all $A \in \cA$, the functors $-\otimes A: \cA \to \cA$ and $A \otimes -: \cA \to \cA$ are exact then the Grothendieck group of $\cA$, $\K{\cA}$, inherits also a ring structure by $[A] \cdot [B] = [A \otimes B]$ for $A, B \in \cA$. In this case, it is customary to call $\K{\cA}$ the \emph{Grothendick ring} or the \emph{$K$-theory ring}.

\begin{rmk}
The requirement that $-\otimes A$ and $A \otimes -$ are exact functors is usually refered to as $\otimes$ being biexact. This property is needed since, if we have a short exact sequence
$$
	0 \longrightarrow B \longrightarrow A \longrightarrow C \longrightarrow 0,
$$
then, for $D \in \cA$, $([B] + [C]) \cdot [D] = [A] \cdot [D]$ if and only if
$$
	0 \longrightarrow B \otimes D \longrightarrow A \otimes D \longrightarrow C \otimes D \longrightarrow 0
$$
is exact, that is, if $- \otimes D$ is exact. In the case that $\otimes$ was not biexact, we could consider, instead of $\K{\cA}$, the abelian group generated by the semi-group $(\Obj{\cA}, \oplus)$ without quotienting by short exact sequence. In that situation, we can still endow such a group with a ring structure but it does not descend to a ring structure on $\K{\cA}$.
\end{rmk}

\begin{ex}\label{ex:grothendieck-groups}
\begin{itemize}
	\item Let us consider the subcategory of finite dimensional $k$-vector spaces, $\Vecto{k}$. It is an abelian monoidal category with direct sum given by cartesian product and monoidal structure given by tensor product, and the tensor product is biexact. The dimension of a vector space gives a correspondence between isomorphism classes of objects of $\Vecto{k}$ and $\NN$. Under this identification, direct sum and tensor product correspond to sum and product of natural numbers, respectively. Hence, $\K{(\Vecto{k})} = \ZZ$ as a ring. Actually, the construction is exactly the same as when we construct $\ZZ$ from the semi-ring $\NN$.
	\item The category of rational sheaves on a topological space $X$, $\Sh{X}$, is an abelian monoidal category with tensor product. Even more, since exactness can be checked stalkwise and the rational vector spaces are flat, $\otimes$ is biexact. Hence, the Grothendieck group $\K{\Sh{X}}$ has a natural ring structure. If we consider the subcategory $\textbf{Vb}(X)$ of rational vector bundles, then its Grothendieck ring $\K{\textbf{Vb}(X)}$ is called the (rational) K-theory of $X$ and it is customary to denote it $K^0(X)$ or just $K(X)$.
\end{itemize}
\end{ex}

\begin{rmk}\label{rmk:KVar}
The category of algebraic varieties over a field $k$, $\Var{k}$, is not an abelian category (it is not even triangulated). However, analogously to the abelian case, we can construct its Grothendieck semi-group as the abelian semi-group generated by the isomorphism classes of algebraic varieties with the relation that $[X] = [Y] + [U]$ if $X = Y \sqcup U$ with $Y \subseteq X$ a closed subvariety. Moreover, this semi-group can be extended to obtain an abelian group, denoted $\K{\Var{k}}$. In that case, the cartesian product of varieties also endows $\K{\Var{k}}$ with a ring structure. For this reason, $\K{\Var{k}}$ is known as the Grothendieck ring of algebraic varieties.

In this context, if $R$ is any ring, a ring homomorphism $E: \K{\Var{k}} \to R$ is called a \emph{generalized Euler characteristic}. Observe that such a homomorphism is uniquely determined by a map $\hat{E}: \Obj{\Var{k}} \to R$ such that $\hat{E}(X) = \hat{E}(Y) + \hat{E}(U)$ for any decomposition $X = Y \sqcup U$ with $Y \subseteq X$ closed and such that $\hat{E}(X \times Y)=\hat{E}(X)\cdot \hat{E}(Y)$ for all varieties $X, Y$.
\end{rmk}

Using the algebraic tool of Grothedieck rings, we can create an important class of $\cC$-algebras (cf.\ Remark \ref{rmk:notation-grothendieck-six}) from simpler data, in the spirit of Grothendieck six operators.

\begin{prop}\label{prop:k-theory-algebra}
Let $\cC$ be a cartesian monoidal category. Let $H$ be a functor
$
	H: \cC \to \Cat^{op} \times \Cat
$
such that, for all $a \in \cC$, $H(a) = (\cH_a, \cH_a)$ for some category $\cH_a$. Given $f: a \to b$ in $\cC$, let us denote $H(f) = (f^*, f_!)$, with $f^*: \cH_b \to \cH_a$ and $f_!: \cH_a \to \cH_b$. Suppose that
\begin{enumerate}
	\item For all $a \in \cA$, $\cH_a$ a triangulated monoidal category with biexact monoidal structure.
	\item For every morphism $f$ of $\cC$, the functor $f^*$ is monoidal.
	\item \label{prop:k-theory-algebra:beck-chevalley} For any pullback diagram in $\cC$
\[
\begin{displaystyle}
   \xymatrix
   {
	d \ar[r]^{\;g'}\ar[d]_{f'} & a_1\ar[d]^{f} \\
	a_2 \ar[r]_{g} & b
   }
\end{displaystyle}
\]
we have a natural isomorphism $g^* \circ f_! \cong f'_! \circ g'^*$. That is, the following diagram commutes up to isomorphism
		\[
\begin{displaystyle}
   \xymatrix
   {
  	 {\cH_d} \ar[d]_{f'_!} & {\cH_{a_1}} \ar[l]_{\;\;\;g'^*} \ar[d]^{f_!} \\
  	 {\cH_{a_2}} & {\cH_{b}} \ar[l]^{\;\;\,g^*}
   }
\end{displaystyle}
\]
	\item\label{prop:k-theory-algebra:ext-prod} For any $f: a \to b$, $g: a' \to b'$ and any $A \in \cH_a$ and $A' \in \cH_{a'}$, we have
	$$
		(f_! \otimes g_!)(p_1^*A \otimes p_2^*A') = \rho_1^* f_!A \otimes \rho_2^* f_! A',
	$$
where $p_1: a \times a' \to a$, $p_2: a \times a' \to a'$, $\rho_1: b \times b' \to b$ and $\rho_2: b \times b' \to b'$ are the projections.
\end{enumerate}

In that situation, $H$ induces a $\cC$-algebra, denoted by $\K{\cH}$.
\begin{proof}
The complete definition of $\K{\cH} = (A,B)$ is the following. The contravariant functor $A: \cC \to \Rng$ assigns, to $a \in \cC$, the Grothendieck ring $\K{\cH}_a$ and, to a morphism $f: a \to b$, the induced map $f^* = \K{f^*}: \K{\cH}_b \to \K{\cH}_a$ which, since $f^*$ is monoidal, is a ring homomorphism.

On the other hand, the functor $B: \cC \to \Mod{\K{\cH}_\star}$ assigns, to $f: a \to b$, the induced homomorphism of abelian groups $f_! = \K{f_!}: \K{\cH}_a \to \K{\cH}_b$ of Proposition \ref{prop:induced-functor-K}. Observe that $f_!$ preserves the exterior product induced by $A$ by condition \ref{prop:k-theory-algebra:ext-prod}. Exactly by the same reason, $f_!$
is a homomorphism of $\K{\cH}_\star$-modules since the module structure can be given in terms of the exterior product (see Remark \ref{rmk:module-exterior-prod}). Finally, $\K{\cH}$ satisfies the Beck-Chevalley condition directly by item \ref{prop:k-theory-algebra:beck-chevalley}.
\end{proof}
\end{prop}

\begin{ex}
A typical case in which this result can be applied is the following. Suppose that, for every $a \in \cC$, we have an abelian category $\cA_a$ and, for every morphism $f: a \to b$ in $\cC$, we have induced morphisms $f^*: \cA_b \to \cA_a$ and $f_!: \cA_a \to \cA_b$. Moreover, suppose that $f^*$ is exact monoidal and that $f_!$ is left exact with bounded derived functor. In that case, we set $H(a) = (\Db{\cA_a},\Db{\cA_a})$ for $a \in \cC$ and, for $f: a \to b$ in $\cC$, we set $H(f)=(f^*, Rf_!)$ as morphisms between the derived categories. If $(f^*, Rf_!)$ have the Beck-Chevalley condition and $Rf_!$ preserves the external product, then $\K{\cH}$ is a $\cC$-algebra with object $\K{\cH_a}=\K{\Db{\cA_a}} = \K{\cA_a}$ for $a \in \cC$. This was the version of the proposition applied in Example \ref{ex:algebra-sheaves}.
\end{ex}

\section{Hodge structures}
\label{sec:hodge-structures}

In this section, we will discuss briefly the fundamentals of Hodge theory. It is a beautiful topic with many links with different areas of mathematics but its main application is as a subtle invariant that traces the complex structure of an algebraic variety. For a more extensive introduction to Hodge theory and its applications, see \cite{Peters-Steenbrink:2008} or \cite{Zein-Trang:2014}.

\subsection{Pure Hodge structures}

The idea of Hodge theory is that the cohomology of a complex algebraic variety carries naturally an extra linear structure called a Hodge structure. In order to explain them, in this section we will focus on the simplest case of Hodge structures, the so-called pure Hodge structures. These are precisely the Hodge structures that are present on projective varieties.

\begin{defn}\label{defnPureHodgeStruc}
Let $H$ be a rational vector space of finite dimension. A \emph{pure (rational) Hodge structure} of weight $k$ on $H$, or just a Hodge structure of weight $k$, is a direct sum decomposition of $H_\CC = H \otimes_\QQ \CC$ as
$$
	H_\CC = \bigoplus\limits_{p+q=k} H^{p,q},
$$
such that $\overline{H^{p,q}}=H^{q,p}$ for all $p,q$. A \emph{morphism} of Hodge structures of weight $k$ is a $\QQ$-linear map $f: H \to H'$ such that, for all $p,q \in \ZZ$ with $p+q=k$ we have that $f_\CC: H^{p,q} \to H'^{p,q}$, where $f_\CC$ is the $\CC$-linear extension of $f$. For $k \in \ZZ$, the category of pure Hodge structures of weight $k$ and morphism of Hodges structures will be denoted by $\HS{k}$. 
\end{defn}

\begin{ex}
Given $m \in \ZZ$, we define the \emph{Tate structures}, $\QQ(m)$, as the pure Hodge structure whose underlying rational vector space is $(2\pi i)^m \QQ \subseteq \CC$ with a single-piece decomposition $\QQ(m) = \QQ(m)^{-m,-m}$ that converts it into a pure Hodge structure of weight $-2m$. For short, we will denote $\QQ_0 = \QQ(0)$ (or just $1 = \QQ(0)$ when it does not cause confusion), the Tate structure of weight zero. From Tate's Hodge structures, if $H$ is another pure Hodge structure of weight $k$ then $H(m) = H \otimes \QQ(m)$ (see below for the definition of tensor product) is a pure Hodge structure of weight $k-2m$, called the \emph{Tate twist} of $H$.
\end{ex}

We can use the general contructions of linear algebra to build new Hodge structures. Suppose that $H$ and $H'$ have pure Hodge structures of weights $k$ and $k'$, respectively.
\begin{itemize}
	\item There is a Hodge structure of weight $k+k'$ on $H \otimes_\QQ H'$ by defining
	$$
		(H \otimes H')_\CC^{p,q} = \bigoplus_{\substack{p_1 + p_2=p \\ q_1+q_2=q}} H^{p_1,q_1} \otimes H^{p_2,q_2}.
	$$
	\item There is a Hodge structure of weight $-k$ on the dual vector space $H^*$ by taking
	$$
		{H_\CC^*}^{p,q} = \left(H^{-p, -q}\right)^*.
	$$
	\item More generally, $\Hom_\QQ(H, H')$ have a pure Hodge structure of weight $k'-k$ by defining
	$$
		Hom(H, H')_\CC^{p,q} = \bigoplus_{\substack{p_2 - p_1=p \\ q_2-q_1=q}} Hom_\CC(H^{p_1,q_1}, H^{p_2,q_2}).
	$$
\end{itemize}

\begin{rmk}
With this definition, $\HS{k}$ can be endowed with a monoidal structure by means of the tensor product. In this way, the forgetful functor $\HS{k} \hookrightarrow \Vect{\QQ}$ that sends a Hodge structure into its underlying vector space is an additive monoidal functor.
\end{rmk}

The origin of pure Hodge structures can be traced back to a classical theorem of Hodge about K\"ahler manifolds. 

\begin{thm}[Hodge]\label{HodgePureStructureKahler}
Let $M$ be a compact K\"ahler manifold. For all $k \geq 0$, the cohomology $H^k(M;\QQ)$ carries a pure Hodge structure of weight $k$ given by the Hodge decomposition
$$
	H^k(M; \CC) = \bigoplus\limits_{p+q=k} H^{p,q}(M),
$$
where $H^{p,q}(M)$ denotes the Dolbeaut cohomology of $M$. Furthermore, if $\Kahc$ is the full subcategory of $\Diffc$ of compact K\"ahler manifolds, then, for all $k$, the cohomology functor $H^k: \Kahc \to \Vect{\QQ}$ factorizes through $\HS{k}$ by means of a monoidal functor $\Kahc \to \HS{k}$
	\[
\begin{displaystyle}
   \xymatrix
   {
   	\Kahc \ar[rr]^{H^k} \ar@{--{>}}[dr]& & \Vect{\QQ} \\
   	& \HS{k} \ar@{^{(}-{>}}[ru]&
   }
\end{displaystyle}   
\]
\end{thm}

The proof of this theorem is one of the greatest advances of classical Hodge theory. The proof is very beautiful since it links hard analysis techniques with subtle algebraic facts. Roughly speaking, the idea is that, in a general compact manifold, every cohomology class of $H^k(M;\CC)$ contains one and only one harmonic representative (with respect to the Laplace-Beltrami operator $\Delta_d = dd^\dag + d^\dag d$). In this way, $H^k(M;\CC)$ can be naturally identified with the space of harmonic $k$-forms. Analogously, if $M$ is a complex manifold, the Dolbeault cohomology $H^{p,q}(M)$ can be naturally identified with the space of $\Delta_{\overline{\partial}}$-harmonic $k$-forms, where $\Delta_{\overline{\partial}} = \overline{\partial}\overline{\partial}^\dag+\overline{\partial}^\dag\overline{\partial}$ is the (anti)-holomorphic Laplace-Beltrami operator. However, in the K\"ahler case, magic happens and the two operators are related in a very simple way by $\Delta_d = 2\Delta_{\overline{\partial}}$. This implies that every harmonic $k$-form can be decomposed into its $(p,q)$-pieces and these pieces are still harmonic. This decomposition gives the desired result.

\begin{ex}
Let $M$ be a compact Kähler manifold of real dimension $2n$. Observe that the integration map in top cohomology
$$\begin{matrix}
	\int: & H^{2n}(M; \QQ) & \to & \QQ\\
	 & [\omega] & \mapsto & \displaystyle{\frac{1}{(2\pi i)^n}\int_M \omega}\\
\end{matrix}
$$
cannot respect the Hodge structures on $H^{2n}(M; \QQ)_\CC=H^{n,n}(M)$ and $\QQ_\CC=\QQ^{0,0}$, since they have different weights. However, if we twist the grading in $\QQ$, then the map $\int: H^{2n}(M; \QQ)  \to \QQ(-n)$ do respect the Hodge structures. Moreover, using this shifted integration map, we have that the bilinear form used in the Poincar\'e duality
$$
	H^{k}(M; \QQ) \otimes_\QQ H^{2n-k}(M; \QQ) \stackrel{\cap}{\to} H^{2n}(M; \QQ) \stackrel{\int}{\to} \QQ(-n)
$$
respects the Hodge structures.
\end{ex}

Generalizing this concept to the general framework, we obtain the so-called polarization on a pure Hodge structure.

\begin{defn}
Let $H$ be a pure Hodge structure of weight $k$. A \emph{polarization} on $H$ is a bilinear form
$$
	Q: H \otimes_\QQ H \to \QQ(-k)
$$
preserving the Hodge structures and such that
\begin{itemize}
	\item $Q$ is symmetric for $k$ even and antisymmetric for $k$ odd.
	\item With respect to its complexification, $Q_\CC: H_\CC \otimes_\CC H_\CC \to \CC$, the spaces $H^{p,q}$ and $H^{p',q'}$ are orthogonal for $p \neq p'$ or $q \neq q'$.
	\item For all $p,q$, the quadratic form $Q^{p,q}: H^{p,q} \otimes H^{p,q} \to \CC$ given by $Q^{p,q}(x,y)=i^{p-q}Q_\CC(x, \overline{y})$ is positive-definite.
\end{itemize}
A pure Hodge structure that admits a polarization is said to be \emph{polarizable}.
\end{defn}

\begin{defn}
The category of polarizable pure Hodge structures, $\MHSq$, is the category whose objects are finite sums of polarizable pure Hodge structures (maybe of different weights) and whose morphism are finite sums of morphisms of pure Hodge structures. It is an abelian category (see \cite{arapura-kang:2006} and \cite{DeligneII:1971}).
\end{defn}

There is an alternative way of giving a pure Hodge structure which will be more useful for defining mixed Hodge structures. 
Given a pure Hodge structure $H$ of weight $k$, we define the subspaces
$$
	F_p H_\CC = \bigoplus_{r \geq p} H^{r, k-r}.
$$
Observe that the $F_p$ form a finite decreasing filtration
$$
	H_\CC \supseteq \ldots \supseteq F_{p-1}H_\CC \supseteq F_p H_\CC \supseteq F_{p+1} H_\CC \supseteq \ldots \supseteq \left\{0\right\}
$$
called the associated \emph{Hodge filtration}. Moreover, the fact that $H^{q,p} = \overline{H^{p,q}}$ translates as $F_p H_\CC \oplus \overline{F_{k-p+1} H_\CC} = H_\CC$ for all $p \in \ZZ$. The original decomposition can be recovered from the filtration by $H^{p,q} = F_p H_\CC \cap \overline{F_q H_\CC}$.
Therefore, we can reformulate the property of having a Hodge structure of weight $k$.

\begin{defn}[Equivalent to \ref{defnPureHodgeStruc}]
Given a finite dimensional rational vector space $H$, a {pure Hodge structure} on $H$ of weight $k$ is a decreasing filtration of $H_\CC$
$$
	H_\CC \supseteq \ldots \supseteq F_{p-1}H_\CC \supseteq F_p H_\CC \supseteq F_{p+1} H_\CC \supseteq \ldots \supseteq \left\{0\right\}
$$
such that, for all $p \in \ZZ$, $F_p H_\CC \oplus \overline{F_{k-p+1} H_\CC} = H_\CC$.
\end{defn}

\subsection{Mixed Hodge structures}
\label{sec:mixed-hodge-structures}

With a view towards the Weil conjectures, in a serie of articles published between 1971 and 1974 (\cite{DeligneI:1971}, \cite{DeligneII:1971} and \cite{DeligneIII:1971}), Deligne extended the notion of pure Hodge structures into a larger category, rather abstract and artificial, known as mixed Hodge structures. With this general definition, he could prove that the cohomology ring of any complex algebraic variety is naturally endowed with a mixed Hodge structure.

\begin{defn}
Let $H$ be a finite dimensional rational vector space. A \emph{mixed Hodge structure} on $H$ is pair of filtrations.
\begin{itemize}
	\item A finite increasing filtration of $H$, $W^\bullet H$, called the \emph{weight filtration}
$$
	0 \subseteq \ldots \subseteq W^{k-1} H \subseteq W^k H \subseteq W^{k+1} H \subseteq \ldots \subseteq H.
$$
	\item A finite decreasing filtration of the $\CC$-vector space $H_\CC$, $F_\bullet H_\CC$, called the \emph{Hodge filtration}
$$
	H_\CC \supseteq \ldots \supseteq F_{p-1} H_\CC \supseteq F_p H_\CC \supseteq F_{p+1} H_\CC \supseteq \ldots \supseteq 0.
$$
\end{itemize}
They have to satisfy that, for every $k \in \ZZ$, the induced filtration of $F_\bullet$ on $\left(\Gr{W}{k} H\right) \otimes_\QQ \CC$ gives a pure Hodge structure of weight $k$ on $\Gr{W}{k} H$. A mixed Hodge structure is said to be \emph{graded polarizable} if each of the pieces $\Gr{W}{k} H$ are polarizable.

Given two mixed Hodge structures $H$ and $H'$, a homomorphism $f: H \to H'$ is called a \emph{morphism of mixed Hodge structures} if $f: H \to H'$ is a filtered morphism with respect to $W^\bullet H$ and $f_\CC: H_\CC \to H'_\CC$ is a filtered morphism with respect to $F_\bullet H_\CC$. With these definitions, the category of (graded polarizable) mixed Hodge structures will be denoted by $\MHS{\QQ}$.
\end{defn}

\begin{ex}\label{ex:MixedHodgeGradingPure}
\begin{itemize}
	\item Let $H$ be a pure Hodge structure of weight $k$, induced by a decreasing filtration $F_\bullet$ of $H_\CC$. Then, $H$ also has a mixed Hodge structure by taking the Hodge filtration as $F_\bullet H_\CC$ and the weight filtration as $W^s H=0$ for $s<k$ and $W^s H =H$ for $s \geq k$.
	\item A mixed Hodge structure $H$ is said to be \emph{balanced} or of \emph{Hodge-Tate} type if $\Gr{F}{p}\Gr{p+q}{W}H = 0$ for $p \neq q$. In that case, for all $p$, $\Gr{2p}{W}H$ has a single step decomposition as $\Gr{2p}{W}H_\CC = \left(\Gr{2p}{W}H\right)^{p,p}$. 
\end{itemize}
\end{ex}

One of the most importants results about the category of mixed Hodge structures is that it is well behaved with respect to kernels and cokernels. The proof of the following theorem can be found in \cite{Cattani-ElZein-Griffiths} or \cite{DeligneI:1971}.

\begin{thm}[Deligne]
The category of mixed Hodge structures, $\MHS{\QQ}$, is an abelian category.
\end{thm}

\begin{rmk} 
In analogy with the pure case, given mixed Hodge structures $H$ and $H'$, the direct sum can be endowed with a mixed Hodge structure just by taking direct sum of the filtrations. Analogous considerations can be done to equip tensor products, dual spaces and homomorphisms of mixed Hodge structures with a mixed Hodge structure. However, the definition is not as straightforward and can be checked in \cite{Peters-Steenbrink:2008}, Example 3.2.2. With these definitions, the forgetful functor $\MHS{\QQ} \to \Vect{\QQ}$ is an additive monoidal functor.
\end{rmk}

The main theorem in this bussiness, and the reason for the introduction of of mixed Hodge structures is the following result, whose proof can be found in \cite{Cattani-ElZein-Griffiths} or \cite{DeligneII:1971}-\cite{DeligneIII:1971}.

\begin{thm}[Deligne]\label{existenceMHS}
The cohomology ring of any complex algebraic variety carries a natural mixed Hodge structure. Furthermore, for any $k \geq 0$, the cohomology functor $H^k: \CVar \to \Vect{\QQ}$ factorizes through $\MHS{\QQ}$ via a monoidal functor $\CVar \to \MHS{\QQ}$.
	\[
\begin{displaystyle}
   \xymatrix
   {
     \CVar \ar[rr]^{H^k} \ar@{--{>}}[rd] && \Vect{\QQ} \\
     & \MHS{\QQ} \ar@{^{(}->}[ru] &
   }
\end{displaystyle}   
\]
\end{thm}

\begin{rmk}
\begin{itemize}
	\item For a projective variety, the mixed Hodge structure that its cohomology has as complex algebraic variety agrees with the pure Hodge structure given as compact K\"ahler manifold. This is not a trivial result at all and was proven by Deligne in \cite{DeligneII:1971}.
	\item The same theorem also holds if we consider compactly supported cohomology $X \mapsto H^k_c(X; \QQ)$.
\end{itemize}
\end{rmk}

Now, let us analize the Grothendieck group of mixed Hodge structures, $\K{\MHS{\QQ}}$. Suppose that $H \in \MHS{\QQ}$ and that $l$ is the smallest index with $W^lH = H$. In that case, we have a short exact sequence
$$
	0 \longrightarrow W^{l-1}H \longrightarrow H \longrightarrow \Gr{l}{W}H \longrightarrow 0,
$$
so, in $\K{\MHS{\QQ}}$, we have $[H] = [\Gr{l}{W}H] + [W^{l-1}H]$. Analogously, we have a short exact sequence $0 \to W^{l-2}H \to W^{l-1}H \to \Gr{l-1}{W}H \to 0$, so $[W^{l-1}H] = [\Gr{l-1}{W}H] + [W^{l-2}H]$. Hence, proceding inductively, we have that the image of $H$ in $\K{\MHS{\QQ}}$ decomposes as $[H] = \sum_{k} \left[\Gr{k}{W}H\right]$, where the objects $\Gr{k}{W}H$ are polarizable pure Hodge structures of weight $k$. For this reason, we have that $\K{\MHS{\QQ}} = \K{\MHSq}$ (for further details, see \cite{arapura-kang:2006}).

A very important property of mixed Hodge structures on the compactly supported cohomology is that the long exact sequence of a stratification is also a sequence of mixed Hodge structures. The proof of this result can be found in \cite{Peters-Steenbrink:2008}.

\begin{prop}\label{prop:longExactSeqHodgeStructure}
Let $X$ be a complex algebraic variety and decompose it as $X = Y \sqcup U$ where $Y \subseteq X$ a closed subvariety and $U \subseteq X$ is an open set. Then, we have a long exact sequence of mixed Hodge structures in compactly supported cohomology
\[
\begin{displaystyle}
\xymatrix{
             0\ar[r] & H^0_c(U) \ar[r] & H^0_c(X) \ar[r] & H^0_c(Y)
                \ar@{-{>}} `r/8pt[d] `/10pt[l] `^dl[ll]|{\delta} `^r/3pt[dll] [dll] \\
             & H^1_c(U) \ar[r] & H^1_c(X) \ar[r] & H^1_c(Y)
                \ar@{--{>}} `r/8pt[d] `/10pt[l] `^dl[ll] `^d[dlll] `^r[ddlll] [ddll]\\
             &&&\\
             & H^n_c(U) \ar[r] & H^n_c(X) \ar[r] & H^n_c(Y)\ar[r] & \ldots
    }
\end{displaystyle}   
\]
\end{prop}

Now, consider the map $\alpha: \Obj{\CVar} \to \K{\MHS{\QQ}}=\K{\MHSq}$ that, to any complex algebraic variety $X$, assigns the Euler characteristic of the mixed Hodge structure on its compactly supported cohomology $\chi\left(H_c^\bullet(X;\QQ)\right) \in \K{\MHSq}$. By the previous proposition, if we have a decomposition $X = Y \sqcup U$ with $Y \subseteq X$ a closed subvariety, then $\chi\left(H_c^\bullet(X;\QQ)\right) = \chi\left(H_c^\bullet(Y;\QQ)\right) + \chi\left(H_c^\bullet(U;\QQ)\right)$. Moreover, $\alpha$ sends cartesian products into tensor products. Hence, by Remark \ref{rmk:KVar}, $\alpha$ gives rise to a generalized Euler characteristic $\K{\alpha}: \K{\CVar} \to \K{\MHSq}$. Given $[X] \in \K{\CVar}$, we will denote its image by $\coh{X;\QQ} \in \K{\MHSq}$ or just by $\coh{X}$.

A step forward, let us define the map $\epsilon:\Obj{\MHSq} \to \ZZ[u^{\pm 1}, v^{\pm 1}]$ that, to a (polarizable) Hodge structure $H$ of weight $k$, assigns the polynomial $\epsilon(H)=\sum_{p + q = k} \dim_\CC H^{p,q}\, u^{p}v^q$. Observe that this map is additive since, if we have a short exact sequence $0 \to H' \to H \to H'' \to 0$ of Hodge structures, then, as the morphisms preserve the decomposition, we also have short exact sequences $0 \to H'^{p,q} \to H^{p,q} \to H''^{p,q} \to 0$ for all $p,q \in \ZZ$. In particular, $\dim_\CC H^{p,q} = \dim_\CC H'^{p,q} + \dim_\CC H''^{p,q}$, which shows that $\epsilon$ is additive. Hence, $\epsilon$ descends to the Grothendieck ring.

\begin{defn}\label{defn:Deligne-Hodge-pol}
The ring homomorphism $\RDelHod = \K{\epsilon}: \K{\MHSq} \to \ZZ[u^{\pm 1}, v^{\pm 1}]$ is called the \emph{Deligne-Hodge polynomial} or the \emph{$E$-polynomial}. Given a complex algebraic variety $X$, we will short $\DelHod{X} = \DelHod{\coh{X;\QQ}}$ and we will call it the \emph{Deligne-Hodge polynomial} of $X$.
\end{defn}

The computation of these Deligne-Hodge polynomials for a special kind of complex algebraic varieties called character varieties is, precisely, the goal of this thesis.

\begin{rmk}
Unraveling the definitions, for a complex algebraic variety $X$, its $E$-polynomial is given by
\begin{align*}
	\DelHod{X} &= \sum_k(-1)^k \DelHod{H_c^k(X;\QQ)} = \sum_{k,l}(-1)^k \DelHod{\Gr{l}{W}H_c^k(X;\QQ)} \\
	&= \sum_{\substack{k,l,p}}(-1)^k \dim_\CC\left(\Gr{F}{p}\Gr{l}{W}H_c^k(X;\QQ)\right)\, u^pv^{l-p}.
\end{align*}
It is customary to denote $H_c^{k;p,q}(X) = \Gr{F}{p}\Gr{p+q}{W}H_c^k(X;\QQ)$ and $h_c^{k;p,q}(X) = \dim_\CC H_c^{k;p,q}(X)$ and to call them the (compactly supported) \emph{mixed Hodge pieces} and \emph{mixed Hodge numbers} of $X$, respectively. In that case, the Deligne-Hodge polynomial can be written down as
$$
	\DelHod{X} = \sum_{k, p, q}(-1)^k h_c^{k;p,q}(X)\,u^pv^q.
$$
This is the classical definition of the $E$-polynomial as in \cite{Hausel-Rodriguez-Villegas:2008} or \cite{LMN}.
\end{rmk}

In the upcoming arguments, the $(-1)$-Tate structure of weight $2$, $\QQ(-1)$, is going to play a very important role. For this reason, we will short $q = [\QQ(-1)] \in \K{\MHSq}$. Observe that, by classical Hodge theory, we know that $\coh{\PP^1} = 1 + q$. Hence, as $\coh{\star} = 1$ and $[\PP^1] = [\CC] + [\star]$, we have that $\coh{\CC} = \coh{\PP^1} - \coh{\star} = q$.
Moreover, observe that the Deligne-Hodge polynomial of $q$ is $\DelHod{q}=uv \in \ZZ[u^{\pm 1}, v^{\pm 1}]$, which justifies to also denote $q=uv$.

\section{Mixed Hodge modules}
\label{sec:mixed-hodge-modules}

In this section, we will describe a large generalization of mixed Hodge structures called mixed Hodge modules. These are a collection of abelian categories for every complex algebraic variety with functors between them induced by regular morphisms.

The idea of Hodge modules is that, while Hodge structures are, in some sense, pointwise, the category of mixed Hodge module on $X$, $\MHM{X}$, captures the idea of a varying mixed Hodge structure along $X$. However, the construction is not as easy since the na\"ive solution of taking sheaves of Hodge modules is not good enough for our purposes. In fact, the construction of this category is extremely involved an requiere the use of many hard concepts from singularity theory, sheaf theory and category theory. 

For this reason, in this section we will just sketch the most important milestones in the construction of $\MHM{X}$ and we will skip lots of details. For a complete description of mixed Hodge modules with plenty of details, see the original reference \cite{Saito:1990} (also \cite{Saito:1986}). There are also very readable surveys about the construction that are very recomendable like \cite{Saito:1989}, \cite{Saito:2017}, \cite{Schnell:2014} and \cite{Schurmann:2011}. The exposition of this section does not, in any way, intend to be original. It is absolutely limited to rephrase Saito's construction and to reinterpret it in the way in which it is most useful for our purposes.

\subsection{Riemann-Hilbert correspondence}
\label{sec:riemann-hilbert}

In this section, we will briefly review some notions of $\cD$-modules and perverse sheaves and we will show that they are strongly related by the celebrated Riemann-Hilbert correspondence. These constructions will be useful in Section \ref{sec:saito-mixed-hodge-modules} since they are the building blocks of mixed Hodge modules.

\subsubsection*{$\cD$-modules}

Let $X$ be a complex algebraic variety and let $\cO_X^{an}$ be its sheaf of holomorphic functions (i.e.\ its structure sheaf as analytic variety). The sheaf of \emph{differential operators}, $(\cD_X, F^\bullet)$, is an increasing filtered subsheaf of $\End_\CC(\cO_X^{an})$ defined recursively as follows.
\begin{itemize}
	\item $F^0\cD_X = \cO_X^{an}$, seen $\cO_X^{an} \subseteq \End_\CC(\cO_X^{an})$ via de multiplication operator i.e.\ $f(g) = fg$ for $f,g \in \cO_X^{an}$.
	\item Given $T \in \End_\CC(\cO_X^{an})$, we declare $T \in F^k\cD_X$ if and only if $[T,f] \in F^{k-1}\cD_X$ for all $f \in \cO_X^{an}$. Here, the commutator is the usual commutator of operators.
\end{itemize}

From this filtration, we define $\cD_X = \bigcup_{k\geq 0} F^k\cD_X$. 
As its name suggests, a $\cD_X$-module (aka.\ a $\cD$-module) is just a sheaf of $\cD_X$-modules on the left. Given $Z \subseteq X$ a subvariety of $X$, a $\cD_X$-module $M$ is said to have \emph{strict support on $Z$} if every non-zero subobject or quotient of $M$ has support equal to $Z$.

\begin{rmk}
\begin{itemize}
	\item By construction, $F^k\cD_XF^l\cD_X \subseteq F^{k+l}\cD_X$.
	\item Suppose that $X$ is smooth of complex dimension $n$ and $(z_1, \ldots, z_{n})$ are local (holomorphic) coordinates in an open set $U \subseteq X$. In that case, any $T \in F^k\cD_X(U)$ can be uniquely written as
	$$
		T = \sum_{0 \leq i_1 + \ldots + i_{n} \leq k} f_{i_1, \ldots, i_{n}} \left(\frac{\partial}{\partial z_1}\right)^{i_1}\cdots \left(\frac{\partial}{\partial z_{n}}\right)^{i_{n}},
	$$
for some $f_{i_1, \ldots, i_{n}} \in \cO_X^{an}(U)$. Equivalently, $\cD_X$  is generated, as a ring, by $\cO_X^{an}$ and the holomorphic tangent sheaf $\cT X$. In particular, $F^1\cD_X = \cO_X^{an} \oplus \cT X$.
\end{itemize}
\end{rmk} 

As we mentioned in the previous remark, if $X$ is smooth, $\cD_X$ is generated by the structure sheaf and the tangent sheaf. For this reason, one can think that a $\cD_X$-module is nothing but a $\cO_X^{an}$-module with a well-behaved action of $\cT X$. This is exactly what the following proposition states, whose proof can be found in \cite{Peters-Steenbrink:2008}, Lemma 13.31.

\begin{prop}
Let $M$ be a $\cO_X^{an}$-module and suppose that $M$ has an action of $\cT X$ such that, for any $f,g \in \cO_X^{an}$, $v, w \in \cT X$ and $m \in M$
\begin{itemize}
	\item (Linearity): $(fv + gw) \cdot m = (fv \cdot m + gw\cdot m)$.
	\item (Leibniz rule): $v \cdot (fm)=fX \cdot m + v(f) \cdot m$.
	\item (Compatibility with bracket): $[v, w] \cdot m = v \cdot w \cdot m - w \cdot v \cdot m$.
\end{itemize}
Then, the action of $\cT X$ lifts to an action of $\cD_X$ in such a way that $M$ becomes a $\cD_X$-module.
\end{prop}

\begin{rmk}
The previous hypothesis are equivalent to have a morphism of Lie algebra sheaves $\cT X \to \cH om(M,M)$.
\end{rmk}

\begin{cor}
Let $\cV$ be a holomorphic vector bundle, seen as a locally free $\cO_X^{an}$-module. Then, a $\cD_X$-module structure on $E$ is equivalent to an itegrable holomorphic connection on $\cV$.
\begin{proof}
An integrable connection $\nabla: \cV \to \cV \otimes \mathbf{\Omega}^1$, where $\mathbf{\Omega}^\bullet$ is the sheaf of holomorphic forms, gives a compatible action by defining $v \cdot \psi = \nabla_v \psi$ for $v \in \cT M$ and $\psi \in \cV$. Moreover, this construction can be clearly reversed.
\end{proof}
\end{cor}

\begin{rmk}
So far, we have considered only $\cD_X$-modules on the left. Analogously, we can also consider right $\cD_X$-modules and, actually, there is an equivalence of categories between left and right $\cD_X$-modules. The point for such an equivalence is to observe that the canonical sheaf $K_X$ has a natural compatible action on the right of $\cT X$ by means of the Lie derivative. In this way, $K_X$ is naturally a right $\cD_X$ module and, thus, given a left $\cD_X$-module $M$ we can turn it into a right $\cD_X$-module by taking $M^{\textrm{rgt}} =K_X \otimes_{\cO_X^{an}} M$. Analogously, if $M$ is a right $\cD_X$-module, then $M^{\textrm{lft}} = \cH om(K_X, M)$ is a left $\cD_X$-module. In particular, we define the \emph{dualizing module} as the left $\cD_X$-module $\cD_X^* = \left(\cD_X^{\textrm{rgt}}\right)^{\textrm{lft}} = \cH om_{\cO_X^{an}}(K_X, K_X \otimes_{\cO_X^{an}} \cD_X)$.
\end{rmk}

The theory of $\cD_X$-modules is compatible with respect to morphisms. Suppose that $f: X \to Y$ is an holomorphic map between complex varieties. From this morphism, given a $\cD_Y$-module $N$ on $Y$, we can define consider the pullback sheaf $f^{-1}N$ on $X$. However, it has no natural $\cD_X$-module structure. If we want to promote it to a $\cD_X$-module, we have to consider $f^*N = \cO_X^{an} \otimes_{f^{-1}\cO_Y^{an}} f^{-1}N$ which is the right pullback in between the categories of $\cD$-modules $f^*: \Mod{\cD_Y} \to \Mod{\cD_X}$. Furthermore, it can be extended to the derived category to obtain a functor $f^!: \Db{\left(\Mod{\cD_Y}\right)} \to \Db{\left(\Mod{\cD_X}\right)}$.  Analogously, the usual direct image morphism $f_*: \Sh{X} \to \Sh{Y}$ can be refined in order to obtain a corresponding morphism $f_+: \Db{\left(\Mod{\cD_X}\right)} \to \Db{\left(\Mod{\cD_Y}\right)}$. For a precise definition of all these functors, see \cite{Peters-Steenbrink:2008}, Section 13.3.4.

Moreover, given a $\cD_X$-module $M$, using the dualizing module we can consider its dual module as the left $\cD_X$-module $M^* = \cH om_{\cD_X}(M, \cD_X^*)$. Moreover, such a dualizing procedure can be promoted to the bounded category to obtain a \emph{dualization operator} $\DD_X: \Db{(\Mod{\cD_X})} \to \Db{(\Mod{\cD_X})}$ which is given as $\DD_X M^\bullet = \cE xt_{\cD_X}(M^\bullet, \cD_X^*)[\dim_\CC X]$, where $\cE xt = R(\cH om)$. From this dualization operator, we define $f^+ = \DD_X \circ f^! \circ \DD_Y$ and $f_! = \DD_Y \circ f_+ \circ \DD_X$ so we have a commutative diagram
\[
\begin{displaystyle}
   \xymatrix
   {
	\Db{\left(\Mod{\cD_X}\right)} \ar@/^1pc/[r]^{f_!} \ar@{<->}[dd]_{\DD_X} & \Db{\left(\Mod{\cD_Y}\right)} \ar@/^1pc/[l]^{f^+} \ar@{<->}[dd]^{\DD_Y} \\
	&\\
	\Db{\left(\Mod{\cD_X}\right)} \ar@/^1pc/[r]^{f_+} & \Db{\left(\Mod{\cD_Y}\right)} \ar@/^1pc/[l]^{f^!}
   }
\end{displaystyle}
\]

\begin{rmk}
If one checks the precise definitions of these morphisms, there are some misterious shifts by the dimensions of $X$ and $Y$, as for the case of $\DD_X$. There is nothing deep inside and they are only considered in order to fit nicely in the Riemann-Hilbert correspondence.
\end{rmk}

\begin{defn}
Let $M$ be a $\cD_X$-module and let $F^\bullet M$ be an increasing filtration. It is called a \emph{good filtration} if
\begin{itemize}
	\item $F^k\cD_X \cdot F^p M \subseteq F^{k+p}M$, with equality for $p$ large enough.
	\item Every point of $X$ has a neighbourhood on which $F^pM = 0$ for $p$ small enough (i.e.\ the filtration is locally bounded below).
	\item Each $F^pM$ is a coherent $\cO_X^{an}$-module.
\end{itemize}
\end{defn}

\begin{ex}
The natural filtration $F^\bullet$ of $\cD_X$ by the order is a good filtration.
\end{ex}

As we explain below, there are very good criterions to decide whether a $\cD_X$-module admits or not a good filtration. For the proof of these results, see \cite{Peters-Steenbrink:2008}, Section 13.4.2.

\begin{prop}
Let $X$ be a complex algebraic variety and let $M$ be a $\cD_X$-module.
\begin{itemize}
	\item A filtration $F^\bullet$ of $M$ is a good filtration if and only if $\Gr{F}{\bullet}M$ is coherent as $\Gr{F}{\bullet}\cD_X$-module.
	\item $M$ admits a locally good filtration (i.e.\ a good filtration on an open set) if and only if $M$ is coherent as $\cD_X$-module.
\end{itemize}
\end{prop}

For a while, let us focus on the order filtration $F^\bullet$ of $\cD_X$. Notice that, as $[F^k\cD_X, F^l\cD_X] \subseteq F^{k+l-1}\cD_X$, each of the graded pieces $\Gr{F}{k}\cD_X = F^k\cD_X/F^{k+1}\cD_X$ is a commutative $\cO_X^{an}$-algebra. In particular, if $X$ is smooth, as $\cD_X$ is generated by $\cT X$ and $\cO_X^{an}$, we have that $\Gr{F}{\bullet}\cD_X$ is isomorphic to the symmetric algebra $\Sym_{\cO_X^{an}}\cT X$. On the other hand, there is a natural identification $\Sym_{\cO_X^{an}}\cT X \cong \pi_*\left(\cO_{\cT^*X}\right)$ where $\pi: \cT^*X \to X$ is the holomorphic cotangent bundle.

As an application, if $(M, F^\bullet)$ is a well-filtered $\cD_X$-module, then $\Gr{F}{\bullet}M$ has a coherent $\pi_*\left(\cO_{\cT^*X}\right)$-module structure. In particular, we can consider its annihilator ideal $\Ann \Gr{F}{\bullet}M \subseteq \pi_*\left(\cO_{\cT^*X}\right)$ and to denote by $I(M, F^\bullet)$ the radical ideal of $\Ann \Gr{F}{\bullet}M$. It can be proven (see \cite{Peters-Steenbrink:2008}, Proposition 13.44) that the ideal $I(M,F^\bullet)$ actually does not depend of the chosen good filtration. For this reason, if we have a coherent $\cD_X$-module then, locally, taking any local good filtration $F^\bullet$ we have a well-defined ideal $I(M,F^\bullet)$. Moreover, since they are uniquely defined, they can be glued together along different open sets to give rise to a global ideal that we will denote $I(M)$.

Observe that, for any $x \in X$, this construction gives an ideal $I(M)_x \subseteq \pi_*\left(\cO_{\cT^*X}\right)_x$ so its zero set $Z\left(I(M)_x\right)$ can be seen as a subvariety of $T_x^*X$. Pasting together these zero sets, we obtain a subvariety
$$
\ChV M = \bigcup_{x \in X} Z\left(I(M)_x\right) \subseteq \cT^*X,
$$
called the \emph{characteristic variety} of $M$. For a more precise definition, see \cite{Peters-Steenbrink:2008}, Section 13.4.2.

\begin{ex}
Seen $\cD_X$ as a $\cD_X$-module, we have that $\Gr{F}{\bullet}\cD_X = \pi_*\left(\cO_{\cT^*X}\right)$. In particular, its annihilator is trivial so $\ChV \cD_X = \cT^*X$.
\end{ex}

The character variety is a very special subvariety of the symplectic manifold $\cT^*X$, as shown in the following result whose proof can be found in \cite{Malgrange:1979}. Recall that, given a symplectic manifold $(M, \omega)$ a submanifold $S \subseteq M$ is called involutive if $T_x S \supseteq T_xS^\perp$ for all $x \in M$, where the orthogonality is measured with respect to the symplectic form $\omega$.

\begin{thm}[Bernstein's inequality]
Let $M$ be a coherent $\cD_X$-module. Its characteristic variety $\ChV M \subseteq \cT^*X$ is an involutive submanifold. In particular, $\dim_\RR \ChV M \geq \dim_\RR X$
\end{thm}

\begin{defn}
Let $M$ be a coherent $\cD_X$-module. It is said to be \emph{holonomic} if $\Ch M \subseteq \cT^*X$ is a lagrangian submanifold or, in other words, if $\dim_\RR \ChV M = \dim_\RR X$. 
\end{defn}

\begin{ex}
It can be shown that, for any $\cD_X$-module that is coherent as $\cO_X^{an}$-module, its characteristic variety is precisely the zero section of the cotangent bundle so, in particular, it is holonomic.
\end{ex}

Another important property of $\cD$-modules is the so-called regularity. Let us consider $\Delta^* \subseteq \CC$ the punctured disc and let $M$ be an holonomic $\cD_{\Delta^*}$-module. As $M$ is coherent as a $\cD_X$-module, this structure is equivalent to see $M$ as a complex vector bundle with a holomorphic connection $\nabla: M \to M \otimes \mathbf{\Omega}^1$. We will say that $M$ has a \emph{regular singularity} at $0$ if there exists an extension bundle $\tilde{M}$ of $M$ to the whole disc $\Delta = \Delta^* \cup \left\{0\right\}$ and a logarithmic connection $\tilde{\nabla}$ extending $\nabla$. Recall that a logarithmic connection is a connection such that the entries of the matrix of connection $1$-forms are meromorphic forms with a single pole at $0 \in \Delta$.

More generally, given an irreducible complex curve $C$ and a good compactification $C \subseteq \overline{C}$ (i.e.\ a compactification whose difference is a finite number of points, see \cite{Peters-Steenbrink:2008}, Definition 4.1), an holonomic $\cD_C$ module is said to have \emph{regular singularities} if it has, locally, regular singularities on each of the points of $\overline{C}-C$ as above.

\begin{defn}\label{defn:regular-holonomic-D-mod}
Let $X$ be a complex algebraic variety and let $M$ be an holonomic $\cD_X$-module. We will say that $M$ is \emph{regular holonomic} if, for all holomorphic maps $f: C \to X$ with $C$ a smooth irreducible curve, we have that cohomology of $f^!M \in \Mod{\cD_C}$ has regular singularities.
\end{defn}

As we will see in the following section, $\cD$-modules are strongly related with another object called perverse sheaves, via an extremely important correspondence called the Riemann-Hilbert correspondence. A notion needed for this correspondence is the so-called de Rham complex. 

\begin{defn}
Suppose that $X$ is a smooth complex algebraic variety and let $M^\bullet$ be a complex of $\cD_X$-modules. The \emph{de Rham} complex of $M$, denoted $\DR(M)$ is the complex of (complex) sheaves
$$
	\DR(M) = \left(\mathbf{\Omega}^\bullet \otimes_{\cO_X^{an}} M^\bullet \right)[\dim_\RR X],
$$
where $\mathbf{\Omega}^\bullet$ is the complex of sheaves of holomorphic differential forms and the tensor product is the one of complexes (see \cite{Weibel}, 2.7.1).
\end{defn}

\begin{ex}
\begin{itemize}
	\item Suppose that $M^\bullet = M$ is a single module. In that case, if we set $d = \dim_\RR X$, we have that $\DR(M)^k = \mathbf{\Omega}^{k+d} \otimes_{\cO_X^{an}} M$.
	\item If $\cV$ is a vector bundle with an integrable holomorphic connection $\nabla$, then the de Rham complex $\DR(\cV)$ is precisely the complex of holomorphic forms with values in $\cV$ (shifted), $\mathbf{\Omega}^\bullet(\cV)$ and the diferential is the covariant derivative $d_\nabla: \mathbf{\Omega}^{k}(\cV) \to \mathbf{\Omega}^{k+1}(\cV)$.
\end{itemize}
\end{ex}

\subsubsection*{Perverse sheaves}

Let $X, Y$ be topological spaces with finite cohomological dimension. Given a continuous morphism $f: X \to Y$, as explained in Example \ref{ex:algebra-sheaves}, we have induced maps on the categories of bounded derived complexes of (rational) sheaves
$$
	Rf_*, Rf_!: \Db{\Sh{X}} \to \Db{\Sh{Y}}, \hspace{2cm} f^*: \Db{\Sh{Y}} \to \Db{\Sh{X}}.
$$
Recall that $f_*$ is given by $f_*\cF(U)= \cF(f^{-1}(U))$, for $\cF \in \Sh{X}$ and $U \subseteq Y$ open, and that $f^*$ is right exact so $Rf^*=f^*$. Notice that $Rf_*$ and $f^*$ mutually adjoint. In the same spirit, the adjoint of $Rf_!$ can be defined and it is denoted by $f^!: \Db{\Sh{Y}} \to \Db{\Sh{X}}$ and it is called the \emph{extra-ordinary pullback}.

\begin{rmk}
It is customary to denote the functor $f^*: \Sh{Y} \to \Sh{X}$ by $f^{-1}$ and to reserve the notation $f^*$ for its extension to modules, as for $\cD$-modules. However, in order to simplify notation and to get in touch with the Grothendieck's six operators setting (see Remark \ref{rmk:notation-grothendieck-six}), we will not make this difference here. 
\end{rmk}

Also, as for $\cD$-modules, there is a duality operator called the \emph{Verdier duality} and denoted $\VerD_X$. Its definition is much more involved than for $\cD$-modules so, here, we just sketch it. For a complete definition, see \cite{Dimca:2004}. Let $\cF^\bullet$ be a complex of rational sheaves. There is a special soft resolution of $\cF^\bullet$ called the $c$-Godement resolution, $\cF^\bullet \to {_{c}\cC_{Gd}^\bullet(\cF^\bullet)}$, which is just the truncation of the usual Godement resolution (see \cite{Peters-Steenbrink:2008}, Example 13.2). From this resolution, given an open set $U \subseteq X$, we define $\VerD_X\cF^\bullet(U) = \Hom_\QQ^\bullet(\Gamma_c(U,\, _{c}\cC_{Gd}^\bullet(\cF^\bullet)), \QQ)$, where $\Hom_\QQ^\bullet$ denotes the complex of graded morphisms of chains as in \cite{Peters-Steenbrink:2008}, Section A.2.1. Verdier duality induces a functor $\VerD_X: \Db{\Sh{X}} \to \Db{\Sh{X}}$. As for the case of $\cD$-modules, using this duality we have a commutative diagram
\[
\begin{displaystyle}
   \xymatrix
   {
	\Db{\Sh{X}} \ar@/^1pc/[r]^{Rf_!} \ar@{<->}[dd]_{\VerD_X} & \Db{\Sh{Y}} \ar@/^1pc/[l]^{f^*} \ar@{<->}[dd]^{\VerD_Y} \\
	&\\
	\Db{\Sh{X}} \ar@/^1pc/[r]^{Rf_*} & \Db{\Sh{Y}} \ar@/^1pc/[l]^{f^!}
   }
\end{displaystyle}
\]

Now, we want to consider a special case of rational sheaves. Suppose that $X$ is an analytic space of pure dimension $n$. Recall that a \emph{Whitney complex analytic stratification} consists of a filtration
$$
	X = X_n \supseteq X_{n-1} \supseteq \ldots \supseteq X_0
$$
such that each $X_k$ is a closed analytic subspace of $X$ and the connected components of $X_k - X_{k-1}$ are $k$-dimensional complex manifolds.

\begin{defn}
Let $X$ be a complex algebraic variety, let $R$ be ring and let $\cF$ be a sheaf of $R$-modules on $X$. $\cF$ is called \emph{constructible} if:
\begin{itemize}
	\item There exists a Whitney stratification of $X$ such that the restriction of $\cF$ to each stratum is locally constant.
	\item For all $x \in X$, the stalk $\cF_x$ is a finitely generated $R$-module.
\end{itemize}
A complex of sheaves $K=\cF^\bullet$ is called \emph{cohomologically constructible} if there exists a Whitney stratification such that $H^k(\cF^\bullet)$ is constructible for all $k \in \ZZ$.
\end{defn}

\begin{defn}
Let $K \in \Db{\Sh{X}}$ be a complex of rational sheaves on $X$. Then, $K$ is called a \emph{perverse sheaf} if $K$ is cohomologically constructible with respect to a stratification
$$
X = X_n \supseteq X_{n-1} \supseteq \ldots \supseteq X_0
$$
such that, for all $0 \leq k \leq n$:
\begin{itemize}
	\item For all $x \in X_k - X_{k-1}$ and $j > -k$, we have $H^j(K)_x =0$.
	\item For all $x \in X_k - X_{k-1}$ and $-j > -k$, we have $H^j_c(K)_x = 0$. Here, $H^j_c(K)_x = \lim\limits_{\substack{\longleftarrow\\x \in U}} H_c^k(U, K)$.
\end{itemize}
\end{defn}

The importance of perverse sheaves is the following fact, whose proof can be found in \cite{Peters-Steenbrink:2008}, Lemma 13.22.

\begin{prop}
The full subcategory of $\Db{\Sh{X}}$ of perverse sheaves, denoted $\Perv{X}$, is an abelian category.
\end{prop}

\begin{rmk}
\begin{itemize}
	\item A perverse sheaf is not a sheaf neither perverse. The reason of this terminology is historical.
	\item Given a morphism of analytic spaces $f: X \to Y$, the functors $f^*$ and $f^!$ preserve the condition of being cohomologically constructible. Moreover, if $f$ is a regular morphism of complex algebraic varieties, then $Rf_*$ and $Rf_!$ also preserve cohomologically constructible complexes. However, none of these morphisms preserve, in general, perversity.
	\item So far, we have focused in the case of rational sheaves since they are the most useful for the definition of mixed Hodge modules. However, thay can be defined for any subfield $k \subseteq \CC$ without further differences.
\end{itemize}
\end{rmk}

With all these preliminaries, we are in disposition of stating the Riemann-Hilbert correspondence. This is an extremely important correspondence that shows that $\cD$-modules and perverse sheaves are the two faces of the same coin. The proof of the result can be found in \cite{Borel:1987}. In order to state it properly, let us denote:
\begin{itemize}
	\item $\RegHol{X}$ is the full subcategory of $\Mod{\cD_X}$ of regular holonomic $\cD$-modules (see Definition \ref{defn:regular-holonomic-D-mod}).
	\item $\Drh{X}$ is the full subcategory of $\Db{\left(\Mod{\cD}\right)}$ of bounded complexes of coherent $\cD_X$-modules whose cohomology sheaves are regular holonomic.
	\item Given a subfield $k$ of $\CC$, $\Dcs{X}{k}$ is the full subcategory of $\Db{\Sh{X,k}}$ of bounded complexes of sheaves of $k$-vector spaces which are cohomologically constructibles.
\end{itemize}

\begin{thm}[Riemann-Hilbert correspondence]
Let $X$ be a smooth complex algebraic variety.
\begin{itemize}
	\item The de Rham complex gives rise to an equivalence of categories $\DR{}: \Drh{X} \to \Dcs{X}{\CC}$. Restricting this functor, we also have an equivalence $\DR{}: \RegHol{X} \to \Perv{X; \CC}$.
	\item Let $f: X \to Y$ be a morphism of smooth complex algebraic varieties. The induced maps from $f$ fit in the following commutative diagrams
	\[
\begin{displaystyle}
   \xymatrix
   {
	\Drh{Y} \ar[r]^{\DR{}} \ar[d]_{f^!}& \Dcs{Y}{\CC} \ar[d]^{f^!}\\
	\Drh{X} \ar[r]_{\DR{}} & \Dcs{X}{\CC}
   }
\hspace{1cm}
   \xymatrix
   {
	\Drh{X} \ar[r]^{\DR{}} \ar[d]_{f_+}& \Dcs{X}{\CC} \ar[d]^{Rf_*}\\
	\Drh{Y} \ar[r]_{\DR{}} & \Dcs{Y}{\CC}
   }
\end{displaystyle}
\]
	\item Under the de Rham functor, the Verdier dual and the dual of $\cD$-modules correspond, in the sense that the following diagram commutes
	\[
\begin{displaystyle}
   \xymatrix
   {
	\Drh{X} \ar[r]^{\DR{}} \ar@{<->}[d]_{\DD_X}& \Dcs{X}{\CC} \ar@{<->}[d]^{\VerD_X}\\
	\Drh{X} \ar[r]_{\DR{}} & \Dcs{X}{\CC}
   }
\end{displaystyle}
\]
\end{itemize}
\end{thm}

\begin{rmk}
Since the dual of $\cD$-modules and the Verdier duality commute, we also have analogous diagrams for the pairs $(f_!, Rf_!)$ and $(f^+,f^*)$.
\end{rmk}

\begin{defn}
Let $X$ be a complex algebraic variety. A \emph{filtered $\cD_X$-module with $\QQ$-structure} is a tuple $(M, F^\bullet, K, \alpha)$, where $(M, F^\bullet)$ is a good-filtration of $\cD_X$-modules, $K$ is a perverse sheaf and $\alpha$ is an isomorphism $\DR{(M)} \cong K \otimes_{{\QQ}} \CC$.
\end{defn}

\subsection{Saito's theory of mixed Hodge modules}
\label{sec:saito-mixed-hodge-modules}

Mixed Hodge modules are the most important algebraic tool used in this thesis and they will be ubiquous from now on.  In this section, we will sketch some of the most important steps in the construction of the category of mixed Hodge modules on an algebraic variety by means of $\cD$-modules and perverse sheaves.

\subsubsection*{Nearby and vanishing cycles functors}

A very important tool for the construction Saito's mixed Hodge modules are the so-called nearby and vanishing cycles functors.  The starting point in Saito's construction is to consider the category of pairs of filtered $\cD$-modules with a $\QQ$-structure. However, this category is too wide for being well behaved under direct and inverse images, so we need to extract an appropiate subcategory, which we will call the category of (pure) Hodge modules. For that purpose, the nearby and vanishing cycles functors will be the key. For this reason, in this section we will sketch the construction of these functors. For a more complete introduction, see \cite{Wu:2017}.

Let $X$ be a complex manifold, let $\Delta \subseteq \CC$ be the unit disc and suppose that we have an holomorphic functor $f: X \to \Delta$ such that $f$ is submersive on $\Delta^* = \Delta - \left\{0\right\}$ (these are the so-called $1$-parameter degenerations). In that case, the fiber over $0$, $X_0 = f^{-1}(0)$, is the singular fiber with an inclusion $i: X_0 \hookrightarrow X$. Let us consider the universal covering $e: \HH \to \Delta^*$, where $\HH$ us the upper half plane and $e(z) = e^{2\pi i z}$. With these spaces, we build the pullback diagram
\[
\begin{displaystyle}
   \xymatrix
   {
	X_\infty \ar[r]^{k}\ar[d] & X\ar[d]^{f} \\
	\HH \ar[r]_{e} & \Delta^*
   }
\end{displaystyle}
\]
The space $X_\infty$ can be understood as the generic fiber of $f$, since all the fibers are homotopic.

\begin{defn}
Let $f: X \to \Delta$ be a $1$-parameter degeneration. The \emph{nearby cycles} functor, $\psi_f: \Dcs{X}{\CC} \to \Dcs{X_0}{\CC}$ is the functor that, for any $K \in \Dcs{X}{\CC}$, gives
$$
	\psi_f\,K = i^{*}Rk_*\left(k^{*}K\right).
$$
Moreover, the \emph{vanishing cycles} functor, $\phi_f: \Dcs{X}{\CC} \to \Dcs{X_0}{\CC}$ is the functor that, to any $K \in \Dcs{X}{\CC}$, assigns the cone of the morphism $i^{*}K \to \psi_f\,K$ (see \cite{Weibel}). 
\end{defn}

\begin{rmk}
\begin{itemize}
	\item The morphism $i^{*}K \to \psi_f\,K$ is constructed as follows. By adjointness of $k^{*}$ and $Rk_*$, we have a morphism $K \to Rk_*(k^{*}\,K)$ and the desired map is obtained by applying $i^{*}$.
	\item By construction, there is a map $\can: \psi_f\,K \to \phi_f\,K$. Another map in the other way around, $\var: \phi_f\,K \to \psi_f\,K$, can also be constructed.
	\item It can be proven (\cite{Dimca:2004}, Proposition 4.2.2) that, for all $x \in X_0$, there is a natural isomorphism $(H^k\psi_f\,K)_x \cong H^k(F_x, K)$, where $F_x = f^{-1}(t) \cap B_\epsilon(x)$ for $\epsilon > 0$ small enough is the local Milnor fiber around $x$.
	\item If $K$ is a (complex) perverse sheaf, it can be proven \cite{Schnell:2014} that $\psi_f\,K[-1]$ and $\phi_f\,K[-1]$ are both perverse sheaves. It is customary to denote them $^{p}\psi_f\,K$ and $^{p}\phi_f\,K$ respectively.
\end{itemize}
\end{rmk}

Suppose that $K \in \Perv{X; \CC}$ with vanishing cycle complex $^{p}\psi_f\,K$. Inside a small ball around $x \in X_0$, $R^kf_*\,K$ is a local system and, thus, we have a monodromy action $\pi_1(\Delta^*) \to \Aut H^k(F_x, K)$ (see \cite{Peters-Steenbrink:2008}, Appendix B.3). In particular, the positive oriented loop around $0$ induces an automorphism $H^k(F_x, K)$ and, by the third item of the previous remark, also an automorphism of $^{p}\psi_f\,K$, denoted $T: {^{p}\psi_f}\,K \to {^{p}\psi_f}\,K$ and called the \emph{monodromy operator}.

As the category of perverse sheaves is an abelian category, we can decompose in generalized eigenspaces of $T$
$$
	^{p}\psi_f\,K = \bigoplus_{\lambda \in \CC^*} {^{p}\psi_{f, \lambda}}\,K, 
$$
where $\psi_{f, \lambda}\,K = \ker\left(T-\lambda I\right)^l$ for $l$ large enough. In particular, in the unipotent part $\psi_{f, 1}\,K = \bigcup_{l} \ker\left(T-I\right)^l$ we can consider the \emph{unipotent monodromy operator}
$$
	N = -\frac{1}{2 \pi i}\,\log T = \frac{1}{2 \pi i} \sum_{l = 1}^\infty \frac{(I-T)^k}{k}.
$$
It can be proven that, actually, $N = \var \circ \can$.

By the Riemann-Hilbert correspondence, these morphisms must have a counterpart in the side of $\cD$-modules. In order to state it, given a $\cD$-module $M$, we consider a special filtration $\cV^\bullet M$ of $M$ called the \emph{Kashiwara-Malgrange} filtration. We will not construct it here but a complete description can be checked in \cite{Peters-Steenbrink:2008}, Section 14.2.2 or \cite{Schnell:2014}, Section 8. However, roughtly speaking, suppose that we have local coordinates $(t, z_1, \ldots, z_{n-1})$ on $X$ such that $X_0$ corresponds to $t=0$. Take $\partial_t$ to be the vector field associated to the coordinate $t$ so that $[\partial_t, t]=1$. In that case, the filtration $\cV^\bullet M$ is the unique increasing filtration such that the action of $t$ on $\cV^\bullet M$ descends the degree by one (with equality for small enough degree), the action of $\partial_t$ ascends it by one and the eigenvalue of the so-called Fuchs field $t \partial_t$ on $\Gr{\cV}{k}$ have real part in $(k-1,k]$.

In that case, given such a Kashiwara-Malgrange filtration $\cV^\bullet$ on $M$ and $\alpha \in \CC$, we define $M_\alpha$ as the generalized eigenspaces of $t\partial_t$ of eigenvalue $\alpha$
 on $\Gr{\cV}{\lceil \re{\alpha} \rceil}M$, where $\lceil \re{\alpha} \rceil$ denotes the smallest integer greatest or equal than $\re{\alpha}$. Then, setting $\lambda = e^{2 \pi i \alpha}$, we have isomorphisms
$$
	\DR{\left(M_\alpha\right)} \cong \,^{p}\psi_{f, \lambda}\DR{\left(M\right)} \hspace{0.2cm}(-1 \leq \re{\alpha} < 0),
\hspace{0.7cm}
	\DR{\left(M_\alpha\right)} \cong \,^{p}\phi_{f, \lambda}\DR{\left(M\right)} \hspace{0.2cm}(-1 < \re{\alpha} \leq 0).
$$
Under this correspondence, the morphisms $\can$ and $\var$ are assigned to the actions of $\partial_t: M_{-1} \to M_0$ and $t: M_0 \to M_{-1}$. In particular, the unipotent monodromy operator $N$ corresponds to $\partial_t t = t \partial_t + 1$.

\begin{rmk}
To be precise, in the general case of a function $f: X \to \CC$, we have to consider the graph of $f$, $\Gamma f: X \to X \times \CC$. In that situation, we now have the local coordinate $t$ as desired on $X \times \CC$ and we have to substitute $M$ by $(\Gamma f)_+M$ (see \cite{Schnell:2014}, Section 8). However, for simplicity, we will omit this technical point. 
\end{rmk}

Now, suppose that $(M, F^\bullet M, K, \alpha)$ is a $\cD$-module with $\QQ$-structure and suppose that all the eigenvalues of $^{p}\psi_f\,K$ are roots of unit. This implies that the eigenvalues of $t \partial_t$ are rational numbers. In that situation, from the Kashiwara-Malgrange filtration, $\cV^\bullet M$, we can define a rational filtration $V^\bullet M$ called the \emph{V}-filtration. For $\alpha \in \QQ$, we consider the direct sum
$$
	\bigoplus_{\lceil \alpha\rceil - 1 < \beta \leq \alpha} M_\beta \subseteq \Gr{\cV}{\lceil \alpha\rceil}M.
$$
Then, we set $V^\alpha M \subseteq \cV^{\lceil \alpha \rceil}M$ to be the preimage of such a sum under the projection $\cV^{\lceil \alpha \rceil}M \to \Gr{\cV}{\lceil \alpha\rceil}M$.

\begin{defn}\label{defn:quasi-unipotent}
Let $X$ be a complex algebraic variety and let $(M, F^\bullet M, K, \alpha)$ be a filtered $\cD_X$-module with $\QQ$-structure and with $V$-filtration $V^\bullet M$. Given an holomorphic function, $f:X \to \CC$, we will say that it is \emph{quasi-unipotent along $f=0$} if:
\begin{itemize}
	\item All the eigenvalue of $^{p}\psi_f\,K$ are roots of unit.
	\item $t: F^pV^\alpha M \to F^pV^{\alpha -1}M$ is surjective for $\alpha <0$.
	\item $\partial_t: F^p\Gr{V}{\alpha}M \to F^{p+1}\Gr{V}{\alpha+1}M$ is surjective for $\alpha > 1$.
\end{itemize}
On the other hand, it is said to be \emph{regular along $f=0$} if $F^\bullet \Gr{V}{\alpha}M$ is a good filtration for every $-1 \leq \alpha \leq 0$.
\end{defn}

\subsubsection*{Definition of pure Hodge modules}

With the previous notions at hand, we are in disposition of defining pure Hodge modules of weight $w \in \ZZ$ on a complex smooth algebraic variety $X$.
The basic objects of this category are regular holonomic filtered $\cD$-modules with $\QQ$-structure, $(M, F^\bullet M, K, \alpha)$. Moreover, we will require that $(M, F^\bullet)$ is quasi-unipotent and regular along $f = 0$ for every locally defined holomorphic function $f : U \subseteq X \to \CC$ and that it admits a decomposition by strict support. For short, we will denote such a triple just by $M$.

The category of \emph{pure Hodge modules of weight $w$ on $X$}, $\HM{w}{X}$, is a full subcategory of the previous objects. It is defined recursively from two auxiliar categories 
$$
	\textrm{HM}^w_Z(X)  = \left\{\txt{Hodge modules of weight\\$w$ with support on $Z \subseteq X$}\right\}, \hspace{1.5cm} \textrm{HM}^w_{< d}(X) = \bigoplus_{\substack{Z \subseteq X \\ \dim Z < d}} \textrm{HM}^w_Z(X).  
$$
The base case is $Z = \left\{x\right\}$ for $x \in X$ with inclusion $i: \left\{x\right\} \hookrightarrow X$. In that case, we set $\textrm{HM}^w_{\left\{x\right\}}(X) = i_*\PHS{w}{\QQ}$. For the induction case, if $Z \subseteq X$ is a subvariety of dimension $d$, we say that $M \in \textrm{HM}^w_Z(X)$ if and only if, for every holomorphic function $f: U \subseteq X \to \CC$ that does not vanish identically on $Z \cap U$, we have that
$$
	\Gr{k}{W(N)}\phi_fM \in \textrm{HM}^{w + k -1}_{< d}(X).
$$
Here, $W(N)$ is the monodromy filtration associated to the unipotent operator $N$ on the nearby cycles (see Remark \ref{rmk:monodromy-filtration}). From these categories, we define
$$
	\HM{w}{X} = \bigoplus_{Z \subseteq X} \textrm{HM}^w_Z(X).
$$

\begin{rmk}
Analogously to the case of variations of Hodge structures (see Section \ref{sec:monodromy-as-mhm}), it can be defined a notion of polarization of a pure Hodge module. Roughly speaking, on a complex algebraic variety $X$ of complex dimension $n$, a polarization of a pure Hodge modules $(M, F^\bullet M, K, \alpha)$ of weight $w$ is a non-degenerate morphism $K(w) \to R\cH om(K, \underline{\QQ}_X(n))[2n]$ that is compatible with the filtrations and the strict support decomposition and recovers the usual notion of polarization of Hodge structure for $\dim\,Z=0$.
\end{rmk}

\subsubsection*{Definition of mixed Hodge modules}

With all the background developed in the previous section, in this section we finally define mixed Hodge modules. The idea of these modules is that, as for mixed Hodge structures, the pure case is not good enough since it loses some important properties. For that reason, we have to consider a refinement of the definition of the pure case. However, in the case of mixed Hodge modules, the definition is extremely involved by its own. In particular, it is really hard to decide whether an object is a mixed Hodge module or not. Nevertheless, at the end of the day, we will find a very well-behaved abelian category.

The threatening aspect of this definition should not scare us. The accurate point is that this category can (and should) be used as a black box with some opperations defined. Actually, in Section \ref{sec:monodromy-as-mhm} we will show that some very tangible objects can be automatically interpreted as mixed Hodge modules and this will give us a large amount of elements of this category. In this way, we will never need to check the very definition of mixed Hodge modules.

For simplicity, along this section we will suppose that $X$ is a smooth complex algebraic variety. The singular case can also be treated similarly by embedding the singular variety into smooth manifolds. As in the pure case, the definition of mixed Hodge modules is recursive by nature. The basic building block is the category $\HMW{X}$ of (graded polarizable) \emph{weakly mixed Hodge modules}. It is constructed in analogy with the case of mixed Hodge structures as the full subcategory of pairs $(M, W^\bullet M)$ where $M$ is a filtered regular holonomic $\cD_X$-module with an increasing filtration $W^\bullet$ such that $\Gr{W}{k}M \in \HM{k}{X}$ for all $k \in \ZZ$ and it is polarizable.

It can be proven that, if $(M, W^\bullet M) \in \HMW{X}$, then $M$ is quasi-unipotent and regular along $f$ for any non-constant holomorphic function $f: X \to \CC$. However, this category is not good enough for our purposes and we need to refine it. In the mixed case, the key propery is the so-called admissibility.

\begin{defn}
Let $(M, W^\bullet M) \in \HMW{X}$ be a weakly mixed Hodge module and let $f: X \to \CC$ be a non-constant holomorphic function. Define the limit filtrations for the nearby and vanishing functors $L^k\,\psi_f(M) = \psi_f\left(W^{k+1} M\right)$ and $L^k\phi_{f,1}(M) = \phi_{f,1}(W^k M)$. We will say that it is \emph{admissible along $f=0$} if:
\begin{itemize}
	\item The weight filtration $W^\bullet$, the natural filtration $F^\bullet$ and the $V$-filtration $V^\bullet$ of $M$ are compatible in the sense that, when computing the graded complex, the order of the filtrations is indiferent.
	\item The relative monodromy filtrations $W^\bullet(N, L^\bullet \psi_f(M))$ and $W^\bullet(N, L^\bullet \phi_{f,1}(M))$ for the action of the unipotent monodromy operator $N$ exists.
\end{itemize}
\end{defn}

\begin{rmk}\label{rmk:monodromy-filtration}
The relative monodromy filtration is a special type of filtration constructed from simpler data. The idea is that, given a filtered object and a nilpotent filtered morphism, the relative monodromy filtration is the unique filtration that respects these data in a controlled way (for a precise definition, see \cite{Saito:1990} Section 1). In this case, the role of the filtration of played by $L^\bullet$ and the morphism is the unipotent monodromy $N$. In the case that $L^\bullet$ is trivial, we recover the notion of an (absolute) monodromy filtration, as for pure Hodge modules.
\end{rmk}

With this definition at hand, we can finally define what a mixed Hodge module is. Let $(M, W^\bullet M)$ be a graded polarizable weakly mixed Hodge module.
We will say that it is a (graded polarizable) mixed Hodge module if, for every non-constant holomorphic function, $f: U \to \CC$, defined on an open set $U \subseteq X$ such that $f^{-1}(0)$ does not contain any irreducible component of $U \cap \supp M$, it is admisible along $f = 0$ and both $(\psi_f(M), W^\bullet(N, L^\bullet \psi_f(M)))$ and $(\phi_{f,1}(M), W^\bullet(N, L^\bullet \phi_{f,1}(M)))$ are mixed Hodge modules. Observe that, as in the pure case, this definition makes sense since the nearby and the vanishing cycles have support of strictly smaller dimension than $X$.

\begin{defn} Let $X$ be a complex algebraic variety.
The full subcategory of $\HMW{X}$ of (graded polarizable) \emph{mixed Hodge modules on $X$} will be denoted by $\MHM{X}$.
\end{defn}

\begin{rmk}
The notation $\MHM{X}$ is not standard at all. In the papers of Saito, such a category is usually denoted by $\textrm{HMH}^p(X)$. However, due to its omnipresence along this thesis, we will shorten it.
\end{rmk}

The category $\MHM{X}$ is an abelian category, as proved in \cite{Saito:1986} and \cite{Saito:1990}. By construction, we have functors
$$
	\rat_X: \MHM{X} \to \Perv{X}, \hspace{1cm} \dmod_X: \MHM{X} \to \RegHol{X},
$$
that just project to the underlying perverse sheaf and $\cD_X$-module, respectively. They can be extended to the (bounded) derived categories as functors
$$
	\rat_X: \Db{\MHM{X}} \to \Dcs{X}{\QQ}, \hspace{1cm} \dmod_X: \Db{\MHM{X}} \to \Drh{X}.
$$
Moreover, given a regular morphism $f: X \to Y$, we can define the functors
$$
	f_*, f_!: \Db{\MHM{X}} \to \Db{\MHM{Y}}, \hspace{1cm} f^*, f^!: \Db{\MHM{Y}} \to \Db{\MHM{X}},
$$
by lifing the corresponding functors $Rf_*, Rf_!, f^*$ and $f^!$ on the level of constructible sheaves (and their counterparts on $\cD$-modules via Riemann-Hilbert correspondence). The fact that they send mixed Hodge modules into mixed Hodge modules is not trivial at all and can be checked in \cite{Saito:1990}.

Finally, the tensor and external product of constructibles complexes (and complexes of $\cD$-modules) also lift to give bifunctors
$$
	\otimes: \Db{\MHM{X}} \times \Db{\MHM{X}} \to \Db{\MHM{X}} \hspace{1cm} \boxtimes: \Db{\MHM{X}} \times \Db{\MHM{Y}} \to \Db{\MHM{X \times Y}}.
$$

\begin{rmk}
As in the case of constructible complexes, the external product can be defined in terms of the usual tensor product by
	$$
		M^\bullet \boxtimes N^\bullet = \pi_1^* M^\bullet \otimes \pi_2^* N^\bullet,
	$$
for $M^\bullet \in \Db{\MHM{X}}$, $N^\bullet \in \Db{\MHM{Y}}$ and $\pi_1: X \times Y \to X$, $\pi_2: X \times Y \to Y$ the corresponding projections. See also Section \ref{sec:C-algebras}.
\end{rmk}

A very important feature of these induced functors is that they behave in a functorial way, as the following result shows. The proof of this claim is a compendium of Section 4.2 of \cite{Schurmann:2011} and Proposition 4.3.2 and Section 4.4 (in particular 4.4.3) of \cite{Saito:1990}.

\begin{thm}[Saito]\label{thm:funct-induced-functors-mhm}
\begin{itemize}
	\item Let $f: X \to Y$ be a regular morphism. The morphisms $f_!, f_*: \Db{\MHM{X}} \to \Db{\MHM{Y}}$ and $f^!: \Db{\MHM{Y}} \to \Db{\MHM{X}}$ commute with the external product. Moreover, the morphism $f^*: \Db{\MHM{Y}} \to \Db{\MHM{X}}$ commutes with the tensor product.
	\item The induced functors commute with composition. More explicitly, let $f: X \to Y$ and $g: Y \to Z$ regular morphisms of complex algebraic varieties, then
$$
	(g \circ f)_* = g_* \circ f_* \hspace{1cm} (g \circ f)_! = g_! \circ f_! \hspace{1cm} (g \circ f)^* = f^* \circ g^* \hspace{1cm} (g \circ f)^! = f^! \circ g^!
$$
	\item Suppose that we have a cartesian square of complex algebraic varieties (i.e.\ a pullback diagram in $\CVar$)
\[
\begin{displaystyle}
   \xymatrix
   {
   		W \ar[r]^{g'} \ar[d]_{f'} & X \ar[d]^{f} \\
   		Y \ar[r]_{g} & Z
   }
\end{displaystyle}   
\]
Then we have a natural isomorphism of functors $g^* \circ f_! \cong f'_! \circ (g')^*$.
	\item For a singleton we have a natural isomorphism of categories $\MHM{\star} = \MHS{\QQ}$.
	\item Given a complex algebraic variety $X$, let us denote the unit of the monoidal structure on $\MHM{X}$ by $\underline{\QQ}_X$. Then, if $c_X: X \to \star$ is the projection onto a point, under the previous isomorphism we have that, for all $k$, $H^k\left((c_X)_!\underline{\QQ}_X\right)$ is the natural mixed Hodge structure on $H_c^k(X;\QQ)$. Analogously, $H^k\left((c_X)_*\underline{\QQ}_X\right) = H^k(X;\QQ)$.
\end{itemize}
\end{thm}

With this theorem at hand, let us define the functor $M: \CVar \to \Cat^{op} \times \Cat$ that, for $X \in \CVar$, assigns $M(X) = (\Db{\MHM{X}}, \Db{\MHM{X}})$, and, for a regular morphism $f: X \to Y$, it assigns $M(f) = (f^*, f_!)$. By the previous theorem, the functor $M$ satisfies the hypothesis of Proposition \ref{prop:k-theory-algebra}. This proves the following result, which is the main theorem of this chapter.

\begin{thm}\label{thm:mixed-hodge-modules-cvar-alg}
The Grothendieck rings of mixed Hodge modules, $\KM{}$, forms a $\CVar$-algebra that assigns $\KM{X}$ for $X \in \CVar$, and pullback $f^*$ and pushout $f_!$ for a regular morphism $f:X \to Y$.
\end{thm}

\begin{rmk}
\begin{itemize}
	\item To be precise, the $\CVar$-algebra assigns, to $X \in \CVar$, the ring $\K{\left(\Db{\MHM{X}}\right)}$. However, as mentioned in \ref{sec:grothendieck-groups}, this ring is naturally isomorphic to $\KM{X}$.
	\item For the case of a singleton $\star \in \CVar$, this $\CVar$-algebra assigns $\KM{\star} = \K{\MHS{\QQ}} = \K{\MHSq}$, where the first equality follows from the equality $\MHM{\star} = \MHS{\QQ}$ and the second one from Section \ref{sec:mixed-hodge-structures}.
	\item Let $\underline{\QQ}_X \in \KM{X}$ be the unit of the ring structure on $\KM{X}$. By Theorem \ref{thm:funct-induced-functors-mhm}, the measure of $\underline{\QQ}_X$ is $\mu(\underline{\QQ}_X) = (c_X)_!\underline{\QQ}_X = \coh{X;\QQ}$, where $\coh{X;\QQ} \in \K{\MHSq}$ denotes the $K$-theory image of the (mixed) Hodge structure on $H_c^\bullet(X; \QQ)$.
\end{itemize}
\end{rmk}

Another very useful property of mixed Hodge modules is that they behave well with respect to stratification of the underlying space.

\begin{prop}\label{prop:strat-mixed-hodge-modules}
Let $X$ be a complex algebraic variety and suppose that we have a decomposition $X = Y \sqcup U$ where $i: Y \hookrightarrow X$ is a closed subvariety and $j: U \hookrightarrow X$ is a Zariski open subset. Then, the ring homomorphism
	$$
		i_! + j_!: \KM{Y} \oplus \KM{U} \to \KM{X}
	$$
is an isomorphism.
\begin{proof}
It is easy to check that $i_! + j_!$ is injective since it has no kernel. This is due to the fact that the supports of elements of $\KM{Y}$ and $\KM{U}$ do not intersect so the image would vanish if and only if each of the initial pieces vanish.

On the other hand, by Saito's result (4.4.1) of \cite{Saito:1990}, there is a distinguised triangle in $\Db{\MHM{X}}$
$$
	j_!j^* \longrightarrow 1_{\MHM{X}} \longrightarrow i_!i^* \stackrel{[1]}{\longrightarrow} \cdots.
$$
Hence, at the level of Grothendieck rings, we have an equality of morphisms $1_{\KM{X}} = i_!i^* + j_!j^*$. In particular, it shows that $i_! + j_!$ is surjective since, for every $M \in \KM{X}$, we have $M = i_!(i^*M) + j_!(j^*M)$ with $i^*M \in \KM{Y}$ and $j^*M \in \KM{U}$. This finishe the proof.
\end{proof}
\end{prop}

\begin{rmk}
Given a closed subvariety $Y \subseteq X$ and $M \in \KM{X}$, for short we will denote $M|_{Y} = i^*M$. Observe that, by the previous proposition, if $X = Y \sqcup U$ we have that $M = i_! \left(M|_Y\right) + j_!\left(M|_U\right)$ for all $M \in \KM{X}$.
\end{rmk}

\section{Monodromy as mixed Hodge module}
\label{sec:monodromy-as-mhm}

Despite the abstract definition of mixed Hodge modules, it is still possible to identify some especific subcategories. That will be very useful since, as we will show, this allows us to formulate monodromy actions as mixed Hodge modules. This is the key for putting in context the so-called Hodge monodromy representation.

\begin{defn}
Let $B$ a complex manifold. A \emph{variation of mixed Hodge structures} is a triple $(\mathbb{V}, W^\bullet, \cF^\bullet)$ composed by:
\begin{itemize}
	\item A local system $\mathbb{V}$ of $\QQ$-vector spaces on $B$ i.e.\ a locally constant sheaf.
	\item A finite increasing filtration $W^\bullet$ of $\mathbb{V}$ by local subsystems, called the \emph{weight filtration}.
	\item A finite decreasing filtration $\cF^\bullet$ of the holomorphic vector bundle $\cV = \mathbb{V} \otimes_{\underline{\QQ}_B} \cO_B^{an}$ by holomorphic subbundles, called the \emph{Hodge filtration}.
\end{itemize}
such that:
\begin{enumerate}
	\item For any $b \in B$, the stalks $W^\bullet_b, \cF^\bullet_b$ induce a mixed Hodge structure on $\mathbb{V}_b$.
	\item (Griffiths' transversality condition) If $\nabla: \cV \to \cV \otimes_{\cO_B^{an}} \mathbf{\Omega}^1_B$ is the induced connection, then, for any $p$
	$$
		\nabla(\cF^p) \subseteq \cF^{p-1} \otimes_{\cO_B^{an}} \mathbf{\Omega}_B^1
	$$
\end{enumerate}
\end{defn}

\begin{rmk}
\label{rmk:variations-hodge-struct}
	\begin{itemize}
		\item If the weight filtration $W^\bullet$ is constant and with a single step, the variation is called a variation of (pure) Hodge structures. In that case, the jumping degree is called the weight of the variation.
		\item The induced connection $\nabla$ on $\cV$ is the so-called \emph{Gauss-Manin connection} which is an integrable holomorphic connection such that $\mathbb{V} = \ker \nabla$. For a detailed description see \cite{Peters-Steenbrink:2008}, 10.24.
		\item The Griffiths' transversality condition appears naturally in the following context. Let $f: X \to B$ be a proper morphism between complex manifolds of maximal rank and such that $X$ is bimeromorphic to a compact K\"ahler manifold. Then, if $\underline{\QQ}_X$ denotes the constant sheaf on $X$, then the direct images sheaves $R^kf_*\underline{\QQ}_X$ have a natural structure of variations of pure Hodge structures such that, for any $b \in B$, the Hodge structure on the stalk over $b$ corresponds to the usual pure Hodge structure on $(R^kf_*\underline{\QQ}_X)_b = H^k(f^{-1}(b); \QQ)$. In this case, the Gauss-Manin connection automatically satisfies Griffiths' transversality condition (see \cite{Peters-Steenbrink:2008}, Corollary 10.31).
		\item In general, let $f: X \to B$ be a regular morphism between complex algebraic varieties with $B$ smooth which is locally trivial in the analytic topology. We will say that $f$ is a \emph{nice fibration}. In this case, the direct images sheaves $R^kf_!\underline{\QQ}_X$ have a natural structure of variations of mixed Hodge structures such that, for any $b \in B$, we have an isomorphism $(R^kf_!\underline{\QQ}_X)_b \cong H_c^k(f^{-1}(b); \QQ)$ as mixed Hodge structures (see \cite{Schurmann:2011}, Example 3.11). This construction is the so-called \emph{geometric variation of mixed Hodge structures}. Analogous considerations hold for $R^kf_*\underline{\QQ}_X$, whose stalks are $(R^kf_*\underline{\QQ}_X)_b \cong H^k(f^{-1}(b); \QQ)$ as mixed Hodge structures.
	\end{itemize}
\end{rmk}

\begin{defn}
A \emph{polarization} of a variation $\mathbb{V}$ of pure Hodge structures of weight $k$ over $B$ is a morphism $Q: \mathbb{V} \otimes \mathbb{V} \to \underline{\QQ(-k)}_B$ that induces a polarization of pure Hodge structures on each stalk. A variation of mixed Hodge structures is said to be \emph{graded-polarizable} if the induced variations of pure Hodge structures $\Gr{k}{W} \mathbb{V}$ are polarizable.
\end{defn}

Now, let $B$ a smooth irreducible algebraic variety of dimension $d$ and suppose that, seen as a complex manifold, we have a variation of mixed Hodge structures $(\mathbb{V}, W^\bullet, \cF^\bullet)$ on $B$. Using the Gauss-Manin connection, $\nabla$, induced by the local system $\mathbb{V}$, $\cV$ becomes a regular holonomic $\cD_B$-module and Griffiths's transversality condicion says that $\cF^\bullet$ is a good filtration. Moreover, taking the shifted de Rham complex then, by the Riemann-Hilbert correspondence (see Section \ref{sec:riemann-hilbert}), $\DR{(\cV)}$ is a perverse sheaf on $B$ with an isomorphism $\alpha: \DR{(\cV)} \to \mathbb{V}[d] \otimes_\QQ \CC$, where $d = \dim_\CC B$. Shifting $W^\bullet$ to $W'^k(\mathbb{V}[d]) = W^{k-d}\mathbb{V}[d]$ we build a tuple $M_\mathbb{V}= ((\cV, F^\bullet, \mathbb{V}[d], \alpha), W'^\bullet)$ as in the definition of a mixed Hodge module. 

In general, this tuple $M_{\mathbb{V}}$ does not define a mixed Hodge module if we do not requiere more conditions to $\mathbb{V}$. First, suppose that there exists a compact algebraic variety $\overline{B} \supseteq B$ such that $D = \overline{B}-B$ has smooth irreducible components and looks, locally, like the crossing of coordinate hyperplanes. This is sometimes called a \emph{good compactification of $B$} and $D$ is said to be a \emph{simple normal crosssing divisor}. The existence of such good compactifications is shown in \cite{Hironaka:1964}. 

In that situation, a local system $\mathbb{V}$ is said to have quasi-unipotent monodromy at infinity if the monodromy of the loops around $D$ is quasi-unipotent (i.e.\ some power is unipotent). Moreover, there is a set of technical conditions on a variation of mixed Hodge structures, called \emph{admissibility}, as described in \cite{Steenbrink-Zucker:1985} and \cite{Kashiwara:1986}. We will not need an explicit formulation of these conditions but they should be thought of as the incarnation of the corresponding conditions for mixed Hodge modules in $M_{\mathbb{V}}$.

\begin{rmk}
To be honest, historically the process went in the other way around. Actually, as Saito explained, he came up the the right definition of admissibility by generalizing the concept for variations of Hodge structures.
\end{rmk}

\begin{defn}
A variation of mixed Hodge structures is said to be \emph{good} if it is admissible and it has quasi-unipotent monodromy at infinity. The full subcategory of graded polarizable good variations of mixed Hodge structures on $B$ is denoted by $\goodVMHS{B}$, which is an abelian category stable under pullbacks and tensor products.
\end{defn}

\begin{ex}
In general, to check whether a variation of mixed Hodge structures is good may be quite hard. However, a large class of such a variations are automatically good, as the following results show:
\begin{itemize}
	\item (Schmid). A polarizable variation of pure Hodge structures is admisible (see \cite{Schmid:1973}). Hence, it is good if and only if it has quasi-unipotent monodromy at infinity.
	\item (Steenbrink-Zucker). Geometric variations of mixed Hodge structures (i.e.\ those coming from nice morphisms, see Remark \ref{rmk:variations-hodge-struct}, fourth item) are good. For a proof, see \cite{Steenbrink-Zucker:1985} and \cite{Kashiwara:1986}.
\end{itemize}
\end{ex}

%Recall that a variation of mixed Hodge structures $\mathbb{V}$ on a $d$-dimensional smooth variety $B$ is said to be good if it satisfies some technical conditions called \emph{admissibility} as described in [REF Steenbrink, Zucker, citado en Schurmann] and if it has quasi-unipotent monodromy at infinity with respect to a good compactification of $B$. We will denote the full subcategory of good variations of mixed Hodge structures on $B$ by $\goodVMHS{B}$, which is an abelian category stable under pullbacks and tensor products.

\begin{thm}[Saito]\label{thm:saito-smooth-mhm}
Let $B$ a complex algebraic variety. If $\mathbb{V}$ is a good variation of mixed Hodge structures on $B$, then $M_{\mathbb{V}}$ is a mixed Hodge module. Moreover, this correspondence is an isomorphism of categories
$$
	\goodVMHS{B} \stackrel{\cong}{\to} \MHM{B}^{sm},
$$
where $\MHM{B}^{sm} \subseteq \MHM{B}$ is the full subcategory of smooth mixed Hodge modules, that is, mixed Hodge modules $M \in \MHM{B}$ such that $\rat_B\,M[-\dim B]$ is a local system.
\end{thm}

\begin{rmk}
Actually, this theorem also identifies the unit of the monoidal structure on $\MHM{B}$: it is precisely the mixed Hodge module associated to the trivial variation of pure Hodge modules of weight zero, $\underline{\QQ}_B$, on $B$ that has the constant $\QQ$ sheaf as local system and a single step Hodge filtration jumping at $p = 0$. In order to get in touch with the notation of \cite{LMN}, we will also denote this unit as $T_B$ (or just $T$ if $B$ is understood) when we want to emphasize its monodromic nature. Finally, if we want to focus on the monoidal structure of $\MHM{B}$, such a unit will be simply denoted by $1 \in \MHM{B}$.
\end{rmk}

\begin{cor}\label{cor:hodge-monodromy-as-monodromy}
Let $f: X \to B$ be a nice morphism and denote the associated variations of mixed Hodge structures as $\mathbb{V}_f^k = R^kf_! \underline{\QQ}_X$. Then, in $\K{\MHM{B}}$, we have the equality
$$
	 f_!\,\underline{\QQ}_X = \sum_k (-1)^k M_{\mathbb{V}_f^k}
$$
In particular, if $f$ has trivial monodromy with fiber $F$, then $
	f_!\,\underline{\QQ}_X = \coh{F;\QQ}\,\underline{\QQ}_B$.
\begin{proof}
Let $\underline{\QQ}_X$ be the trivial variation of mixed Hodge modules on $X$. By definition, in the Grothendieck ring of variations of mixed Hodge structures on $B$ we have that
$$
	\K{f_!}\,\underline{\QQ}_X = \chi\left(Rf_!\,\underline{\QQ}_X \right)= \sum_k (-1)^k\,\mathbb{V}_f^k.
$$
In particular, since all the $\mathbb{V}_f^k$ are good, $\K{f}_!\,\underline{\QQ}_X \in \K{\left(\goodVMHS{B}\right)}$ and the first part follows by passing to $\K{\MHM{B}}$ via Theorem \ref{thm:saito-smooth-mhm}. For the second part, just observe that, if $f$ has trivial monodromy, then $\mathbb{V}_f^k = H_c^k(F;\QQ)\,\underline{\QQ}_B$ as variations of mixed Hodge structures and, thus, as mixed Hodge modules when passing to $\MHM{B}$ via Theorem \ref{thm:saito-smooth-mhm}.
\end{proof}
\end{cor}

\begin{cor}[Logares-Mu\~noz-Newstead]\label{cor:mhs-trivial-monodromy}
If $f:X \to B$ is a nice morphism with trivial monodromy and fiber $F$, then
$$
	[H_c^\bullet(X;\QQ)] = [H_c^\bullet(F;\QQ)] [H_c^\bullet(B;\QQ)].
$$
In particular, $\DelHod{X} = \DelHod{F}\DelHod{B}$.
\begin{proof}
Just apply $(c_B)_!$ to both sides of the second part of Corollary \ref{cor:hodge-monodromy-as-monodromy}, where $c_B: B \to \star$ in the projection onto a point, and use Theorem \ref{thm:funct-induced-functors-mhm}. 
\end{proof}
\end{cor}

\begin{rmk}\label{rmk:trivial-monodromy}
As explained in \cite{LMN}, Remark 2.5, there are some cases in which we can automatically know that a nice fibration $f$ has trivial monodromy, as the following:
\begin{itemize}
	\item If $f$ is locally trivial in the Zariski topology.
	\item If $f$ is a principal $G$-bundle with $G$ a connected algebraic group.
	\item If the fiber is $F = \PP^N$ for some $N$.
\end{itemize}
\end{rmk}

Let $B$ be an algebraic smooth variety and let $\rho: \pi_1(B) \to \Aut(V)$ be a representation where $V$ is a polarized pure Hodge structure preserved by $\rho$. Then, the local system associated to $\rho$ (see \cite{Peters-Steenbrink:2008}, Appendix B.3.1), call it $\mathbb{V}_\rho$, has a natural structure of polarizable variation of pure Hodge structures just by considering the corresponding subbundles for the Hodge filtration. If $\rho$ is also quasi-unipotent at infinity, then $\mathbb{V}_\rho$ is a good variation and, thus, it defines a mixed Hodge module that, for simplicity, we will denote also by $\rho \in \MHM{B}$. 

\begin{ex}
\begin{itemize}
	\item If $V$ is a graded-polarizable mixed Hodge structure and $\rho: \pi_1(B) \to \Aut(V)$ is a representation preserving it with quasi-unipotent monodromy at infinity, then, the associated variation of mixed Hodge structures $\mathbb{V}_\rho$ is graded-polarizable. It might be not admissible by its own but, in $K$-theory, it decomposes as a sum of polarizable variations of pure Hodge structures. Hence, in any case, we can associate to $\rho$ the sum of the associated mixed Hodge modules of each of these summands. The resulting mixed Hodge module is also denoted by $\rho \in \K{\MHM{B}}$.
	\item If $f: B' \to B$ is a regular morphism and $\rho: \pi_1(B) \to \Aut(V)$ is a representation preserving a graded polarizable mixed Hodge structure $V$, then $f^*\mathbb{V}_\rho = \mathbb{V}_{\rho \circ f}$ as variations of mixed Hodge structures. Hence, as elements of $\K{\MHM{B'}}$, we have $f^*\rho = (\rho \circ f)$ which justifies the notational ambiguity with pullbacks of representations.
	\item Let $\rho: \pi_1(B) \to \Aut(V)$ be a representation of graded-polarizable mixed Hodge structures with finite image. Then, automatically, $\rho$ is quasi-unipotent so, in particular, it is quasi-unipotent at infinity. Hence, it defines a mixed Hodge module.
	\item Let $f: X \to B$ be a nice fibration with fiber $F$ and monodromy representations $h_k: \pi_1(B) \to \Aut\left(H^k_c(F;\QQ)\right)$. If the representations $h_k$ preserve the mixed Hodge structure and have quasi-unipotent monodromy at infinity then, as variations of mixed Hodge modules, we have $R^kf_! \underline{\QQ}_X = \mathbb{V}_{h_k}$. So, by Corollary \ref{cor:hodge-monodromy-as-monodromy}, we have that, in $\K{\MHM{B}}$
	$$
		f_!\,\underline{\QQ}_X = \sum_k (-1)^k h_k
	$$
	\item In particular, if $f: X \to B$ is a finite covering of degree $d$, then $H_c^0(F;\QQ) = \QQ_0^d$ and $H_c^k(F;\QQ) = 0$ for $k>0$ as Hodge structures. In that case, the only non-trivial monodromy is $h = h_0: \pi_1(B) \to \Aut\mkern-4mu\left(\QQ_0^d\right)$ which preserves the Hodge structure and has finite image. Hence, we have that $f_!\,\underline{\QQ}_X = h$.
\end{itemize}
\end{ex}

\begin{defn}
Given a nice morphism $f: X \to B$, we define the \emph{Hodge monodromy representation} of $X$ on $B$ by
$$
	\RMc{f}{X}{B}=f_!\,\underline{\QQ}_X \in \K{\MHM{B}}.
$$
When the morphism $f$ it clear from the context, we will just write $\RM{X}{B}$.
\end{defn}

In order to get in touch with the definition of Hodge monodromy representation of \cite{LMN}, recall that $q = \QQ(-1)$ denotes the $(-1)$-Tate structure of weight $2$. In that case, $\K{\MHM{B}}$ has a natural $\ZZ[q]$-module structure inherited from the one as $\K{\MHSq}$-module. 

Now, suppose that we have a representation of mixed Hodge modules $\rho: \pi_1(B) \to \Aut(V)$, where $V$ is a balanced mixed Hodge structure (see \ref{ex:MixedHodgeGradingPure}, second item). In that case, in $\KM{B}$, we have that $\rho = \sum_p \rho^{p,p}\,q^p$ where $\rho^{p,p}: \pi_1(B) \to \Aut(V^{p,p})$ is the restriction of $\rho$ to the pieces of $V$ which are $V^{p,p}$ of degree $2p$ and concentrated at $(p,p)$.

In particular, for the monodromy representations, $h_k: \pi_1(B) \to \Aut(H_c^k(F;\QQ))$, of a nice fibration $X \to B$ with fiber $F$ of balanced type, we have $h_k = \sum_{p} h_{k}^{p,p} q^p$ where $h_k^{p,p}: \pi_1(B) \to \Aut(H_c^{k;p,p}(F;\QQ))$ are the restrictions. Therefore, adding up we have that
$$
	\RM{X}{B} = \sum_k (-1)^k\,h_k = \sum_{k,p}(-1)^k\, h_k^{p,p}(X)\,q^p,
$$
which can be seen as an element of the representation ring of $\pi_1(B)$ tensorized with $\ZZ[q^{\pm 1}]$. This is precisely the definition of \cite{LMN} of Hodge monodromy representations. In particular, reinterpreting Hodge monodromy representations as mixed Hodge modules, all the computations of \cite{LMN}, \cite{MM} and \cite{Martinez:2017} can be seen as calculations of mixed Hodge modules.

\subsection{Hodge monodromy representations of covering spaces}
\label{sec:Hodge-monodromy-covering-spaces}

Suppose that $X$ and $B$ are smooth complex varieties and that $\pi: X \to B$ is a regular morphism which is a covering space with finite fiber $F$ and degree $d$. In that case, the monodromy of $\pi$ is given by path lifting. 

More precisely, there is an action of $\pi_1(B)$ on $F$ given as follows. Let $\gamma \in \pi_1(B)$ and $x \in F$ so there exists an unique path lift $\tilde{\gamma}_x: [0,1] \to X$ of $\gamma$ such that $\tilde{\gamma}_x(0)=x$. Then, we set $\gamma \cdot x = \tilde{\gamma}_x(1) \in F$.
Since $F$ is a finite set of $d$ points, $H_c^k(F; \QQ)=0$ for $k > 0$ and $H_c^0(F;\QQ)=\QQ^d$. Hence, the only non-trivial action induced on cohomology $h = h_0: \pi_1(X) \to \GL{}(H_c^0(F;\QQ)) = \GL{d}(\QQ)$ is the monodromy action so $\RM{X}{B} = h$.

Using this description, some important examples of monodromy can be computed.

\begin{ex}
Let $X = \left\{(x,y) \in \CC^* \times \CC^* \,|\, y^2 = x\right\} \cong \CC^*$, $B = \CC^*$ and consider the projection $\pi: X \to B$ given by $\pi(x,y) = x$. It is a double covering with non-trivial monodromy. In order to compute this action, consider the fiber over $1 \in \CC^*$, $F = \pi^{-1}(1)=\left\{(1,1), (1,-1)\right\}$. Given the generator loop $\gamma(t) = e^{2\pi i t}$, for $0 \leq t \leq 1$, the lifts are loops $\tilde{\gamma} = \pm e^{\pi i t}$. Hence, if $\tilde{\gamma}(0) = (1,1)$ then $\tilde{\gamma}(1) = (1,-1)$ and viceversa. Therefore, on cohomology, the monodromy action $h: \pi_1(\CC^*) \to \GL{}(H_c^0(F;\QQ)) = \GL{2}(\QQ)$ is given by
$$
	\gamma \mapsto \begin{pmatrix}
	0 & 1\\
	1 & 0\\
\end{pmatrix} \sim \begin{pmatrix}
	1 & 0\\
	0 & -1\\
\end{pmatrix}.
$$
Hence, $h = T + S$, where $T$ is the trivial action (i.e.\ the unit of $\KM{B}$) and $S: \pi_1(\CC^* ) \to \QQ^* = \GL{1}(\QQ)$ is the action $\gamma \mapsto -1$. Thus, the Hodge monodromy representation of $\pi$ is $\RM{X}{B} = T + S$.

Given $M \in \KM{B}$, we will short $\mu(M) = c_!M \in \K{\MHSq}$, where $c: B \to \star$ is the projection onto a point. Observe that $\mu(B)= \mu(T_B)$ by construction. In particular, since $\coh{X} = \coh{\CC^*} = q-1$ and $\mu{(T)} = \coh{\CC^*}=q-1$, we must have $\mu(S) = 0$ because it holds
$$
	q-1 = \coh{X} = \mu(\RM{X}{B}) = \mu(T) + \mu(S) = q-1 + \mu(S). 
$$

In general, let us fix $p_0, \ldots, p_s \in \CC$ with $p_i \neq \pm p_j$ for $0 \leq i,j \leq s$ and consider $B = \CC - \left\{p_0, \ldots, p_s\right\}$. We have the collection of representations $S_{p_i}: \pi_1(B) \to \QQ^*$, for $0 \leq i \leq s$, that are given by $S_{p_i}(\gamma_{p_j}) = 1$ if $i \neq j$ and $S_{p_i}(\gamma_{p_i}) = -1$, where $\gamma_{p_j}$ is a small loop around $p_j$

The importance of these representations appears when we consider $X = \left\{(x,y) \in B \times \CC \,|\, y^2 = x-p_0\right\}$ with projection $\pi: X \to B$, $\pi(x,y)=x$. An analogous analysis to the previous one shows that the Hodge monodromy representation of $\pi$ is $\RM{X}{B} = T + S_{p_0}$. However, in this case $X \cong \CC^* - \left\{\pm \sqrt{p_1}, \ldots, \pm \sqrt{p_s}\right\}$ so $\coh{X} = q - 2s -1$. Thus, it holds
$$
	q -2s-1 = \coh{X} = \mu(\RM{X}{B}) = \mu(T) + \mu(S_{p_0}) = q-(s+1) + \mu(S_{p_0}).
$$
Hence, $\mu(S_{p_0}) = -s$. This agrees with discussion after Theorem 6 of \cite{MM}. Actually, a modification of this argument with $X = \left\{y^2 = (x-p_i)(x-p_j)\right\}$ shows that, also, $\mu(S_{p_{i}} \otimes S_{p_{j}}) = -s$ for any $i \neq j$.
\end{ex}

\begin{ex}
Let pick a point $\lambda \in \CC - \left\{0,1\right\}$, set $B = \CC - \left\{0, 1, \lambda\right\}$ and consider the variety $X = \left\{(x, y) \in B \times \CC^*\,|\,y^2 = x(x-1)(x-\lambda)\right\}$ with a projection $\pi: X \to B$, $\pi(x,y)=x$. As above, we have that $\RM{X}{B} = T + S_0 \otimes S_1 \otimes S_{\lambda}$. Moreover, $X$ is a cubic curve minus four points (the three removed values of $x$ and the point at infinity), so $\coh{X} = -3 -u -v + q$, where $u$ (resp.\ $v$) is the $1$-dimensional pure Hodge structure of degree $1$ concentrated at $(1,0)$ (resp.\ $(0,1)$). Since $\mu(T) = \coh{B} = q-3$, we have that $\mu(S_0 \otimes S_1 \otimes S_{\lambda}) = -u-v$. This is important since $-u-v \in \K{\MHSq}$ is not of balanced type.
\end{ex}

\begin{ex}\label{ex:monodromy-plane-minus-points}
Let $X = \CC^* - \left\{\pm 1\right\}$ and $B = \CC - \left\{\pm 2\right\}$ with the morphism $t: X \to B$ given by $t(\lambda) = \lambda + \lambda^{-1}$ for $\lambda \in \CC^* - \left\{\pm 1\right\}$. Again, the morphism $t$ is a double cover and the monodromy action of $\pi_1(\CC - \left\{\pm 2\right\})$ on the fiber $F$ is given by path lifting. The lifts of a small loop around $2$, $\gamma_{2}(t) = 2 + r e^{2\pi it}$, for $0 \leq t \leq 1$ and $r >0$ small, are the paths
$$
	\tilde{\gamma}_{2}(t) = \frac{2 + r e^{2\pi it} \pm r^{1/2} e^{\pi i t}\sqrt{4+r e^{it}}}{2} \approx 1 \pm r^{1/2}e^{\pi i t} + \frac{r}{2} e^{2\pi it}.
$$
Hence, the endpoint of $\tilde{\gamma}_{2}(t)$ does not agree with the initial point, so $\gamma_2$ acts on $F$ by interchanging its elements. Analogously, a small loop around $-2$, $\gamma_{-2}$, also interchanges the elements. Thus, on cohomology, the monodromy action $h: \pi_1(\CC - \left\{\pm 2\right\}) \to \GL{2}(\QQ)$ is given by
$$
	\gamma_2 \mapsto \begin{pmatrix}
	0 & 1\\
	1 & 0\\
\end{pmatrix} \sim \begin{pmatrix}
	1 & 0\\
	0 & -1\\
\end{pmatrix}, \hspace{2cm} \gamma_{-2} \mapsto \begin{pmatrix}
	0 & 1\\
	1 & 0\\
\end{pmatrix} \sim \begin{pmatrix}
	1 & 0\\
	0 & -1\\
\end{pmatrix}.
$$
Therefore, $\RM{\CC^* - \left\{\pm 1\right\}}{\CC - \left\{\pm 2\right\}} = h = T + S_2 \otimes S_{-2}$, where $T$ is the trivial action and $S_{\pm 2}: \pi_1(\CC - \left\{\pm 2\right\}) \to \QQ^* = \GL{1}(\QQ)$ are the actions $\gamma_{\pm 2} \mapsto -1$ and $\gamma_{\mp 2} \mapsto 1$ respectively. Observe that, as explained above, $\mu(S_2 \otimes S_{-2}) = -1$ so
\begin{align*}
	\coh{\CC^* - \left\{\pm 1\right\}} &= \mu(\RM{\CC^* - \left\{\pm 1\right\}}{\CC - \left\{\pm 2\right\}}) \\
	&= \mu(T) + \mu(S_2 \otimes S_{-2}) = q- 2 -1 = q-3,
\end{align*}
as we already knew.
\end{ex}

\subsection{Equivariant Hodge monodromy representations}
\label{sec:equivariant-hodge-mono}

In order to finish this chapter, let us look more closely to a special type of coverings coming from group actions. Let $X$ be an algebraic variety with a free action of $\ZZ_2$ on it. We define its equivariant Hodge structures as the elements of $\K{\MHSq}$ given by
$$
	\coh{X;\QQ}^+ = \coh{X/\ZZ_2;\QQ}, \hspace{1cm} \coh{X;\QQ}^- = \coh{X;\QQ} - \coh{X;\QQ}^+.
$$
Moreover, if $X \to Z$ is a nice fibration equivariant for the $\ZZ_2$-action, we denote $\RM{X}{Z}^+ = \RM{X/\ZZ_2}{Z}$ and $\RM{X}{Z}^- = \RM{X}{Z} - \RM{X}{Z}^+$.

Observe that the action of $\ZZ_2$ on $X$ induces an action of $\ZZ_2 = \pi_1(X/\ZZ_2)/\pi_1(X)$ on the cohomology of $X$. In that case, $\coh{X;\QQ}^+$ can also be seen as the fixed part of the cohomology of $X$ by this action (see Proposition 4.3 of \cite{Florentino-Silva:2017}). Analogously, $\RM{X}{Z}^+$ can be seen as the invariant part of the monodromy representation by the action of $\ZZ_2$.

Now, let us consider a nice fibration $f: X \to B$ with fiber $F$ and trivial monodromy. In addition, suppose that there is a free action of $\ZZ_2$ on $X$ and $B$ such that $f$ is equivariant, so it descends to a regular morphism $\tilde{f}: X/\ZZ_2 \to B/\ZZ_2$. This means that they fit in a diagram of fibrations
\[
\begin{displaystyle}
   \xymatrix
   {
	F \ar[r] & X \ar[r]^f \ar[d] & B \ar[d]\\
	F \ar[r] & X/\ZZ_2 \ar[r]_{\tilde{f}} & B/\ZZ_2
   }
\end{displaystyle}
\]
where the vertical arrows are the quotient maps and the upper fibration has trivial monodromy. In \cite{LMN}, Proposition 2.6, it is proven that, on $\K{\MHSq}$, we have an equality
$$
	\coh{X;\QQ}^+ = \coh{B;\QQ}^+ \coh{F;\QQ}^+ + \coh{B;\QQ}^- \coh{F;\QQ}^-.
$$
Here, $\coh{F;\QQ}^+$ denotes the invariant part of the mixed Hodge structure on the cohomology of $F$ under the action of $\ZZ_2$.

\begin{ex}\label{ex:equiv-product}
Proceding recursively with this formula, we have that, if $\ZZ_2$ acts on $X^n$ with a simultaneous action, then
$$
    \coh{X^n}^+ = \frac{1}{2}\left[\coh{X}^n + \left(\coh{X}^+ - \coh{X}^-\right)^n\right],
$$
$$
    \coh{X^n}^- = \frac{1}{2}\left[\coh{X}^n - \left(\coh{X}^+ - \coh{X}^-\right)^n\right].
$$
\end{ex}

The same argument can be adapted to consider the case of Hodge monodromy representations. Suppose that, in addition to previous diagram, we also have a nice morphism $\pi: B \to Z$ which is invariant for the $\ZZ_2$ action. From this map, we also obtain natural morphisms $X \to Z$, $X/\ZZ_2 \to Z$ and $B/\ZZ_2 \to Z$. In that case, we claim that
$$
	\RM{X}{Z}^+ = \RM{B}{Z}^+ \coh{F;\QQ}^+ + \RM{B}{Z}^- \coh{F;\QQ}^-.
$$
In order to check it, observe that $\RM{X}{Z}^+ = \RM{X/\ZZ_2}{Z}$ is precisely the invariant part of the monodromy of $X \to Z$ by the action of $\ZZ_2$. Hence, as $X \to B$ has trivial monodromy, we have $f_! \underline{\QQ}_X = \coh{F;\QQ} \underline{\QQ}_B$ and, thus, we can compute
\begin{align*}
	\RM{X}{Z}^+ &= \left[(\pi \circ f)_! \underline{\QQ}_X\right]^+ = \left[\pi_! \left(f_!\underline{\QQ}_X\right)\right]^+ = \left[\coh{F;\QQ} \boxtimes \pi_!\underline{\QQ}_B\right]^+ = \left[\coh{F;\QQ} \boxtimes \RM{B}{Z}\right]^+ \\
	&= \coh{F;\QQ}^+ \RM{B}{Z}^+ + \coh{F;\QQ}^- \RM{B}{Z}^-.
\end{align*}

\begin{rmk}
This kind of arguments deeply reminds to equivariant cohomology. Actually, we expect that, following the ideas of \cite{Florentino-Silva:2017}, equivariant Hodge representations might be formulated in great generality for equivariant mixed Hodge modules. It is a future work to explore. However, in this thesis, we will not need this general case but only for $\ZZ_2$ actions.
\end{rmk}

%% file: Chapters/RepresentationVarieties.tex
% Chapter 3

\chapter{TQFTs and Representation Varieties} % Main chapter title

\label{chap:representation} % For referencing the chapter elsewhere, use \ref{Chapter1} 

\lhead{Chapter 3. \emph{TQFTs and Representation varieties}} % This is for the header on each page - perhaps a shortened title

%----------------------------------------------------------------------------------------

\section{Group representations as algebraic varieties}
\label{sec:intro-representation}

In this chapter, we are going to focus on the main objects of study of this thesis, the so-called representation varieties. For this reason, in this section we will review the definition of representation and character varieties and we will describe some variants as parabolic representation varieties.

Along this section, we will work over a ground algebraically closed field $k$.

\begin{defn}
Let $\Gamma$ be a finitely generated group and $G$ an algebraic group. The set of representations of $\Gamma$ into $G$
$$
	\Rep{G}(\Gamma)=\Hom(\Gamma, G)
$$
is called the \emph{representation variety}. 
\end{defn}

As its name suggests, $\Rep{G}(\Gamma)$ has a natural algebraic structure. In order to define it, let $\Gamma = \langle \gamma_1, \ldots, \gamma_r\;|\; R_\alpha(\gamma_1, \ldots, \gamma_r)=1\rangle$ be a presentation of $\Gamma$ with finitely many generators, where $R_\alpha$ are the relations (possibly infinitely many). In that case, we define the injective map $\psi: \Hom(\Gamma, G) \to  G^r$ given by $\psi(\rho)=(\rho(\gamma_1), \ldots, \rho(\gamma_r))$. The image of $\psi$ is the algebraic subvariety of $G^r$ 
$$
	\img{\psi} = \left\{(g_1, \ldots, g_r) \in G^r\;\right|R_\alpha(g_1, \ldots, g_r)=1\left.\right\}.
$$
Hence, we can impose an algebraic structure on $\Hom(\Gamma, G)$ by declaring that $\psi$ is a regular isomorphism over its image. Observe that this algebraic structure does not depend on the chosen presentation.

\subsection{Character varieties}
\label{sec:char-var}
The representation variety $\Rep{G}(\Gamma)$ has a natural action of $G$ by conjugation i.e.\ $g \cdot \rho (\gamma) = g\rho(\gamma) g^{-1}$ for $g \in G$, $\rho \in \Rep{G}(\Gamma)$ and $\gamma \in \Gamma$. Recall that two representations $\rho, \rho'$ are said to be isomorphic if and only if $\rho' = g \cdot \rho$ for some $g \in G$.

\begin{defn}
Let $\Gamma$ be a finitely generated group and $G$ an algebraic reductive group. The Geometric Invariant Theory quotient
$$
	\Char{G}(\Gamma) = \Rep{G}(\Gamma) \sslash G,
$$
is called the \emph{character variety}.
\end{defn}

\begin{rmk}
The Geometric Invariant Theory quotient, also known as GIT quotient, is a kind of quotient that makes sense for algebraic varieties. The problem is that, in general, the orbit space is no longer an algebraic variety, so we have to substitute the quotient space by an algebraic variety that behaves as expected for quotient spaces. This can be done if the acting group $G$ is reductive and the unique solution is the GIT quotient. We will review Geometric Invariant Theory in Section \ref{sec:review-git}.
\end{rmk}

\begin{ex}\label{ex:gamma-for-repr}
\begin{itemize}
	\item Take $\Gamma = F_n$, the free group with $n$ generators. In that case, for short, we will denote the associated representation variety by $\Xf{n}(G) = \Rep{G}(F_n) = G^n$ or, when the group $G$ is understood, just by $\Xf{n}$. The corresponding character variety is $\Char{G}(F_n) = \Xf{n}(G) \sslash G = G^n \sslash G$. The importance of this case comes from the fact that, if $\Gamma$ is any finitely generated group with $n$ generators, then the epimorphism $F_n \to \Gamma$ gives an inclusion $\Rep{G}(\Gamma) \subseteq \Xf{n}(G)$.
	\item Let $M$ be a compact manifold with fundamental group $\Gamma = \pi_1(M)$. The fundamental group of such a manifold is finitely generated since a compact manifold has the homotopy type of a finite CW-complex. Hence, we can form its representation variety, that we will shorten $\Rep{G}(M) = \Rep{G}(\pi_1(M))$. The corresponding character variety is called the character variety of $M$ and it is denoted by $\Char{G}(M) = \Rep{G}(M) \sslash G$.
	\item As a particular case, take $\Gamma = \pi_1(\Sigma_g)$, the fundamental group of the genus $g$ compact surface. In that case, we will denote the associated representation variety by $\Xs{g}(G) = \Rep{G}(\Sigma_g)$ or even just by $\Xs{g}$ when the group $G$ is understood from the context. Since the standard presentation of $\Gamma$ is
$$
	\pi_1(\Sigma_g) = \left\langle \alpha_1, \beta_1 \ldots, \alpha_{g}, \beta_g\;\;\left|\;\; \prod_{i=1}^g [\alpha_i, \beta_i] = 1 \right.\right\rangle,
$$
then $\Xs{g}(G) \subseteq \Xf{2g}(G)$. Actually, $\Xs{g}(G)$ is given by tuples of $2g$ elements of $G$ satisfying the relation of $\pi_1(\Sigma_g)$. Explicitly
$$
	\Rep{G}(\Sigma_g) = \left\{(A_1, B_1 \ldots, A_{g}, B_g) \in G^{2g}\;\;\left|\;\; \prod_{i=1}^g [A_i, B_i] = 1 \right.\right\}.
$$
The GIT quotient of this representation variety under the conjugacy action will be denoted by $\Chars{g}(G)$ or just $\Chars{g}$. These representation varieties will be studied in detail in Sections \ref{almost-TQFT-strategy} and \ref{sec:sl2-repr-var} and the corresponding character varieties in Sections \ref{sec:surface-groups-nonpar} and \ref{sec:surface-groups-par}.
\end{itemize}
\end{ex}

\subsection{Parabolic representation varieties}
\label{sec:parabolic-repr-var}

A step further in the construction of character varieties can be done by considering an extra structure on them, called a parabolic structure. As we will see in \ref{sec:non-abelian-hodge}, these parabolic structures allow us to extend the non-abelian Hodge theory. In this way, we obtain further correspondences between parabolic character varieties, moduli spaces of parabolic Higgs bundles and moduli spaces of logarithmic flat connections.

As above, let $\Gamma$ be a finitely generated group and let $G$ be an algebraic group. A \emph{parabolic structure} on $\Rep{G}(\Gamma)$, $Q$, is a finite set of pairs $(\gamma, \lambda)$, where $\gamma \in \Gamma$ and $\lambda \subseteq G$ is a locally closed subset which is closed under conjugation. Given such a parabolic structure, we define the \emph{parabolic representation variety}, $\Rep{G}(\Gamma, Q)$, as the subset of $\Rep{G}(\Gamma)$
$$
	\Rep{G}(\Gamma, Q) = \left\{\rho \in \Rep{G}(\Gamma)\,\,|\,\, \rho(\gamma) \in \lambda\,\,\textrm{for all } (\gamma, \lambda) \in Q\right\}.
$$
As in the non-parabolic case, $\Rep{G}(\Gamma, Q)$ has a natural algebraic variety structure given as follows. Suppose that $Q = \left\{(\gamma_1, \lambda_1), \ldots, (\gamma_s, \lambda_s)\right\}$. Then, we choose a finite set of generators $S$ of $\Gamma$ that contains all the $\gamma_{i}$, say $S = \left\{\eta_1, \ldots, \eta_r, \gamma_{1}, \ldots, \gamma_{s}\right\}$.
Using $S$, we can identify $\Rep{G}(\Gamma)$ with a closed subvariety of $G^{r + s}$ (see Section \ref{sec:rep-var}). In that case, we also have a natural identification $\Rep{G}(\Gamma, Q) = \Rep{G}(\Gamma) \cap \left(G^r \times \lambda_1 \times \ldots \times \lambda_s\right)$. We impose in $\Rep{G}(\Gamma, Q)$ the algebraic structure inherited from this identification and, as before, such a structure does not depend on $S$.

The conjugacy action of $G$ on $\Rep{G}(\Gamma)$ restricts to an action on $\Rep{G}(\Gamma, Q)$ since the subsets $\lambda_i$ are closed under conjugation. The GIT quotient of the representation variety by this action,
$$
	\Char{G}(\Gamma, Q) = \Rep{G}(\Gamma, Q) \sslash G,
$$
is called the \emph{parabolic character variety}.

\begin{rmk}
As a particular choice for the subset $\lambda \subseteq G$, we can choose the conjugacy classes of an element $h \in G$, denoted $[h]$. Observe that $[h] \subseteq G$ is locally closed since, by \cite{Newstead:1978} Lemma 3.7, it is an open subset of its Zariski closure, $\overline{[h]}$. \end{rmk}

\begin{ex}\label{ex:par-struc-repr}
\begin{itemize}
	\item The simplest example can be obtained by taking $\Gamma = F_{n+s}$ and the parabolic structure $Q = \left\{(\gamma_{1}, \lambda_1), \ldots, (\gamma_s, \lambda_s)\right\}$, where $\gamma_1, \ldots, \gamma_s \in F_{r+s}$ is an independent set. The corresponding representation variety is, thus, $\Rep{G}(F_{r+s}, Q) = G^r \times \lambda_1 \times \ldots \times \lambda_s$.	
	\item Let $\Sigma = \Sigma_g - \left\{p_1, \ldots, p_s\right\}$ with $p_i \in \Sigma_g$ distinct points, called the punctures or the marked points. In that case, we have a presentation of the fundamental group of $\Sigma$ given by
$$
	\pi_1(\Sigma) = \left\langle \alpha_1, \beta_1 \ldots, \alpha_{g}, \beta_g, \gamma_1, \ldots, \gamma_s\;\;\left|\;\; \prod_{i=1}^g [\alpha_i, \beta_i]\prod_{j=1}^s \gamma_s = 1 \right.\right\rangle,
$$
where the $\gamma_i$ are the positive oriented simple loops around the punctures. As parabolic structure, we will take $Q = \left\{(\gamma_1, \lambda_1), \ldots, (\gamma_s, \lambda_s)\right\}$. The corresponding parabolic representation variety is
$$
	\Rep{G}(\pi_1(\Sigma), Q) = \left\{(A_1, B_1 \ldots, A_{g}, B_g, C_1, \ldots, C_s) \in G^{2g+s}\;\;\left|\;\; \begin{matrix}
	\displaystyle \prod_{i=1}^g [A_i, B_i]\prod_{j=1}^s C_s = 1 \\
	C_j \in \lambda_j
\end{matrix}	
	\right.\right\}.
$$
Observe that the epimorphism $F_{2g +s} \to \pi_1(\Sigma)$ gives an inclusion $\Rep{G}(\pi_1(\Sigma), Q) \subseteq \Rep{G}(F_{2g+s}, Q)$.
\end{itemize}
\end{ex}

Let $M$ be a compact differentiable manifold and let $S \subseteq M$ be a codimension $2$ closed connected submanifold with a co-orientation (i.e.\ an orientation of its normal bundle). Embed the normal bundle as a small tubular neighbourhood $U \subseteq M$ around $S$. Fixed $s \in S$, consider an oriented local trivialization $\psi: V \times \RR^2 \to U$ of the normal bundle around an open neighbourhood $V \subseteq S$ of $s$. In that situation, the loop $\gamma(t) = \psi(s, (\cos{t}, \sin{t})) \in \pi_1(M - S)$ is called the \emph{positive meridian} around $s$.

\begin{figure}[h]
	\begin{center}
	\includegraphics[scale=0.25]{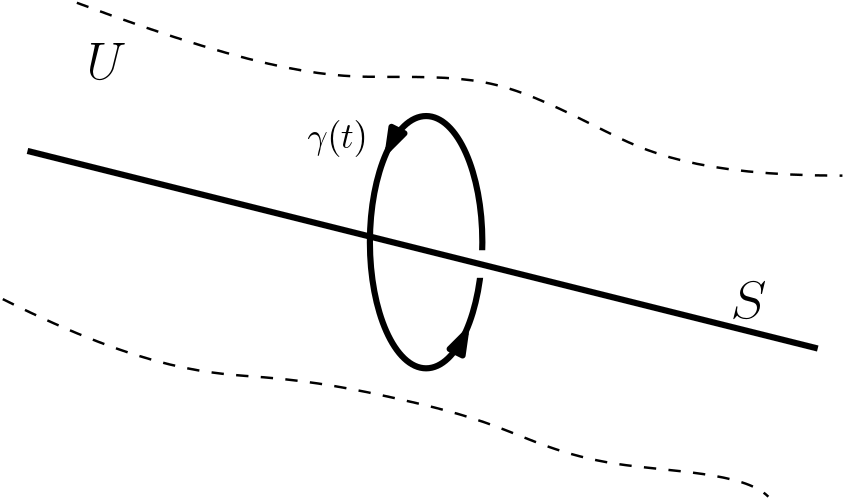}
	\end{center}
\vspace{-0.5cm}
\end{figure}

\begin{rmk}
\begin{itemize}
	\item Two positive meridians around $s, s' \in S$ are conjugate to each other on $\pi_1(M-S)$. For that, choose any path in $S$ between $s$ and $s'$ and move it slighly outwards $S$ in order to obtain a path outsite $S$ very close to the original. By means of such a path, both meridians become conjugate.
	\item If $S$ is not connected, then the conjugacy classes of meridians are in correspondence with the connected components of $S$.
	\item A loop $\gamma \in \pi_1(M-S)$ is a generalized knot. In that sense, the kernel of $\pi_1(M - S) \to \pi_1(M)$ are the loops in $M-S$ that `surround' $S$. It can be proven (see \cite{Smith:1978}) that this kernel is the smallest normal subgroup containing all the meridians around the connected components of $S$. Hence, the meridians capture all the information of the loops surrounding $S$.
\end{itemize}
\end{rmk}

Take $\Lambda$ to be a collection of locally closed subsets of $G$ that are invariant under conjugation. Suppose that $Q=\left\{(S_1, \lambda_1), \ldots, (S_s, \lambda_s)\right\}$ is a parabolic structure over $\Lambda$, in the sense of Section \ref{sec:some-extra-struct}. Let us denote $S = \bigcup_i S_i$. To this parabolic structure, we can build the parabolic structure on $\Rep{G}(\pi_1\left(M - S\right))$, also denoted by $Q$, $Q = \left\{(\gamma_{1,1}, \lambda_{1}), \ldots, (\gamma_{1,m_1}, \lambda_{1}), \ldots, (\gamma_{s, 1}, \lambda_s), \ldots, (\gamma_{s, m_s}, \lambda_s)\right\}$, where $\gamma_{i, 1}, \ldots, \gamma_{i, m_i}$ is a generating set of positive meridians around the connected components of $S_i$. As for the non-parabolic case, we will shorten the corresponding parabolic character variety by $\Rep{G}(M, Q) = \Rep{G}(\pi_1\left(M-S\right), Q)$.

\begin{rmk}
Suppose that $M=\Sigma_g$ is a closed oriented surface. A parabolic structure is given by $Q=\left\{(p_1, \lambda_1), \ldots, (p_s, \lambda_s)\right\}$, with $p_i \in \Sigma_g$ points with a preferred orientation of a small disk around them (see Example \ref{ex:parab-struct}). In that case, the meridian of $p_i$ is given by a small loop encycling $p_i$ positively with respect to the orientation of the small disk around it. Therefore, the associated parabolic structure of representation variety is
$$
	\Rep{G}(\Sigma_g, Q) = \left\{(A_1, B_1 \ldots, A_{g}, B_g, C_1, \ldots, C_s) \in G^{2g+s}\;\;\left|\;\; \begin{matrix}
	\displaystyle \prod_{i=1}^g [A_i, B_i]\prod_{j=1}^s C_s = 1 \\
	C_j^{\epsilon_j} \in \lambda_j
\end{matrix}	
	\right.\right\},	
$$
where $\epsilon_j = 1$ if the orientation of the disk around $p_j$ agrees with the global orientation and $\epsilon_j = -1$ if it does not. Notice that, they agree with the ones of Example \ref{ex:par-struc-repr}.
\end{rmk}

\subsection{A brief on non-abelian Hodge theory}
\label{sec:non-abelian-hodge}

The importance of character varieties comes from the strong relation of these spaces with the moduli spaces of Higgs bundles and the moduli spaces of flat connections. It is a very deep and active area of reseach, known as \emph{non-abelian Hodge theory}. In this section, we will sketch the fundamental results of the theory, maybe in a rather informal way. For a detailed account on this wide subject, see \cite{SimpsonNonAbelian} or \cite{Hausel:2005}.

As the name sugests, the starting point of this theory is the following interplay between algebro-geometric objects. Since every compact Riemann surface $X$ is a compact Kähler manifold, classical Hodge theory (see Theorem \ref{HodgePureStructureKahler}) gives us a decomposition $H^1(X; \CC) \cong H^{1,0}(X) \oplus H^{0,1}(X)$.
Considering Dolbeault cohomology as a sheaf cohomology, if $\mathbf{\Omega}^p$ denote the sheaf of holomorphic $p$-forms, we have isomorphisms
$$
	H^{1,0}(X) \cong H^0(X; \mathbf{\Omega}^1), \hspace{1.5cm} H^{0,1}(X) \cong H^1(X; \mathbf{\Omega}^0) = H^1(X; \mathcal{O}_X).
$$
Each element of $H^1(X, \mathcal{O}_X)$ determines a degree zero line bundle via its image under the exponential sequence $H^1(X; \mathcal{O}_X) \to H^1(X; \mathcal{O}_X^*) \to H^2(X; \ZZ)$. Therefore, by means of these isomorphisms, the previous decomposition can be reinterpreted as that it is equivalent to have a degree $1$ cohomology class as to have a pair of a degree zero algebraic line bundle and a holomorphic $1$-form.

On the other hand, if $H_B^1(X; \CC)$ is the singular cohomology of $X$ with coefficients in $\CC$ (also called Betti cohomology), then de Rham theorem give us an isomorphism $H_B^1(X; \CC) \cong H_{DR}(X; \CC)$. Moreover, by Hurewicz theorem, we have that $H_B^1(X; \CC) \cong \frac{\pi_1(X)}{[\pi_1(X), \pi_1(X)]} \otimes_\ZZ \CC \cong \Hom(\pi_1(X), \CC)$.

Philosophically, non-abelian Hodge theory translates this abelian framework to the more general setting of moduli spaces. Along this section, $G$ will be a complex reductive linear group and $X$ a compact Riemann surface. In the sense of non-abelian Hodge theory, the non-abelian analogues of the previous ones are the following:
\begin{itemize}
	\item De Rham cohomology: Its non-abelian analogous is the moduli space of flat $G$-connections, $\mathcal{M}_{DR}(X, G)$, i.e.\ the space that parametrices flat connections on $X$.
	\item Dolbeault cohomology: Its non-abelian analogous is the moduli space of $G$-Higgs bundles, $\mathcal{M}_{Dol}(X, G)$ i.e.\ the space that parametrices Higgs bundles on $X$.
	\item Betti cohomology: Its non-abelian analogous is character variety, $\Char{G}(X)$, also denoted in this context $\mathcal{M}_{B}(X, G)$.
\end{itemize}

The first step in this non-abelian Hodge theory is to understand the relation between flat $G$ connections and representations $\pi_1(X) \to G$. Recall from \cite{Peters-Steenbrink:2008}, Appendix B.3, that every $G$-local system (i.e.\ a local system with transition functions on $G$) induce a representation $\pi_1(X) \to G$, just by considering its action over a point. Now, notice that a $G$-local system is the same as a $G$-locally constant sheaf and the later is equivalent to a $G$-vector bundle equipped with a flat connection. Hence, to each flat $G$-vector bundle on $X$ we can associate a representation of $\pi_1(X)$ onto $G$, the so-called monodromy representation.
 
Therefore, putting together these correspondences, we have a mapping from flat $G$-vector bundles, modulo gauge equivalence, to representations $\pi_1(X) \to G$ modulo conjugation. The fact that this identification which respects the topology of the moduli spaces is the content of the so-called {Riemann-Hilbert correspondence}, whose proof can be found in \cite{SimpsonII}.

\begin{thm}[Riemann-Hilbert correspondece]
The moduli spaces of $G$-flat connections $\mathcal{M}_{DR}(X, G)$ and the character variety $\Char{G}(X)$ are real analytically isomorphic. 
\end{thm}

\begin{rmk}
The coincidence between the name of this result and the one in the theory of $\cD$-modules and perverse sheaves is not accidental. Both theorems show the same principle: it is equivalent to have a vector bundle with an integrable connection (i.e.\ a $\cD$-module) that some kind of algebraic data as perverse sheaves or Higgs bundles.
\end{rmk}

Even more, we can twist these spaces in order to obtain a more general Riemann-Hilbert correspondence. For this purpose, let $G$ be a linear group with a faithful representation of rank $n$ and pick a point $p \in X$. For $d \in \NN$, we consider the parabolic structure $Q_d = \left\{\left(p, \left[e^{\frac{2\pi i d}{n}}\Id\right]\right)\right\}$ on $X$. The associated character variety
$$
	\Char{G}(X,Q_d) = \left\{(A_1, \ldots, A_g, B_1, \ldots, B_g) \in G^{2g} \,|\, \prod_{k=1}^g[A_k, B_k] = e^{\frac{2\pi i d}{n}}\Id \right\}\sslash G
$$
is usually called the twisted character variety and it is denoted $\mathcal{M}_{DR}^d(X, G)$ in the context of non-abelian Hodge theory.

\begin{thm}[Riemann-Hilbert correspondece, twisted case]
The moduli spaces of logarithmic $G$-flat connections on $X$ with a single pole at $p$ with residue $-\frac{d}{n}\Id$, $\mathcal{M}_{DR}^d(X, G)$, and the twisted character variety, $\mathcal{M}_{B}^d(X, G)$, are real analytically isomorphic. 
\end{thm}

Now, we can also focus our attention to the moduli space of $G$-Higgs bundles or, more restrictive, of holomorphic vector bundles. In this context, the starting point of the theory was a result of Narasimahan and Seshadri that relates poly-stable holomorphic vector bundles of rank $n$ and unitary character varieties. In the original proof in \cite{NarasimhanSeshadri}, they used only algebraic methods to prove the theorem. However, Donaldson, in a later paper \cite{Donaldson:1983}, gave a new proof of this theorem using gauge-theoretical methods. This iniciated the study of this kind of theorems from the point of view of Higgs bundles.

\begin{thm}[Narasimhan-Seshadri]
The moduli space of poly-stable holomorphic vector bundles of rank $n$ and degree $d$, $\mathcal{M}_{\textrm{VB}}(M, n, d)$, is homeomorphic to the twisted character variety $\mathcal{M}_{B}^d(X, U(n))$.
\end{thm}

Hitchin, in \cite{Hitchin}, proved a generalization of this theorem for the case of $SU(2)$-Higgs bundles using a totally different proof based on gauge theory. Later, the combined work of Donaldson, Corlette and Simpson in \cite{Donaldson:1987}, \cite{Corlette:1988} and \cite{Simpson:1992}, among others, proved the following version.

\begin{thm}[Hitchin-Kobayashi correspondence]
Let $G \subseteq \GL{n}(\CC)$ be a reductive Lie group. The moduli space of poly-stable $G$-vector bundles of rank $n$ and degree $d$, $\mathcal{M}_{Dol}^d(M, G)$, is real analitically isomorphic to the twisted character variety, $\Char{G}(X, Q_d) = \mathcal{M}_{B}^d(X, G)$.
\end{thm}

Therefore, with this result, the relation between moduli spaces is the following.
\[
\begin{displaystyle}
   \xymatrix
   {
  	 & \mathcal{M}_{B}^d(X, G) \ar@{{<}-{>}}[rd]^{}\ar@{{<}-{>}}[ld]_{} & \\
   	 \mathcal{M}_{DR}^d(X, G) \ar@{{<}-{>}}[rr]_{} & & \mathcal{M}_{Dol}^d(X, G) \\
   }
\end{displaystyle}   
\]
\begin{rmk}
These equivalences are real analytic but they are far from being holomorphic (complex analytic). Hence, these three different points of view allow us to introduce a very special structure on these spaces, known as a hyperkähler structure, which consists of three compatible Kähler structures.
\end{rmk}

Finally, in the context of parabolic $G$-Higgs bundles, Mehta and Seshadri in \cite{Mehta-Seshadri:1980}, first, and later Simpson in \cite{Simpson:parabolic}, proved the following non-abelian Hodge correspondence for parabolic Higgs bundles.

\begin{thm}
Let $G \subseteq \GL{n}(\CC)$ reductive. Let us choose parabolic points $p_1, \ldots, p_s \in X$ and let us fix conjugacy classes $\lambda_1, \ldots, \lambda_s \subseteq G$ of semisimple elements. Denoting the effective Weil divisor $D=p_1 + \ldots + p_s$ and $Q = \left\{(p_i, \lambda_i)\right\}_{i=1}^s$ we have:
\begin{itemize}
	\item The moduli space of parabolic $G$-Higgs bundles of parabolic degree zero with parabolic structures $\alpha$ on $D$ (see \cite{Boden-Yokogawa}, Definition 2.1), $\mathcal{M}_{Dol}^\alpha(X, G)$, is real analitically isomorphic to the parabolic character variety, $\Char{G}(X, Q) = \mathcal{M}_{B}^Q(X, G)$.
	\item The moduli space of logarithmic flat $G$-vector bundles with poles in $D$, $\mathcal{M}_{DR}^D(X, G)$, is real analytically isomorphic to the parabolic character variety, $\Char{G}(X, Q) = \mathcal{M}_{B}^Q(X, G)$.
\[
\begin{displaystyle}
   \xymatrix
   {
  	 & \mathcal{M}_{B}(X, Q, G) \ar@{{<}-{>}}[rd]\ar@{{<}-{>}}[ld]& \\
   	 \mathcal{M}_{DR}^D(X, G) \ar@{{<}-{>}}[rr] & & \mathcal{M}_{Dol}^\alpha(X, G) \\
   }
\end{displaystyle}   
\]
\end{itemize}
\end{thm}

\begin{rmk}
For the definition of a parabolic Higgs bundle, see \cite{Boden-Yokogawa}.
\end{rmk}

In particular, take $G=\SL{2}(\CC)$, $X$ an elliptic curve and two punctures. Using the Nahm transform (see \cite{Jardim:2004}, \cite{Jardim}, \cite{JardimThesis}, \cite{Jardim:2001} and \cite{Biquard-Jardim}), it can be proven that $\mathcal{M}_{Dol}^\alpha(X, G)$ is diffeomorphic to the moduli space of doubly periodic instantons (i.e.\ anti self-dual solutions to Yang-Mills equations which are invariant in two directions). This exceptional equivalence is one of the reasons for studying $\SL{2}(\CC)$-character varieties in Section \ref{sec:sl2-repr-var}.

\section{A TQFT for representation varieties}
\label{sec:TQFT-for-repr}

This section is the heart of this PhD thesis. Here, we are going to construct a lax monoidal TQFT (of pairs) that computes the Hodge structures of representations varieties of any complex group. For that purpose, we will use the construction technique described in Theorem \ref{thm:physical-constr-sheaf}.

As the category of fields for this construction, we take the category of complex algebraic varieties, $\CVar$. In this setting, the role of the field theory will be played by a sort of generalized representation varieties that consider not only representations of the fundamental group but also of the fundamental groupoid. For the quantisation, we will use the $\CVar$-algebra $\KM{}$ of Grothendieck rings of mixed Hodge modules constructed in Theorem \ref{thm:mixed-hodge-modules-cvar-alg}.

Recall that a \emph{groupoid}, $\cG$, is a category in which all the morphisms are invertible. We will denote by $\Grpd$ the category of small groupoids with natural transformations as morphisms. In particular, a group can be seen as a groupoid with a single object. Given a groupoid $\cG$, we will say that $\cG$ is \emph{finitely generated} if $\Obj{\cG}$ is finite and, for any object $a$ of $\cG$, $\Hom_{\cG}(a,a)$ (usually denoted $\cG_a$, the vertex group of $a$) is a finitely generated group. We will denote by $\Grpdo$ the full subcategory of $\Grpd$ of finitely generated groupoids.

\begin{rmk}
Given a groupoid $\cG$, two objects $a,b \in \cG$ are said to be connected if $\Hom_{\cG}(a,b)$ is not empty. In particular, this means that $\cG_a$ and $\cG_b$ are isomorphic groups so, in order to check whether $\cG$ is finitely generated, it is enough to check it on a point of every connected component.
\end{rmk}

Recall from \cite{Brown}, Chapter 6, that given a topological space $X$ and a subset $A \subseteq X$, we can define the \emph{fundamental groupoid of $X$ with respect to $A$}, $\Pi(X, A)$, as the category whose objects are the points in $A$ and, given $a, b \in A$, $\Hom(a,b)$ is the set of homotopy classes (with fixed endpoints) of paths in $X$ between $a$ and $b$. It is a straighforward check that this category is actually a groupoid and only depends on the homotopy type of the pair $(X,A)$. In particular, if $A = \left\{x_0\right\}$, $\Pi(X, A)$ has a single object whose automorphism group is $\pi_1(X, x_0)$, the fundamental group of $X$ based on $x_0$. For convenience, if $A$ is any set, not necessarily a subset of $X$, we will denote $\Pi(X, A) = \Pi(X, X \cap A)$ and we declare that $\Pi(\emptyset, \emptyset)$ is the singleton category.

\begin{rmk}
Let $M$ be a compact connected manifold (possibly with boundary) and let $A \subseteq X$ be finite. As we mentioned above, for any $a \in A$, $\Pi(M,A)_a = \pi_1(M, a)$. But a compact connected manifold has finitely generated fundamental group (c.f.\ \ref{ex:gamma-for-repr}, second item). Hence, in that case $\Pi(M,A)$ is a finitely generated groupoid.
\end{rmk}

As for the fundamental group, if $f: (X_1, A_1) \to (X_2, A_2)$ is a continuous map of pairs (i.e.\ $f(A_1) \subseteq A_2$) then we have an induced groupoid homomorphism $f_*: \Pi(X_1, A_1) \to \Pi(X_2, A_2)$. Hence, the fundamental groupoid of a pair gives a functor $\Pi: \Embpc \to \Grpdo$. 

\begin{figure}[h]
	\begin{center}
	\includegraphics[scale=0.35]{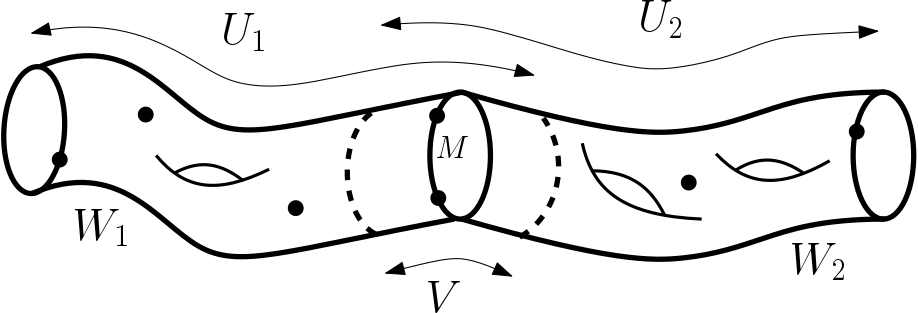}
	\caption{}
	\label{img:seifert-kampen}
	\end{center}
\vspace{-0.6cm}
\end{figure}

The functor $\Pi$ send the initial object of $\Embpc$, namely $(\emptyset, \emptyset)$, into the initial object of $\Grpdo$ since, by definition $\Pi(\emptyset, \emptyset) = 1$ is the unit group. Moreover, $\Pi$ sends gluing pushouts in $\Embpc$ to pushouts. In order to see that, let $(M, A) \hookrightarrow (\partial W_1^+, A) \subseteq (W_1, A_1)$ and $(M,A) \hookrightarrow (\partial W_2^+, A) \subseteq (W_2, A_2)$ be a gluing pushout. Let $V \subseteq W_1 \cup_M W_2$ be an open bicollar around $M$ such that $V \cap (A_1 \cup A_2) = A$. Set $U_1 = W_1 \cup V$ and $U_2 = W_2 \cup V$ (see Figure \ref{img:seifert-kampen}).

By construction, $\left\{U_1, U_2\right\}$ is an open covering of $W_1 \cup_M W_2$ such that $(U_1, A_1)$ is homotopically equivalent to $(W_1, A_1)$, $(U_2, A_2)$ is homotopically equivalent to $(W_2, A_2)$ and $(U_1 \cap U_2, (A_1 \cup A_2) \cap V) = (V, A)$ which is homotopically equivalent to $(M, A)$. Therefore, by Seifert-van Kampen theorem for fundamental groupoids (see \cite{Brown}, \cite{Brown:1967} and \cite{Higgins}) we have a pushout diagram induced by inclusions
	\[
\begin{displaystyle}
   \xymatrix
   {	\Pi(V, A)=\Pi(M, A) \ar[r] \ar[d] & \Pi(U_1, A_1) = \Pi(W_1, A_1) \ar[d] \\
   		\Pi(U_2, A_2) = \Pi(W_2, A_2) \ar[r]& \Pi(W_1 \cup_M W_2, A_1 \cup A_2)
   }
\end{displaystyle}   
\]
Hence $\Pi(W_1 \cup_M W_2, A_1 \cup A_2)$ is a pushout, as claimed.

Now, let $G$ be an algebraic group. Seeing $G$ as a (small) groupoid, we can consider the functor $\Hom_{\Grpd}(-, G): \Grpd \to \Sets$. Moreover, if $\cG$ is finitely generated, then $\Hom(\cG, G)=\Hom_{\Grpd}(\cG, G)$ has a natural structure of complex algebraic variety. To see that, pick a set $J = \left\{a_1, \ldots, a_s\right\}$ of objects of $\cG$ such that every connected component of $\cG$ contains exactly one element of $J$. Moreover, for any object $a$ of $\cG$, pick a morphism $f_a: a \to a_i$ where $a_i$ is the object of $J$ in the connected component of $a$. Hence, if $\rho: \cG \to G$ is a groupoid homomorphism, it is uniquely determined by the group representations $\rho_i: \cG_{a_i} \to G$ for $a_i \in J$ together with elements $g_a$ corresponding to the morphisms $f_a$ for any object $a$. Since the elements $g_a$ can be chosen without any restriction, if $\cG$ has $n$ objects, we have
$$
	\Hom(\cG, G) \cong \Hom(\cG_{a_1},G) \times \ldots \times \Hom(\cG_{a_s}, G) \times G^{n-s},
$$
and each of these factors has a natural algebraic structure as representation variety. This endows $\Hom(\cG, G)$ with an algebraic structure which can be shown not to depend on the choices.

\begin{defn}
Let $(X, A)$ be a pair of topological spaces such that $\Pi(X,A)$ is finitely generated (e.g.\ $X$ is a compact manifold and $A$ is finite). The variety
$$
	\Rep{G}(X,A) = \Hom(\Pi(X, A), G)
$$
is called the \emph{representation variety} of $(X, A)$ into $G$.
\end{defn}

\begin{rmk}
\begin{itemize}
	\item In particular, if $A$ is a singleton, we recover the usual representation varieties.
	\item If we drop out the requirement of $\cG$ being finitely generated, we can still endow $\Hom(\cG, G)$ with a scheme structure following the same lines (see, for example \cite{Nakamoto}). However, in general, this scheme is no longer of finite type. For this reason, in the definition of $\Bordp{n}$, we demand the subset $A \subseteq W$ to be finite.
\end{itemize}
\end{rmk}

Therefore, we can promote this functor to a contravariant functor $\Hom(-,G): \Grpdo \to \Var{k}$. Recall that $\Hom(-,G)$ sends colimits of $\Grpd$ into limits of $\Var{k}$ so, in particular, sends pushouts into pullback. Hence, the functor $\Rep{G} = \Hom(-,G) \circ \Pi: \Embpc \to \Grpdo$ satisfies the Seifert-van Kampen property. Therefore, by Proposition \ref{prop:field-theory}, we obtain an associated field theory
$$
	\Fld{\Rep{G}}: \Bordp{n} \to \Span{\Var{k}}.
$$
Observe that, since $\Pi$ and $\Hom(-,G)$ are both strict monoidal functors, $\Fld{\Rep{G}}$ is strict monoidal.

\begin{rmk}
The key point in the proof that $\Rep{G} = \Hom(-,G) \circ \Pi: \Embpc \to \Grpdo$ sends gluing pushouts into pullbacks was precisely the Seifert-van Kampen theorem. That is the reason for the name Seifert-van Kampen property.
\end{rmk}

Now, suppose that the ground field is $k = \CC$. By Theorem \ref{thm:mixed-hodge-modules-cvar-alg}, the $K$-theory of mixed Hodge modules, $\KM{}$, is a $\CVar$-algebra. From this, Theorem \ref{thm:physical-constr-sheaf} implies that, for any $n \geq 1$, there exists a lax monoidal TQFT
$$
	Z_{\Rep{G}, \KM{}}: \Bordp{n} \to \Modt{\K{\MHSq}}.
$$
For short, we will denote $\Zs{G} = Z_{\Rep{G}, \KM{}}$. Using the properties of this TQFT as explained in Remark \ref{rmk:expression-TQFT}, we have obtained the following result.

\begin{thm}\label{thm:existence-s-LTQFT}
There exists a lax monoidal TQFT of pairs, $\Zs{G}: \Bordp{n} \to \Modt{\K{\MHSq}}$ such that, for any $n$-dimensional connected closed manifold $W$ and any non-empty finite subset $A \subseteq W$, we have
$$
	\Zs{G}(W, A)\left(1\right) = \coh{\Rep{G}(W); \QQ} \otimes \coh{G; \QQ}^{|A|-1}
$$ 
where $1 = \QQ_0 \in \K{\MHSq}$ is the unit Hodge structure. 
\end{thm}

\begin{rmk}\label{rmk:Deligne-Hodge-soft-TQFT}
\begin{itemize}
	\item In particular, taking the Deligne-Hodge polynomial of the classes of Hodge structures above, we obtain that
$$
	\DelHod{\Rep{G}(W)} = \frac{\DelHod{\Zs{G}(W, A)(1)}}{\DelHod{G}^{|A|-1}}.
$$
Observe that $\DelHod{G}$ is usually known for the standard groups. For this reason, we say that $\Zs{G}$ computes the Deligne-Hodge polynomials of representation varieties.
	\item Consider bordisms $(W, A), (W, A'): (M_1, A_1) \to (M_2, A_2)$ with $A' \subseteq A$. This defines a $2$-cell $(W, A') \Rightarrow (W, A)$ in $\Bordp{n}$ that, under the field theory, becomes the $2$-morphism of spans
		\[
\begin{displaystyle}
   \xymatrix
   {	& \Rep{G}(W, A) = \Rep{G}(W, A') \times G^{|A|-|A'|} \ar[dd]^{\alpha} \ar[ld]_{i_1} \ar[dr]^{i_2}& \\
   	\Rep{G}(M_1, A_1)  & & \Rep{G}(M_2, A_2) \\
   	& \Rep{G}(W, A') \ar[ur]_{j_2} \ar[ul]^{j_1} &
   }
\end{displaystyle}   
\]
where $\alpha: \Rep{G}(W, A') \times G^{|A|-|A'|} \to \Rep{G}(W, A')$ is the projection onto the first component and $i_1, i_2, j_1, j_2$ are the restrictions of the representation varieties to the boundaries. In this way, we have $\Zs{G}(W, A') = (j_2)_!(j_1)^*$ and $\Zs{G}(W, A) = (i_2)_!(i_1)^* = (j_2)_!\alpha_!\alpha^*(j_1)^*$. Thus, $\alpha_!\alpha^*$ is the twist defining the $2$-cell $\Zs{G}(W,A') \Rightarrow \Zs{G}(W,A)$.
	\item As $\K{\MHM{}}$ is a $\CVar$-algebra, in particular, it can be used to obtain a pre-quantisation $\sQtm{\K{\MHM{}}}: \CVar \to \Bim{\K{\MHSq}}$, as described in Section \ref{sec:soft-TQFT}. With this datum, we can also construct a soft-TQFT, $\mathscr{Z}_{\Rep{G}, \K{\MHM{}}} = \sQtm{\K{\MHM{}}} \circ \cF_{\Rep{G}}: \Bord{n} \to \Bim{\K{\MHSq}}
$, that also computes Hodge structures of representation varieties.
	\item This result was proven for the first time in \cite{GPLM-2017}, but with a remarkable lower degree of generality. In a parallel work, Ben-Zvi, Gunningham and Nadler in \cite{Ben-Zvi-Gunningham-Nadler} have constructed a $2$-dimensional Extended Topological Quantum Field computing the Dolbeault cohomology of the character stack. For that purpose, instead of generalized representation varieties, they use the stack of $G$-local systems. Moreover, in that case, the construction of the functor is not explicit but it is done by means of Lurie's cobordisms hypothesis. Recall from Remark \ref{rmk:cobord-hypoth}, that this result classifies Extended TQFTs in terms of a special type of objects, the so-called fully dualizable objects (see \cite{Lurie-eTQFT:2009} for further information). It is an interesting future work to study in which way the construction of \cite{Ben-Zvi-Gunningham-Nadler} is related with the one introduced in this thesis.
\end{itemize}
\end{rmk}

\subsection{Almost-TQFT description and computational strategy}
\label{almost-TQFT-strategy}

From the lax monoidal TQFT of pairs described in the previous section, we can obtain an almost-TQFT, $\mathfrak{Z}_G: \Tubp{n} \to \Mod{\K{\MHSq}}$, following the recipe of Section \ref{sec:almost-TQFT-tubes}. That almost-TQFT allows the computation of Hodge structures of representation varieties. To be precise, from the previous construction we obtain the following result.

\begin{cor}\label{cor:almost-TQFT}
There exists an almost-TQFT of pairs, $\mathfrak{Z}_G: \Tubp{n} \to \Mod{\K{\MHSq}}$, such that, for any $n$-dimensional connected closed orientable manifold $W$ and any non-empty finite set $A \subseteq W$, we have
$$
	\mathfrak{Z}_G(W,A)(1) = \coh{\Rep{G}(W); \QQ)} \otimes \coh{G}^{|A|-1}
$$
where $1 \in \K{\MHSq}$ is the unit of the ring.
\end{cor}

\begin{figure}[h]
	\begin{center}
	\includegraphics[scale=0.5]{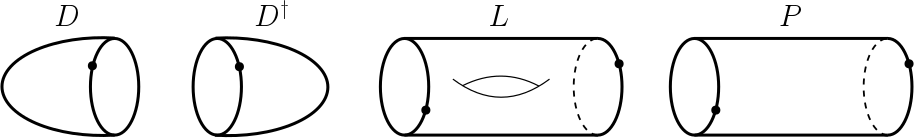}
	\caption{}
	\label{img:gen-tubes-pairs}
	\end{center}
\vspace{-0.8cm}
\end{figure}

In the case $n=2$ (surfaces), we can describe explicitly this almost-TQFT following the lines of Section \ref{sec:almost-TQFT-tubes}. In this case, since the objects of $\Bordp{2}$ are equipped with a finite configuration of points, we have to substitute the generators of Section \ref{sec:almost-TQFT-tubes} by the set $\Delta = \left\{D, D^\dag, L, P\right\}$ depicted in Figure \ref{img:gen-tubes-pairs}.

Recall that $\pi_1(S^1)=\ZZ$ so, for any point $\star \in S^1$, we have $\Rep{G}(S^1, \star)=\Hom(\ZZ, G)=G$. Therefore, on objects, the almost-TQFT assigns
$$
	\mathfrak{Z}_G(\emptyset) = \K{\MHM{1}} = \K{\MHSq}, \hspace{1cm} \mathfrak{Z}_G(S^1, \star) = \K{\MHM{G}}.
$$

Regarding morphisms, observe that the disc is simply connected so $\Rep{G}(D) = \Rep{G}(D^\dag) = 1$. Therefore, their images under the field theory, $\Fld{\Rep{G}}$, are
$$
\Fld{\Rep{G}}(D) = \left[1 \longleftarrow 1 \stackrel{i}{\longrightarrow} G\right],
	 \hspace{1cm} \Fld{\Rep{G}}\left(D^\dag\right) = \left[G \stackrel{\,i}{\longleftarrow} 1 \longrightarrow 1\right],
$$
where $i: 1 \to G$ is the inclusion map. Thus, their images under $\mathfrak{Z}_G$ are
$$
	\mathfrak{Z}_G(D) = i_!: \K{\MHSq}=\K{\MHM{1}} \to \K{\MHM{G}} \hspace{1cm} \mathfrak{Z}_G(D^\dag) = i^*: \K{\MHM{G}} \to \K{\MHM{1}} = \K{\MHSq}.
$$

\begin{figure}[h]
	\begin{center}
	\includegraphics[scale=0.2]{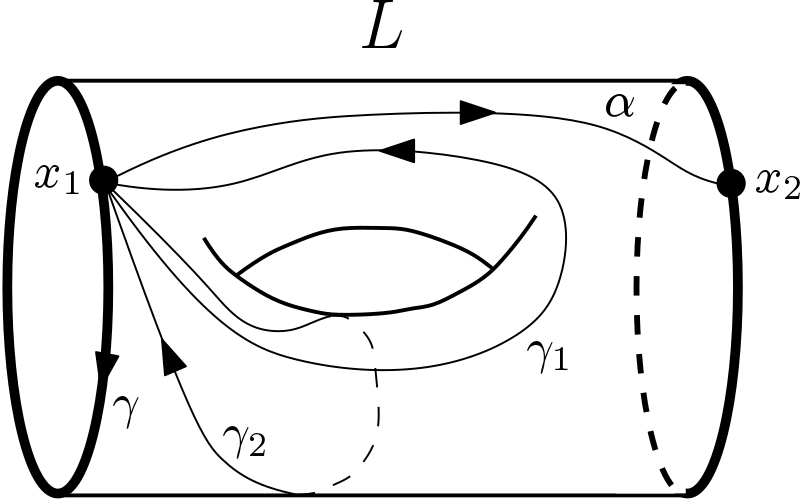}
	\caption{}
	\label{img:paths-L}
	\end{center}
\vspace{-0.6cm}
\end{figure}

For the holed torus $L: (S^1,\star) \to (S^1, \star)$ the situation is a bit more complicated. Let $L=(L_0, A)$ with $A = \left\{x_1, x_2\right\}$ the set of marked points of $L$, where $x_1$ lies in the ingoing boundary and $x_2$ lies in the outgoing boundary. Recall that $L_0$ is homotopically equivalent to a bouquet of three circles so its fundamental group is the free group with three generators. Thus, we can take $\gamma, \gamma_1, \gamma_2$ as the set of generators of $\Pi(L)_{x_1} = \pi_1(L_0, x_1)$ depicted in Figure \ref{img:paths-L} and $\alpha$ the shown path between $x_1$ and $x_2$.

With this description, $\gamma$ is a generator of $\pi_1(S^1, x_1)$ and $\alpha\gamma[\gamma_1, \gamma_2]\alpha^{-1}$ is a generator of $\pi_1(S^1, x_2)$, where $[\gamma_1, \gamma_2] = \gamma_1\gamma_2\gamma_1^{-1}\gamma_2^{-1}$ is the group commutator. Hence, since $\Rep{G}(L) =  \Hom(\Pi(L_0, A), G) = G^4$ we have that the associated field theory, $\Fld{\Rep{G}}(L)$, is the span
$$
\begin{matrix}
	G & \stackrel{p}{\longleftarrow} & G^4 & \stackrel{q}{\longrightarrow} & G \\
	g & \mapsfrom & (g, g_1, g_2, h) & \mapsto & hg[g_1,g_2]h^{-1}
\end{matrix}
$$
where $g, g_1, g_2$ and $h$ are the images of $\gamma, \gamma_1, \gamma_2$ and $\alpha$, respectively. Hence, we obtain that
$$
	\mathfrak{Z}_G(L): \K{\MHM{G}} \stackrel{p^*}{\longrightarrow} \K{\MHM{G^4}} \stackrel{q_!}{\longrightarrow} \K{\MHM{G}}.
$$

\begin{figure}[h]
	\begin{center}
	\includegraphics[scale=0.2]{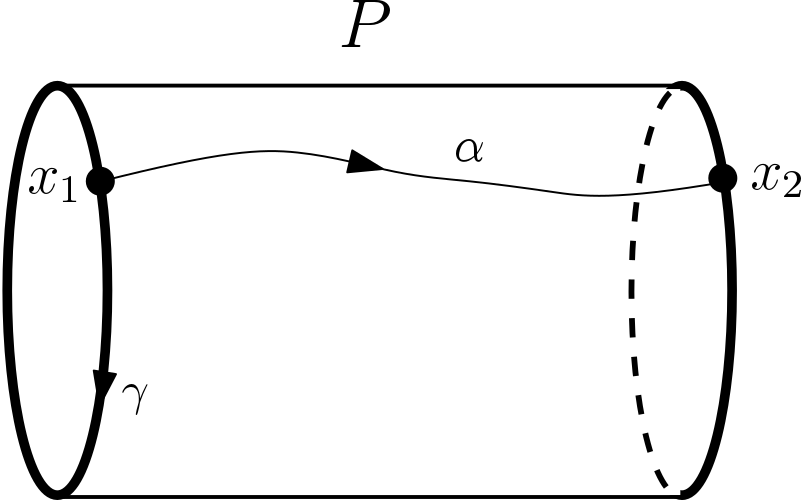}
	\caption{}
	\label{img:paths-P}
	\end{center}
\vspace{-0.6cm}
\end{figure}

For the morphism $P$, let $P=(S^1 \times [0,1], A)$ where $A = \left\{x_1, x_2\right\}$ with $x_1, x_2$ the ingoing and outgoing boundary points, respectively. Since $\pi_1(S^1 \times [0,1]) = \ZZ$, the fundamental groupoid $\Pi(P)$ has two vertices isomorphic to $\ZZ$ and, thus, $\Rep{G}(P) = \Hom(\Pi(P),G)=G^2$. Let $\gamma$ be a generator of $\pi_1(S^1, x_1)$ and $\alpha$ a path between $x_1$ and $x_2$, as depicted in Figure \ref{img:paths-P}.

Since $\alpha\gamma\alpha^{-1}$ is a generator of $\pi_1(S^1, x_2)$ we obtain that $\Fld{\Rep{G}}(P)$ is the span
$$
\begin{matrix}
	G & \stackrel{u}{\longleftarrow} & G^2 & \stackrel{v}{\longrightarrow} & G \\
	g & \mapsfrom & (g, h) & \mapsto & hgh^{-1}
\end{matrix}
$$
Hence, we have that
$$
	\mathfrak{Z}_G(P): \K{\MHM{G}} \stackrel{u^*}{\longrightarrow} \K{\MHM{G^3}} \stackrel{v_!}{\longrightarrow} \K{\MHM{G}}.
$$

Now, let $\Sigma_g$ be the closed orientable surface of genus $g$. Choose any $g+1$ points on $\Sigma_g$ so we have a decomposition of the bordism $\Sigma_g: \emptyset \to \emptyset$ as $\Sigma_g = D^\dag \circ L^g \circ D$.

\begin{figure}[h]
	\begin{center}
	\includegraphics[scale=0.4]{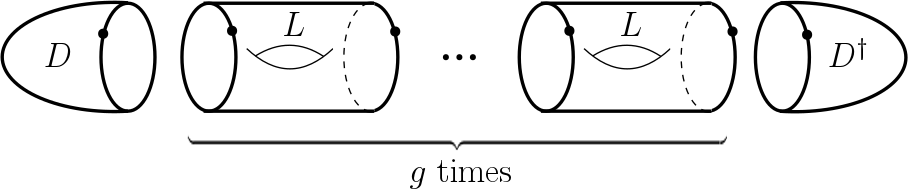}
	\end{center}
\vspace{-0.6cm}
\end{figure}

By Corollary \ref{cor:almost-TQFT}, if $\coh{G} \in \K{\MHSq}$ is not a zero divisor, we have that
$$
	\coh{\Rep{G}(\Sigma_g)} = \frac{1}{\coh{G}^g}\,{\mathfrak{Z}_G(D^\dag) \circ \mathfrak{Z}_G(L)^g \circ \mathfrak{Z}_G(D) (1)} = \frac{1}{\coh{G}^g}\,{i^* \circ (q_!p^*)^g \circ i_! (1)}.
$$
Therefore, in order to compute the Hodge structure of a representation variety, we just need to compute three linear maps
$$
	i_!: \K{\MHSq} \to \KM{G} \hspace{0.5cm} i^*: \KM{G} \to \K{\MHSq} \hspace{0.5cm} q_!p^*: \KM{G} \to \KM{G} 
$$
Actually, as we will show in \ref{sec:discs-tubes}, $i_!$ and $i^*$ are quite easy to work with. Therefore, the most important part is to understand the linear map $q_!p^*: \KM{G} \to \KM{G}$.

\begin{rmk}
These kind of computations were carried out in the paper \cite{MM} for $G=\SL{2}(\CC)$ and in \cite{Martinez:2017} for $G=\PGL{2}(\CC)$. There, the underlying calculations hid the almost-TQFT described here. However, since this machinery was not available at that time, all the computations were done by performing an explicit pasting. For a more precise description of the interpretation of these results in the context of TQFTs, see Section \ref{sec:genus-tube}.
\end{rmk}

\subsection{Parabolic case}
\label{sec:parabolic-case-TQFT}

The previous construction can be easily extended to the parabolic case. Following the construction above, it is just necessary to adapt the geometrisation functor to the parabolic context.

As above, let us fix an algebraic group $G$ and, as parabolic data, we choose for $\Lambda$ a collection of locally closed subsets of $G$ which are closed under conjugacy. Consider a compact manifold $M$, a finite subset $A \subseteq M$ and a parabolic structure $Q = \left\{(S_1, \lambda_1), \ldots, (S_s, \lambda_s)\right\}$ on $M$ (see Section \ref{sec:TQFT-over-sheaf} and Section \ref{sec:parabolic-repr-var}). Denote $S = \bigcup_i S_i$. In analogy with Section \ref{sec:intro-representation}, the representation variety of $\Pi(M,A)$ into $G$ with parabolic structure $Q$ is
$$
	\Rep{G}(M, A, Q) = \left\{\rho: \Pi(M-S,A) \to G\,\left|\,\,\, \begin{matrix}
	\rho(\gamma_i) \in \lambda_i \textrm{ for all }\gamma_i \\
	\textrm{ positive meridian around }S_i\\
\end{matrix}
 \right.\right\}.
$$
Observe that positive meridians with different basepoints are conjugated, so the notion of positive meridian is well defined even for fundamental groupoids. Actually, if $M = M_1 \sqcup \ldots \sqcup M_s$ is the decomposition of $M$ into connected components, we have an identification
$$
	\Rep{G}(M, A, Q) = \Rep{G}(M_1, Q) \times \ldots \times \Rep{G}(M_s, Q) \times G^{|A|-s}
$$
Hence, $\Rep{G}(M, A, Q)$ has also a natural algebraic structure. Using these varieties, we define the functor $\Rep{G}: \Embpparc{\Lambda} \to \Var{k}$ that assigns, to $(M,A) \in \Embpparc{\Lambda}$, the variety $\Rep{G}(M,A)$.

As in the non-parabolic case, $\Rep{G}$ has the Seifert-van Kampen property. Hence, for $k = \CC$, this functor together with the $\CVar$-algebra $\KM{}$ gives us, by means of Theorem \ref{thm:physical-constr-sheaf}, a lax monoidal TQFT
$$
	\Zs{G} = Z_{\Rep{G}, \KM{}}: \Bordppar{n}{\Lambda} \to \Modt{\K{\MHSq}}.
$$
Again, by the description of Remark \ref{rmk:expression-TQFT}, this TQFT satisfies that, for any closed connected $n$-dimensional manifold $W$, any non-empty finite subset $A\subseteq W$ and any parabolic structure $Q$, we have
$$
	\Zs{G}(W, A, Q)(1) = \coh{\Rep{G}(W, Q); \QQ} \otimes \coh{G; \QQ}^{|A|-1}.
$$

Furthermore, following the construction of Section \ref{sec:almost-TQFT-tubes}, we can also modify this result to obtain an almost-TQFT, $\mathfrak{Z}_G: \Tubppar{n}{\Lambda} \to \Mod{\K{\MHSq}}$, where $\CTubpp{n}(\Lambda)$ is the category of tubes of pairs with parabolic data $\Lambda$. 

\begin{figure}[h]
	\begin{center}
	\includegraphics[scale=0.2]{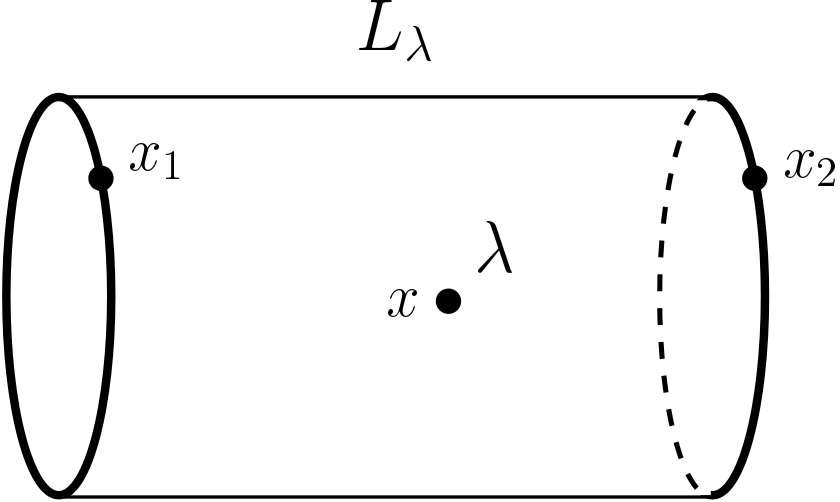}
	\caption{}
	\label{img:parabolic-tube}
	\end{center}
\vspace{-0.6cm}
\end{figure}

In the particular case $n=2$, recall from Example \ref{ex:parab-struct} that a parabolic structure is given by a bunch of points on the surface. Hence, in order to obtain a set of generators of $\CTubpp{n}(\Lambda)$, it is enough to consider the elements of the set of generators $\Delta$ previouly defined with no parabolic structure and to add a collection of tubes $L_\lambda = (S^1 \times [0,1], \left\{x_1, x_2\right\}, \left\{(x, \lambda)\right\})$ for $\lambda \in \Lambda$. Here, $x$ is any interior point of $S^1 \times [0,1]$ with parabolic structure $\lambda$ and $x_1, x_2$ are points in the ingoing and outgoing boundaries, repectively, as depicted in Figure \ref{img:parabolic-tube}.

In this case, observe that $\pi_1((S^1 \times [0,1]) - \left\{x\right\})$ is the free group with two generators and that the fundamental groupoid of $L_\lambda$ has two vertices. These generators can be taken to be around the incoming boundary and around the marked point so $\Rep{G}(L_\lambda) = G^2 \times \lambda$. Thus, the image under the field theory of $L_\lambda$ is the span
$$
\begin{matrix}
	G & \stackrel{r}{\longleftarrow} & G^2 \times \lambda & \stackrel{s}{\longrightarrow} & G \\
	g & \mapsfrom & (g, h, \xi) & \mapsto & hg\xi h^{-1}
\end{matrix}
$$
Hence, the image $\mathfrak{Z}_G(L_\lambda)$ is the morphism
$$
	\mathfrak{Z}_G(L_\lambda): \K{\MHM{G}} \stackrel{r^*}{\longrightarrow} \K{\MHM{G^2 \times \lambda}} \stackrel{s_!}{\longrightarrow} \K{\MHM{G}}.
$$
With this description, as an extension of Corollary \ref{cor:almost-TQFT}, we have proven the following result.

\begin{thm}\label{thm:almost-tqft-parabolic}
Let $\Sigma_g$ be a closed oriented surface of genus $g$ and $Q$ a parabolic structure on $\Sigma_g$ with $s$ marked points with data $\lambda_1, \ldots, \lambda_s \in \Lambda$. Then, 
$$
	\coh{\Rep{G}(\Sigma_g, Q)} = \frac{1}{\coh{G}^{g+s}}{\mathfrak{Z}_G(D^\dag) \circ \mathfrak{Z}_G(L_{\lambda_s}) \circ \ldots \circ \mathfrak{Z}_G(L_{\lambda_1}) \circ \mathfrak{Z}_G(L)^g \circ \mathfrak{Z}_G(D) (1)}.
$$
\end{thm}

\section{Other TQFTs for representation varieties}
\label{sec:other-tqft-repr}

In this section, we will describe several variants of the previous TQFT more focused on character varieties than on representation varieties. As we will see, these solutions are not as elegant as the previous one. However, with a view towards explicit computations, some of them are more suitable so they will be used in the upcoming sections.

\subsection{TQFT for character varieties}
\label{sec:TQFT-char-var}

Let $G$ be an algebraic reductive group. In Section \ref{sec:TQFT-for-repr}, we built the contravariant functor $\Rep{G}: \Embpc \to \Var{k}$. However, as for representation of the fundamental group, for a pair $(M, A) \in \Embpc$, the associated representation variety $\Rep{G}(M,A) = \Hom(\Pi(M,A),G)$ has a natural action of $G$ by conjugation. That is, given $g \in G$, $\rho: \Pi(M, A) \to G$ and a path $\gamma: a \to a'$ the action is given by
	$$
		(g \cdot \rho)(\gamma) = g \rho(\gamma) g^{-1}.
	$$
Therefore, we can also consider the associated character variety $\Char{G}(M,A) = \Rep{G}(M, A) \sslash G$. 

Even more, suppose that $f: (M_1, A_1) \to (M_2, A_2)$ is a continuous map with $f(A_1) \subseteq A_2$. By composing the induced morphism $\Rep{G}(f): \Rep{G}(M_2, A_2) \to \Rep{G}(M_1, A_1)$ with the quotient map of $\Rep{G}(M_1, A_1)$, we obtain a $G$-equivariant morphism $\Rep{G}(M_2, A_2) \to \Char{G}(M_1, A_1)$. Hence, by the universal property of GIT quotients, it can be promoted to a morphism $\Char{G}(f): \Char{G}(M_2, A_2) \to \Char{G}(M_1, A_1)$. Therefore, as in Section \ref{sec:TQFT-for-repr}, this description gives us a contravariant functor $\Char{G}: \Embpc \to \Var{k}$. However, this functor has not the Seifert-van Kampen property since, in general, if $(W_1, A_1)$ and $(W_2, A_2)$ are compact manifolds to be pasted along a common boundary $(M, A)$, we have that $\Char{G}(W_2 \circ W_1, A_1 \cup A_2)$ is not the pullback of $\Char{G}(W_1, A_1)$ and $\Char{G}(W_2, A_2)$ along $\Char{G}(M, A)$.

In any case, $\Char{G}: \Embpc \to \Var{k}$ does send the initial object into the final object. Hence, by Remarks \ref{rmk:field-theory-very-lax} and \ref{rmk:field-theory-very-lax-sheaf}, the associated field theory $\Fld{\Char{G}}: \Bordp{n} \to \Span{\Var{k}}$ is a lax $2$-functor. In particular, if $k = \CC$, using the $\CVar$-algebra $\KM{}$, we obtain a very lax Topological Quantum Field Theory
$$
	Z_{\Char{G}, \KM{}}: \Bordp{n} \to \Modt{\K{\MHSq}}.
$$
The importance of this TQFT is that, if $W$ is a closed connected $n$-dimensional manifold and $A \subseteq W$ is a finite subset, then, by construction
$$
	Z_{\Char{G}, \KM{}}(W, A)\left(1\right) = \coh{\Char{G}(W, A); \QQ}.
$$

\begin{rmk}
\begin{itemize}
	\item Despite that $\Rep{G}(W, A) = \Rep{G}(W) \times G^{|A|-1}$, the same decomposition no longer holds for character varieties since $G$ acts simultaneously on the components. For this reason, the later cohomology no longer decomposes as in Theorem \ref{thm:existence-s-LTQFT}
	\item Since our main aim is to compute Deligne-Hodge polynomials of character varieties, this TQFT seems more suited than the one in Section \ref{sec:TQFT-for-repr}. However, its very lax condition, together with the involved structure of the intermediate varieties of the spans, make this TQFT not too much useful for explicit computations.
	\item The construction above can be extended to the parabolic case without major modifications, giving rise to a very lax TQFT $Z_{\Char{G},\KM{}}: \Bordppar{n}{\Lambda} \to \Modt{\K{\MHSq}}$.
\end{itemize}
\end{rmk}

\subsection{Piecewise regular varieties}
\label{sec:piecewise-regular-varieties}

In this section, we extend slightly the notion of an algebraic variety in order to allow topological spaces which are not algebraic varieties, but disjoint union of algebraic varieties, as in the case of orbit spaces. This will allow us to considier a variation of the TQFT above that, instead of GIT quotients, uses usual quotients as orbit spaces

\begin{defn}
Let $k$ be an algebraically closed field. We define the \emph{category of piecewise varieties over $k$}, $\PVar{k}$ as the category given by:
\begin{itemize}
	\item Objects: The objects of $\PVar{k}$ is the Grothendieck semi-ring of algebraic varieties over $k$ (see Remark \ref{rmk:KVar}). This is the semi-ring generated by the symbols $[X]$, for $X$ an algebraic variety, with the relation $[X] = [U] + [Y]$ if $X = U \sqcup Y$ with $U \subseteq X$ open and $Y \subseteq X$ closed. It is a semi-ring with disjoint union as sum and cartesian product as product. For further information, see \cite{kuber:2015}.
	\item Morphisms: Given $[X], [Y] \in \PVar{k}$ a morphism $f: [X] \to [Y]$ is given by a function $f: X \to Y$ that decomposes as a disjoint union of regular maps. More precisely, it must exist decompositions $[Y] = \sum_i [Y_i]$ and $[X] = \sum_{i,j} [X_{ij}]$ into algebraic varieties such that $f = \bigsqcup\limits_{i,j} f_{ij}$ with $f_{ij}: X_{ij} \to Y_i$ a regular morphism.
	\item Composition is given by the usual composition of maps. Observe that, given $f: [X] \to [Y]$ and $g: [Y] \to [Z]$, by Chevalley theorem \cite{EGAI}, we can find a common decomposition $Y = \sqcup_{i,j} Y_{ij}$ such that
	$$
	f = \bigsqcup_{i,j,k} f_{ijk}: \bigsqcup_{i,j,k} X_{ijk} \to \bigsqcup_i Y_{ij}, \hspace{1cm} g = \bigsqcup_{i,j} g_{ij}: \bigsqcup_{i,j} Y_{ij} \to \bigsqcup_i Z_i
$$
decompose as regular maps. Hence, $g \circ f = \bigsqcup\limits_{i,j,k} \left(g_{ij} \circ f_{ijk}\right)$ is a decomposition into regular maps of $g \circ f$.
\end{itemize}
\end{defn}

\begin{rmk}\label{rmk:forgetful-piecewise-var}
\begin{itemize}
	\item There is a functor $\Var{k} \to \PVar{k}$ that sends $X \mapsto [X]$ and analogously for regular maps. On the other way around, there is a forgetful functor $\PVar{k} \to \Sets$ that recovers the underlying set of a piecewise variety.
	\item The functor $\Var{k} \to \PVar{k}$ is not an embedding of categories. For example, let $X = \left\{y^2 = x^3\right\}$ a cuspidal cubic plane curve. Observe that, removing the origin in $X$, we have a decomposition $[X] = [X - \star] + [\star] = [\mathbb{A}^1 - \star] + [\star] = [\mathbb{A}^1]$. Hence, the images of $X$ and $\mathbb{A}^1$ under this functor agree, even though they are not isomorphic.
\end{itemize}
\end{rmk}
	
\begin{ex}
The affine line with a double point at the origin $X$ is not an algebraic variety. However, it is a piecewise variety since we can decompose $[X] = [\mathbb{A}^1] + [\star]$ with both pieces locally closed subsets.
\end{ex}

\begin{rmk} There is a slightly more abstract (and precise) way of defining $\PVar{k}$ that we sketch below:
\begin{itemize}
	\item Let $\cU$ be the collection of finite sets of algebraic varieties. Given $\left\{X_1, \ldots, X_r\right\} \in \cU$, we declare the binary relation $\left\{X_1, \ldots, X_r\right\} \prec \left\{X_1, \ldots, \hat{X_i},\ldots, X_r, U_i, Y_i\right\}$ where $X_i = U_i \sqcup Y_i$ with $Y_i \subseteq X_i$ closed and $\hat{X_i}$ denotes the omision of $X_i$. If $\leq$ is the partial order obtained by taking the transitive closure of $\prec$, then $(\cU, \leq)$ is a partially ordered set (also an abstract rewriting system, see \cite{Leeuwen:1990}).
	\item Given $X, Y \in \cU$, we define the relation $X \sim Y$ if the open intervals $[X, \infty) \cap [Y, \infty) \neq \emptyset$. As $\leq$ is a confluent order, $\sim$ is an equivalence relation. The class of $X \in \cU$, $[X]$, can be interpreted as all the possible partitions of $X$ into algebraic varieties.
	\item In other terms, the interval $[X, \infty)$ is a (principal) filter on $(\cU, \leq)$ and $[X]$ is the associated ultrafilter (see \cite{Marker:2002}). If $[X]$ has a least element, it is also a principal filter, and such a least element is understood as a minimal decomposition of $X$.
	\item Suppose that, for $X' = \left\{X_i\right\} \in [X]$ and $Y' = \left\{Y_i\right\} \in [Y]$ big enough we have a set of regular maps $\sqcup_i f_i: X_i \to Y_i$. By Chevalley theorem, given any other representatives bigger than $X'$ and $Y'$, we can refine them to obtain another representation of $f$ by restriction. Hence, $f$ does not depend on the chosen representatives.
\end{itemize}
\end{rmk}

Let $G$ be an algebraic group acting on a variety $X$. Maybe after restricting to the open subset of semi-stable points, we can suppose that a (good) GIT quotient $\pi: X \to X \sslash G$ exists (see Remark \ref{sec:review-git}). On $X$, we find the subset $X^1 \subseteq X$ of poly-stable points whose orbits have maximum dimension (see Remark \ref{rmk:polystable-points}). It is an open subset by an adaptation of Proposition 3.13 of \cite{Newstead:1978} and it is non-empty since the set of poly-stable points is so. On the poly-stable points, the GIT quotient is an orbit space so we have a $G$-invariant regular map $\pi_1: X_1 \to X_1 / G$.
Now, let $Y = X - X_1$. As $Y$ is closed on $X$ and the action of $G$ restricts to an action on $Y$ and we can repeat the argument to obtain a regular $G$-invariant map $\pi_2: X_2 \subseteq Y \to X_2 / G$, where $X_2 \subseteq Y$ is the set of poly-stable points of $Y$ of maximum dimension. Repeating this procedure, we obtain a stratification $X = X_1 \sqcup \ldots \sqcup X_r$ and a set of regular maps $\pi_i: X_i \to X_i/G$ where each $X_i / G$ has a natural algebraic structure.

\begin{defn}\label{defn:piecewise-quotient}
Let $G$ be a reductive group and let $X$ be an algebraic variety. With the decomposition above, the \emph{piecewise quotient} of $X$ by $G$, denoted by $[X / G]$ is the object of $\PVar{k}$ given by
$$
	[X/G] = [X_1/G] + \ldots + [X_r / G].
$$
We also have a piecewise quotient morphism $\pi = \sqcup_i \pi_i: [X] \to [X/G]$ where $\pi_i: X_i \to X_i/G$ are the projections on the orbit space of each of the strata.
\end{defn}

\begin{rmk}
Actually, since all the strata considered are made of poly-stable points over the previous statum, no orbits are identified under the quotient. Hence, each orbit of $X$ is contained in exactly one of the strata $X_1, \ldots, X_r$. Therefore, the underlying set of $[X/G]$ is, precisely, $X/G$ so the previous construction can be restated as that the orbit space of $X$ by $G$ has a piecewise variety structure.
\end{rmk}

The piecewise quotients have similar properties than GIT quotients but in the category $\PVar{k}$. In particular, they satisfy a similar universal property.

\begin{prop}
Let $G$ be a reductive group acting on a variety $X$ and let $f: [X] \to [Y]$ be a morphism of piecewise varietes which is $G$-invariant. Then, there exists a unique piecewise morphism $\tilde{f}: [X / G] \to [Y]$ such that $\tilde{f} \circ \pi = f$.
	\[
\begin{displaystyle}
   \xymatrix
   {	[X] \ar[r]^f \ar[d]_{\pi} & [Y] \\
   		[X/G] \ar@{--{>}}[ru]_{\tilde{f}}
   }
\end{displaystyle}   
\]
\begin{proof}
Suppose that $f = \sqcup \, f_{jk}: \sqcup X_{jk}' \to \sqcup Y_j$ is a decomposition of $f$ into regular morphisms. We refine this stratification of $X$ by taking $X_{ijk} = X_i \cap X_{jk}'$, where the statum $X_i$ are the ones considered in Definition \ref{defn:piecewise-quotient}.
Since $f$ is $G$-invariant, $f_{ijk} = f|_{X_{ijk}}: X_{ijk} \to Y_j$ is a $G$ invariant morphism of algebraic varieties. Thus, it descends to a regular morphism $\tilde{f}_{ijk}: X_{ijk}/G \to Y_j$. Putting them together, we obtain a map
$$
	\bigsqcup_{i,j,k} \tilde{f}_{ijk}: \bigsqcup_{i,j,k} X_{ijk}/G \to \bigsqcup_{j} Y_j.
$$
Since $X_i/G = \bigsqcup_{j,k} X_{ijk}/G$, we have a decomposition $[X/G] = \sum_{i,j,k} [X_{ijk}/G]$ so this map gives rise to a piecewise morphism $\tilde{f}: [X/G] \to [Y]$. By construction, $\tilde{f} \circ \pi = f$ and it is unique.
\end{proof}
\end{prop}

\begin{defn}
Let $G$ be an algebraic group. Given a pair of topological spaces $(X, A)$ such that $\Pi(X,A)$ is finitely generated, we define the \emph{piecewise character variety} as the object of $\PVar{k}$
$$
	\CharW{G}(X, A) = [\Rep{G}(X,A) / G].
$$ 
Here, $G$ acts on $\Rep{G}(X,A)$ by conjugation (see Section \ref{sec:char-var}).
\end{defn}

As an application, we can construct a TQFT focused on piecewise character varieties, in analogy with the one of Section \ref{sec:TQFT-char-var} for genuine character varieties. For that, given a reductive algebraic group $G$, we can define a new geometrisation functor $\CharW{G}: \Embpc \to \PVar{k}$ by replacing the GIT quotient of Section \ref{sec:TQFT-char-var} by a piecewise quotient. With respect to morphisms, given an embedding $f: (M_1, A_1) \to (M_2, A_2)$, the induced map to the quotient $\Rep{G}(M_2, A_2) \to \Rep{G}(M_1, A_1) \to \CharW{G}(M_1, A_1)$ is a $G$-invariant piecewise morphism. Hence, it descends to a piecewise map $\CharW{G}(f): \CharW{G}(M_2, A_2) \to \CharW{G}(M_1, A_1)$.
However, as in Section \ref{sec:TQFT-char-var}, $\CharW{G}$ no longer sends gluing pushouts into pullbacks but it sends the initial object into the final object of $\PVar{k}$ (which is still $\star$). Hence, it gives a monoidal lax $2$-functor $\Fld{\CharW{G}}: \Bordp{n} \to \Span{\PVar{k}}$.

With respect to the quantization part, we have the following fundamental result.

\begin{prop}
The $\CVar$-algebra $\KM{}$ extends to a $\PVarC$-algebra.
\begin{proof}
For an object $[X] = X_1 + \ldots + X_r \in \PVarC$, we define $\KM{[X]} = \bigoplus_i \KM{X_i}$. This definition does not depend on the chosen representative of $[X]$. Indeed, if $X_i = X_i^1 \sqcup X_i^2$, with $X_i^1 \subseteq X_i$ closed, and $i_k: X_i^k \to X_i$ are the inclusions for $k = 1,2$, then, by Proposition \ref{prop:strat-mixed-hodge-modules}, the map $(i_1)_! + (i_2)_! : \KM{X_i^1} \oplus \KM{X_i^2} \to \KM{X_i}$ is an isomorphism.

For a morphism $f: [X] \to [Y]$, given as a set of regular maps $f_{i,j}: X_{i,j} \to Y_i$, we define
	$$
		f^* = \bigoplus_i \left[\bigoplus_jf_{i,j}^*: \KM{Y_i} \to \bigoplus_j \KM{X_{i,j}}\right], \hspace{0.5cm}  f_! = \bigoplus_i \left[\sum_j(f_{i,j})_!: \bigoplus_j \KM{X_{i,j}} \to \KM{Y_i}\right].
	$$
Again, they do not depend on the chosen representatives. In order to prove it, let us suppose that $f:X \to Y$ is a morphism between algebraic varieties. If $X = X_1 \sqcup X_2$, $i_k: X_k \to X$ are the inclusions and $f_k = f|_{X_k}: X_k \to X$ are the corresponding restrictions, we have commutative diagrams

\vspace{-0.6cm}
\begin{minipage}[t]{0.5\textwidth}
	\[
\begin{displaystyle}
   \xymatrix
   {	
   		\KM{X_1} \oplus \KM{X_2} \ar[rr]^{\;\;\;\;\;\;\;\;\;(f_1)_! + (f_2)_!} \ar[d]_{(i_1)_! + (i_2)_!} && \KM{Y} \\
   		\KM{X} \ar[rru]_{f_!}&&
   }
\end{displaystyle}   
\]
\end{minipage}
\begin{minipage}[t]{0.5\textwidth}
	\[
\begin{displaystyle}
   \xymatrix
   {	
   		\KM{Y} \ar[rr]^{f_1^* \oplus f_2^* \;\;\;\;\;\;\;\;\;\;} \ar[rrd]_{f^*} && \KM{X_1} \oplus \KM{X_2} \ar[d]^{(i_1)_! + (i_2)_!} \\
   		&& \KM{X}
   }
\end{displaystyle}   
\]
\end{minipage}

The commutativity of the first diagram is immediate and the commutativity of the second one follows from the fact that $(i_1)_!f_1^* + (i_2)_!f_2^* = (i_1)_!(i_1)^*f^* + (i_2)_!(i_2)^*f^* = f^*$ since $(i_1)_!(i_1)^* + (i_2)_!(i_2)^* = 1_{\KM{X}}$ as shown in the proof of Proposition \ref{prop:strat-mixed-hodge-modules}. On the other hand, suppose that we have a decomposition $Y = Y_1 \sqcup Y_2$ with inclusions $j_k: Y_k \to Y$ and we set $X_k = f^{-1}(Y_k)$ with inclusions $i_k: X_k \to X$. Setting $f_k = f|_{X_k}: X_k \to Y_k$, we have commutative diagrams

\vspace{-1cm}
\begin{minipage}[t]{0.5\textwidth}
	\[
\begin{displaystyle}
   \xymatrix
   {	
   		\KM{X_1} \oplus \KM{X_2} \ar[rr]^{\;\;(f_1)_! \oplus (f_2)_!} \ar[d]_{(i_1)_! + (i_2)_!} && \KM{Y_1} \oplus \KM{Y_2} \ar[d]^{(j_1)_! + (j_2)_!} \\
   		\KM{X} \ar[rr]_{f_!} && \KM{Y}
   }
\end{displaystyle}   
\]
\end{minipage}
\begin{minipage}[t]{0.5\textwidth}
	\[
\begin{displaystyle}
   \xymatrix
   {	
   		\KM{Y_1} \oplus \KM{Y_2} \ar[rr]^{\;\;f_1^* \oplus f_2^*} \ar[d]_{(j_1)_! + (j_2)_!} &&  \KM{X_1} \oplus \KM{X_2} \ar[d]^{(i_1)_! + (i_2)_!} \\
   		\KM{Y} \ar[rr]_{f^*} && \KM{X}
   }
\end{displaystyle}   
\]
\end{minipage}

The first one follows from the fact that $[(j_1)_! + (j_2)_!] \circ [(f_1)_! \oplus (f_2)_!] = (j_1 \circ f_1)_! + (j_2 \circ f_2)_! = (f \circ i_1)_! + (f \circ i_2)_! = f_! \circ [(i_1)_! + (i_2)_!]$. For the second one, observe that $[(i_1)_! + (i_2)_!] \circ (f_1^* \oplus f_2^*) = (i_1)_!f_1^* + (i_2)_!f_2^* = f^* \circ [(j_1)_! + (j_2)_!]$, where the last equality follows from the Beck-Chevalley property of $\KM{}$ and the fact that, for $k = 1,2$, the square
\[
\begin{displaystyle}
   \xymatrix
   {	
  		X_k \ar[r]^{f_k} \ar[d]_{i_k} & Y_k \ar[d]^{j_k} \\
  		X \ar[r]_{f} & Y
   }
\end{displaystyle}   
\]
is a pullback.
\end{proof}
\end{prop}

With this results at hand, if $G$ is a complex group, Theorem \ref{thm:physical-constr-sheaf} gives the following result.

\begin{cor}
The assignment
$
	Z_{\CharW{G}, \KM{}}: \Bordp{n} \to \Modt{\K{\MHSq}}
$
is a very lax Topological Quantum Field Theory.
\end{cor}

\subsection{Geometric and reduced TQFT}
\label{sec:geometric-reduced-TQFT}

In this section, we will describe an intermediate step between the TQFTs for character varieties introduced in Section \ref{sec:TQFT-char-var} and the one for piecewise character varieties presented in the previous section. Thanks to its mixed nature, this TQFT will be suited for explicit calculations of mixed Hodge structures of representation varieties. Actually, as it will be clear from the computations in Sections \ref{sec:tube-J+} and \ref{sec:genus-tube}, this procedure underlies the calculations of \cite{LMN}, \cite{MM:2016} and \cite{MM}. The method of this section is based on the reduction of a TQFT of Section \ref{sec:reduction-TQFT}.

From a complex algebraic group $G$, we can consider the representation variety geometrisation $\Rep{G}: \Embpc \to \CVar$ as in Section \ref{sec:TQFT-for-repr}. After the functor $\CVar \to \PVarC$ of Remark \ref{rmk:forgetful-piecewise-var}, we may consider the reduction $\pi$ that assigns, to $(M, A) \in \Bordp{n}$, the piecewise character variety $\CharW{G}(M, A) = [\Rep{G}(M, A) / G]$ together with the piecewise quotient morphism $\pi_{(M, A)}: \Rep{G}(M, A) \to \CharW{G}(M, A)$. As explained in Section \ref{sec:reduction-TQFT}, this reduction can be used to modify the standard TQFT, $\Zs{G}$ of Theorem \ref{thm:existence-s-LTQFT}, in order to obtain a very lax TQFT $\Zg{G} = (\Zs{G})_{\pi}: \Bordp{n} \to \Modt{\K{\MHSq}}$. 

Let us we denote by $\cV_{(M, A)}$ the submodule of $\KM{\CharW{G}(M, A)}$ generated by the images $\Zg{G}(W_r) \circ \ldots \circ \Zg{G}(W_1)(1) \in \cV_{(M, A)}$ for strict tubes $W_r \circ \ldots W_1: \emptyset \to (M, A)$, as in Section \ref{sec:reduction-TQFT}. If the morphism $\eta_{(M, A)} = (\pi_{(M, A)})_! \circ (\pi_{(M, A)})^*: (\cV_{(M, A)})_0 \to (\cV_{(M, A)})_0$ is invertible as a $\Ko{\MHSq}$-module homomorphism for all $(M, A) \in \Bordp{n}$, then the morphisms $\cZg{G}(W,A) = \Zg{G}(W,A) \circ \eta_{(M, A)}$ give rise to a strict TQFT, the $\pi$-reduction of $\Zs{G}$ (Corollary \ref{cor:existence-reduction-TQFT}). Here, $(\cV_{(M, A)})_0$ stands for the localization $S^{-1}\cV_{(M, A)}$ by an appropiate multiplicative system $S \subseteq \K{\MHSq}$ (for example, $S$ can be the set of non-zero divisors of the ring $\K{\MHSq}$).

Moreover, this construction can be extended to the parabolic case without major modifications. Suppose that we have a set of parabolic data $\Lambda$ of conjugacy invariant subvarieties of $G$. Then, given $(M, A, Q) \in \Bordppar{n}{\Lambda}$, as $\Char{G}(M, A,Q)$ has a natural $G$-action, we can also consider the piecewise quotient $\CharW{G}(M, A,Q)= [\Rep{G}(M, A,Q) / G]$. Hence, we can extend the reduction $\pi$ to the parabolic context giving rise to a prerreduction, $\Zg{G}: \Bordppar{n}{\Lambda} \to \Modt{\K{\MHSq}}$. If the morphisms $\eta_{(M, A, Q)}$ are also invertible after some localization, then we can also define the reduction $\cZg{G}: \Bordppar{n}{\Lambda} \to \Modt{\Ko{\MHSq}}$.

\begin{defn}
Let $G$ be a complex algebraic group and $\Lambda$ a collection of conjugacy closed subvarieties of $G$. The Topological Quantum Field Theories
$$
	\Zg{G}: \Bordppar{n}{\Lambda} \to \Modt{\K{\MHSq}}, \hspace{1cm} \cZg{G}: \Bordppar{n}{\Lambda} \to \Modt{\Ko{\MHSq}},
$$
are called the \emph{geometric TQFT} and the \emph{reduced TQFT} for representation varieties, respectively. 
\end{defn}
Such a TQFTs will be extensively used along Section \ref{sec:sl2-repr-var} since, as we will see, they allow us to compute Hodge structures on representation varieties in a simple way.

\begin{ex}
In the case of surfaces ($n=2$), the morphisms of $\Zg{G}$ can be explicitly written down from the description of Section \ref{almost-TQFT-strategy}. Consider the set of generators $\Delta = \left\{D, D^\dag, L, P\right\}$ of $\Bordp{2}$ as depicted in Figure \ref{img:gen-tubes-pairs}. Since, by definition, $\Zg{G}(W,A) = \pi_! \circ Z(W, A) \circ \pi^*$ for any bordism $(W, A)$, then we have that
\begin{align*}
	\Zg{G}(D) &= (\pi \circ i)_!: \K{\MHSq} \to \KM{[G/G]}, \hspace{1cm} \Zg{G}(D^\dag) =  (\pi \circ i)^*: \KM{[G/G]} \to \K{\MHSq}, \\
	\Zg{G}(L) &= \hat{q}_! \circ \hat{p}^*: \KM{[G/G]} \to \KM{[G/G]}, \hspace{1.1cm} \Zg{G}(P) = \hat{v}_! \circ \hat{u}^*: \KM{[G/G]} \to \KM{[G/G]},
\end{align*}
where $\hat{q} = \pi \circ q, \hat{p} = \pi \circ p: G^4 \to [G/G]$ and $\hat{v} = \pi \circ l, \hat{u} = \pi \circ t: G^2 \to [G/G]$. Observe that $\pi_{\emptyset}$ is the identity map since $\Rep{G}(\emptyset) = \CharW{G}(\emptyset) = 1$ and that $\Zg{G}(S^1,\star)=\Rep{G}(S^1, \star) = [G/G]$.

If we also consider parabolic structures with parabolic data $\Lambda$, then we have to add to the set of generators the tube $L_{\lambda}: (S^1, \star) \to (S^1, \star)$, for $\lambda \in \Lambda$, as depicted in Figure \ref{img:parabolic-tube}. For that tube, the recipe above says that
$$
	\Zg{G}(L_{\lambda}) = \hat{s}_! \circ \hat{r}^*: \KM{[G/G]} \to \KM{[G/G]},
$$
where $\hat{r} = \pi \circ r, \hat{s} = \pi \circ s: G^2 \times \lambda \to [G/G]$. 

Finally, as mentioned in Example \ref{ex:reduction-TQFT-surfaces}, the only relevant submodule is $\cV = \cV_{(S^1, \star)} \subseteq \KM{\CharW{G}(S^1, \star)} =  \KM{[G/G]}$, so we can restrict the previous maps to $\cV$. Moreover, if the morphism $\eta = \eta_{(S^1, \star)}: \cV_0 \to \cV_0$ is invertible, the reduction $\cZg{G}$ can also be constructed. 
\end{ex}

\section{$\SL{2}(\CC)$-representation varieties}
\label{sec:sl2-repr-var}

In this section, we will focus on the case $G = \SL{2}(\CC)$ and $n = 2$ (surfaces). As an application of the previous TQFTs, we will compute the image in $K$-theory of the Hodge structures of some parabolic representations varieties for arbitrary genus. For that, we will give explicit expressions of the homomorphisms of Section \ref{sec:geometric-reduced-TQFT}.

First of all, let us describe some general properties of the group $\SL{2}(\CC)$ that will be used in the following computations. 
Recall that $\SL{2}(\CC)$ is the group of volume preserving $\CC$-linear automorphisms of $\CC^2$ or, equivalently, the set of $2 \times 2$ matrices with determinant $1$. As a variety, it is an affine variety given as a hypersurface of $\CC^4$. In particular, its (complex) dimension is $3$. Moreover, the action of $\SL{2}(\CC)$ on itself by left multiplication becomes it into a homogeneous variety so it is smooth. For short, from now on, we will denote $\SL{2} = \SL{2}(\CC)$.

The first observation is that, with respect to the action of $\SL{2}$ on itself by conjugation, we have that the GIT quotient map is the usual trace map $\tr: \SL{2} \to \CC = \SL{2} \sslash \SL{2}$. As $\SL{2}$ is affine, all the points are semi-stable for this action. However, there are no stable points for the action since every representation $\rho: \ZZ \to \SL{2}$ is reducible (see Corollary \ref{cor:stability-repr-var} with $\Gamma = \ZZ$).

The unique points which are not poly-stable are the non-diagonalizable matrices of $\SL{2}$. These matrices are conjugated to the Jordan type matrices
$$
	J_+ = \begin{pmatrix}
	1 & 1\\
	0 & 1\\
\end{pmatrix}, \hspace{1cm}
	J_- = \begin{pmatrix}
	-1 & 1\\
	0 & -1\\
\end{pmatrix}.
$$
The stabilizer of $J_\pm$ is the set
$$
	\Stab{J_\pm} = \left\{\left.\begin{pmatrix}
	\pm 1 & \alpha\\
	0 & \pm 1\\
\end{pmatrix}\,\right|\, \alpha \in \CC \right\} \cong \CC.
$$
Hence, their orbits $[J_\pm]$ are isomorphic to
$\SL{2} / \Stab{J_{\pm}} = \SL{2} / \CC
$.

With respect to the poly-stable points, the points $\Id, -\Id \in \SL{2}$ are poly-stable with zero dimensional orbits (indeed, the action is trivial on them). The remaining set is $D = \SL{2} - \left\{\pm \Id\right\} - [J_{\pm}] = \left\{A \in \SL{2}\,|\, \tr(A) \neq \pm 2\right\}$, which are the poly-stable points of maximal dimension. Every element $A \in D$ is conjugated to
$$
	D_\lambda = \begin{pmatrix}
	\lambda & 0\\
	0 & \lambda^{-1}\\
\end{pmatrix},
$$
where $\lambda \in \CC - \left\{0, \pm 1\right\}$ satisfies $\lambda + \lambda^{-1} = \tr{A}$. In particular, $D_\lambda \sim D_{\lambda^{-1}}$, where $\sim$ stands for `conjugate to', and these are the unique diagonal elements mutually conjugated. Actually, $D_{\lambda^{-1}} = P_0 D_\lambda P_0^{-1}$ where
$$
	P_0 = \begin{pmatrix}
	0 & -1\\
	1 & 0\\
\end{pmatrix}.
$$
The stabilizer of $D_\lambda$ is the set
$$
	\Stab{D_\lambda} = \left\{\left.\begin{pmatrix}
	\mu & 0\\
	0 & \mu^{-1}\\
\end{pmatrix}\,\right|\, \mu \in \CC^* \right\} \cong \CC^*,
$$
where $\CC^* = \CC - \left\{0\right\}$. Hence, $[D_\lambda] = \SL{2}/\CC^*$ and $D$ is isomorphic to
$$
	D = \frac{\left(\CC^* - \left\{\pm 1\right\}\right) \times \SL{2}/\CC^*}{\ZZ_2}.
$$
Here, $\ZZ_2$ acts on $\left(\CC^* - \left\{\pm 1\right\}\right) \times \SL{2}/\CC^*$ by $-1 \cdot (\lambda, P) = (\lambda^{-1}, P_0PP_0^{-1})$ for $\lambda \in \CC^* - \left\{\pm 1\right\}$ and $P \in \SL{2}/\CC^*$. The isomorphism sends $(\lambda, P) \mapsto PD_\lambda P^{-1}$. Observe that, in particular, the quotient map given by the trace $\tr: D \to \CC - \left\{\pm 2\right\} = (\CC^* - \left\{\pm 1\right\}) / \ZZ_2$ has non-trivial monodromy of order $2$.

From this description, we obtain a decomposition of $\SL{2}$ by
\begin{equation} \label{eqn:dec-sl2}
	\SL{2} = \left\{\Id\right\} \sqcup \left\{-\Id\right\} \sqcup [J_+] \sqcup [J_-] \sqcup D.
\end{equation}
The stratum $D$ are the poly-stable points of maximal dimension. On $\SL{2}-D$, the poly-stable points of maximal dimension are $\pm \Id$ and, on $[J_\pm] = \SL{2}-D-\left\{\pm \Id\right\}$, all the points are poly-stable (on itself). Hence, this decomposition corresponds, stratum by stratum, with the piecewise quotient
\begin{equation*}
	[\SL{2} / \SL{2}] = \left\{\Id\right\} + \left\{-\Id\right\} + \left\{[J_+]\right\} + \left\{[J_-]\right\} + (\CC-\left\{\pm 2\right\}).
\end{equation*}
In particular, using the description of the geometric TQFT given in Section \ref{sec:geometric-reduced-TQFT}, we obtain that
$$
	\Zg{G}(S^1, \star) = \KM{[\SL{2} / \SL{2}]} = \KM{\Id} \oplus \KM{-\Id} \oplus \KM{[J_+]} \oplus \KM{[J_-]} \oplus \KM{\CC - \left\{\pm 2\right\}},
$$
where the first four summands are mixed Hodge modules over a point and, thus, they are naturally isomorphic to $\K{\MHSq}$.

\subsection{Hodge theory on $\SL{2}(\CC)$}
\label{sec:hodge-theory-sl2}

Thanks to the explicit expression of $\SL{2}$, one can compute the $K$-theory image of the mixed Hodge structure on its cohomology. Notice that we have a surjective regular map $\pi: \SL{2} \to \CC^2-\left\{(0,0)\right\}$ given by $\pi(A) = A(1,0)^t$. For any $(a,b) \in \CC^2-\left\{(0,0)\right\}$, the fiber of this morphism is the affine line $\left\{(x,y)\,|\,ay - bx = 1\right\}$ so $\pi^{-1}(v) \cong \CC$. Hence, we have a fibration
$$
	\CC \longrightarrow \SL{2} \stackrel{\pi}{\longrightarrow} \CC^2 - \left\{(0,0)\right\}.
$$
Moreover, this fibration is locally trivial in the Zariski topology so, in particular, it has trivial monodromy. Hence, by Corollary \ref{cor:mhs-trivial-monodromy} we have
$$
\coh{\SL{2}} = \coh{\CC} \cdot \coh{\CC^2- \left\{(0,0)\right\}} = q(q^2-1) = q^3-q,
$$
where $q = \coh{\CC} \in \K{\MHSq}$, see Remark \ref{sec:mixed-hodge-structures}.

\begin{rmk}
On the finite field of $q$ elements, $\FF_q$, let us consider the group $\SL{2}(\FF_q)$. Such a group has $q^3-q$ elements, precisely because, fixed the image of $(1,0)$ (for which we have $q^2-1$ possibilities) we can assign to $(0,1)$ only $q$ different vectors. Hence, $|\SL{2}(\FF_q)| = q^3-q$ so, by Katz' theorem (see \cite{Hausel-Rodriguez-Villegas:2008}), the Deligne-Hodge polynomial of $\SL{2}(\CC)$ is balanced (see Example \ref{ex:MixedHodgeGradingPure}) and it is $\DelHod{\SL{2}(\CC)} = q^3-q$, in agreement with the previous computation. Actually, the proof is essentially the same (a counting matter) which ilustrates the principle that Katz' theorem and mixed Hodge modules are the two faces of the same coin.
\end{rmk}

From this calculation, the Hodge structures of the strata of $\SL{2}$ of decomposition (\ref{eqn:dec-sl2}) can also be computed. For the orbits of the Jordan type elements, the argument is simple since the fibration
$$
	\CC = \Stab{J_\pm} \longrightarrow \left\{J_\pm \right\} \times \SL{2} {\longrightarrow} [J_\pm]
$$
is locally trivial in the Zariski topology. Hence, we have
$$
	\coh{[J_\pm]} = \frac{\coh{\SL{2}}}{\coh{\CC}} = \frac{q^3-q}{q} = q^2-1.
$$
Analogously, for the orbit of $[D_\lambda]$ we have a locally trivial fibration $\CC^* \to D_\lambda \times \SL{2} \to [D_\lambda]$ so $\coh{[D_\lambda]} = (q^3-q)/(q-1) = q^2+q$.

\begin{ex}\label{ex:Hodge-monodromy-diagonal-class}
Another way of computing such Hodge structure is to notice that, fixed $\lambda \in \CC^*-\left\{\pm 1\right\}$, an element of the conjugacy class $[D_\lambda]$ is determined by the directions of its two different eigenvectors. Such a possible directions are $(\PP^1 \times \PP^1)-\Delta$, where $\Delta \subseteq \PP^1 \times \PP^1$ is the diagonal, which is isomorphic to $\PP^1$. Hence $\coh{[D_\lambda]} = \coh{\PP^1 \times \PP^1} - \coh{\PP^1} = (q+1)^2 - (q+1) = q^2+q$, in agreement with the previous result.

Moreover, with this description, consider the map $\rho: (\PP^1 \times \PP^1) - \Delta \to \PP^2$ given by $\rho([a:b],[c:d]) = [ac: bd : ad+bc]$. This map is $\ZZ_2$-invariant so it descends to an injective morphism $\tilde{\rho}: [(\PP^1 \times \PP^1)-\Delta]/{\ZZ_2} \to \PP^2$ whose image is the whole $\PP^2$ minus the conic $C=\left\{x_2^2 = 4 x_0x_1\right\} \subseteq \PP^2$. Hence, $[(\PP^1 \times \PP^1)-\Delta]/{\ZZ_2} \cong \PP^2 - C$ so, in the notation of Section \ref{sec:equivariant-hodge-mono},
$$
	\coh{\SL{2}/\CC^*}^+ = \coh{(\PP^1 \times \PP^1)-\Delta}^+ = (q^2+q+1)-(q+1)=q^2.
$$
Therefore, $\coh{\SL{2}/\CC^*}^- = \coh{\SL{2}/\CC^*} - \coh{\SL{2}/\CC^*}^+ = q$.
\end{ex}

The Hodge structure of $D$ cannot be obtained with the same kind of argument. However, using the stratification (\ref{eqn:dec-sl2}), we have $D = \SL{2} - \left\{\pm \Id\right\} - [J_\pm]$ so
$$
	\coh{D} = \coh{\SL{2}} - \coh{\pm \Id} - \coh{[J_\pm]} = (q^3-q) - 2 -2(q^2-1) = q^3-2q^2-q.
$$
Observe that $\coh{\CC-\left\{\pm 2\right\}} \coh{[D_\lambda]} = q^3-q^2-2q \neq \coh{D}$, so the quotient map $\tr: D \to \CC-\left\{\pm 2\right\}$ no longer has trivial monodromy. Thus, if we want to compute $\coh{D}$ by means for this map, we need to control the monodromy.

For that purpose, mixed Hodge modules coming from variations of Hodge structures (see Section \ref{sec:monodromy-as-mhm}) are the key ingredients. The point is that we can cover the fibration $\tr: D \to \CC-\left\{\pm 2\right\}$ by another fibration with trivial monodromy. Consider the diagram of fibrations
\[
\begin{displaystyle}
   \xymatrix
   {	
  		\SL{2}/\CC^* \ar[r] & (\CC^*-\left\{\pm 1\right\}) \times \SL{2}/\CC^* \ar[r] \ar[d] & \CC^*-\left\{\pm 1\right\} \ar[d]^{t} \\
   		\SL{2}/\CC^* \ar[r] & D \ar[r]_{\tr} & \CC-\left\{\pm 2\right\} \\
   }
\end{displaystyle}   
\]
Here, $t: \CC^*-\left\{\pm 1\right\} \to \CC-\left\{\pm 2\right\}$ is the morphism $t(\lambda)=\lambda + \lambda^{-1}$ and the middle vertical arrow is the quotient by the action of $\ZZ_2$ on $(\CC^*-\left\{\pm 1\right\}) \times \SL{2}/\CC^*$ given by $-1 \cdot (\lambda, P)= (\lambda^{-1}, P_0PP_0^{-1})$, for $\lambda \in \CC^*-\left\{\pm 1\right\}$ and $P \in \SL{2}/\CC^*$. Since the upper fibration has trivial monodromy, by the results of Section \ref{sec:equivariant-hodge-mono}, we have that
\begin{align*}
	\RM{D}{\CC-\left\{\pm 2\right\}} &= \coh{\SL{2}/\CC^*}^+ \RM{\CC-\left\{\pm 2\right\}}{\CC-\left\{\pm 2\right\}} \\&
	+ \coh{\SL{2}/\CC^*}^- \left(\RM{\CC^*-\left\{\pm 1\right\}}{\CC-\left\{\pm 2\right\}} -\RM{\CC-\left\{\pm 2\right\}}{\CC-\left\{\pm 2\right\}} \right) \\
	&= q^2 T + q\left(T + S_2 \otimes S_{-2} - T\right) = q^2T + qS_2 \otimes S_{-2}.
\end{align*}
Here, we have used that $\RM{\CC^*-\left\{\pm 1\right\}}{\CC-\left\{\pm 2\right\}} = T + S_2 \otimes S_{-2}$ by Example \ref{ex:monodromy-plane-minus-points} and that $\coh{\SL{2}/\CC^*}^+ = q^2$ by Example \ref{ex:Hodge-monodromy-diagonal-class}. Observe that this agrees with the previous computation since
$$
	\coh{D}=\mu(\RM{D}{\CC-\left\{\pm 2\right\}}) = q^2 \mu(T) + q \mu(S_2 \otimes S_{-2}) = q^2(q-2) - q = q^3-2q^2-q,
$$
where $\mu = c_!: \KM{\CC - \left\{\pm 2\right\}} \to \K{\MHSq}$ is the pushout via the projection $c: \CC - \left\{\pm 2\right\} \to \star$.

\subsubsection*{The core submodule}

Recall that the piecewise quotient of $\SL{2}$ under conjugation decomposes as $[\SL{2}/\SL{2}] = \left\{\Id\right\} + \left\{-\Id\right\} + \left\{J_+\right\} + \left\{J_-\right\} + \left(\CC - \left\{\pm 2\right\}\right)$, so $
	\KM{[\SL{2} / \SL{2}]} = \KM{\Id} \oplus \KM{-\Id} \oplus \KM{J_+} \oplus \KM{J_-} \oplus \KM{\CC - \left\{\pm 2\right\}}$. 
For short, in the following we will denote $\Bt = \CC - \left\{\pm 2\right\}$ and we will call it the \emph{space of traces}, as it is the image of the quotient map $\tr: D \to \CC - \left\{\pm 2\right\}$. 

Let $T_1 \in \KM{\Id}$, $T_{-1} \in \KM{-\Id}$, $T_+ \in \KM{J_+}$ and $T_- \in \KM{J_-}$ be the units of these rings (recall that all are isomorphic to $\KM{\star}=\K{\MHSq}$). On $\KM{\Bt}$, we consider the unit $T_{\Bt}$ and $S_2, S_{-2} \in \KM{\Bt}$ the one dimensional representations of $\pi_1(\CC - \left\{\pm 2\right\})$ that are non-trivial on small loops $\gamma_2, \gamma_{-2}$ around $2$ and $-2$, respectively, as in Example \ref{ex:monodromy-plane-minus-points}.

\begin{defn}
The set $\cS = \left\{T_1, T_{-1}, T_+, T_-, T_{\Bt}, S_2, S_{-2}, S_2 \otimes S_{-2}\right\} \subseteq \KM{[\SL{2} / \SL{2}]}$ is called the set of \emph{core elements}. The submodule of $\KM{[\SL{2} / \SL{2}]}$ generated by $\cS$, $\cW$, is called the \emph{core submodule}.
\end{defn}

\begin{rmk}
The submodule $\cW$ will be very important in the incoming computations since we will show that $\cW$ is the submodule generated by the elements $\Zg{\SL{2}}(L)^g \circ \Zg{\SL{2}}(D)(1)$ for $g \geq 0$, that is $\cW = \cV_{(S^1, \star)}$ using the notations of Section \ref{sec:reduction-TQFT}.
\end{rmk}

Now, consider the quotient map $\tr: \SL{2} \to [\SL{2}/\SL{2}]$. From this map, we can build the endomorphism $\eta = \tr_! \circ \tr^*: \KM{[\SL{2}/\SL{2}]} \to \KM{[\SL{2}/\SL{2}]}$. This morphism will be useful in the upcoming sections due to its role in Corollary \ref{cor:existence-reduction-TQFT}.

\begin{prop}\label{prop:eta-morph}
The core submodule $\cW$ is invariant for the morphism $\eta$. Actually, we have
$$
\begin{tabular}{ccc}
	$\eta(T_{\pm 1}) = T_{\pm 1}$ &\hspace{0.5cm}\null& $\eta(T_{\pm}) = (q^2-1)T_{\pm}$\\
	$\eta(T_{\Bt}) = q^2 T_{\Bt} + q S_2 \otimes S_{-2}$ && $\eta(S_2 \otimes S_{-2}) = q T_{\Bt} + q^2 S_2 \otimes S_{-2}$\\
	$\eta(S_2) = q^2 S_2 + q S_{-2}$ && $\eta(S_{-2}) = q S_2 + q^2 S_{-2}$
\end{tabular}
$$
\begin{proof}
First of all, observe that, if $i: X \hookrightarrow [\SL{2} / \SL{2}]$ is the inclusion of a locally closed subset and $\bar{X} = \tr^{-1}(X) \subseteq \SL{2}$, they fit in a pullback diagram
\[
\begin{displaystyle}
   \xymatrix
   {	
	\bar{X} \ar[r]^{\tr|_{\bar{X}}} \ar@{^{(}-{>}}[d]& X \ar@{^{(}-{>}}[d]^i \\
	\SL{2} \ar[r]_{\tr\;\;\;\;}& [\SL{2}/\SL{2}]
   }
\end{displaystyle}   
\]
Thus, $\tr_! \circ \tr^* \circ i_! = (\tr|_{\bar{X}})_! \circ (\tr|_{\bar{X}})^*$ that is, we can compute the image of $\eta$ within $\bar{X}$. For this reason, since $\tr|_{\left\{\pm 1\right\}}: \left\{\pm 1\right\} \to \left\{\pm 1\right\}$ is the identity map, we have $\eta(T_{\pm 1}) = T_{\pm 1}$. Analogously, since $(\tr|_{[J_\pm]})^*$ is a ring homomorphism, it sends units into units so, for $T_{\pm}$, we have that $\eta(T_{\pm}) = (\tr|_{[J_\pm]})_!\underline{\QQ}_{[J_\pm]}T_\pm = \coh{[J_+]} = (q^2-1)T_{\pm}$.

For the first identity of the second row, we have $(\tr|_D)_!(\tr|_D)^*T_{\Bt} = (\tr|_D)_!T_D = \RM{D}{\Bt}$ and the result follows from the computation of Section \ref{sec:hodge-theory-sl2}. For $S_2$, consider the auxiliar variety $Y = \left\{(t,y) \in \Bt \times \CC^*\,|\,y^2 = t-2\right\}$, for which $\RM{Y}{\Bt} = T_{\Bt} + S_2$ (see Section \ref{sec:Hodge-monodromy-covering-spaces}), and the variety $Y' = \left\{(A,y) \in D \times \CC^*\,|\,y^2 = \tr(A)-2\right\}$. They fit in a commutative diagram
	\[
\begin{displaystyle}
   \xymatrix
   {	
   		& Y'\ar[d] \ar[r] \ar[dl]& Y \ar[d] \\
   		\Bt & D \ar[r]_{\tr|_D} \ar[l]^{\;\;\tr|_D} & \Bt
   }
\end{displaystyle}   
\]
where the rightmost square is a pullback. Hence, we have that $\eta(T_{\Bt}) + \eta(S_2) = (\tr|_D)_! (\tr|_D)^*(T_{\Bt} + S_2) = \RM{Y'}{\Bt}$ and, thus, $\eta(S_2) = \RM{Y'}{\Bt} - q^2T_{\Bt}-qS_2 \otimes S_{-2}$ by the computation for $T_{\Bt}$.

In order to compute $\RM{Y'}{\Bt}$ observe that the fiber of the morphism $Y' \to \Bt$ over a point $t$ is $[D_{\lambda}] \times \left\{\pm \sqrt{t-2}\right\}$, which is isomorphic to $\left[(\PP^1 \times \PP^1)-\Delta\right] \times \left\{\pm \sqrt{t-2}\right\}$. If we compactify the space $Y'$ fiberwise, we obtain a variety $\overline{Y}'$ and a morphism $\overline{Y}' \to \Bt$ with fiber, $\overline{F}$, equal to two copies of $\PP^1 \times \PP^1$. Therefore, the Hodge structure on the cohomology of $\overline{F}$ is pure and its unique non-vanishing pieces are $H_c^{0;0,0}(\overline{F})= \QQ \oplus \QQ$, $H_c^{2;1,1}(\overline{F})= \QQ^2 \oplus \QQ^2$ and $H_c^{4;2,2}(\overline{F})= \QQ \oplus \QQ$.  Hence, as mixed Hodge structures, $H_c^0(\overline{F})= 2 \QQ_0$, $H_c^2(\overline{F})= 4q$ and $H_c^4(\overline{F}) = 2q^2$. Working explicitly with the description above, we obtain that the monodromy action of $\pi_1(\Bt)$ on $H_c^k(\overline{F})$ is given as in the following table.

\begin{center}
\begin{tabular}{|c|c|c|c|}
\hline
	 & $k=0$ & $k=2$ & $k=4$\\
\hline
	$\gamma_2$ & $\begin{pmatrix}
	0 & 1\\
	1 & 0\\
\end{pmatrix} \sim \begin{pmatrix}
	1 & 0\\
	0 & -1\\
\end{pmatrix}$
 & $\begin{pmatrix}
	0 & 0 & 0 & 1\\
	0 & 0 & 1 & 0\\
	0 & 1 & 0 & 0\\
	1 & 0 & 0 & 0\\
\end{pmatrix} \sim \begin{pmatrix}
	1 & 0 & 0 & 0\\
	0 & -1 & 0 & 0\\
	0 & 0 & -1 & 0\\
	0 & 0 & 0 & 1\\
\end{pmatrix}$ & $\begin{pmatrix}
	0 & 1\\
	1 & 0\\
\end{pmatrix} \sim \begin{pmatrix}
	1 & 0\\
	0 & -1\\
\end{pmatrix}$\\
\hline
	$\gamma_{-2}$ & $\begin{pmatrix}
	1 & 0\\
	0 & 1\\
\end{pmatrix}$ & $\begin{pmatrix}
	0 & 1 & 0 & 0\\
	1 & 0 & 0 & 0\\
	0 & 0 & 0 & 1\\
	0 & 0 & 1 & 0\\
\end{pmatrix} \sim \begin{pmatrix}
	1 & 0 & 0 & 0\\
	0 & 1 & 0 & 0\\
	0 & 0 & -1 & 0\\
	0 & 0 & 0 & -1\\
\end{pmatrix}$ & $\begin{pmatrix}
	1 & 0\\
	0 & 1\\
\end{pmatrix}$\\
\hline
\end{tabular}
\end{center}
Therefore, we have $\RM{\overline{Y}'}{\Bt} = T_{\Bt} + S_2 + q(T_{\Bt} + S_2 + S_{-2} + S_2 \otimes S_{-2}) + q^2(T_{\Bt} + S_2)$.

The projection of the complement $\overline{Y}' - Y' \to \Bt$ has fiber $\Delta \times \left\{\pm \sqrt{t-2}\right\} \cong \PP^1 \times \left\{\pm \sqrt{t-2}\right\}$. Hence, for the monodromy action we have that $\gamma_{-2}$ acts trivially and $\gamma_2$ interchanges the two copies of $\PP^1$, so $\RM{\overline{Y}' - Y'}{\Bt} = T_{\Bt} + S_2 + q(T_{\Bt} + S_2)$. Thus, we obtain that
$$
	\RM{Y'}{\Bt} = \RM{\overline{Y}'}{\Bt}-\RM{\overline{Y}' - Y'}{\Bt} = q(S_{-2} + S_{2} \otimes S_{-2}) + q^2(T_{\Bt} + S_2).
$$
From this, it follows that $\eta(S_2) = q^2S_2 + qS_{-2}$, as claimed. The calculations for $S_{-2}$ and $S_2 \otimes S_{-2}$ are analogous.
\end{proof}
\end{prop}

The morphism $\eta: \cW \to \cW$ is not invertible since, on $\K{\MHSq}$, the elements $q-1$ and $q+1$ have no inverse. We can solve this problem by considering the localization of $\K{\MHSq}$ over the multiplicative system generated by $q-1$ and $q+1$, denoted $\Ko{\MHSq}$. Extending the localization to $\cW$, we obtain the $\Ko{\MHSq}$-module $\cW_0$.

\begin{cor}\label{cor:eta-morph}
The matrix of the $\K{\MHSq}$-module endomorphism $\eta: \cW \to \cW$ in the set of generators $\cS$ of core elements is
$$
\begin{pmatrix}
	1 & 0 & 0 & 0 & 0 & 0 & 0 & 0\\
	0 & 1 & 0 & 0 & 0 & 0& 0 & 0\\
	0 & 0 & q^2-1 & 0 & 0 & 0& 0 & 0\\
	0 & 0 & 0 & q^2-1 & 0 & 0& 0 & 0\\
	0 & 0 & 0 & 0 & q^2 & 0 & 0 & q\\
	0 & 0 & 0 & 0 & 0 & q^2 & q & 0\\
	0 & 0 & 0 & 0 & 0 & q & q^2 & 0\\
	0 & 0 & 0 & 0 & q & 0 & 0 & q^2\\
\end{pmatrix}
$$
In particular, it is an isomorphism $\eta: \cW_0 \to \cW_0$.
\end{cor}

Let us consider the parabolic data $\Lambda = \left\{\left\{1\right\}, \left\{-1\right\}, [J_+], [J_-]\right\}$, let $\Zs{\SL{2}}: \Bordppar{n}{\Lambda} \to \Modt{\K{\MHSq}}$ be the associated TQFT as done in Section \ref{sec:parabolic-case-TQFT} and let $\Zg{\SL{2}}: \Bordppar{n}{\Lambda} \to \Modt{\K{\MHSq}}$ be the geometric TQFT of Section \ref{sec:geometric-reduced-TQFT}. As a result of the upcoming computations of Sections \ref{sec:tube-J+} and \ref{sec:genus-tube}, we will show that $\cW = \cV_{(S^1, \star)}$. Hence, applying Corollary \ref{cor:existence-reduction-TQFT}, we obtain the following result.

\begin{thm}
Let $\Zs{\SL{2}}: \Bordppar{n}{\Lambda} \to \Modt{\Ko{\MHSq}}$ be the TQFT for $\SL{2}(\CC)$-representation varieties and let $\Zg{\SL{2}}: \Bordppar{n}{\Lambda} \to \Modt{\Ko{\MHSq}}$ be the geometric TQFT. There exists an almost-TQFT, $\cZg{\SL{2}}: \Tub{2}(\Lambda) \to \Mod{\Ko{\MHSq}}$, such that:
\begin{itemize}
	\item The image of the circle is $\cZg{\SL{2}}(S^1, \star) = \cW \subseteq \KMo{[\SL{2}/\SL{2}]}$.
	\item For a strict tube $(W, A, Q): (S^1, \star) \to (S^1, \star)$ it assigns $\cZg{\SL{2}}(W, A, Q) = \Zg{\SL{2}}(W, A, Q) \circ \eta^{-1}: \cW_0 \to \cW_0$.
	\item For a closed surface $W: \emptyset \to \emptyset$, it gives $
		\cZg{\SL{2}}(W, A, Q)(1) = \coh{\Rep{\SL{2}}(W, Q)} \otimes (q^3-q)^{|A|-1}$.
\end{itemize}
\end{thm}

\subsection{Discs and first tubes}
\label{sec:discs-tubes}

As a warm lap, in this section we will compute the morphisms $\Zg{\SL{2}}(L_\lambda)$, for $\lambda = \left\{\Id\right\}$ and $\lambda = \left\{-\Id\right\}$, where $L_\lambda: (S^1, \star) \to (S^1, \star)$ denotes the trivial tube with a marked point with the class $\lambda$. Moreover, we will show that these morphisms preserve the submodule $\cW$.

First, let us focus on the case $\lambda = \left\{\Id\right\}$. The image under the field theory of $\Zg{\SL{2}}$ of the tube $L_{\Id}$ is the span
$$
\begin{matrix}
	[\SL{2}/\SL{2}] & \stackrel{\tr \circ \pi_1}{\longleftarrow} & \SL{2}^2 & \stackrel{\tr \circ \pi_1}{\longrightarrow} & [\SL{2}/\SL{2}]
\end{matrix},
$$
where $\pi_1: \SL{2}^2 \to \SL{2}$ denotes the projection onto the first component. 

Therefore, $\Zg{\SL{2}}(L_{\Id}) = \tr_! \circ (\pi_1)_! \circ (\pi_1)^* \circ \tr^* = (q^3-q)\eta: \KM{[\SL{2}/\SL{2}]} \to \KM{[\SL{2}/\SL{2}]}$, where we have used that $(\pi_1)_! \circ (\pi_1)^*: \KM{\SL{2}} \to \KM{\SL{2}}$ is just multiplication by $\coh{\SL{2}} = q^3-q$. Thus, by Proposition \ref{prop:eta-morph}, we have that $\Zg{\SL{2}}(L_{\Id})(\cW) \subseteq \cW$ and, moreover, its matrix in the generators $\cS$ is the one given in Corollary \ref{cor:eta-morph}. In consequence, the reducted TQFT does $\cZg{\SL{2}}(L_{\Id}) = \Zg{\SL{2}}(L_{\Id}) \circ \eta^{-1} = (q^3-q)1_{\KMo{[\SL{2}/\SL{2}]}}$.

For the case $\lambda = \left\{-\Id\right\}$, the computation is similar. In that case, the image under the field theory is the span
$$
\begin{matrix}
	[\SL{2}/\SL{2}] & \stackrel{\tr \circ \pi_1 }{\longleftarrow} & \SL{2}^2 & \stackrel{-\tr \circ \pi_1}{\longrightarrow} & [\SL{2}/\SL{2}]
\end{matrix}.
$$
Hence, we have $\Zg{\SL{2}}(L_{-\Id}) = (-\tr)_! \circ (\pi_1)_! \circ (\pi_1)^* \circ \tr^* = (q^3-q)\sigma_! \circ \eta$, where $\sigma: [\SL{2}/\SL{2}] \to [\SL{2}/\SL{2}]$ is the reflection that sends $\sigma([A]) = [-A]$ for $A \in \SL{2}$. It is straighforward to check that
$$
	\sigma_! (T_{\pm 1}) = T_{\mp 1} \hspace{0.9cm} \sigma_! (T_{\pm}) = T_{\mp} \hspace{0.9cm} \sigma_! (T_{\Bt}) = T_{\Bt} \hspace{0.9cm} \sigma_! (S_{\pm 2}) = S_{\mp 2} \hspace{0.9cm} \sigma_! (S_2 \otimes S_{-2}) = S_2 \otimes S_{-2}. 
$$
In particular, $\sigma_!(\cW) \subseteq \cW$ and $\cZg{\SL{2}}(L_{-\Id}) = \Zg{\SL{2}}(L_{-\Id}) \circ \eta^{-1} = (q^3-q)\sigma_!$. Thus, the matrix of $\cZg{\SL{2}}(L_{-\Id})$ in the generators $\cS$ of $\cW$ is
$$
	\Zg{\SL{2}}(L_{-\Id})|_V = (q^3-q)\begin{pmatrix}
	0 & 1 & 0 & 0 & 0 & 0 & 0 & 0\\
	1 & 0 & 0 & 0 & 0 & 0 & 0 & 0\\
	0 & 0 & 0 & 1 & 0 & 0 & 0 & 0\\
	0 & 0 & 1 & 0 & 0 & 0 & 0 & 0\\
	0 & 0 & 0 & 0 & 1 & 0 & 0 & 0\\
	0 & 0 & 0 & 0 & 0 & 0 & 1 & 0\\
	0 & 0 & 0 & 0 & 0 & 1 & 0 & 0\\
	0 & 0 & 0 & 0 & 0 & 0 & 0 & 1\\
	\end{pmatrix}.
$$

Finally, with respect to the discs $D: \emptyset \to (S^1,\star)$ and $D^\dag: (S^1,\star) \to \emptyset$ recall that $\Zg{\SL{2}}(D) = i_!$ and $\Zg{\SL{2}}(D^\dag) = i^*$, where $i: \left\{\Id\right\} \hookrightarrow [\SL{2}/\SL{2}]$ is the inclusion. The image of the first one is very easy to identify since, by definition $\Zg{\SL{2}}(D)(1) = i_!(1)=T_1$. For the second one, observe that $i^*T_1=1$ but $i^*T_{-1} = i^*T_\pm = i^*T_{\Bt} = i^*S_{\pm 2} = 0$, so $\Zg{\SL{2}}(D^\dag): \cW \to \K{\MHSq}$ is just the projection onto $T_1$. Since $\eta$ fixes $T_1$, we also have $\Zg{\SL{2}}(D)= \cZg{\SL{2}}(D)$ and $\Zg{\SL{2}}(D^\dag)=\cZg{\SL{2}}(D^\dag)$.

\subsection{The tube with a Jordan type marked point}
\label{sec:tube-J+}

In this section, we will focus on the computation of the image of the tube $L_{[J_+]}: (S^1, \star) \to (S^1, \star)$. Again, the field theory for $\Zg{\SL{2}}$ on $L_{[J_+]}$ is the span
$$
\begin{matrix}
	[\SL{2}/\SL{2}] & \stackrel{\hat{r}}{\longleftarrow} & \SL{2}^2 \times [J_+] & \stackrel{\hat{s}}{\longrightarrow} & [\SL{2}/\SL{2}]\\
	\tr(A) & \mapsfrom & (A, B, C) & \mapsto & \tr(BACB^{-1}) = \tr(AC)
\end{matrix}
$$
so $\Zg{\SL{2}}(L_{[J_+]}) = \hat{s}_! \circ \hat{r}^*: \KM{[\SL{2}/\SL{2}]} \to \KM{[\SL{2}/\SL{2}]}$.

Notice that, if $i: X \hookrightarrow [\SL{2} / \SL{2}]$ is an inclusion of a locally closed subset, then we have a commutative diagram
\[
\begin{displaystyle}
   \xymatrix
   {	
	\hat{r}^{-1}(X) \ar[r]^{\;\;\;\hat{r}|_{\hat{r}^{-1}(X)}} \ar@{^{(}-{>}}[d]& X \ar@{^{(}-{>}}[d]^i \\
	\SL{2}^2 \times [J_+] \ar[r]_{\hat{r}\;}& [\SL{2}/\SL{2}]
   }
\end{displaystyle}   
\]
Hence, for any $M \in \KM{X}$, $\hat{s}_!\circ \hat{r}^* \circ i_!M = \hat{s}_!\circ (\hat{r}|_{\hat{r}^{-1}(X)})^*M$. Moreover, if we decompose $[\SL{2}/\SL{2}] = \sum_k Y_k$ with inclusions $j_k: Y_k \hookrightarrow [\SL{2}/\SL{2}]$, observe that, by the proof of Proposition \ref{prop:strat-mixed-hodge-modules}, we have $\sum_k (j_k)_!(j_k)^* = 1$ so $\hat{s}_! = \sum_k (j_k)_!(j_k)^*\hat{s}_! = \sum_k (j_k)_!(\hat{s}|_{\hat{s}^{-1}(Y_k)})_!$. Therefore,
$$
	\Zg{\SL{2}}(L_{[J_+]}) (i_!M) = \hat{s}_!\circ \hat{r}^* \circ i_!M = \sum_k (j_k)_! \circ (\hat{s}|_{\hat{r}^{-1}(X) \cap \hat{s}^{-1}(Y_k)})_! \circ (\hat{r}|_{\hat{r}^{-1}(X) \cap \hat{s}^{-1}(Y_k)})^* M
$$
Notice that we can decompose
\begin{align*}
	\hat{r}^{-1}(X) \cap \hat{s}^{-1}(Y) &= \left\{(A, B, C) \in \hat{r}^{-1}(X) \times [J_+] \,\,|\,\, \tr(AC) \in Y\right\} \\
	&\cong \left\{A \in \tr^{-1}(X) \,\,|\,\, \tr(AJ_+) \in Y\right\} \times \left(\SL{2}/\CC\right) \times \SL{2},
\end{align*}
where the last isomorphism is given by sending $(A,P,B) \in \tr^{-1}(X) \times \left(\SL{2}/\CC\right) \times \SL{2}$, in the variety downstairs, to $(PAP^{-1}, B, PJ_+P^{-1})$, in the variety upstairs. For short, we will denote $Z_{X,Y} = \hat{r}^{-1}(X) \cap \hat{s}^{-1}(Y)$ and $Z_{X,Y}^0 = \left\{A \in \tr^{-1}(X) \,\,|\,\, \tr(AJ_+) \in Y\right\}$.
Using as stratification the decomposition $[\SL{2}/\SL{2}] = \left\{\pm 1\right\} + \left\{[J_\pm]\right\} + \Bt$, we can compute:
\begin{itemize}
	\item For $T_{1}$, we have $\hat{r}^{-1}(\Id) = Z_{\Id, [J_+]} = [J_+] \times \SL{2}$ and $\hat{s}$ is just the projection onto a point $\hat{s}: [J_+] \times \SL{2} \to \left\{[J_+]\right\}$. Hence, $\Zg{\SL{2}}(L_{[J_+]})(T_1) = (q^2-1)(q^3-q) T_+$.
	\item For $T_{-1}$, since $[-J_+] = [J_-]$, we have $\hat{r}^{-1}(-\Id) = Z_{-\Id, [J_-]} = [J_-] \times \SL{2}$. Again, $\hat{s}$ is just the projection $[J_-] \times \SL{2} \to \left\{[J_-]\right\}$ so $\Zg{\SL{2}}(L_{[J_+]})(T_{-1}) = (q^2-1)(q^3-q) T_-$.
	\item For $T_+$, the situation becomes more involved since we have a non-trivial decomposition of $\hat{r}^{-1}([J_+])$. We analize each stratum separately:
	\begin{itemize}
		\item For $\Id$, we have $Z_{[J_+], \Id} = [J_+] \times \SL{2}$ since all the elements of $Z_{[J_+], \Id}$ have the form $(A^{-1}, B, A)$. Thus, $\hat{s}: [J_+]  \times \SL{2} \to \left\{\Id\right\}$ is the projection onto a point so $\Zg{\SL{2}}(L_{[J_+]})(T_{+})|_{\Id} = (q^2-1)(q^3-q)T_1$.
		\item For $-\Id$, observe that $Z_{[J_+], -\Id} = \emptyset$ so $\Zg{\SL{2}}(L_{[J_+]})(T_{+})|_{-\Id} = 0$.
		\item For $[J_+]$, a straightforward check shows that
		$$
			Z_{[J_+], [J_+]}^0 = \left\{\left.\begin{pmatrix}
	1 & b\\
	0 & 1\\
\end{pmatrix}\,\right|\,\, b \in \CC-\left\{0, -1\right\}
\right\}.
		$$
	Hence, $Z_{[J_+], [J_+]} = \left(\CC-\left\{0, -1\right\}\right) \times \left(\SL{2}/\CC\right)  \times \SL{2}$ and $\hat{s}$ is a projection onto the singleton $\left\{[J_+]\right\}$. Thus, $\Zg{\SL{2}}(L_{[J_+]})(T_{+})|_{[J_+]} = (q-2)(q^2-1)(q^3-q)T_+$.
		\item For $[J_-]$, we have that 
		$$
			Z_{[J_+], [J_-]}^0 = \left\{\left.\begin{pmatrix}
	1+a & \frac{a^2}{2}\\
	-4 & 1-a\\
\end{pmatrix}\,\right|\,\, a \in \CC
\right\}.
		$$
		Therefore, $Z_{[J_+], [J_-]} = \CC \times \left(\SL{2}/\CC\right)  \times \SL{2}$ so $\Zg{\SL{2}}(L_{[J_+]})(T_{+})|_{[J_-]} = q(q^2-1)(q^3-q)T_-$.
		\item For $\Bt$, we have that
$$
		Z_{[J_+], \Bt}^0 = \left\{\left.\begin{pmatrix}
	1+a & -\frac{a^2}{c}\\
	c & 1-a\\
\end{pmatrix}\,\right|\,\, a \in \CC, c \in \CC^* - \left\{-4\right\}
\right\},
		$$
with the projection $Z_{[J_+], \Bt}^0 \to \Bt$ given by $A \mapsto c + 2$. Hence $\hat{s}: Z_{[J_+], \Bt} = Z_{[J_+], \Bt}^0 \times \left(\SL{2}/\CC\right)  \times \SL{2} \to \Bt$ applies $\hat{s}(A, B, P) = c + 2$. In particular, $\hat{s}$ has trivial monodromy there, so $\Zg{\SL{2}}(L_{[J_+]})(T_+)|_{\Bt} = \RM{Z_{[J_+], \Bt}}{\Bt} = q(q^2-1)(q^3-q)T_{\Bt}$.
	\end{itemize}
Summarizing, we have
$$
	\Zg{\SL{2}}(L_{[J_+]})(T_{+}) = (q^2-1)(q^3-q)\left[T_1 + (q-2)T_+ + qT_- + qT_{\Bt}\right].
$$
	\item For $T_-$, the calculation is very similar to the one of $T_+$, since $Z_{[J_-], Y}^0 = Z_{[J_+], -Y}^0$ for every stratum $Y$. Hence, we have:
		\begin{itemize}
		\item For $\Id$, $Z_{[J_-], \Id} = \emptyset$, so $\Zg{\SL{2}}(L_{[J_+]})(T_{-})|_{\Id} = 0$. Analogously, for $-\Id$ we have $Z_{[J_-], -\Id} = -[J_+] \times \SL{2}$, so $\Zg{\SL{2}}(L_{[J_+]})(T_{-})|_{-\Id} = (q^2-1)(q^3-q)T_{-1}$.
		\item For $[J_+]$, $Z_{[J_-], [J_+]}^0 = Z_{[J_+], [J_-]}^0 \cong \CC$, so $\Zg{\SL{2}}(L_{[J_+]})(T_{-})|_{[J_+]} = q(q^2-1)(q^3-q)T_+$.
		\item For $[J_-]$, we have $Z_{[J_-], [J_-]}^0 = Z_{[J_+], [J_+]}^0 \cong \CC - \left\{0, -1\right\}$. Hence, $\Zg{\SL{2}}(L_{[J_+]})(T_{-})|_{[J_-]} = (q-2)(q^2-1)(q^3-q)T_-$.
		\item For $\Bt$, since $-\Bt=\Bt$, we have $Z_{[J_-], \Bt}^0 = Z_{[J_+], -\Bt}^0 \cong \CC \times \left(\CC^*-\left\{-4\right\}\right)$. Hence, $\hat{s}$ also has trivial monodromy on $Z_{[J_-], \Bt}$ so $\Zg{\SL{2}}(L_{[J_+]})(T_-)|_{\Bt} = q(q^2-1)(q^3-q)T_{\Bt}$.
	\end{itemize}
Summarizing, this calculation says that
$$
	\Zg{\SL{2}}(L_{[J_+]})(T_{-}) = (q^2-1)(q^3-q)\left[T_{-1} + qT_+ + (q-2)T_- + qT_{\Bt}\right].
$$
	\item For $T_{\Bt}$, observe that we have $Z_{\Bt, \pm \Id} = \emptyset$ so $\Zg{\SL{2}}(L_{[J_+]})(T_{\Bt})$ has no components on $\pm \Id$. In this section, we will compute the image of $\Zg{\SL{2}}(L_{[J_+]})(T_{\Bt})$ on $[J_\pm]$. Its image on $\Bt$ is much harder and it will be described later. For $[J_\pm]$ observe that $Z_{\Bt, [J_\pm]}^0 \cong Z_{[J_\pm], \Bt}^0 \cong \CC \times (\CC^* - \left\{-4\right\})$. Hence, $\Zg{\SL{2}}(L_{[J_+]})(T_{\Bt})|_{[J_\pm]} = q(q-2)(q^2-1)(q^3-q)T_\pm$.
	\item For $S_2$, $\Zg{\SL{2}}(L_{[J_+]})(S_2)$ has no components on $\pm \Id$. Again, we will only focus on its image on $[J_\pm]$. We consider the auxiliar variety
$$
	Y_2 = \left\{(t,y) \in \Bt \times \CC^*\,\,|\,\, y^2 = t-2\right\}
$$
with nice fibration $Y_2 \to \Bt$ given by $(t,y) \mapsto t$ whose Hodge monodromy representation is $\RM{Y_2}{\Bt} = T_{\Bt} + S_2$. In that case, we have a commutative diagram
\[
\begin{displaystyle}
   \xymatrix
   {	
	& \hat{Y}_2\ar[d]\ar[ld] \ar[r]& Y_2 \ar[d]\\
	[J_+] & Z_{\Bt, [J_+]} \ar[l]^{\hat{s}} \ar[r]_{\;\;\;\;\hat{r}} & \Bt 
   }
\end{displaystyle}   
\]
where $\hat{Y}_2$ is the pullback of $Z_{\Bt, [J_+]}$ and $Y_2$ over $\Bt$. Computing directly, we obtain that
$$
	\hat{Y}_2 = \left\{(t, x, y, P, B) \in \Bt \times \CC \times \CC^* \times \left(\SL{2}/\CC\right) \times \SL{2}\,|\, y^2=t-2\right\} \cong Y_2 \times \CC \times \left(\SL{2}/\CC\right) \times \SL{2}.
$$
Hence, $\RM{\hat{Y}_2}{[J_+]} = q(q^2-1)(q^3-q)\coh{Y_2}T_+ = q(q^2-1)(q^3-q)(q-3)T_+$. Since
$$
	\RM{\hat{Y}_2}{[J_+]} = \hat{s}_!\hat{r}^* \RM{Y_2}{\Bt} = \Zg{\SL{2}}(L_{[J_+]})(T_{\Bt})|_{[J_+]} + \hat{s}_!\hat{r}^* S_2,
$$
	we obtain that $\Zg{\SL{2}}(L_{[J_+]})(S_2)|_{[J_+]} =\hat{s}_!\hat{r}^* S_2 = -q(q^2-1)(q^3-q)T_+$. Furthermore, as $Z_{\Bt,[J_+]} \cong Z_{\Bt,[J_-]}$, we also have $\Zg{\SL{2}}(L_{[J_+]})(S_2)|_{[J_-]} = -q(q^2-1)(q^3-q)T_-$.
\item For $S_{-2}$ the calculation is analogous to the one of $S_2$ so $\Zg{\SL{2}}(L_{[J_+]})(S_{-2})|_{[J_+]} = -q(q^2-1)(q^3-q)T_+$ and $\Zg{\SL{2}}(L_{[J_+]})(S_{-2})|_{[J_-]} = -q(q^2-1)(q^3-q)T_-$. The same holds for $S_2 \otimes S_{-2}$ so $\Zg{\SL{2}}(L_{[J_+]})(S_2 \otimes S_{-2})|_{[J_+]} = -q(q^2-1)(q^3-q)T_+$ and $\Zg{\SL{2}}(L_{[J_+]})(S_2 \otimes S_{-2})|_{[J_-]} = -q(q^2-1)(q^3-q)T_-$.
\end{itemize}

\subsubsection*{The components in $\Bt$}

The calculation of the components of the mixed Hodge modules $T_{\Bt}, S_2, S_{-2}$ and $S_2 \otimes S_{-2}$ on $\Bt$ is harder than the previous ones and it requires a more subtle analysis. For this reason, we compute it separately. First of all, observe that
$$
	Z_{\Bt, \Bt}^0 = \left\{\left.\begin{pmatrix}
	\frac{t}{2} + a & b\\
	t'-t & \frac{t}{2} - a\\
\end{pmatrix}\,\right|\,\, \begin{matrix}
	t, t' \in \Bt, a,b \in \CC\\
	\frac{t^2}{4}-a^2-(t'-t)b = 1
\end{matrix}
\right\}.
$$
The projections onto each copy of $\Bt$ are $\hat{r}(a, b, t, t') = t$ and $\hat{s}(a, b, t, t') = t'$. We decompose $Z_{\Bt, \Bt}^0 = Z_1 \sqcup Z_2$ where $Z_1 = Z_{\Bt, \Bt}^0 \cap \left\{t=t'\right\}$ and $Z_2 = Z_{\Bt, \Bt}^0 \cap \left\{t \neq t'\right\}$.

For $Z_1$, we have the explicit expression $Z_1 = \left\{4a^2 = (t+2)(t-2)\right\} \times \CC$. Hence, by the same argument that in Section \ref{sec:Hodge-monodromy-covering-spaces}, with respect to the projection $(t,a,b) \mapsto t$ we have $\RM{Z_1}{\Bt} = q(T_{\Bt} + S_2 \otimes S_{-2})$. Thus, the contribution of $Z_1$ to $\Zg{\SL{2}}(L_{[J_+]})(T_{\Bt})|_{\Bt}$ is $q(q^2-1)(q^3-q)(T_{\Bt} + S_2 \otimes S_{-2})$.
	
	With respect to $S_2$, we again consider $Y_2 = \left\{y^2 = t-2\right\}$, for which we have a commutative diagram
\[
\begin{displaystyle}
   \xymatrix
   {	
	& \hat{Y}_2\ar[d]\ar[ld] \ar[r]& Y_2 \ar[d]\\
	\Bt & Z_1 \ar[l]^{\hat{s}} \ar[r]_{\hat{r}} & \Bt 
   }
\end{displaystyle}   
\]
where $\hat{Y}_2$ is the pullback of $Z_1$ and $Y_2$ over $\Bt$. Computing directly
$$
	\hat{Y}_2 = \left\{(t, a, y) \in \Bt \times \CC^* \times \CC \times \CC^* \,\left|\, \begin{matrix}4a^2=(t+2)(t-2)\\y^2=t-2\end{matrix}\right.\right\} \times \CC,
$$
with projection $\hat{Y}_2 \to \Bt$ given by $(t, a, y, b) \mapsto t$. The fiber over $t \in \Bt$ is the set of four lines $F=\left\{\left(t, \pm\frac{1}{2} \sqrt{(t+2)(t-2)}, \pm \sqrt{t-2}\right)\right\} \times \CC$ so the monodromy action of $\pi_1(\Bt)$ on $H_c^1(F)=\QQ(-1)^4 = 4q$ is given by
$$
\gamma_2 \mapsto \begin{pmatrix}
	0 & 0 & 0 & 1\\
	0 & 0 & 1 & 0\\
	0 & 1 & 0 & 0\\
	1 & 0 & 0 & 0\\
\end{pmatrix} \sim \begin{pmatrix}
	1 & 0 & 0 & 0\\
	0 & -1 & 0 & 0\\
	0 & 0 & -1 & 0\\
	0 & 0 & 0 & 1\\
\end{pmatrix},
\hspace{0.5cm}
\gamma_{-2} \mapsto \begin{pmatrix}
	0 & 1 & 0 & 0\\
	1 & 0 & 0 & 0\\
	0 & 0 & 0 & 1\\
	0 & 0 & 1 & 0\\
\end{pmatrix} \sim \begin{pmatrix}
	1 & 0 & 0 & 0\\
	0 & 1 & 0 & 0\\
	0 & 0 & -1 & 0\\
	0 & 0 & 0 & -1\\
\end{pmatrix}.
$$
Therefore, $\RM{\hat{Y}_2}{\Bt} = q(T_{\Bt} + S_2 + S_{-2} + S_2 \otimes S_{-2})$. Since $\RM{\hat{Y}_2}{\Bt} = \hat{s}_!\hat{r}^*\RM{Y_2}{\Bt} = \hat{s}_!\hat{r}^*(T_{\Bt} + S_2)= \RM{Z_1}{\Bt} + \hat{s}_!\hat{r}^* S_2$, we finally obtain that $\hat{s}_!\hat{r}^*S_2 = q(S_2 + S_{-2})$. Therefore, the contribution of $Z_1$ to $\Zg{\SL{2}}(L_{[J_+]})(S_2)|_{\Bt}$ is $q(q^2-1)(q^3-q)(S_2 + S_{-2})$. By symmetry, the contribution to $\Zg{\SL{2}}(L_{[J_+]})(S_{-2})|_{\Bt}$ is $q(q^2-1)(q^3-q)(S_2 + S_{-2})$. An analogous computation for $S_2 \otimes S_{-2}$ shows that the contribution of $Z_1$ to $\Zg{\SL{2}}(L_{[J_+]})(S_2 \otimes S_{-2})|_{\Bt}$ is $q(q^2-1)(q^3-q)(T_{\Bt} + S_2 \otimes S_{-2})$.

\begin{figure}[h]
	\begin{center}
	\includegraphics[scale=0.35]{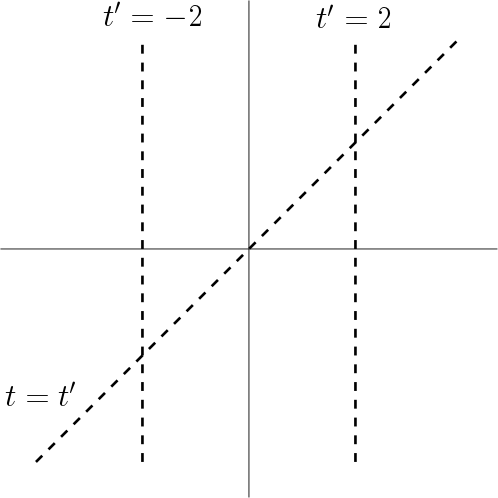}
	\caption{}
	\label{img:Z2}
	\end{center}
\vspace{-0.6cm}
\end{figure}

The most involved stratum is $Z_2$. For it, we have the explicit expression
$$
	Z_2 = \left\{\left.\begin{pmatrix}
	\frac{t}{2} + a & \frac{t^2 - 4a^2 -4}{4(t'-t)}\\
	t'-t & \frac{t}{2} - a\\
\end{pmatrix}\,\right|\,\, \begin{matrix}
	t, t' \in \Bt, t' \neq t\\
	a \in \CC
\end{matrix}
\right\} \cong \CC \times \left[(\CC-\left\{\pm 2\right\})^2- \Delta\right],
$$
where $\Delta \subseteq (\CC-\left\{\pm 2\right\})^2$ is the diagonal. Here, the projections are $\hat{r}(a,t,t')=t$ and $\hat{s}(a,t,t')=t'$. In this way, the second factor of this decomposition of $Z_2$ is as depicted in Figure \ref{img:Z2}. In particular, observe that, under the projection $\hat{s}$, the monodromy is trivial so $\RM{Z_2}{\Bt} = q(q-3)T_{\Bt}$. Hence, 
as $Z_{\Bt, \Bt}= Z_{\Bt, \Bt}^0 \times \SL{2}/\CC \times \SL{2}$, the contribution of $Z_2$ to $\Zg{\SL{2}}(L_{[J_+]})(T_{\Bt})|_{\Bt}$ is $q(q-3)(q^2-1)(q^3-q)T_{\Bt}$.

For $S_2$, again consider the auxiliar variety $Y_2 = \left\{y^2 = t-2\right\}$ for $t \neq \pm 2$. In that case, the pullback of $Y_2$ and $Z_2$ over $\Bt$ is
\begin{align*}
	\hat{Y}_2 &= \left\{(a, t, t', y) \in \CC \times \left[(\CC-\left\{\pm 2\right\})^2- \Delta\right] \times \CC^* \,\left|\, y^2=t-2\right.\right\} \\
	&= \CC \times \left\{(t, t', y) \in \left[(\CC-\left\{\pm 2\right\})^2- \Delta\right] \times \CC^* \,\left|\, y^2=t-2\right.\right\}.
\end{align*}
Denote $X_2 = \left\{(t, t', y) \in \left[(\CC-\left\{\pm 2\right\})^2- \Delta\right] \times \CC^* \,\left|\, y^2=t-2\right.\right\}$. Under the projection $\hat{s}$, the fiber of $X_2$, $F$, is a parabola with five points removed. Compactifying the fibers of $X_2$ we obtain a variety $\overline{X}_2$ whose fiber $\overline{F}$ is $\PP^1$. By Remark \ref{rmk:trivial-monodromy} (see also Remark 2.5 of \cite{LMN}), a nice fibration whose fiber is $\PP^1$ has trivial monodromy. Therefore, $\RM{\overline{X}_2}{\Bt} = (q+1)T_{\Bt}$.

Now, consider the difference $\overline{X}_2 - X_{2}$, whose fiber over $t' \in \Bt$ is the set of six points
$$
	\overline{F} - F = \left\{(2, t', 0), (-2, t', \pm 2i), (t',t', \pm \sqrt{t'-2}), \infty\right\},
$$
where $\infty$ denotes the (unique) point at infinity of the parabola over $t'$. From this expression, we obtain that the monodromy action of $\pi_1(\Bt)$ on the fiber of the covering $\overline{X}_2 - X_{2}$ is given by
$$
\gamma_2 \mapsto \begin{pmatrix}
	1 & 0 & 0 & 0 & 0 & 0\\
	0 & 1 & 0 & 0 & 0 & 0\\
	0 & 0 & 1 & 0 & 0 & 0\\
	0 & 0 & 0 & 0 & 1 & 0\\
	0 & 0 & 0 & 1 & 0 & 0\\
	0 & 0 & 0 & 0 & 0 & 1\\
\end{pmatrix} \sim \begin{pmatrix}
	1 & 0 & 0 & 0 & 0 & 0\\
	0 & 1 & 0 & 0 & 0 & 0\\
	0 & 0 & 1 & 0 & 0 & 0\\
	0 & 0 & 0 & 1 & 0 & 0\\
	0 & 0 & 0 & 0 & -1 & 0\\
	0 & 0 & 0 & 0 & 0 & 1\\
\end{pmatrix} \hspace{1cm}
\gamma_{-2} \mapsto \Id
$$
Therefore, $\RM{\overline{X}_2-X_2}{\Bt} = 5T_{\Bt} + S_2$ so $\RM{X_2}{\Bt} = \RM{\overline{X}_2}{\Bt} - \RM{\overline{X}_2-X_2}{\Bt} = (q+1)T_{\Bt} - 5T_{\Bt} - S_2 = (q-4) T_{\Bt} - S_2$. Hence, since $\hat{Y}_2 = \CC \times X_2$, we finally obtain $\RM{\hat{Y}_2}{\Bt} = q\RM{X_2}{\Bt} = q(q-4)T_{\Bt} - S_2$. With this information at hand, for the restriction to $Z_2$ we have
$$
	\RM{\hat{Y}_2}{\Bt} = (\hat{s}|_{Z_2})_!(\hat{r}|_{Z_2})^*\RM{Y_2}{\Bt} = (\hat{s}|_{Z_2})_!(\hat{r}|_{Z_2})^*(T_{\Bt} + S_2) = \RM{Z_2}{\Bt} + (\hat{s}|_{Z_2})_!(\hat{r}|_{Z_2})^*S_2.
$$
Therefore, the contribution of $Z_2$ to $\Zg{\SL{2}}(L_{[J_+]})(S_2)|_{\Bt}$ is $(q^2-1)(q^3-q)(\hat{s}|_{Z_2})_!(\hat{r}|_{Z_2})^*S_2 = (q^2-1)(q^3-q)\left[-qT_{\Bt}-qS_2\right]$. For $S_{-2}$, the calculation is absolutely analogous by considering $Y_{-2} = \left\{y^2 = t+2\right\}$, so we obtain that the contribution of $Z_2$ to $\Zg{\SL{2}}(L_{[J_+]})(S_{-2})|_{\Bt}$ is $(q^2-1)(q^3-q)\left[-qT_{\Bt}-qS_{-2}\right]$.

For $S_2 \otimes S_{-2}$, the computation has some peculiarities that deserve a brief explanation. Consider the auxiliar variety $Y_{2,-2} = \left\{y^2 = (t-2)(t+2)\right\}$ so, under the projection $(t,y) \mapsto t$, we have $\RM{Y_{2,-2}}{\Bt} = T_{\Bt} + S_2 \otimes S_{-2}$. In this case, the pullback of $Y_{2,-2}$ and $Z_2$ under $\Bt$ is $\hat{Y}_{2,-2} = \CC \times X_{2,-2}$, with
$$
	X_{2,-2} = \left\{(t, t', y) \in \left[(\CC-\left\{\pm 2\right\})^2- \Delta\right] \times \CC^* \,\left|\, y^2=(t-2)(t+2)\right.\right\}.
$$
Hence, this time the fiber $F$ is an hyperbola with four points removed. Again, for the fiberwise compactification of $X_{2,-2}$, $\overline{X}_{2,-2}$, with fiber $\overline{F}$, we have $\RM{\overline{X}_{2,-2}}{\Bt} = (q+1)T_{\Bt}$. The fiber over $t' \in \Bt$ is the set of six points
$$
	\overline{F} - F = \left\{(\pm 2, t', 0), (t',t', \pm \sqrt{(t'-2)(t'+2)}), \pm \infty\right\},
$$
where $\pm \infty$ denotes the two points at infinity of the hyperbola over $t'$. Hence, the monodromy action on $\overline{X}_{2,-2} - X_{2,-2}$ is given by
$$
\gamma_2, \gamma_{-2} \mapsto \begin{pmatrix}
	1 & 0 & 0 & 0 & 0 & 0\\
	0 & 1 & 0 & 0 & 0 & 0\\
	0 & 0 & 1 & 0 & 0 & 0\\
	0 & 0 & 0 & 0 & 1 & 0\\
	0 & 0 & 0 & 1 & 0 & 0\\
	0 & 0 & 0 & 0 & 0 & 1\\
\end{pmatrix} \sim \begin{pmatrix}
	1 & 0 & 0 & 0 & 0 & 0\\
	0 & 1 & 0 & 0 & 0 & 0\\
	0 & 0 & 1 & 0 & 0 & 0\\
	0 & 0 & 0 & 1 & 0 & 0\\
	0 & 0 & 0 & 0 & -1 & 0\\
	0 & 0 & 0 & 0 & 0 & 1\\
\end{pmatrix}
$$
Therefore, $\RM{\overline{X}_{2,-2}}{\Bt} = 5T_{\Bt} + S_2 \otimes S_{-2}$. Repeating the computations above as for $S_2$, this implies that the contribution of $Z_2$ to $\Zg{\SL{2}}(L_{[J_+]})(S_2 \otimes S_{-2})|_{\Bt}$ is $(q^2-1)(q^3-q)\left[-qT_{\Bt}-qS_2 \otimes S_{-2}\right]$.

Summarizing, if we put together the calculations for $Z_1$ and $Z_2$, we finally obtain that
\begin{align*}
 \Zg{\SL{2}}(L_{[J_+]})(T_{\Bt})|_{\Bt} & = (q^2-1)(q^3-q)\left[q (T_{\Bt} + S_2 \otimes S_{-2}) + q(q-3) T_{\Bt}\right] \\
& = (q^2-1)(q^3-q)\left[(q^2-2q)T_{\Bt} + qS_2 \otimes S_{-2}\right],\\
\Zg{\SL{2}}(L_{[J_+]})(S_2)|_{\Bt} &= (q^2-1)(q^3-q)\left[q(S_2 + S_{-2}) - q(T_{\Bt} + S_2)\right] \\
&= (q^2-1)(q^3-q)\left[-qT_{\Bt} + qS_2\right],\\
\Zg{\SL{2}}(L_{[J_+]})(S_{-2})|_{\Bt} &= (q^2-1)(q^3-q)\left[q (S_2 + S_{-2}) - q(T_{\Bt} + S_{-2})\right] \\
&= (q^2-1)(q^3-q)\left[-qT_{\Bt} + qS_{-2}\right],
\\
\Zg{\SL{2}}(L_{[J_+]})(S_2 \otimes S_{-2})|_{\Bt} &= (q^2-1)(q^3-q)\left[q (T_{\Bt} + S_2 \otimes S_{-2}) - q(T_{\Bt} + S_2 \otimes S_{-2})\right] = 0.
\end{align*}
Therefore, we have completed the proof of the following result.
\begin{thm}
The core submodule, $\cW$, is invariant under the morphism $\Zg{\SL{2}}(L_{[J_+]}): \KM{[\SL{2}/\SL{2}]} \to \KM{[\SL{2}/\SL{2}]}$. Moreover, in the set of generators $\cS$ of $\cW$, the matrix of $\Zg{\SL{2}}(L_{[J_+]})$ is

$$
\Zg{\SL{2}}(L_{[J_+]}) = (q^2-1)(q^3-q)\begin{pmatrix}
0 & 0 & 1 & 0 & 0 & 0 & 0 & 0 \\
0 & 0 & 0 & 1 & 0 & 0 & 0 & 0 \\
1 & 0 & q - 2 & q & q^{2} - 2 \, q & -q & -q & -q \\
0 & 1 & q & q - 2 & q^{2} - 2 \, q & -q & -q & -q \\
0 & 0 & q & q & q^{2} - 2 \, q & -q & -q & 0 \\
0 & 0 & 0 & 0 & 0 & 0 & q & 0 \\
0 & 0 & 0 & 0 & 0 & q & 0 & 0 \\
0 & 0 & 0 & 0 & q & 0 & 0 & 0
\end{pmatrix}.
$$
\end{thm}

In particular, multiplying by $\eta^{-1}$ on the right, we obtain the following corollary.

\begin{cor}
In the set of generators $\cS$, the matrix of $\cZg{\SL{2}}(L_{[J_+]}): \cW_0 \to \cW_0$ is
$$
\cZg{\SL{2}}(L_{[J_+]}) = (q^3-q)\begin{pmatrix}
0 & 0 & 1 & 0 & 0 & 0 & 0 & 0 \\
0 & 0 & 0 & 1 & 0 & 0 & 0 & 0 \\
q^{2} - 1 & 0 & q - 2 & q & (q-1)^2 & -q
+ 1 & -q + 1 & -2 \, q + 2 \\
0 & q^{2} - 1 & q & q - 2 & (q-1)^2& -q
+ 1 & -q + 1 & -2 \, q + 2 \\
0 & 0 & q & q & q^{2} - 2 \, q & -q + 1 & -q + 1
& -q + 2 \\
0 & 0 & 0 & 0 & 0 & -1 & q & 0 \\
0 & 0 & 0 & 0 & 0 & q & -1 & 0 \\
0 & 0 & 0 & 0 & q & 0 & 0 & -1
\end{pmatrix}.
$$
\end{cor}

\begin{rmk}
From this computation, it is relatively easy to calculate the image of the tube $L_{[J_-]}: (S^1, \star) \to (S^1, \star)$. The point is that, since $[J_-] = -[J_+]$, the span associated to $L_{[J_-]}$ fits in a commutative diagram
\[
\begin{displaystyle}
   \xymatrix
   {
   	& \SL{2}^2 \times [J_-] \ar[dd]^{\varphi} \ar[rd]\ar[ld] & \\
	\SL{2} & & \SL{2} \\
	& \SL{2}^2 \times [J_+] \ar[ru]_{r} \ar[lu]^{\varsigma \circ s} & \\
   }
\end{displaystyle}   
\]
where $\varphi(A,B,C) = (A,B,-C)$ and $\varsigma: \SL{2} \to \SL{2}$ is given by $\varsigma(A)=-A$. Hence, as $\varphi$ is an isomorphism, by Lemma \ref{lem:iso-C-algebra} we have $\varphi_! \circ \varphi^* = 1$, so $\Zs{\SL{2}}(L_{[J_-]}) = \varsigma_! \circ s_! \circ r^*$. This implies that $\cZg{\SL{2}}(L_{[J_-]}) = \sigma_! \circ \cZg{\SL{2}}(L_{[J_+]})$, where $\sigma: [\SL{2}/\SL{2}] \to [\SL{2}/\SL{2}]$ is the factorization through the quotient of $\varsigma$. Hence, using the description of $\sigma$ of Section \ref{sec:discs-tubes}, we obtain that the matrix of $\cZg{\SL{2}}(L_{[J_-]})$ in $\cW$ is
$$
\cZg{\SL{2}}(L_{[J_-]}) = (q^3-q)\begin{pmatrix}
0 & 0 & 0 & 1 & 0 & 0 & 0 & 0 \\
0 & 0 & 1 & 0 & 0 & 0 & 0 & 0 \\
0 & q^{2} - 1 & q & q - 2 & (q-1)^2 & -q
+ 1 & -q + 1 & -2 \, q + 2 \\
q^{2} - 1 & 0 & q - 2 & q & (q-1)^2 & -q
+ 1 & -q + 1 & -2 \, q + 2 \\
0 & 0 & q & q & q^{2} - 2 \, q & -q + 1 & -q + 1
& -q + 2 \\
0 & 0 & 0 & 0 & 0 & q & -1 & 0 \\
0 & 0 & 0 & 0 & 0 & -1 & q & 0 \\
0 & 0 & 0 & 0 & q & 0 & 0 & -1
\end{pmatrix}
$$ 
\end{rmk}

\begin{rmk}
Observe that, despite that we need to localize in order to invert $\eta$, $\cZg{\SL{2}}(L_{[J_+]})$ is a well-defined map on $\cW$ itself. Moreover, the matrix of $\cZg{\SL{2}}(L_{[J_+]})$ for the generators $\cS$ is diagonalizable. Actually, using any kind of symbolic computation software, like SageMath \cite{sagemath}, we find that $\cZg{\SL{2}}(L_{[J_+]}) = PDP^{-1}$, where
$$ D= (q^3-q)\left(\begin{array}{c|ccc|ccc|c}
q^{2} - 1 & 0 & 0 & 0 & 0 & 0 & 0 & 0 \\
\hline
 0 & -q - 1 & 0 & 0 & 0 & 0 & 0 & 0 \\
 0 & 0 & -q - 1 & 0 & 0 & 0 & 0 & 0 \\
 0 & 0 & 0 & -q - 1 & 0 & 0 & 0 & 0 \\
 \hline
0 & 0 & 0 & 0 & q - 1 & 0 & 0 & 0 \\
0 & 0 & 0 & 0 & 0 & q - 1 & 0 & 0 \\
0 & 0 & 0 & 0 & 0 & 0 & q - 1 & 0 \\
\hline
0 & 0 & 0 & 0 & 0 & 0 & 0 & 0
\end{array}\right)
$$
$$
P = \begin{pmatrix}
1 & 1 & 0 & 0 & 1 & 0 & 0 & 1 \\
1 & 0 & 1 & 0 & 0 & 1 & 0 & 1 \\
q^{2} - 1 & -q - 1 & 0 & 0 & q - 1 & 0 & 0 &
0 \\
q^{2} - 1 & 0 & -q - 1 & 0 & 0 & q - 1 & 0 &
0 \\
q^{2} & \frac{q}{q - 1} & \frac{q}{q - 1} & 0 & 0 &
0 & 1 & 1 \\
0 & 0 & 0 & 1 & \frac{1}{2} \, q & \frac{1}{2} \, q
& \frac{1}{2} \, q - \frac{3}{2} & 0 \\
0 & 0 & 0 & -1 & \frac{1}{2} \, q & \frac{1}{2} \, q
& \frac{1}{2} \, q - \frac{3}{2} & 0 \\
q & -\frac{q}{q - 1} & -\frac{q}{q - 1} & 0 & 0 & 0
& 1 & q
\end{pmatrix}
$$
The suggestive form of these matrices deserves an explanation.
\end{rmk}

\subsection{The genus tube}
\label{sec:genus-tube}

In this section, we discuss the case of the holed torus tube $L: (S^1, \star) \to (S^1, \star)$, see Sections \ref{sec:almost-TQFT-tubes} and \ref{almost-TQFT-strategy} for the precise definition. 
As described there, the standard TQFT induces a morphism $\Zs{\SL{2}}(L): \KM{\SL{2}} \to \KM{\SL{2}}$.
This morphism has been implicitly studied in the articles of \cite{LMN} and \cite{MM} (see also \cite{MM:2016}). However, there, the approach was focused on computing Hodge monodromy representations of representation varieties since the TQFT formalism was not known at that time. Despite of that, as we saw in Section \ref{sec:monodromy-as-mhm}, Hodge monodromy representation has a direct interpretation in terms of $K$-theory classes of mixed Hodge modules. 

Using this interpretation, in this section we will show that the computations of \cite{MM} is nothing but the calculation of the geometric TQFT as above, with some peculiarities. For a while, let us consider the general case of an arbitrary (reductive) group $G$. As explained in \ref{almost-TQFT-strategy}, $\Zs{G}(L)=q_! \circ p^*$ with $p: G^4 \to G$ given by $p(a, g, h, b) = a$ and $q: G^4 \to G$ given by $q(a, g, h, b) = ba[g,h]b^{-1}$ (beware of the change of notation for the arguments). Hence, $\Zs{G}(L)^g = (q_!p^*)^g$. In order to compute it explicitly, recall that we have a pullback diagram
\[
\begin{displaystyle}
   \xymatrix
   {	
	G^7 \ar[r]\ar[d] & G^4 \ar[d]^{q} \\
	G^4 \ar[r]_p & G
   }
\end{displaystyle}   
\]
where the upper arrow is $(a, g, h, g', h', b, b') \mapsto (b^{-1}ab,b^{-1}gb,b^{-1}hb,b)$ and the leftmost arrow is $(a, g, h, g', h', b, b')  \mapsto (a[g,h], g', h', b')$. Repeating the procedure, we have that $\Zs{G}(L)^g = (\beta_g)_! \circ (\alpha_g)^*$ where $\alpha_g, \beta_g: G^{3g+1} \to G$ are given by
\begin{align*}
	\alpha_g(a, g_1, h_1, g_2, \ldots, g_g, h_g, b_1, \ldots, b_g) &= b_g\cdots b_1ab_1^{-1} \cdots b_g^{-1}, \\
	\beta_g(a, g_1, h_1, g_2, \ldots, g_g, h_g, b_1, \ldots, b_g) &= a\prod_{i = 1}^g [g_i,h_i].
\end{align*}
In particular, $\Zs{G}(L)^g \circ \Zs{G}(D)(1) = \RM{G^{3g}}{G}$ under the projection $\gamma_g(g_1, h_1, \ldots, g_g, h_g, b_1, \ldots, b_g) = \beta_g(1, g_1, h_1, \ldots, g_g, h_g, b_1, \ldots, b_g) = \prod_{i\geq 1} [g_i,h_i]$.

Let us come back to the case $G = \SL{2}(\CC)$. The aim of the paper \cite{MM}, as iniciated in \cite{LMN}, is to compute, for $g \geq 0$, the Hodge monodromy $\RM{Y_g/\ZZ_2}{\Bt}$ where
$$
	Y_g = \left\{(A_1, B_1, \ldots, A_g, B_g, \lambda) \in \SL{2}^{2g} \times \left(\CC^*-\left\{\pm 1\right\}\right)\,\left|\,\,\, \prod_{i=1}^g [A_i,B_i] = D_{\lambda} \right.\right\},
$$
and the action of $\ZZ_2$ on $Y_g$ is $-1\cdot(A_i,B_i, \lambda) = (P_0A_iP_0^{-1}, P_0B_iP_0^{-1}, \lambda^{-1})$, where $P_0 \in \SL{2}$ is the matrix of Section \ref{sec:sl2-repr-var} and the projection $w_g: Y_g/\ZZ_2 \to \Bt$ is $w_g(A_i,B_i, \lambda) = \lambda + \lambda^{-1}$. In \cite{MM}, it is proven that $\RM{Y_g/\ZZ_2}{\Bt} \in \cW|_{\Bt} \subseteq \KM{[\SL{2}/\SL{2}]}$ for all $g \geq 0$. Here, $\cW|_{\Bt}$ denotes the submodule of $\cW$ generated by $T_{\Bt}, S_{\pm 2}$ and $S_2 \otimes S_{-2}$. Moreover, they proved that there exists a module homomorphism $M: \cW|_{\Bt} \to \cW|_{\Bt}$ such that, for all $g \geq 1$
$$
	\RM{Y_{g}/\ZZ_2}{\Bt} = M\left(\RM{Y_{g-1}/\ZZ_2}{\Bt}\right).
$$

\begin{rmk}
\begin{itemize}
	\item Strictly speaking, the paper \cite{MM} performs this computation seeing the Hodge monodromy representation as an element of the representation ring of $\Bt = \CC - \left\{\pm 2\right\}$ and focusing on Deligne-Hodge polynomials. However, if we reinterpret the representation ring as mixed Hodge modules (see Section \ref{sec:monodromy-as-mhm}) a complete review of that paper shows that, formally, all the computations are valid as elements of $\KM{\Bt} = \KM{\CC - \left\{\pm 2\right\}}$, seen as a $\K{\MHSq}$-module. For that reason, under this interpretation, \cite{MM} proves the claimed result.
	\item In the paper \cite{MM}, the variety $Y_g$ is denoted by $\overline{X}_4^g$. However, we will not use this notation since it is confusing with ours.
\end{itemize}
\end{rmk}

\begin{lem}
Let $X_{g} = \left\{(A_i, B_i) \in \SL{2}^{2g}\,|\, \prod_{i=1}^g [A_i,B_i]\in D\right\}$ and consider $\omega_g: X_g \to D$ given by $\omega_g(A_i,B_i) =\prod_{i=1}^g [A_i,B_i]$. Then, $\tr^*\RM{Y_g /\ZZ_2}{\Bt} = \RM{X_g}{D}$.
\begin{proof}
We have a pullback diagram
\[
\begin{displaystyle}
   \xymatrix
   {	
	\overline{Y_g/\ZZ_2} \ar[r] \ar[d]_{\omega_g} & Y_g/\ZZ_2 \ar[d]^{w_g} \\
	D \subseteq \SL{2} \ar[r]_{\;\;\;\;\;\tr} & \Bt
   }
\end{displaystyle}   
\]
where the pullback variety is
$$
	\overline{Y_g/\ZZ_2} = \left\{(\overline{(A_i,B_i, \lambda)}, P) \in Y_g/\ZZ_2 \times D \,\left|\,\,\,\, \prod_{i=1}^g [A_i,B_i] = D_\lambda \sim P \right.\right\},
$$
with map $\overline{Y_g/\ZZ_2} \to D$ given by projection over the last variable. Here, $\overline{(A_i,B_i, \lambda)} \in Y_g/\ZZ_2$ denotes the class of $(A_i,B_i, \lambda) \in Y_g$ under the action of $\ZZ_2$. Now, consider the auxiliar varieties
\begin{align*}
	\overline{Y_g/\ZZ_2}' &= \left\{(\overline{A_i,B_i, \lambda, Q}) \in \left(Y_g \times \SL{2}\right)/\ZZ_2 \,\left|\,\,\,\, \prod_{i=1}^g [A_i,B_i] = D_\lambda\right.\right\}, \\
	X_g' &= \left\{(A_i,B_i, \overline{Q}) \in X_g \times \SL{2}/\ZZ_2 \,\left|\,\,\,\, Q^{-1}\prod_{i=1}^g [A_i,B_i]Q = D_\lambda \right.\right\},
\end{align*}
where $\ZZ_2$ acts on $\SL{2}$ by $-1 \cdot Q = QP_0$. These varieties have morphisms $\overline{Y_g/\ZZ_2}' \to \overline{Y_g/\ZZ_2}$ and $X_g' \to X_g$ given by $(\overline{A_i,B_i, \lambda, Q}) \mapsto ((\overline{A_i,B_i, \lambda}), QD_\lambda Q^{-1})$ and $(A_i,B_i, \overline{Q}) \mapsto (A_i,B_i)$, both with trivial monodromy and fiber $\CC^*$.

Moreover, we have an isomorphism $\varphi: \overline{Y_g/\ZZ_2}' \to X_g'$ given by $\varphi(\overline{A_i,B_i, \lambda, Q}) = (QA_iQ^{-1}, QB_iQ^{-1}, Q)$, which fits in a commutative diagram with fibrations as rows
	\[
\begin{displaystyle}
   \xymatrix
   {	
	\CC^* \ar[r] & \overline{Y_g/\ZZ_2}' \ar[r]\ar@{{<}-{>}}[d]^\varphi & \overline{Y_g/\ZZ_2} \\
	\CC^* \ar[r] & X_g' \ar[r] & X_g
   }
\end{displaystyle}   
\]
Therefore, since isomorphisms do not change Hodge monodromy representations, we have
$$
	\tr^*\RM{Y_g/\ZZ_2}{\Bt} = \RM{\overline{Y_g/\ZZ_2}}{D} = \frac{1}{q-1}\RM{\overline{Y_g/\ZZ_2}'}{D} = \frac{1}{q-1}\RM{X_g'}{D} = \RM{X_g}{D},
$$
as we wanted to prove.
\end{proof}
\end{lem}

\begin{cor}\label{prop:formula-M}
For all $g \geq 0$, we have
$$
	\Zg{\SL{2}}(L) \left(\RM{Y_{g-1}/\ZZ_2}{\Bt}\right)|_{\Bt} = (q^3-q)\,\eta \left(\RM{Y_g/\ZZ_2}{\Bt}\right).
$$
\begin{proof}
If we compute the right hand side, we have $\eta \left(\RM{Y_g/\ZZ_2}{\Bt}\right) = \tr_! \circ \tr^* \left(\RM{Y_g/\ZZ_2}{\Bt}\right) = \RM{X_g}{\Bt}$, by the previous Lemma, where the projection is $\tr \circ \omega_g: X_g \to \Bt$.

On the other hand, the left hand side can be rewritten as $\Zg{\SL{2}}(L) \left(\RM{Y_{g-1}/\ZZ_2}{\Bt}\right)|_{\Bt} = \tr_! \circ \Zs{\SL{2}}(L) \circ \tr^* \left(\RM{Y_{g-1}/\ZZ_2}{\Bt}\right)|_{\Bt} = \tr_! \circ \Zs{\SL{2}}(L)\left(\RM{X_{g-1}}{D}\right)|_{\Bt}$. Observe that $\gamma_{g}^{-1}(D) = X_{g} \times \SL{2}^{g}$ so $\RM{X_g}{D} = \frac{1}{(q^3-q)^g}(\gamma_{g}|_{\gamma_{g}^{-1}(D)})_! \underline{\QQ} = \frac{1}{(q^3-q)^g}\left(\Zs{\SL{2}}(L)^g \circ Z(D)(1)\right)|_{D}$. Using that, we finally obtain that
\begin{align*}
	\tr_! \circ \Zs{\SL{2}}(L)\left(\RM{X_g}{D}\right) &= \frac{1}{(q^3-q)^{g-1}}\,\tr_! \circ \Zs{\SL{2}}(L)\left[\Zs{\SL{2}}(L)^{g-1} \circ Z(D)(1)\right]|_{D} \\
	&= \frac{1}{(q^3-q)^{g-1}}\,\tr_! \left[\Zs{\SL{2}}(L)^g \circ Z(D)(1)\right]|_{D} = (q^3-q)\,\tr_! \left(\RM{X_g}{D}\right) \\
	&= (q^3-q)\RM{X_g}{\Bt}.
\end{align*}
This proves the desired equality.
\end{proof}
\end{cor}

\begin{cor}
The morphism $M: \cW|_{\Bt} \to \cW|_{\Bt}$ satisfies that $(q^3-q)M = \eta^{-1} \circ \Zg{\SL{2}}(L)|_{\Bt}$.
\begin{proof}
By the results of \cite{MM}, the set $\left\{\RM{Y_g/\ZZ_2}{\Bt}\right\}_{g \geq 0}$ generate the same submodule as $T_{\Bt}, S_2, S_{-2}$ and $S_2 \otimes S_{-2}$. Therefore, since $\RM{Y_g/\ZZ_2}{\Bt} = M\left(\RM{Y_{g-1}/\ZZ_2}{\Bt}\right)$, from Proposition \ref{prop:formula-M} we obtain that $(q^3-q) \eta \circ M = \Zg{\SL{2}}(L)$ on the submodule generated by $T_{\Bt}, S_2, S_{-2}$ and $S_2 \otimes S_{-2}$.
\end{proof}
\end{cor}

Indeed, in \cite{MM}, Section 9, a larger homomorphism $M: \cW \to \cW$ is defined for which the previous one is just its restriction to $\cW|_{\Bt}$. The explicit expression of $M$ in the set of generators $\cW$ is
$$
{\left( {\fontsize{0.1}{4}\selectfont
\begin{array}{c@{\hspace{2em}}c@{\hspace{2em}}c@{\hspace{2em}}c@{\hspace{2em}}c@{\hspace{2em}}c@{\hspace{2em}}c@{\hspace{2em}}c}
 q^4+4q^3 & q^3-q & q^5 -2q^4-4q^3 & q^5+3q^4   & q^6 -2q^5 -4q^4 & -q^5-4q^4 & 2q^5 -7q^4 -3q^3 & -5q^4-q^3  \\
 -q^2-4q  &       & +2q^2+3q       & -q^3 -3q^2 & +3q^2 +2q       & +4q^2+q  & +7q^2+q         & +5q^2 +q    \\
& & & & & & & \\ 
q^3 -q & q^4 +4q^3 & q^5+3q^4  & q^5-2q^4-4q^3 & q^6-2q^5-4q^4 & 2q^5-7q^4-3q^3 & -q^5-4q^4 & -5q^4-q^3 \\
       & -q^2-4q   & -q^3-3q^2 & +2q^2+3q      & +3q^2+2q      & +7q^2 +q       & +4q^2+q   &  +5q^2+q \\
 & & & & & & & \\
q^3-2q^2 & q^3+3q^2 & q^5+q^4  & q^5 -3q^3 & q^6-2q^5-3q^4 & -q^5+2q^4  & -q^5 -q^4   & -2q^4 -q^3 \\
-3q      &          & +3q^2+3q & -6q^2     & +q^3+3q^2     & -4q^3+3q^2 & -4q^3 +6q^2 & +3q^2 \\
 & & & & & & & \\
q^3+3q^2 & q^3-2q^2 & q^5-3q^3 & q^5+q^4  & q^6-2q^5-3q^4 & -q^5-q^4   & -q^5 +2q^4  & -2q^4 -q^3 \\
         & -3q      & -6q^2    & +3q^2+3q & +q^3+3q^2     & -4q^3+6q^2 & -4q^3 +3q^2 & +3q^2 \\
 & & & & & & & \\
q^3 & q^3 & q^5-3q^3 & q^5-3q^3 & q^6 -2q^5-2q^4 & -q^5-q^4 & -q^5-q^4 & -2q^4 \\
    &     &          &          & +4q^3+q^2      & +2q^3    & +2q^3    &       \\
 & & & & & & & \\
-3q & 3q^2 & 3q^2+3 & -6q^2 & -3q^3+3q^2 & 4q^4-6q^3 +4q^2 & -8q^3+6q^2 & -3q^3+3q^2 \\
& & & & & & & \\
3q^2 & -3q & -6q^2 & 3q^3+3q & -3q^3+3q^2 & -8q^3+6q^2 & 4q^4-6q^3+4q^2 & -3q^3+3q^2 \\
& & & & & & & \\
-1 & -1 & 2q^2 & 2q^2 & -4q^2+2 & -2q^2+q+1 & -2q^2+q+1 & q^4 -2q^2 \\
   &    &      &      &         &           &           & +2q+1 
\end{array}} \right)}.
$$
Analogous (and simpler) calculations can be done as in Corollary \ref{prop:formula-M} in order to show that $(q^3-q) \eta \circ M(T_{\pm 1}) = \Zg{\SL{2}}(L)(T_{\pm 1})$ and $(q^3-q)\eta \circ M(T_{\pm}) = \Zg{\SL{2}}(L)(T_{\pm})$. Hence, $(q^3-q)\eta \circ M = \Zg{\SL{2}}(L)$ on the whole $\cW$. This proves the following result.

\begin{cor}\label{prop:formula-M-TQFT}
The morphism $M: \cW \to \cW$ satisfies that $(q^3-q)M = \eta^{-1} \circ \Zg{\SL{2}}(L)$.
\end{cor}

Corollary \ref{prop:formula-M-TQFT} can be restated as that $M: \cW_0 \to \cW_0$ is (up to a constant) the left $\tr$-reduction of the standard TQFT, as explained in Remark \ref{rmk:left-reduction}. However, we have chosen to consider right $\tr$-reductions in Section \ref{sec:geometric-reduced-TQFT} since they are more natural from the geometric point of view. In any case, from this expression we can recover the geometric TQFT as $\Zg{\SL{2}}(L) = (q^3-q)\eta M$. Hence, in the set of generators $\cS$, the matrix of $\Zg{\SL{2}}(L)$ is 
$$
\scriptsize{
\fontsize{0.8}{10} \selectfont (q^3-q)^2      	
\begin{pmatrix}
q + 4 & 1 & q^{2} - 2 q - 3 & q^{2} + 3  q & \substack{q^{3}
- 2  q^{2}\\ - 3  q - 2} & -q^{2} - 4  q - 1 & 2  q^{2} - 7
 q - 1 & -5  q - 1 \\[0.15cm]
1 & q + 4 & q^{2} + 3  q & q^{2} - 2  q - 3 & \substack{q^{3}
- 2  q^{2}\\ - 3  q - 2} & 2  q^{2} - 7  q - 1 & -q^{2} - 4
 q - 1 & -5  q - 1 \\[0.15cm]
q^{2} - 2  q - 3 & q^{2} + 3  q & \substack{q^{4} + q^{3} \\+ 3  q + 3}
& q^{4} - 3  q^{2} - 6  q & \substack{q^{5} - 2  q^{4} - 3  q^{3}\\
+ q^{2} + 3  q} & \substack{-q^{4} + 2  q^{3}\\ - 4  q^{2} + 3  q} &
\substack{-q^{4} - q^{3}\\ - 4  q^{2} + 6  q} & -2  q^{3} - q^{2} + 3  q
\\[0.15cm]
q^{2} + 3  q & q^{2} - 2  q - 3 & q^{4} - 3  q^{2} - 6 
q & \substack{q^{4} + q^{3} \\+ 3  q + 3} & \substack{q^{5} - 2  q^{4} - 3  q^{3}
\\+ q^{2} + 3  q} & \substack{-q^{4} - q^{3} \\- 4  q^{2} + 6  q} & \substack{-q^{4}
+ 2  q^{3}\\ - 4  q^{2} + 3}  q & -2  q^{3} - q^{2} + 3  q \\[0.15cm]
q^{2} + 1 & q^{2} + 1 & q^{4} - 2  q^{2} & q^{4} - 2 
q^{2} & \substack{q^{5} - 2  q^{4} - q^{3}\\ + 2  q^{2} - 2} & \substack{-q^{4} - q^{3} \\+ q^{2} - q - 1} & \substack{-q^{4} - q^{3} \\+ q^{2} - q - 1} & \substack{-2 
q^{3} \\+ q^{2} - 2  q - 1} \\[0.15cm]
0 & 3  q & 3  q^{2} & -3  q & -3  q^{2} & 4
 q^{3} - 6  q^{2} & -4  q^{2} & -3  q^{2} \\[0.15cm]
3  q & 0 & -3  q & 3  q^{2} & -3  q^{2} & -4
 q^{2} & 4  q^{3} - 6  q^{2} & -3  q^{2} \\[0.15cm]
q & q & q^{3} & q^{3} & \substack{q^{4} - 2  q^{3}\\ - q^{2} - 2
 q} & -q^{3} - q^{2} - q & -q^{3} - q^{2} - q & q^{3} - 2
 q^{2} - q
\end{pmatrix}}.
$$
Analogous calculation can be done to obtain $\cZg{\SL{2}}(L) = (q^3-q)\eta M\eta^{-1}$.

Let $\Sigma_g$ be the closed surface of genus $g \geq 0$. For the parabolic data $\Lambda = \left\{\left\{- \Id\right\},[J_+], [J_-]\right\}$ we consider a parabolic structure on $\Sigma_g$, $Q = \left\{(p_1, [C_1]), \ldots, (p_s, [C_s])\right\}$, where $C_i = J_+, J_-$ or $-\Id$ and $p_1, \ldots, p_s \in \Sigma_g$. Let us denote by $r_+$ be the number of $[J_+]$ in $Q$, $r_-$ the number of $[J_-]$ and $t$ the number of $-\Id$ (so that $r_+ + r_- + t = s$). Set $r = r_+ + r_-$ and $\sigma = (-1)^{r_- + t}$.

\begin{thm}\label{thm:Hodge-repr-sl2}
The $K$-theory image of the cohomology of $\Rep{\SL{2}(\CC)}(\Sigma_g, Q)$ is
\begin{itemize}
	\item If $\sigma = 1$, then
\begin{align*}
	\coh{\Rep{\SL{2}(\CC)}(\Sigma_g, Q)} =&\, {\left(q^2 - 1\right)}^{2g + r - 1} q^{2g - 1} +
\frac{1}{2} \, {\left(q -
1\right)}^{2g + r - 1}q^{2g -
1}(q+1){\left({2^{2g} + q - 3}\right)} 
\\ &+ \frac{\left(-1\right)^{r}}{2} \,
{\left(q + 1\right)}^{2g + r - 1} q^{2g - 1} (q-1){\left({2^{2g} +q -1}\right)}.
\end{align*}
	\item If $\sigma = -1$, then
\begin{align*}
	\coh{\Rep{\SL{2}(\CC)}(\Sigma_g, Q)} = &\, {\left(q - 1\right)}^{2g + r - 1} (q+1)q^{2g - 1}{{\left( {\left(q + 1\right)}^{2 \,
g + r-2}+2^{2g-1}-1\right)} } \\
&+ \left(-1\right)^{r + 1}2^{2g - 1}  {\left(q + 1\right)}^{2g + r
- 1} {\left(q - 1\right)} q^{2g - 1}.
\end{align*}
\end{itemize}
\begin{proof}
By Theorem \ref{thm:almost-tqft-parabolic}, we have that
$$
	\coh{\Rep{\SL{2}(\CC)}(\Sigma_g, Q)} = \frac{1}{(q^3-q)^{g+s}}\cZg{\SL{2}}(D^\dag) \circ \cZg{\SL{2}}(L_{[C_s]}) \circ \ldots \circ \cZg{\SL{2}}(L_{[C_1]}) \circ \cZg{\SL{2}}(L)^g \circ \cZg{\SL{2}}(D)(1).
$$
Moreover observe that, as $\cZg{\SL{2}}$ is an almost-TQFT and the strict tubes commute, all these linear morphisms commute. Hence, we can group the $-\Id$ and $J_-$ tubes together and, using that $\cZg{\SL{2}}(L_{-\Id}) \circ \cZg{\SL{2}}(L_{-\Id}) = (q^3-q)^2 1_{\KM{\SL{2}}}$ and $\cZg{\SL{2}}(L_{-\Id}) \circ \cZg{\SL{2}}(L_{[J_-]}) = (q^3-q)\cZg{\SL{2}}(L_{[J_+]})$, we finally have:
\begin{itemize}
	\item If $\sigma = 1$, then all the $-\Id$ tubes cancel so we have
$$
	\coh{\Rep{\SL{2}(\CC)}(\Sigma_g, Q)} = \frac{1}{(q^3-q)^{g+r}}\cZg{\SL{2}}(D^\dag) \circ \cZg{\SL{2}}(L_{[J_+]})^r \circ \cZg{\SL{2}}(L)^g \circ \cZg{\SL{2}}(D)(1).
$$
Now, observe that, as $\cZg{\SL{2}}(L)$ and $\cZg{\SL{2}}(L_{[J_+]})$ commute, they can be simultaneously diagonalized. Hence, there exists $P, A, B: \cW_0 \to \cW_0$ such that $PAP^{-1} = \cZg{\SL{2}}(L)$ and $PBP^{-1}=\cZg{\SL{2}}(L_{[J_+]})$. Therefore, we have that
$$
	\coh{\Rep{\SL{2}(\CC)}(\Sigma_g, Q)} = \frac{1}{(q^3-q)^{g+r}}\cZg{\SL{2}}(D^\dag) PB^rA^gP \cZg{\SL{2}}(D)(1).
$$
From the matrices of $\cZg{\SL{2}}(L)$ and $\cZg{\SL{2}}(L_{[J_+]})$ for the set of generators $\cS$, the matrices $P,A$ and $B$ can be explicitly calculated with a symbolic computation software (like SageMath \cite{sagemath}). Using that matrices, and recalling that $\cZg{\SL{2}}(D)$ and $\cZg{\SL{2}}(D^\dag)$ are the inclusion and projection onto $T_1$ respectively, the calculation follows.
	\item If $\sigma = -1$, then all the $-\Id$ tubes cancel except one, so we have
$$
	\coh{\Rep{\SL{2}(\CC)}(\Sigma_g, Q)} = \frac{1}{(q^3-q)^{g+r+1}}\cZg{\SL{2}}(D^\dag) \circ \cZg{\SL{2}}(L_{-\Id}) \circ \cZg{\SL{2}}(L_{[J_+]})^r \circ \cZg{\SL{2}}(T)^g \circ \cZg{\SL{2}}(D)(1).
$$
Again, decomposing $PAP^{-1} = \cZg{\SL{2}}(T)$ and $PBP^{-1}=\cZg{\SL{2}}(L_{[J_+]})$, we obtain that
$$
	\coh{\Rep{\SL{2}(\CC)}(\Sigma_g, Q)} = \frac{1}{(q^3-q)^{g+r+1}}\cZg{\SL{2}}(D^\dag) \cZg{\SL{2}}(L_{-\Id}) PB^rA^gP \cZg{\SL{2}}(D)(1).
$$
Hence, using a symbolic computation software, the result follows.
\end{itemize}
\end{proof}
\end{thm}

\begin{rmk}
The homomorphism $\cZg{\SL{2}}(L_{[J_+]}): \cW_0 \to \cW_0$ has non-trivial kernel. In this way, if $A$ is the diagonal matrix associated to $\cZg{\SL{2}}(L_{[J_+]})$, as in the proof of Theorem \ref{thm:Hodge-repr-sl2}, we have that $A^0 \neq \Id$ because of the vanishing eigenvalues. For this reason, we cannot directly set $r=0$ in the previous formula in order to obtain the Hodge structure on non-parabolic representation varieties. However, it may be traced back the extra term lost when we set $r=0$ and to obtain the formula for the non-parabolic case
\begin{align*}
	\coh{\Rep{\SL{2}(\CC)}(\Sigma_g)} &= \left.\coh{\Rep{\SL{2}(\CC)}(\Sigma_g, Q)}\right|_{r=0} + q(q^2-1)^{2g-1} \\
	& = {\left(q^2 - 1\right)}^{2g - 1} q^{2g - 1} +
\frac{1}{2} \, {\left(q -
1\right)}^{2g - 1}q^{2g -
1}(q+1){\left({2^{2g} + q - 3}\right)} 
\\ &+ \frac{1}{2} \,
{\left(q + 1\right)}^{2g + r - 1} q^{2g - 1} (q-1){\left({2^{2g} +q -1}\right)} + q(q^2-1)^{2g-1}.
\end{align*}
Here $\left.\coh{\Rep{\SL{2}(\CC)}(\Sigma_g, Q)}\right|_{r=0}$ denotes the Hodge structure obtained when setting $r=0$ in the formula for $\coh{\Rep{\SL{2}(\CC)}(\Sigma_g, Q)}$ in the case $\sigma = 1$ of Theorem \ref{thm:Hodge-repr-sl2}

Analogous considerations can be done in the case of only one marked point with conjugacy class $[-\Id]$. In that case, if we denote $Q_0 = \left\{(p, [-\Id])\right\}$ with $p \in \Sigma_g$ we have 
\begin{align*}
	\coh{\Rep{\SL{2}(\CC)}(\Sigma_g, Q_0)} &= \left.\coh{\Rep{\SL{2}(\CC)}(\Sigma_g, Q)}\right|_{r=0} + q(q^2-1)^{2g-1} \\ &\, {\left(q - 1\right)}^{2g + r - 1} (q+1)q^{2g - 1}{{\left( {\left(q + 1\right)}^{2 \,
g + r-2}+2^{2g-1}-1\right)} } \\
&+ \left(-1\right)^{r + 1}2^{2g - 1}  {\left(q + 1\right)}^{2g + r
- 1} {\left(q - 1\right)} q^{2g - 1} + q(q^2-1)^{2g-1}.
\end{align*}
Now, $\left.\coh{\Rep{\SL{2}(\CC)}(\Sigma_g, Q)}\right|_{r=0}$ denotes the Hodge structure of the case $\sigma = -1$ of Theorem \ref{thm:Hodge-repr-sl2} after setting $r=0$.
\end{rmk}

\begin{rmk}
When interpreting $q$ as a variable, the previous polynomials are also the Deligne-Hodge polynomials of the parabolic representations varieties with $q=uv$. This is due to the fact that $\DelHod{q} = \DelHod{\coh{\CC}} = uv$.

Under this point of view, the results of Theorem \ref{thm:Hodge-repr-sl2} generalize the previous work on Deligne-Hodge polynomials of parabolic representation varieties to the case of an arbitrary number of marked points. For genus $g=1,2$ and one marked point in the class $[J_\pm]$ (i.e.\ $r_+ = 1$ or $r_-=1$), this Theorem agrees with the calculations of \cite{LMN}, Sections 4.3, 4.4, 11 and 12. For arbitrary genus and zero or one marked points, this result agrees with \cite{MM}, Proposition 11. Finally, for the case $g=1$ and two marked points in the classes $[J_\pm]$, this Theorem agrees with the computations of \cite{LM}, Section 3. 
\end{rmk}

%% file: Chapters/TopGIT.tex
% Chapter 4

\chapter{Topological Geometric Invariant Theory} % Main chapter title

\label{chap:git} % For referencing the chapter elsewhere, use \ref{Chapter1} 

\lhead{Chapter 4. \emph{Topological Geometric Invariant Theory}} % This is for the header on each page - perhaps a shortened title

%----------------------------------------------------------------------------------------

\section{Review of Geometric Invariant Theory}
\label{sec:review-git}
Except when explicitly said, along this section and sections \ref{section:pseudo-quotients} and \ref{section:stratification}, we will work over an arbitraty algebraically closed field $k$.
In this section, we will review some of the most important notions of Geometric Invariant Theory (GIT for short). The definitions are essentially the ones used in Newstead's book \cite{Newstead:1978}.

In order to fix some notation, given an algebraic group $G$ acting algebraically on a variety $X$ we will denote by $Gx$ (or by $[x]$ when $G$ is clear from the context) the orbit of $x \in X$ and by $\overline{Gx}$ (or by $\overline{[x]}$) its Zariski closure. In general, the space of orbits $X/G$ has no structure of algebraic variety so we need to consider more subtle quotients.

Recall that a \emph{categorical quotient} of $X$ by $G$ is a $G$-invariant regular morphism $\pi: X \to Y$ onto some algebraic variety $Y$ such that, for any $G$-invariant regular morphism $f: X \to Z$, with $Z$ and algebraic variety, there exists a unique regular morphism $\tilde{f}: Y \to Z$ such that the following diagram commutes 
\[
\begin{displaystyle}
   \xymatrix
   {
   	X \ar[r]^{f} \ar[d]_{\pi} & Z \\
   	Y \ar@{--{>}}[ru]_{\tilde{f}}&
  	}
\end{displaystyle}
\]
Using this universal property, it follows that the categorical quotient, if exists, is unique up to regular isomorphism.

\begin{rmk}
Along this thesis, we will always work with categorical quotients within the category of algebraic varieties and regular morphisms. However, sometimes in the literature, larger categories are considered.
\end{rmk}

\begin{ex}\label{ex:GIT-affine}
Let $X = Spec(R)$ be an affine algebraic variety, where $R$ is a finitely generated torsion-free $k$-algebra. By considering $R$ as the algebra of regular functions on $X$, the action of $G$ on $X$ induces an action on $R$. By Nagata's theorem (see \cite{Nagata:1963} or Theorem 3.4 of \cite{Newstead:1978}), if $G$ is reductive (i.e.\ if its radical is isomorphic to a torus group) then the ring of $G$-invariant elements of $R$, $R^G$, is finitely generated. Therefore, taking $Y=Spec(R^G)$, the inclusion $R^G \hookrightarrow R$ induces a regular morphism $\pi: X \to Y$ that can be shown to be a categorical quotient.
\end{ex}

The definition of a categorical quotient does not say anything about the geometry of $Y$. Trying to capture these geometric properties, a regular morphism $\pi: X \to Y$ is called a \emph{good quotient} if it satisfies:
\begin{enumerate}[label=$\roman*)$,ref=$\roman*)$]
	\item $\pi$ is $G$-invariant.
	\item $\pi$ is surjective.
	\item For any open set $U \subseteq Y$, $\pi$ induces an isomorphism
	$$
		\pi^*: \cO_Y(U) \stackrel{\cong}{\longrightarrow} \cO_X(\pi^{-1}(U))^G \subseteq \cO_X(\pi^{-1}(U)),
	$$
where $\cO_X$ is the sheaf of regular functions on $X$.
	\item If $W \subseteq X$ is a closed $G$-invariant set, then $\pi(W) \subseteq Y$ is closed.
	\item Given two closed $G$-invariant subsets $W_1, W_2\subseteq X$, $W_1 \cap W_2 = \emptyset$ if and only if $\pi(W_1) \cap \pi(W_2) = \emptyset$.
\end{enumerate}

If $\pi: X \to Y$ is a good quotient, then it is a categorical quotient (see Corollary 3.5.1 of \cite{Newstead:1978}). On the other side, the categorical quotients described in Example \ref{ex:GIT-affine} actually are good quotients (see Theorem 3.5 of \cite{Newstead:1978}).

The last useful standard quotient considered in the literature is the so-called \emph{geometric quotient}. It is a good quotient $\pi: X \to Y$ such that, for any $y \in Y$, $\pi^{-1}(y) = Gy$. This last condition is sometimes referred to as an \emph{orbit space}. Obviously, geometric quotients are categorical so they are unique. Observe that, from the properties of good quotients, a geometric quotient is the same as a good quotient in which the action of $G$ on $X$ is closed (i.e.\ $Gx$ is a closed subset of $X$ for all $x \in X$).

\begin{rmk}
Sometimes in the literature, good quotients are referred to as good categorical quotients and geometric quotients are called good geometric quotients. Moreover, some authors (notably \cite{Newstead:1978}) add to the definition of good quotient the requirement that $\pi$ is affine. The inclusion of this hypothesis is purely technical and is only useful for proving some restrictions to the existence of geometric quotients (see Proposition 3.24 and Remark 3.25 of \cite{Newstead:1978}). In particular, it has no influence in the validity of existence results.
Actually, the addition of this hypothesis is not too restrictive since most times it can be taken for granted (see Propositions 0.7 and 0.8 of \cite{MFK:1994}).
\end{rmk}

In this framework, GIT deals with the problem of existence of these quotients for the action of a reductive group $G$ on an algebraic variety $X$. Suppose also that the action of $G$ is \emph{linearizable}, meaning that there exists a line bundle $L \to X$ with a fiberwise linear action of $G$ compatible with the one on $X$. Observe that, if $L$ is ample, then some tensor power of $L$ gives an embedding $X \hookrightarrow \mathbb{P}^N$ for $N$ large enough. In this setting, the linearization reduces to a linear representation $G \to \GL{N+1}(k)$ such that, when restricted to $X \subseteq \mathbb{P}^N$, gives the original action.

If $L \to X$ is a fixed linearization of the action of $G$, then a point $x \in X$ is called \emph{semi-stable} if there exists a $G$-invariant section $f$ of $L^r$ for some $r>0$ such that $f(x) \neq 0$ and $X_f = \left\{x \in X \,|\, f(x) \neq 0\right\} \subseteq X$ is affine. If the action of $G$ on $X_f$ is, in addition, closed and $\dim Gx = \dim G$, then $x$ is called \emph{stable}. The set of semi-stable and stable points are open subsets of $X$ and we will denote them by $X^{SS}$ and $X^S$, respectively.

With these definitions, the most important result of GIT about the existence of quotients says that if $G$ is a reductive group acting via a linearizable action on a (quasi-projective) variety $X$, then there exists a good quotient on $X^{SS}$ that restricts to a geometric quotient on $X^S$. It is customary to call this good quotient the \emph{GIT quotient} and denote it by $X^{SS} \to X \sslash G$, or $X^{SS} \sslash G$ when we want to emphasize that it is defined only on $X^{SS}$. The proof of this result is just an appropriate gluing of the good quotients constructed in Example \ref{ex:GIT-affine} for an affine covering of $X^{SS}$. Despite that, a priori, the result of this gluing is an algebraic scheme, if we push forward the ample line bundle $L|_{X^{SS}} \to X^{SS}$ to $X \sslash G$ we obtain an ample line bundle there that embedds it into the projective space, turning the scheme into a variety.

\begin{rmk}\label{rmk:polystable-points}
There is another important class of points on $X$, called the \emph{poly-stable points} on $X$ and denoted by $X^{PS}$. These are the points $x \in X$ such that there exists a $G$-invariant section $f$ of $L^r$ for some $r>0$ with $f(x) \neq 0$, $X_f$ affine and the action of $G$ on $X_f$ is closed. Hence, stable points are nothing but poly-stable points with orbits of full dimension (actually, poly-stable points are just called stable on \cite{MFK:1994}). The closure of the orbit of each (semi-stable) point contains exactly one poly-stable point so the orbits of poly-stable points are in natural bijection with the GIT quotient. However, $X^{PS} \subseteq X$ may be not open nor closed, as it happens for $X^{SS}$ and $X^{S}$, so they are less useful for GIT considerations. 
\end{rmk}

\section{Pseudo-quotients and stratification of algebraic quotients}
\label{section:pseudo-quotients}

Let $\pi: X \to Y$ be a good quotient. If $U \subseteq Y$ is an open subset, then the restriction $\pi: \pi^{-1}(U) \to U$ is again a good quotient (see Theorem 3.10 of \cite{Newstead:1978}). However, if we take a closed set $W \subseteq Y$, the restriction $\pi: \pi^{-1}(W) \to W$ could be no longer a good quotient. However, the topological properties of the good quotient remain valid when restricted to $\pi^{-1}(W)$, so it is useful to collect them as a kind of weaker quotient.

\begin{defn}
Let $X$ be an algebraic variety with an action of an algebraic group $G$. A \emph{pseudo-quotient} for the action of $G$ on $X$ is a surjective $G$-invariant regular morphism $\pi: X \to Y$ such that, for any disjoint $G$-invariant closed sets $W_1, W_2 \subseteq X$, $\overline{\pi(W_1)} \cap \overline{\pi(W_2)} = \emptyset$.

\end{defn}

\begin{rmk}\label{rmk:prop-pseudo-quotients}
Suppose that $\pi: X \to Y$ is a pseudo-quotient:
\begin{enumerate}[label=$\roman*)$,ref=$\roman*)$]
	\item Let $x_1,x_2 \in X$. Since $\pi$ is $G$-invariant, $\pi$ maps every point of $\overline{Gx_i}$ into the same point of $Y$ (i.e.\ $\pi(x_i)$). Therefore, since the $\overline{Gx_i}$ are closed $G$-invariant sets, $\overline{Gx_1} \cap \overline{Gx_2} = \emptyset$ if and only if $\pi(x_1) \neq \pi(x_2)$.
	\item \label{rmk:prop-pseudo-quotients-closed-sets} Let $W \subseteq X$ be a $G$-invariant closed set and suppose that $\pi(W)$ were not closed. Then, for any $y \in \overline{\pi(W)} - \pi(W)$, we would have that $\pi^{-1}(y)$ and $W$ are closed $G$-invariant sets so $\left\{y\right\} \cap \overline{\pi(W)} = \emptyset$, which is a contradiction. Thus, the image of any $G$-invariant closed set is closed. In particular, good quotients are pseudo-quotients.
	\item \label{rmk:prop-pseudo-quotients-good} Let $U \subseteq Y$ be an open set and let $\pi^*: \cO_Y(U) \to \cO_X(\pi^{-1}(U))$ be the induced ring morphism. Since $\pi$ is $G$-invariant, this morphism factorizes through the inclusion $\cO_X(\pi^{-1}(U))^G \subseteq \cO_X(\pi^{-1}(U))$ so it defines a ring morphism $\pi^*: \cO_Y(U) \to \cO_X(\pi^{-1}(U))^G$. However, in contrast with good quotients, we requiere no longer this morphism to be an isomorphism. For all these reasons, a pseudo-quotient is a regular map satisfying conditions $i), ii), iv)$ and $v)$ of a good quotient, but maybe failing $iii)$.
	\item \label{rmk:prop-pseudo-quotients-restrictions} If $\pi: X \to Y$ is a pseudo-quotient and $W \subseteq X$ is a closed $G$-invariant set, then the restriction $\pi: W \to \pi(W)$ is also a pseudo-quotient. For open sets, an easy adaptation of Lemma 3.6 of \cite{Newstead:1978} shows that, if $U \subseteq Y$ is an open set, then $\pi: \pi^{-1}(U) \to U$ is a pseudo-quotient.
\end{enumerate}
\end{rmk}

\begin{defn}\label{defn:orbitwise-closed}
Let $X$ be a variety and let $G$ be an algebraic group acting on $X$. A subset $A \subseteq X$ is said to be \emph{orbitwise-closed} if, for any $a \in A$, $\overline{Ga} \subseteq A$. Analogously, $A$ is said to be \emph{completely orbitwise-closed} if both $A$ and $X - A$ are orbitwise-closed.
\end{defn}

\begin{ex} A closed invariant subset is orbitwise-closed. An open orbitwise-closed set is completely orbitwise-closed.
\end{ex}

\begin{lem}\label{lem:completely orbitwise-closed}
Let $G$ be an algebraic group acting on a variety $X$ and let $\pi: X \to Y$ be a pseudo-quotient.
\begin{enumerate}[label=$\roman*)$,ref=$\roman*)$]
	\item\label{lem:completely orbitwise-closed:enum:saturated} If $A \subseteq X$ is completely orbitwise-closed, then it is saturated for $\pi$, i.e.\ $\pi^{-1}(\pi(A))=A$.
	\item\label{lem:completely orbitwise-closed:enum:good} If $U$ is open and orbitwise-closed, then $\pi(U)$ is open. Moreover, $\pi|_U: U \to \pi(U)$ is a pseudo-quotient.
\end{enumerate}
\begin{proof}
For \ref{lem:completely orbitwise-closed:enum:saturated}, trivially, $A \subseteq \pi^{-1}(\pi(A))$. For the other inclusion, if $x \in \pi^{-1}(\pi(A))$ then $\pi(x) = \pi(a)$ for some $a \in A$. Hence, since $\pi$ is a pseudo-quotient, $\overline{Gx} \cap \overline{Ga} \neq \emptyset$ which implies that $x \in A$ since $A$ is completely orbitwise-closed. For \ref{lem:completely orbitwise-closed:enum:good}, observe that, by Remark \ref{rmk:prop-pseudo-quotients} \ref{rmk:prop-pseudo-quotients-closed-sets}, if we set $W = X - U$, then $\pi(W)$ is closed so $Y - \pi(W)$ is open. But, since $U$ is completely orbitwise-closed and $\pi$ is surjective, $Y = \pi(U) \sqcup \pi(W)$, so $\pi(U) = Y - \pi(W)$ is open. Finally, the fact that $\pi|_U: U \to \pi(U)$ is a pseudo-quotient follows from Remark \ref{rmk:prop-pseudo-quotients} \ref{rmk:prop-pseudo-quotients-restrictions} since $U$ is saturated.
\end{proof}
\end{lem}

\begin{ex}\label{rmk:uniqueness-pseudo-quotients}
In contrast with categorical quotients, pseudo-quotients may not be unique. As an example, let us take $X = \mathbb{A}^2$ and $G = k^*$ acting by $\lambda \cdot (x,y)=(\lambda x, \lambda^{-1}y)$, for $\lambda \in k^*$ and $(x,y) \in \mathbb{A}^2$. Since the ring of $G$-invariant regular functions on $X$ is $\cO_X(X)^G = k[xy]$, standard GIT results show that the inclusion $k[xy] \hookrightarrow k[x,y] =  \cO_X(X)$ induces a good quotient $\pi: X \to \mathbb{A}^1$. Now, let $C = \left\{y^2 = x^3\right\} \subseteq \mathbb{A}^2$ be the nodal cubic curve. The map $\alpha: \mathbb{A}^1 \to C$, $\alpha(t) = (t^2, t^3)$, is a regular bijective morphism so $\alpha \circ \pi: X \to C$ is a pseudo-quotient for $X$. However, $C$ is not isomorphic to $\mathbb{A}^1$ since $C$ is not normal.
\end{ex}

\subsection{Uniqueness of pseudo-quotients}
\label{sec:uniqueness-pseudo-quot}

The previous Example \ref{rmk:uniqueness-pseudo-quotients} is general in the way that, if $\pi: X \to Y$ is a pseudo-quotient and $\alpha: Y \to Y'$ is any regular bijective morphism, then $\alpha \circ \pi: X \to Y'$ is also a pseudo-quotient. However, this non-uniqueness of pseudo-quotients is, in some sense, the typical case. We are going to devote the rest of this section to study this kind of uniqueness properties. The first result in this direction is the following.

\begin{prop}\label{prop:pseudo-quotient-regular-bijective}
Let $X$ be an algebraic variety acted on by an algebraic group $G$. Suppose that $\pi: X \to Y$ and $\pi': X \to Y'$ are pseudo-quotients and that $\pi$ is a categorical quotient. Then, there exists a regular bijective morphism $\alpha: Y \to Y'$.
\begin{proof}
By definition, the map $\pi': X \to Y'$ is $G$-invariant so, using the categorical property of $Y$, it defines a regular map $\alpha: Y \to Y'$ such that $\pi'=\alpha \circ \pi$. The surjectivity of $\alpha$ follows from the one of $\pi'$. For the injectivity, suppose that $\alpha(y)=\alpha(y')$ for some $y,y' \in Y$. Then, there exists $x,x' \in X$ such that $y = \pi(x)$ and $y' = \pi(x')$ so $\pi'(x)=\alpha(\pi(x))=\alpha(\pi(x')) = \pi'(x')$. Thus, since $\pi'$ is a pseudo-quotient, $\overline{Gx} \cap \overline{Gx'} \neq \emptyset$ so $y=\pi(x) = \pi(x') = y'$.
\end{proof}
\end{prop}

If there exists a pseudo-quotient $\pi: X \to Y$ which is also categorical, we will say that $X$ admits a categorical pseudo-quotient. Observe that, if that occurs, any categorical quotient of $X$ is automatically a pseudo-quotient. This happens, for example, if $G$ is a reductive group acting linearly on $X$ and all the points of $X$ are semi-stable for the action.

In order to obtain stronger results about uniqueness of pseudo-quotients, for the rest of the section we will suppose that $k$ is an algebraically closed field of characteristic zero.

\begin{rmk}\label{rmk:properties-regular-bijective}
\begin{enumerate}[label=$\roman*)$,ref=$\roman*)$]
	\item \label{rmk:properties-regular-bijective:polynomial} If there exists a regular bijective morphism $\alpha: X \to Y$ between algebraic varieties, then $X$ and $Y$ define the same element in the $K$-theory of algebraic varieties (see Remark \ref{rmk:KVar}). This is a consequence of the fact that, in characteristic zero, every dominant injective regular morphism is birational (essentially, because every field extension is separable and, thus, the degree of $\alpha$ is the degree of the field extension $K(X)/K(Y)$). Thus, there exists proper subvarieties $X' \subseteq X$ and $Y' \subseteq Y$ such that $X-X'$ and $Y-Y'$ are isomorphic.  Therefore, arguing inductively on $X'$ and $Y'$, there exists a stratification of $X$ and another one for $Y$ with isomorphic pieces by pairs. Hence, in the $K$-theory of algebraic varieties, $X$ and $Y$ define the same object.
	\item\label{rmk:properties-regular-bijective:isomorphism} If $\alpha: X \to Y$ is a regular bijective morphism and $Y$ is normal then, indeed, $\alpha$ is an isomorphism. To check that, observe that, as we argued above, $\alpha$ is a birational equivalence. Then, by Zariski's main theorem (see \cite{EGAII}), $\alpha$ is a regular isomorphism of $X$ with some open subvariety $V \subseteq Y$. But, since $\alpha$ is surjective, $V = Y$, proving that $\alpha$ was an isomorphism. 
	In particular, if $\pi: X \to Y'$ is a pseudo-quotient with $Y'$ normal and the action on $X$ admits a categorical pseudo-quotient $X \to Y$ then $Y \cong Y'$ and $\pi$ is categorical too.
\end{enumerate}
\end{rmk}

In the general case, it can happen that the categorical quotient does not exist. Even in this case, in characteristic zero it is possible to compare pseudo-quotients. The key point is the following proposition that adapts Proposition 0.2 of \cite{MFK:1994} to the context of pseudo-quotients. Also, compare it with Remark \ref{rmk:properties-regular-bijective} above.

\begin{prop}\label{prop:pseudo-quotient-is-good}
Let $\pi: X \to Y$ be a pseudo-quotient for the action of some algebraic group $G$ on $X$. If $X$ is irreducible and $Y$ is normal, then $\pi$ is a good quotient.
\begin{proof}
As we mentioned in Remark \ref{rmk:prop-pseudo-quotients} \ref{rmk:prop-pseudo-quotients-good}, it is enough to prove that, for any open set $U \subseteq Y$, the induced ring morphism $\pi^*: \cO_Y(U) \to \cO_X(\pi^{-1}(U))^G \subseteq \cO_X(\pi^{-1}(U))$ is an isomorphism. Recall that, since $\pi$ is $G$-invariant, $\pi^*\left(\cO_Y(U)\right) \subseteq \cO_X(\pi^{-1}(U))^G$. Moreover, since $\pi$ is surjective, $\pi^*$ is injective.

In order to prove the surjectivity of $\pi^*$, let $f: \pi^{-1}(U) \to k=\mathbb{A}^1$ be a regular $G$-invariant function. We want to show that $f$ can be lifted to a regular morphism $\tilde{f}: U \to \mathbb{A}^1$ such that $f= \tilde{f} \circ \pi$. The only possible candidate is $\tilde{f}(y)=f(x)$ for any $x \in \pi^{-1}(y)$, which is well-defined since $\pi$ is a pseudo-quotient.
%The candidate is unique since we must define such map by $\tilde{f}(y)=f(x)$ for any $x \in \pi^{-1}(y)$. Since $\pi$ is a pseudo-quotient, $\tilde{f}$ is well-defined.
Observe that $\tilde{f}$ is continuous in the Zariski topology. To check that, let $W \subseteq \mathbb{A}^1$ a closed set. Then, the set $f^{-1}(W)$ is a closed $G$-invariant set so, by Remark \ref{rmk:prop-pseudo-quotients} \ref{rmk:prop-pseudo-quotients-closed-sets}, $\pi(f^{-1}(W)) \subseteq U$ is closed. But, by construction, $\pi(f^{-1}(W)) = \tilde{f}^{-1}(W)$, proving the continuity of $\tilde{f}$.
%Since $f$ is continuous, $f^{-1}(W) \subseteq \pi^{-1}(U)$ is closed. Moreover, it is $G$-invariant since $f$ also is. Thus, since $\pi$ is a pseudo-quotient, $\pi(f^{-1}(W)) \subseteq U$ is closed. But, by construction, $\pi(f^{-1}(W)) = \tilde{f}^{-1}(W)$, proving that $\tilde{f}$ is continuous.

Therefore, it is enough to prove that $\tilde{f}: U \to \mathbb{A}^1$ is regular. For this purpose, construct the morphism $\pi'= f \times \pi : \pi^{-1}(U) \to \mathbb{A}^1 \times U$ and let $U' \subseteq \mathbb{A}^1 \times U$ be the closure of $\pi'(\pi^{-1}(U))$. Denoting by $p_1: \mathbb{A}^1 \times U \to \mathbb{A}^1$ and $p_2: \mathbb{A}^1 \times U \to U$ the first and second projections respectively and $\omega = p_2|_{U'}$ we have a commutative diagram
\[
\begin{displaystyle}
   \xymatrix
   {
  	 \pi^{-1}(U) \ar[rr]^{\pi'}\ar[rd]_{\pi} \ar@/^2.5pc/[rrrr]^{f} && U' \ar@{^{(}-{>}}[r] \ar[ld]^{\omega} & \mathbb{A}^1 \times U \ar[r]^{p_1} & \mathbb{A}^1 \\
  	 & U &&&
   }
\end{displaystyle}   
\]
Observe that, in order to finish the proof, it is enough to prove that $\omega$ is an isomorphism, since, in that case $\tilde{f} = p_1 \circ \omega^{-1}$ would be regular. First of all, $\omega$ is surjective since $\pi$ is. Moreover, $\omega$ is injective on $\pi'(\pi^{-1}(U)) \subseteq U'$ since, for any $x,x' \in \pi^{-1}(U)$, if $\omega(\pi'(x)) = \omega(\pi'(x'))$ then $\pi(x) = \pi(x')$ which happens if and only if $\overline{Gx} \cap \overline{Gx'} \neq \emptyset$. In that case, since $f$ is $G$-invariant, $f(x)=f(x')$ and, thus, $\pi'(x)=(f(x), \pi(x)) = (f(x'), \pi(x')) = \pi'(x')$.

Therefore, in order to conclude that $\omega$ is bijective, it is enough to prove that $\pi'$ is surjective. Let $A = U' - \pi'(\pi^{-1}(U))$ which can be described as the set of points $z \in U'$ such that, for any $x \in \pi^{-1}(U)$ with $\pi(x)=\omega(z)$, it holds $\pi'(x) \neq z$. Using the continuity of $\tilde{f}$,  we can rewrite $A$ as the open set
$$
	A = \left\{(l, y) \in U' \,|\, \tilde{f}(y) \neq l\right\}.
$$
Thus, if $A$ was non-empty, then, since $\pi'(\pi^{-1}(U)) \subseteq U'$ is dense, $\pi'(\pi^{-1}(U)) \cap A \neq \emptyset$, which is impossible.

Hence, $\pi'$ is surjective and, thus, $\omega: U' \to U$ is a regular bijective morphism. But, since the characteristic of $k$ is zero, as mentioned in Remark \ref{rmk:properties-regular-bijective}, Zariski's main theorem implies that $\omega$ is an isomorphism, as desired.
\end{proof}
\end{prop}

\begin{cor}\label{cor:equality-epol-pseudo-quotients}
Let $X$ be an algebraic variety with an action of an algebraic group $G$. For any pseudo-quotients $\pi_1: X \to Y_1$ and $\pi_2: X \to Y_2$, the varieties $Y_1$ and $Y_2$ define the same element in the $K$-theory of algebraic varieties.
\begin{proof}
First of all, observe that, restricting to the irreducible components of $X$ if necessary, we can suppose that $X$ is irreducible. Let $Y^{'}_i \subseteq Y_i$ be the open subset of normal points of $Y_i$ for $i = 1,2$. Since the $\pi_i^{-1}(Y_i^{'}) \subseteq X$ are open orbitwise-closed sets, $U = \pi_1^{-1}(Y_1^{'}) \cap \pi_2^{-1}(Y_2^{'}) \subseteq X$ is an open orbitwise-closed set. Therefore, the restrictions $\pi_i|_U: U \to \pi_i(U) \subseteq Y_i^{'}$ are pseudo-quotients onto normal varieties so, by Proposition \ref{prop:pseudo-quotient-is-good}, they are good quotients. In particular, they are categorical quotients so, by uniqueness, $\pi_1(U)$ is isomorphic to $\pi_2(U)$.

Therefore, proceding inductively on $X-U$, we find a stratification $Y_1 = Z_1 \sqcup \ldots \sqcup Z_s$ and $Y_2 = \hat{Z}_1 \sqcup \ldots \sqcup \hat{Z}_s$ such that $Z_j$ is isomorphic to $\hat{Z}_j$ for $j = 1, \ldots, s$, so they define the same object in the $K$-theory of algebraic varieties.
\end{proof}
\end{cor}

\subsection{Stratification techniques}
\label{section:stratification}

In this section, we will show that the previous results can be used to reconstruct the quotient of an algebraic variety from the corresponding quotients of the pieces of a suitable decomposition. This kind of arguments will be very useful in the following computations. Except when explicitly noted, the arguments of this section works on any algebraically closed ground field $k$.

\begin{thm}\label{prop:decomposition-quotient}
Let $X$ be an algebraic variety with an action of an algebraic group $G$ with a pseudo-quotient $\pi: X \to \overline{X}$. Suppose that we have a decomposition $X = Y \sqcup U$ where $Y$ is a closed subvariety and $U$ is an open orbitwise-closed subvariety. Then, we have a decomposition
$$
	\overline{X} = \pi(Y) \sqcup \pi(U)
$$
where $\pi(U) \subseteq \overline{X}$ is open, $\pi(Y) \subseteq \overline{Y}$ is closed and the maps $\pi|_U: U \to \pi(U)$ and $\pi|_Y: Y \to \pi(Y)$ are pseudo-quotients.

Furthermore, if $k$ has characteristic zero then, for any pseudo-quotients $Y \to {\overline{Y}}$ and $U \to {\overline{U}}$, we have that $[\overline{X}] = [\overline{Y}] + [\overline{U}]$ in the $K$-theory of algebraic varieties.

\begin{proof}
The decomposition $\overline{X} = \pi(Y) \sqcup \pi(U)$ and properties of $\pi(Y)$ and $\pi(U)$ follows immediately from the surjectivity of $\pi$ together with Remark \ref{rmk:prop-pseudo-quotients} \ref{rmk:prop-pseudo-quotients-restrictions} and Lemma \ref{lem:completely orbitwise-closed} \ref{lem:completely orbitwise-closed:enum:good}. For the last part, use Corollary \ref{cor:equality-epol-pseudo-quotients}.
\end{proof}
\end{thm}

\begin{rmk} Along this remark, with the notations of Theorem \ref{prop:decomposition-quotient}, suppose that $\pi: X \to \overline{X}$ is good.
\begin{itemize}
	\item In that case, $\pi|_U: U \to \pi(U)$ is also good. However, a priori nothing more can be said about the closed part $\pi|_Y: Y \to \pi(Y)$.
	\item Furthermore, if $\pi(Y)$ is normal and $k$ has characteristic zero, then, for any categorical quotients $Y \to {\overline{Y}}$ and $U \to {\overline{U}}$ we have $\overline{X} = \overline{Y} \sqcup \overline{U}$. It follows immediately from the uniqueness of categorical quotients together with the observation that the pseudo-quotient $\pi|_Y: Y \to \pi(Y)$ is good by Proposition \ref{prop:pseudo-quotient-is-good}.
\end{itemize}
\end{rmk}

\begin{ex}
The hypothesis that $\pi: X \to \overline{X}$ is a pseudo-quotient is needed in Theorem \ref{prop:decomposition-quotient}, even if $k = \CC$. Consider $G = \CC$ and $X = \mathbb{A}^2$ with the action $\lambda \cdot (x,y)=(x, y + \lambda x)$, for $\lambda \in \CC$ and $(x,y) \in \mathbb{A}^2$. Recall that $G$ is not reductive, so classical GIT theory does not garantee that a good quotient for the action exists. Actually, the map $\pi: X \to \overline{X} = \mathbb{A}^1$, $\pi(x,y)=x$, is a categorical quotient but is not a pseudo-quotient so it is not good.

Let us take $U = \left\{x \neq 0\right\} \subseteq X$ and $Y = \left\{x=0\right\} \subseteq X$. Observe that $U$ is orbitwise-closed (actually, the full action is closed), the restriction $\pi|_U: U \to \overline{U} = \pi(U) = \mathbb{A}^1-\left\{0\right\}$ is a good quotient and $\pi(Y) = \left\{0\right\}$. On the other hand, the identity map $Y \to Y$ is a categorical pseudo-quotient for the trivial action of $G$ on $Y$ but $[\overline{X}] \neq [Y] + [\overline{U}]$ (it can be checked using the Deligne-Hodge polynomial) so the conclusion of Theorem \ref{prop:decomposition-quotient} fails.
\end{ex}

Another useful application of pseudo-quotients arises when we we find a subset of an algebraic variety that concentrates the action. In this case, we can obtain much information of the pseudo-quotients of the whole variety from the ones of these special subvarieties.

\begin{prop}\label{prop:core}
Let $X$ be an algebraic variety with an action of an algebraic group $G$. Suppose that there exists a subvariety $Y \subseteq X$ and an algebraic subgroup $H \subseteq G$ such that:

\begin{enumerate}[label=\roman*),ref=$\roman*)$]
	\item\label{prop:core:enum:0} $Y$ is orbitwise-closed for the action of $H$.
	\item\label{prop:core:enum:1} For any $x \in X$, $\overline{Gx} \cap Y \neq \emptyset$.
	\item\label{prop:core:enum:2} For any $W_1, W_2 \subseteq Y$ closed (in $Y$) $H$-invariant subsets, we have $W_1 \cap W_2 \neq \emptyset$ if and only if $\overline{GW_1} \cap \overline{GW_2} \neq \emptyset$.
\end{enumerate}
Such a pair $(Y,H)$ is called a \emph{core}. Suppose that there exists a pseudo-quotient $\pi: X \to \overline{X}$ for the action of $G$ on $X$. Then, $\pi$ restricts to a pseudo-quotient $\pi|_{Y}: Y \to \overline{X}$ for the action of $H$ on $Y$.
\begin{proof}
For the surjectivity of $\pi|_Y$, take $\overline{x} \in \overline{X}$. Since $\pi$ is surjective, $\overline{x} = \pi(x)$ for some $x \in X$ and, by hypothesis \ref{prop:core:enum:1}, there exists $y \in Y$ such that $\overline{Gy} \cap \overline{Gx} \neq \emptyset$. Thus, $\pi(y) = \pi(x) = \overline{x}$.

Now, let $W_1, W_2 \subseteq Y$ be two disjoint closed $H$-invariant subsets. If $\overline{\pi(W_1)} \cap \overline{\pi(W_2)} \neq \emptyset$, then $\pi^{-1}\left(\overline{\pi(W_1)}\right) \cap \pi^{-1}\left(\overline{\pi(W_2)}\right) \neq \emptyset$. But we claim that $\pi^{-1}\left(\overline{\pi(W_i)}\right) = \overline{GW_i}$, so this is impossible. In order to check that, observe that the inclusion $\overline{GW_i} \subseteq \pi^{-1}\left(\overline{\pi(W_i)}\right)$ is trivial. For the other inclusion, since $\overline{GW_i}$ is a closed $G$-invariant set, then $\pi(\overline{GW_i})$ is a closed subset containing $\pi(W_i)$ and, thus, $\overline{\pi(W_i)} \subseteq \pi(\overline{GW_i})$ which implies that $\pi^{-1}\left(\overline{\pi(W_i)}\right) \subseteq \pi^{-1}\left(\pi(\overline{GW_i})\right) = \overline{GW_i}$, as desired.
\end{proof}
\end{prop}

As a final remark, observe that the results about uniqueness of pseudo-quotients are good enough to prove equality of a generalized Euler characteristic $\EuChS$ (see Remark \ref{rmk:KVar}) since we have uniqueness of the $K$-theory class (Corollary \ref{cor:equality-epol-pseudo-quotients}). Therefore, rephrasing the results above, we obtain:

\begin{itemize}
	\item Let $X$ be an algebraic variety with an action of an algebraic group $G$. Then, for any pseudo-quotients $\pi: X \to Y$ and $\pi': X \to Y'$ (categorical or not) we have $\EuCh{Y} = \EuCh{Y'}$. It follows from Corollary \ref{cor:equality-epol-pseudo-quotients}.
	\item Let $X$ be an algebraic variety with a decomposition $X = Y \sqcup U$ with $Y \subseteq X$ closed and $U \subseteq X$ an open orbitwise-closed subset. Then, for any pseudo-quotients $X \to \overline{X}$, $Y \to \overline{Y}$ and $U \to \overline{U}$ we have $\EuCh{\overline{X}} = \EuCh{\overline{Y}} + \EuCh{\overline{U}}$. It follows from Theorem \ref{prop:decomposition-quotient}. In particular, if $G$ is a reductive group acting linearly on $X$ and $X = X^{ss}$, we can consider the usual GIT quotients, which are good, so we have an equality of generalized Euler characteristics
$$
	\EuCh{X \sslash G} = \EuCh{Y \sslash G} + \EuCh{U \sslash G}.
$$
	\item If $(Y, H)$ is a core for the action of $G$ on $X$ then, 
for any pseudo-quotient $Y \to \overline{Y}$ we have $\EuCh{\overline{X}} = \EuCh{\overline{Y}}$. It follows from Proposition \ref{prop:core}.
\end{itemize}

We will focus our attention to Deligne-Hodge polynomials, which are a generalized Euler characteristic (Definition \ref{defn:Deligne-Hodge-pol}). Even more, we will need a subtler analysis of how Deligne-Hodge polynomial behaves for geometric quotiens. To do so, have to show that such geometric quotient is actually a principal bundle. This is precisely the content of the so-called Luna's slice theorem, whose proof can be found in \cite{Luna:1973} (see also \cite{Drezet:2004}, Proposition 5.7).

\begin{thm}[Luna's slice theorem]\label{thm:Luna}
Let $X$ be an affine variety with an action of an algebraic reductive group $G$ on it. Let $X_0 \subseteq X$ the set of points where the action of $G$ is free and closed (on $X$). Then, $X_0$ is open and saturated for the GIT quotient $\pi: X \to X \sslash G$ and the restriction $\pi|_{X_0}: X_0 \to \pi(X_0)$ is a principal $G$-bundle of analytic spaces.
\end{thm}

\begin{cor}\label{cor:luna-thm-epol}
%Let $X$ be an affine variety with an action of an algebraic reductive group $G$ and let $\pi: X \to X \sslash G$ be its GIT quotient. Let $U \subseteq X$ be a $G$-invariant open set and suppose that the action of $G$ on $U$ is closed and free. Then, $\pi|_U : U \to \pi(U) = U \sslash G$ is a principal $G$-bundle.
With the notations and hypothesis of Theorem \ref{thm:Luna}, if $U \subseteq X_0$ is a $G$-invariant open set, then $\pi|_U : U \to \pi(U) = U \sslash G$ is a principal $G$-bundle. Moreover, if $G$ is connected, we have $
	\DelHod{U \sslash G} = \DelHod{U}/\DelHod{G}$.
\begin{proof}
By Luna's slice theorem \ref{thm:Luna}, $\pi|_{X_0}: X_0 \to \pi(X_0)$ is a principal $G$-bundle. Since $U$ is an orbitwise-closed open set, Lemma \ref{lem:completely orbitwise-closed} \ref{lem:completely orbitwise-closed:enum:saturated} implies that $U$ is saturated and $\pi(U)$ is open. In that setting, the restriction $\pi|_U: U \to \pi(U)$ is also a principal $G$-bundle. For the last part, just use Remark 2.5 of \cite{LMN} to observe that $\pi|_U$ has trivial monodromy, so $\DelHod{U} = \DelHod{G} \DelHod{U \sslash G}$.
\end{proof}
\end{cor}

\section{GIT properties of character varieties}

\label{sec:rep-var}

In this section, we will review some well-known results about the structure of representation varieties that will be also useful to fix the notation. We will work on an arbitrary algebraically closed field $k$.

Recall from Section \ref{sec:intro-representation} that, for a finitely generated group $\Gamma$ and an algebraic group $G$, the representation variety of $\Gamma$ into $G$ is denoted by $\Rep{G}(\Gamma) = \Hom(\Gamma, G)$. The GIT quotient of the representation variety is $\Char{G}(\Gamma) = \Rep{G}(\Gamma) \sslash G$, the so-called character variety as explained in Section \ref{sec:parabolic-repr-var}. The aim of this chapter is to study the properties of this quotient and to show how to compute the Deligne-Hodge polynomial of $\Char{G}(\Gamma)$ from the Hodge structures of $\Rep{G}(\Gamma)$ that we computed in Theorem \ref{thm:Hodge-repr-sl2}.

We will focus on the cases of Example \ref{ex:gamma-for-repr}, that is $\Gamma = F_n$, the free group with $n$ generators, for which the representation variety will be shorten $\Xf{n}(G) = \Rep{G}(F_n) = G^n$ or just $\Xf{n}$. The other important case will be for $\Gamma = \pi_1(\Sigma_g)$, the fundamental group of the genus $g$ compact surface. In that case, we will shorten the associated representation variety by $\Xs{g}(G) = \Rep{G}(\pi_1(\Sigma_g))$, or even just $\Xs{g}$, and the corresponding character variety by $\cR_g(G)$ or just $\cR_g$.

\subsection{Stability of representation varieties}

Before studying the general case of representation varieties, let us consider a simpler case. Fixed an integer $m \geq 1$, consider the action of $\GL{m}(k)$ on itself by conjugation, that is $P \cdot X = PXP^{-1}$ for $P, X \in \GL{m}(k)$. This action is linearizable since, if $\mathbb{M}_m$ is the vector space of $m \times m$ matrices, then, for any $P \in \GL{m}(k)$, the map $\mathbb{A}^1_k \oplus \mathbb{M}_m \to \mathbb{A}^1_k \oplus \mathbb{M}_m $ given by $\lambda + X \mapsto \lambda + P \cdot X = \lambda + PXP^{-1}$ is linear and commutes with homotheties. Hence, we can see $\GL{m}(k)$ as a quasi-affine variety of $\PP(k \oplus \mathbb{M}_m) = \PP_k^{m^2}$ in which the action by commutation is the restriction of a linear one. In other words, $L=\cO_{\PP_k^{m^2}}(1)|_{\GL{m}(k)} \to \GL{m}(k)$ is an ample $G$-line bundle compatible with the conjugacy action.

Actually, the same holds for any linear algebraic group $G \subseteq \GL{m}(k)$ and its action by conjugation. For that, let $i: G \to \GL{m}(k)$ be the embedding of $G$ as linear group. Since conjugation is just twice product on the group, $i$ is a $G$-equivariant map for the conjugacy action. Hence, if $L \to \GL{m}(k)$ is the linearization of the action above, then $L' = i^*L \to G$ is a linearization for the action of $G$ by conjugation. Moreover, since $L$ is ample, $L'$ is too.

\begin{prop}
Let $\Gamma$ be a finitely generated group and let $G$ be a linear algebraic group. Then, the action of $G$ on $\Rep{G}(\Gamma)$ is linearizable via an ample line bundle.
\begin{proof}
As we mentioned above, $\Rep{G}(\Gamma) \subseteq \Xf{n}(G) = G^n$ for some $n > 0$, so it is enough to prove it in the free case.
In that situation, just take the line bundle $\bigotimes_{i=1}^n \pi_i^* L' \to \Xf{n}(G)$ where $\pi_i: G^n \to G$ is the $i$-th projection and $L' \to G$ is the line bundle constructed above.
%Just embedd $\Xf{n} = G^n \to \PP^m \times \ldots \PP^m \hookrightarrow \PP^N$ for $N$ large enough and observe that the conjugation action on $\PP^N$ is compatible with the Segre embedding.
%Since the action of $G$ on itself is linearizable via an ample line bundle, we have an action of $G$ on $\mathbb{P}^m$ and a $G$-equivariant embedding $f: G \hookrightarrow \mathbb{P}^m$. With this embedding, we construct $g = f \times \ldots \times f: \Xf{n} = G^n \hookrightarrow \mathbb{P}^m \times \ldots \times \mathbb{P}^m \subseteq \mathbb{P}^N_k$, for $N$ large enough, with an induced action of $G$ on $\mathbb{P}^N$ such that $g$ is $G$-equivariant. Using the procedure of Remark \ref{rmk:linearization-abstract}, this is the same as a linearization of the action of $G$ on $X$ via an ample bundle.
\end{proof}
\end{prop}

Recall that a linear representation $\rho: \Gamma \to \mathrm{GL}(V)$, where $V$ is a finite dimensional $k$-vector space, is said to be reducible if there exists a proper $\Gamma$-invariant subspace of $V$. If $G$ is a linear group, we can see $G \subseteq \GL{m}(k)$ for $m$ large enough, so it also make sense to speak about reducible $G$-representations. With this definition, we have a decomposition
$$
	\Rep{G}(\Gamma) = \Xred{\mathfrak{X}}_G(\Gamma) \sqcup \Xirred{\mathfrak{X}}_G(\Gamma),
$$
where $\Xred{\mathfrak{X}}_G(\Gamma) \subseteq \Rep{G}(\Gamma)$ is the closed subvariety of reducible representations and $\Xirred{\mathfrak{X}}_G(\Gamma)$ the open set of irreducible ones.

\begin{prop}
Let $\mathfrak{X}_{\GL{m}(k)}(\Gamma)^S$ and $\mathfrak{X}_{\GL{m}(k)}(\Gamma)^{SS}$ be the set of stable and semi-stable points, respectively, of $\mathfrak{X}_{\GL{m}(k)}(\Gamma)$ under the action of $\GL{m}(k)$ by conjugation. Then, we have $\mathfrak{X}_{\GL{m}(k)}(\Gamma)^{SS} = \mathfrak{X}_{\GL{m}(k)}(\Gamma)$ and $\mathfrak{X}_{\GL{m}(k)}(\Gamma)^{S} = \Xirred{\mathfrak{X}}_{\GL{m}(k)}(\Gamma)$.
\begin{proof}

We will use the Hilbert-Mumford criterion of stability (see \cite{Newstead:1978}, Theorem 4.9). Recall that this criterion says that $A \in \mathfrak{X}_{\GL{m}(k)}(\Gamma)^S$ (resp. $A \in \mathfrak{X}_{\GL{m}(k)}(\Gamma)^{SS}$) if and only if $\mu(A, \lambda) > 0$ (resp. $\geq 0$) for all $1$-parameter subgroups $\lambda: k^* \to \GL{m}(k)$, where $\mu(A, \lambda)$ is the minimum $\alpha \in \ZZ$ such that $\lim\limits_{t \to 0} t^\alpha \lambda(t) \cdot A$ exists.

Let us prove that all the elements of $\mathfrak{X}_{\GL{m}(k)}(\Gamma)$ are semi-stable. Given a $1$-parameter subgroup $\lambda: k^* \to \GL{m}(k)$, pick a regular function $\alpha: k^* \to k^*$ such that $\alpha(t)^m = \det{\lambda(t)}$ for all $t \in k^*$. In this case, $\overline{\lambda}(t) = \alpha(t)^{-1}\lambda(t) \in \SL{m}(k)$ for all $t \in k^*$, so it is a $1$-parameter subgroup of $\SL{m}(k)$. By \cite{Newstead:1978}, Theorem 4.11, there exists $P \in \SL{m}(k)$ such that
$$
	P^{-1}\overline{\lambda}(t)P = \begin{pmatrix} t^{s_1} & 0 & \cdots & 0 \\ 0 & t^{s_2} & \cdots & 0 \\ \vdots & & \ddots &  \vdots \\ 0 & \ldots && t^{s_m}\end{pmatrix} =: D_{s_1, \ldots, s_m}(t)
$$
for some $s_1 \geq s_2 \ldots \geq s_m \in \RR$, not all zero, and $\sum_i s_i = 0$. Now, given $A \in \mathfrak{X}_{\GL{m}(k)}(\Gamma)$, we have that
$$
	\lambda(t) \cdot A = \lambda(t) A \lambda(t)^{-1} = \overline{\lambda}(t) A \overline{\lambda}(t)^{-1} = PD_{s_1, \ldots, s_m}(t)P^{-1}APD_{s_1, \ldots, s_m}(t)^{-1}P^{-1}.
$$
Hence, without loss of generality, we can suppose that $\lambda(t) = D_{s_1, \ldots, s_m}(t)$. Let $A = (A_1, \ldots, A_n) \in \mathfrak{X}_{\GL{m}(k)}(\Gamma)$, where $A_l = \left(a_{i,j}^l\right)$ for $1 \leq l \leq n$. A direct computation shows that the $(i,j)$-entry of $\lambda(t) \cdot A_l$ is $a_{i,j}^l t^{s_i - s_j}$. In particular, all the elements in the diagonal and under it are multiplied by $t^\beta$ with $\beta \leq 0$. Since at least one of them have to be non-zero, $\lim\limits_{t \to 0} t^{\alpha} \lambda(t) \cdot A$ cannot exist if $\alpha < 0$. Hence, $\mu(A, \lambda) \geq 0$ for all $1$-parameter subgroups, proving that all the points are semi-stable.

For the stable points, suppose that $A = (A_1, \ldots, A_n)$ is reducible with a proper invariant subspace of dimension $k < m$. In that case, maybe after conjugation, we can get that the entries $a_{i,j}^l = 0$ for $i > k$ and $j \leq k$. Taking $\lambda(t) = D_{s_1, \ldots, s_m}(t)$ with $s_i = 1/k$ for $i \leq k$, and $s_i = -1 /(m-k)$ for $i > k$, we have that all the non-zero entries of $\lambda(t)\cdot A_l$ are multiplied by $t^\beta$ with $\beta \geq 0$. Hence, $\lim\limits_{t \to 0} \lambda(t) \cdot A$ exists so $\mu(A, \lambda)=0$ and $A$ cannot be stable.
 
On the other hand, if $A \in \Xf{n}$ is not stable, then $\lim\limits_{t \to 0}\lambda(t) \cdot A$ exists for some $1$-parameter subgroup $\lambda(t)$. As we explained above, maybe after conjugate $A$, we can suppose that $\lambda(t) = D_{s_1, \ldots, s_m}(t)$. In that case, $\lambda(t) \cdot A$ cannot have non-zero entries with a factor $t^\beta$ with $\beta < 0$, which is possible if and only if there is $1 \leq k<m$ such that $s_1 = s_2 = \ldots = s_k > s_{k+1} = \ldots = s_{m}$ and all the entries $a_{i,j}^l$ with $i > k$ and $j \leq k$ of $A$ vanish, that is, if and only if $A$ has a proper invariant subspace.
\end{proof}
\end{prop}

\begin{cor}\label{cor:stability-repr-var}
Let $\Gamma$ a finitely generated group and $G$ a linear algebraic group. Then all the points of $\mathfrak{X}_{G}(\Gamma)$ are semi-stable for the action of $G$ by conjugation and $\mathfrak{X}_{G}(\Gamma)^{S} = \Xirred{\mathfrak{X}}_{G}(\Gamma)$
\begin{proof}
For the first part, observe that the inclusion $G \subseteq \GL{m}(k)$ induces a natural inclusion $\mathfrak{X}_{G}(\Gamma) \subseteq \mathfrak{X}_{\GL{m}(k)}(\Gamma)$. Hence, since all the points of $\mathfrak{X}_{\GL{m}(k)}(\Gamma)$ are semi-stable for the conjugacy action of $\GL{m}(k)$, the same holds for $\mathfrak{X}_{G}(\Gamma)$. For the stable points, just use the same argument than above taking into account that the reducibility of a representation does not change after a conjugation.
\end{proof}
\end{cor}

By definition, the action of $G$ on the stable locus is closed. However, we cannot expect that the action of $G$ there was free because all the elements in the center of $G$, $G^0 \subseteq G$, act trivially. However, it turns out that, along the irreducible representations, the action of $G/G^0$ is, indeed, free.

\begin{prop}\label{prop:action-free-on-irred}
Suppose that $k$ is an algebraically closed field. The action of $G/G^0$ on the irreducible representations $\Xirred{\mathfrak{X}_{G}}(\Gamma)$ is free. In particular $\Xirred{\mathfrak{X}_{G}}(\Gamma) \to \Xirred{\mathfrak{X}_{G}}(\Gamma) \sslash (G/G^0)$ is a free geometric quotient.
\begin{proof}
Suppose that $G \subseteq \GL{m}(k)$ for $m$ large enough. Let $A = (A_1, \ldots, A_n) \in \Xirred{\mathfrak{X}_{G}}(\Gamma)$ and suppose that $P \in G$ satisfies $P \cdot A = A$, that is, $PA_iP^{-1}=A_i$ for all $i$ or, equivalently, $PA_i = A_iP$, i.e.\ $A_i$ and $P$ commutes. We are going to prove that there exists a basis $v_1, \ldots, v_m \in k^m$ of eigenvectors of $P$ with the same eigenvalue. In that case, $P$ would be a multiple of the identity so $P \in G^0$.

To do so, let $v_1 \in k^m-\left\{0\right\}$ be an eigenvector of $P$ with eigenvalue $\lambda \in k$. Since $A$ is irreducible, the matrices $A_i$ cannot have a common eigenvector so $A_jv_1 \not\in \langle v_1\rangle$ for some $1 \leq j \leq n$. Moreover, since $P$ and $A_j$ commutes, $A_jv_1$ is also an eigenvector of $P$ of eigenvalue $\lambda$ and we set $v_2 = A_jv_1$. By induction, suppose that we have build $v_1, \ldots, v_l$ a set of linearly independent eigenvectors of $P$ of eigenvalue $\lambda$ with $l<m$. Since $A$ is irreducible, $A_jv_i \not\in \langle v_1, \ldots, v_l \rangle$ for some $1 \leq i \leq l$ and $1 \leq j \leq n$. Hence, $A_jv_i$ is an eigenvector of $P$ of eigenvalue $\lambda$ linearly independent with $v_1, \ldots, v_m$, so we can set $v_{l+1} = A_jv_i$.
\end{proof}
\end{prop}

\begin{rmk}\label{rmk:triang-irred}
\begin{itemize}
	\item If $G \subseteq \GL{2}(k)$, then a representation $\rho: \Gamma \to G$ is reducible if and only if all the elements of $\rho(\Gamma)$ have a common eigenvector.
	\item By \cite{Borel-1991} Proposition 1.10, every affine group is linear so the previous results automatically holds for affine groups.
\end{itemize}
\end{rmk}

\section{Deligne-Hodge polynomials of $\SL{2}(\CC)$-character varieties}
\label{subsec:reducible-rep}

For the rest of the paper, we will work within the ground field $k = \CC$ and we will focus on the algebraic group $G = \SL{2}(\CC)$. For short, we will just denote by $\SL{2} = \SL{2}(\CC)$ and $\PGL{2}=\PGL{2}(\CC)$.

The center of $\SL{2}$ are the matrices $\pm  \Id$ so the action of $\SL{2}$ by conjugation descends to an action of $\SL{2}/\left\{\pm  \Id\right\} = \PGL{2}$. Hence, given a finitely generated group $\Gamma$, we have a decomposition $\mathfrak{X}_{\SL{2}}(\Gamma) = \Xred{\mathfrak{X}_{\SL{2}}}(\Gamma) \sqcup \Xirred{\mathfrak{X}_{\SL{2}}}(\Gamma)$.
The first stratum of this decomposition is a closed set and the second stratum is an open set with a geometric free action of $\PGL{2}$. Observe that two representations of $\mathfrak{X}_{\SL{2}}(\Gamma)$ are conjugated by an element of $\GL{2}$ if and only if they are so by an element of $\SL{2}$, so the conjugacy classes agree. In particular, $A \in \Rep{\SL{2}}(\Gamma)$ is reducible if and only if it is $\SL{2}$-conjugate to a upper triangular element.

\begin{prop}\label{prop:isotropy-of-matrices}
Let $\Gamma$ be a finitely generated group. Suppose that $A \in \mathfrak{X}_{\SL{2}}(\Gamma)$ has non trivial isotropy for the action of $\PGL{2}$. Then, $A$ is conjugate to an element of the following classes:
$$
\left(\begin{pmatrix} \pm 1 & \alpha_1 \\ 0 & \pm 1\end{pmatrix}, \ldots, \begin{pmatrix}\pm 1 & \alpha_n \\ 0 & \pm 1\end{pmatrix}\right) \;\;\;\txt{or}\;\;\; \left(\begin{pmatrix}\lambda_1 & 0 \\ 0 & \lambda_1^{-1}\end{pmatrix}, \ldots, \begin{pmatrix}\lambda_n & 0 \\ 0 & \lambda_n^{-1}\end{pmatrix}\right),
$$
for $\alpha_i \in \CC$ and $\lambda_i \in \CC^* = \CC - \left\{0\right\}$.
\begin{proof}
	By Proposition \ref{prop:action-free-on-irred}, $A$ must be reducible so, after conjugation, we can suppose that
	$$
		A = (A_1, \ldots, A_n) = \left(\begin{pmatrix} \lambda_1 & \alpha_1 \\ 0 & \lambda_1^{-1}\end{pmatrix}, \ldots, \begin{pmatrix}\lambda_n & \alpha_n \\ 0 & \lambda_n^{-1}\end{pmatrix}\right).
	$$
If all $\lambda_i = \pm 1$, then $A$ belongs to the first class. Hence, we can suppose that $\lambda_j \neq \pm 1$ for some $1 \leq j \leq n$. In that case, conjugating with
$$
	P_0 = \begin{pmatrix}1 & \frac{\alpha_j}{\lambda_j - \lambda_j^{-1}} \\ 0 & 1\end{pmatrix}
$$ we obtain that $\alpha_j=0$ so $A_j$ is a diagonal matrix. We claim that, after such a conjugation, all the $\alpha_i$ must vanish so $A$ is of the second class. Otherwise, since $P_0$ preserves upper triangular matrices, there must be another upper triangular non-diagonal matrix $A_k$. In that case, if $P \in \SL{2}$ fixes $A$, we obtain in particular that $P$ stabilices both a diagonal matrix and a non-diagonal upper triangular matrix. But this is only possible if $P = \pm  \Id$.
\end{proof}
\end{prop}

Given a finitely generated group $\Gamma$, let us denote by $\XP{\mathfrak{X}_{\SL{2}}}(\Gamma)$ the representations that have the form of the first class and not of the second class and let $\XPh{\mathfrak{X}_{\SL{2}}}(\Gamma)$ be the union of the orbits of elements of $\XP{\mathfrak{X}_{\SL{2}}}(\Gamma)$. Analogously, we will denote by $\XD{\mathfrak{X}_{\SL{2}}}(\Gamma)$ the representations that has the form of only the second class and by $\XDh{\mathfrak{X}_{\SL{2}}}(\Gamma)$ their orbits. Finally, $\XI{\mathfrak{X}_{\SL{2}}}(\Gamma)$ will be the representations that belong to both classes. Observe that the action of $\SL{2}$ on the later stratum is trivial. Setting $\XTilde{\mathfrak{X}_{\SL{2}}}(\Gamma) = \Xred{\mathfrak{X}_{\SL{2}}}(\Gamma) - {\XPh{\mathfrak{X}_{\SL{2}}}}(\Gamma) - {\XDh{\mathfrak{X}_{\SL{2}}}}(\Gamma) - \XI{\mathfrak{X}_{\SL{2}}}(\Gamma)$ we obtain a stratification
$$
	\mathfrak{X}_{\SL{2}}(\Gamma) = \Xirred{\mathfrak{X}_{\SL{2}}}(\Gamma) \sqcup {\XPh{\mathfrak{X}_{\SL{2}}}}(\Gamma) \sqcup {\XDh{\mathfrak{X}_{\SL{2}}}}(\Gamma) \sqcup \XI{\mathfrak{X}_{\SL{2}}}(\Gamma) \sqcup \XTilde{\mathfrak{X}_{\SL{2}}}(\Gamma).
$$
Observe that this stratification also decomposes $\mathfrak{X}_{\SL{2}}(\Gamma)$ in terms of the conjugacy classes of the stabilizers for the action of $\SL{2}$. To be precise, the stabilizer of the points of $\XPh{\mathfrak{X}_{\SL{2}}}(\Gamma)$ and $\XDh{\mathfrak{X}_{\SL{2}}}(\Gamma)$ are conjugated to the subgroups of $\SL{2}$
$$
		\Stab\,J_+ = \left\{\begin{pmatrix} \pm 1 & \beta \\ 0 & \pm 1\end{pmatrix},\, \beta \in \CC\right\} \hspace{0.8cm} \Stab\,D_\lambda = \left\{\begin{pmatrix} \mu & 0\\ 0 & \mu^{-1} \end{pmatrix},\, \mu \in \CC^*\right\},
$$
respectively, where $D_\lambda$ is the diagonal matrix with eigenvalues $\lambda^{\pm 1}$. For the stratum $\XI{\mathfrak{X}_{\SL{2}}}(\Gamma)$, the action is trivial and for the strata $\Xirred{\mathfrak{X}_{\SL{2}}}(\Gamma)$ and $\XTilde{\mathfrak{X}_{\SL{2}}}(\Gamma)$ the action of $\PGL{2}$ is free.

\begin{rmk}
The stratification considered in this paper is strongly related to the so-called Luna's stratification for ${\mathfrak{X}_{\SL{2}}}(\Gamma)$ (see  \cite{Luna:1972} and \cite{Luna:1973}, also Section 6.9 of \cite{Popov-Vinberg:1989}). In this case, the Luna strata are $\XPh{\mathfrak{X}_{\SL{2}}}(\Gamma)$, $\XDh{\mathfrak{X}_{\SL{2}}}(\Gamma)$, $\XI{\mathfrak{X}_{\SL{2}}}(\Gamma)$ and $\Xirred{\mathfrak{X}_{\SL{2}}}(\Gamma) \cup \XTilde{\mathfrak{X}_{\SL{2}}}(\Gamma)$.
\end{rmk}

\subsection{Free groups}

Let us fix $\Gamma = F_n$, the free group of $n$ generators, and recall that we abbreviate $\Xf{n} = \mathfrak{X}_{\SL{2}}(F_n)$. In this case, the stratification considered above has the following properties:

\begin{itemize}
	\item $\XPh{\Xf{n}}$. Given $A = (A_1, \ldots, A_n) \in \XPh{\Xf{n}}$, let
$$
    \left(\begin{pmatrix} \epsilon_1 & \alpha_1 \\ 0 & \epsilon_1\end{pmatrix}, \ldots, \begin{pmatrix}\epsilon_n & \alpha_n \\ 0 & \epsilon_n\end{pmatrix}\right)
$$ be the element of $\XP{\Xf{n}}$ conjugate to $A$ with $\epsilon_i = \pm 1$ and $\alpha_i \in \CC$ not all zero. Observe that such an element of $\XP{\Xf{n}}$ is unique up to simultaneous rescalling of the upper triangular components $\alpha_i$. Thus, the $\SL{2}$-orbit of $A$, $[A]$, is the set of reducible representations $(B_1, \ldots, B_n) \in \Xf{n}$ with a double eigenvalue such that, in their upper triangular form,
$$
    \left(\begin{pmatrix} \epsilon_1 & \beta_1 \\ 0 & \epsilon_1\end{pmatrix}, \ldots, \begin{pmatrix} \epsilon_n & \beta_n \\ 0 & \epsilon_n\end{pmatrix}\right)
$$
there exists $\lambda \neq 0$ such that $(\alpha_1, \ldots, \alpha_n) = \lambda(\beta_1, \ldots, \beta_n)$. Then, taking $\lambda \to 0$, we find that the closure of the orbit, $\overline{[A]}$, is precisely set of reducible representations with double eigenvalue such that their upper triangular components satisfy $(\alpha_1, \ldots, \alpha_n) = \lambda(\beta_1, \ldots, \beta_n)$ for some $\lambda \in \CC$. In particular, for $\lambda = 0$ we have that $(\epsilon_1 \Id, \ldots, \epsilon_n \Id) \in \overline{[A]}$.

	The Deligne-Hodge polynomial of this stratum can be computed observing that the tuple $(\alpha_1, \ldots, \alpha_n) \in \CC^n-\left\{0\right\}$ determines the diagonal form up to projectivization. Hence, we obtain a regular fibration
	$$
		\CC^* \longrightarrow \SL{2}/\Stab\,J_+ \times \left\{\pm 1\right\}^n \times \left(\CC^n-\left\{0\right\}\right) \longrightarrow \XPh{\Xf{n}}.
	$$
Observe that this fibration is locally trivial in the Zariski topology so, by \cite{LMN}, this fibration has trivial monodromy and, thus
$$
    \DelHod{\XPh{\Xf{n}}} = \DelHod{\left\{\pm 1\right\}^n} \DelHod{\CC\PP^{n-1}} \DelHod{\SL{2}/\Stab\,J_+} = 2^n(q^2-1)\frac{q^{n} -1}{q-1}.
$$
	\item $\XDh{\Xf{n}}$. Given $A = (A_1, \ldots, A_n) \in \XDh{\Xf{n}}$, let
$$
    \left(\begin{pmatrix} \lambda_1 & 0 \\ 0 & \lambda_1^{-1}\end{pmatrix}, \ldots, \begin{pmatrix}\lambda_n & 0 \\ 0 & \lambda_n^{-1}\end{pmatrix}\right)
$$ be the element of $\XD{\Xf{n}}$ conjugate to $A$ with $\lambda_i \in \CC^*$ and not all equal to $\pm 1$.
	In order to compute the Deligne-Hodge polynomial of this stratum, recall that the diagonal form of an element of $\XDh{\Xf{n}}$ is determined up to permutation of the columns, so we have a double covering
	$$
		\SL{2}/\Stab\,D_\lambda \times \left((\CC^*)^n-\left\{(\pm 1, \ldots, \pm 1)\right\}\right) \longrightarrow \XDh{\Xf{n}}.
	$$
Therefore, we obtain that
$$
	\XDh{\Xf{n}} = \frac{\SL{2}/\Stab\,D_\lambda \times \big((\CC^*)^n-\left\{(\pm 1, \ldots, \pm 1)\right\}\big)}{\ZZ_2}.
$$
Using Example \ref{ex:equiv-product}, its Deligne-Hodge polynomial is
$$
	\DelHod{\XDh{\Xf{n}}} = \frac{q^3-q}{2} \left((q-1)^{n-1} + (q+1)^{n-1}\right) - 2^nq^2.
$$

	\item $\XI{\Xf{n}}$. Here, the action of $\SL{2}$ is trivial so, in particular, the action is closed.
	This stratum consists of $2^n$ points, so $\DelHod{\XI{\Xf{n}}} = 2^n$.
	\item $\XTilde{\Xf{n}}$. This is the set of reducible representations, not completely reducible, such that some matrix has no double eigenvalue. Given $A \in \XTilde{\Xf{n}}$,  $A$ is conjugate to an element of the form
	$$
		\left(\begin{pmatrix} \lambda_1 & \alpha_1 \\ 0 & \lambda_1^{-1}\end{pmatrix}, \ldots, \begin{pmatrix}\lambda_n & \alpha_n \\ 0 & \lambda_n^{-1}\end{pmatrix}\right),
	$$
where the vector $(\alpha_1, \ldots, \alpha_n)$ is determined up to the action of $U = \CC^* \times \CC$ by
$$
	(\mu, a)\cdot (\alpha_1, \ldots, \alpha_n) = \left(\alpha_1\mu^2 - a\mu\left(\lambda_1-\lambda_1^{-1}\right), \ldots, \alpha_n\mu^2 - a\mu\left(\lambda_n-\lambda_n^{-1}\right)\right)
$$
Observe that, since $A \not\in \Xf{n}^D$, we must have $(\mu, a)\cdot(\alpha_1, \ldots, \alpha_n) \neq (0, \ldots, 0)$ for all $(\mu,a) \in U$. This vanishing can happen if and only if $(\alpha_1, \ldots, \alpha_n) = \lambda(\lambda_1 - \lambda_1^{-1}, \ldots, \lambda_n - \lambda_n^{-1})$ for some $\lambda \in \CC$. Therefore, the allowed values for $(\alpha_1, \ldots, \alpha_n)$ lie in $\CC^n - l$, where $l \subseteq \CC^n$ is the line spanned by $(\lambda_1 - \lambda_1^{-1}, \ldots, \lambda_n - \lambda_n^{-1})$. Thus, we have regular fibration
$$
	U \longrightarrow \PGL{2} \times \Omega
 \stackrel{\pi}{\longrightarrow} \XTilde{\Xf{n}},
$$ 
where $\Omega \cong \big((\CC^*)^n - \left\{\pm 1, \ldots, \pm 1\right\}\big) \times \big(\CC^n - l\big)$ is the set of allowed values for the eigenvalues and the antidiagonal component.
Hence, the Deligne-Hodge polynomial of this stratum is
$$
	\DelHod{\XTilde{\Xf{n}}} = \frac{q^3-q}{(q-1)q} \left((q-1)^n - 2^n\right)\left(q^n-q\right).
$$
\end{itemize}

From this analysis, it is possible to describe the GIT quotient $\Xf{n} \sslash \SL{2}$ quite explicitly. First of all, recall that we have a decomposition $\Xf{n} = \Xred{\Xf{n}} \sqcup \Xirred{\Xf{n}}$ with $\Xirred{\Xf{n}}$ an open subvariety. Since $\SL{2}$ is a reductive affine group, by the results of Section \ref{section:stratification}, we have
$$
	 \DelHod{\Xf{n} \sslash \SL{2}} = \DelHod{\Xred{\Xf{n}} \sslash \SL{2}} + \DelHod{\Xirred{\Xf{n}} \sslash \SL{2}}.
$$
The action of $\PGL{2}$ on $\Xirred{\Xf{n}}$ is closed and free, so the Deligne-Hodge polynomial of the quotient can be computed by Corollary \ref{cor:luna-thm-epol}, obtaining
	$$
		\DelHod{\Xirred{\Xf{n}} \sslash \SL{2}} = \frac{\DelHod{\Xirred{\Xf{n}}}}{\DelHod{\PGL{2}}} = \frac{\DelHod{\Xirred{\Xf{n}}}}{q^3-q}.
	$$
The computation of $\DelHod{\Xirred{\Xf{n}}}$ can be done using the analysis above since $\Xirred{\Xf{n}} = \Xf{n} - \XPh{\Xf{n}} - \XDh{\Xf{n}} - \XI{\Xf{n}} - \XTilde{\Xf{n}}$. Using that $\DelHod{\Xf{n}} = \DelHod{\SL{2}}^n = (q^3-q)^n$, we find that
\begin{align*}
	\DelHod{\Xirred{\Xf{n}}} =&\, 2^{n} q^{2} - \frac{1}{2} \, {\left(q^{3} - q\right)} {\left({\left(q + 1\right)}^{n - 1} + {\left(q - 1\right)}^{n - 1}\right)} \\
	&- {\left(2^{n} q + {\left(q - 1\right)}^{n} q^{n} - {\left(q - 1\right)}^{n} q - 2^{n}\right)} {\left(q + 1\right)} - 2^{n} + {\left(q^{3} - q\right)}^{n}.
\end{align*}

On the other hand, for $\Xred{\Xf{n}}$ the situation is more involved since the action is very far from being closed. Actually, we are going to show that $(\XD{\Xf{n}}, \ZZ_2)$ is a core for the action of $\SL{2}$ on $\Xred{\Xf{n}}$. For that, we need a preliminary result about the behaviour of the orbits.

\begin{prop}\label{prop:intersection-closure-core}
Let $W \subseteq \XD{\Xf{n}}(\SL{2})$ be a closed set. For any $A \in \overline{[W]}$ we have that $\overline{[A]} \cap W \neq \emptyset$.
\begin{proof}

Since $W$ is closed, it is the zero set of some regular functions $f_1, \ldots, f_r: \XD{\Xf{n}} \cong (\CC^*)^n \to \CC$. Since $\cO_{(\CC^*)^n}((\CC^*)^n) = \CC[\lambda_1, \lambda_1^{-1}, \ldots, \lambda_n, \lambda_n^{-1}]$, after multiplying by some polynomials we can suppose that $f_i \in \CC[\lambda_1, \ldots, \lambda_n]$. Now, let us consider the algebraic variety
$$
    \Omega = \left\{(A_i, \lambda_i, v) \in \Xf{n} \times (\CC^*)^n \times \left(\CC^2-\left\{0\right\}\right)\,\left|\\\;\begin{matrix}A_iv = \lambda_i v\\f_i(\lambda_1, \ldots, \lambda_n) = 0\end{matrix}\right.\right\} \subseteq \Xf{n} \times \CC^{n + 2}.
$$
Now, let $p: \Xf{n} \times \CC^{n + 2} \to \Xf{n}$ be the projection onto the first factor, and let $W' = p(W)$. Observe that, if $A = (A_1, \ldots, A_n) \in W'$, then the matrices $A_i$ have a common eigenvector so $A$ is reducible. Hence,
$$
	QAQ^{-1} = \left(\begin{pmatrix}
				\lambda_1 & a_1 \\
				0 & \lambda_1^{-1}
				\end{pmatrix},
				\ldots,
				\begin{pmatrix}
				\lambda_n & a_n \\
				0 & \lambda_n^{-1}
				\end{pmatrix}
				\right)
$$
for some $Q \in \SL{2}$ and satisfying $f_i(\lambda_1, \ldots, \lambda_n)=0$. Therefore, taking $P_m =  \begin{pmatrix}
	m^{-1} & 0\\
	0 & m\\
\end{pmatrix}$ we have
$$
	(P_mQ) \cdot A \to \left(\begin{pmatrix}
	\lambda_1 & 0\\
	0 & \lambda_1^{-1}\\
\end{pmatrix}, \ldots, \begin{pmatrix}
	\lambda_n & 0\\
	0 & \lambda_n^{-1}\\
\end{pmatrix}
\right) \in \XD{\Xf{n}},
$$
for $m \to \infty$. That diagonal element is in the analytic closure of $[A]$ and, thus, also in the Zariski one. Moreover, $f_i(\lambda_1, \ldots, \lambda_n)=0$, so this element belongs to $W$. Therefore, $\overline{[A]} \cap W \neq \emptyset$. Hence, in order to finish the proof, it is enough to show that $\overline{[W]} \subseteq W'$. Trivially, $[W] \subseteq W'$, so it is enough to prove that $W'$ is closed.

To do so, let us consider the projectivization of $\Omega$ to the closed projective set
$$
    \tilde{\Omega} = \left\{\left(A_i, [\mu_0:\mu_1: \ldots:\mu_n], \bar{v}\right) \in \Xf{n} \times \PP^n \times \PP^1\,\left|\\\;\begin{matrix} A_i\mu_0v = \mu_i v\\ \tilde{f}_i(\mu_0, \ldots, \mu_n) = 0\end{matrix}\right.\right\} \subseteq \Xf{n} \times \PP^{n} \times \PP^1,
$$
where $\tilde{f}_i$ denotes the homogenization of the polynomial $f_i$. Any element $\left(A_i, [\mu_0:\mu_1: \ldots:\mu_n], \bar{v}\right) \in \tilde{\Omega}$ must have $\mu_0 \neq 0$ so there are no points at infinity. Hence, we can also write $W' = \rho(\tilde{\Omega})$ where  $\rho: \Xf{n} \times \PP^{n} \times \PP^1 \to \Xf{n}$ is the first projection. Recall that projective spaces are universally closed, so $\rho$ is closed and, thus, $W'=\rho(\tilde{\Omega})$ is closed, as we wanted.
\end{proof}
\end{prop}

\begin{cor}\label{cor:core-diagonal}
Consider the action of $\ZZ_2$ on $\XD{\Xf{n}}$ by permutating the columns. Then $(\XD{\Xf{n}}, \ZZ_2)$ is a core for the action of $\SL{2}$ on $\Xred{\Xf{n}}$. 
\begin{proof}
First, observe that the orbit of an element of $\XD{\Xf{n}}$ by the action of $\ZZ_2$ is finite, so the action is automatically closed and, thus, part \ref{prop:core:enum:0} of Proposition \ref{prop:core} holds. For part \ref{prop:core:enum:1}, observe that, for any $A \in \Xred{\Xf{n}}$, taking $Q$ and $P_m$ as in the proof of Proposition \ref{prop:intersection-closure-core}, we see that $\overline{[A]} \cap \XD{\Xf{n}} \neq \emptyset$. This intersection consists of, at most, two points uniquely determined by the eigenvalues of $A$ so, in particular, they are $\ZZ_2$-equivalent.

In order to prove condition \ref{prop:core:enum:2}, let $W_1, W_2 \subseteq \XD{\Xf{n}}$ be two $\ZZ_2$-invariant disjoint closed subsets and suppose that $A \in \overline{[W_1]} \cap \overline{[W_2]}$. By Proposition \ref{prop:intersection-closure-core}, we have $\overline{[A]} \cap W_1 \neq \emptyset$ and $\overline{[A]} \cap W_2 \neq \emptyset$. However, this cannot happen since the elements of $\overline{[A]} \cap \XD{\Xf{n}}$ are $\ZZ_2$-equivalent and $W_1, W_2$ are $\ZZ_2$-invariant and disjoint. Thus, $\overline{[W_1]} \cap \overline{[W_2]} = \emptyset$, proving condition \ref{prop:core:enum:2}. This proves that $(\XD{\Xf{n}}, \ZZ_2)$ is a core. 
\end{proof}
\end{cor}

\begin{cor} The Deligne-Hodge polynomial of the GIT quotient of the reducible stratum is
$$
	\DelHod{\Xred{\Xf{n}} \sslash \SL{2}} = \frac{1}{2} \left((q-1)^{n} + (q+1)^{n}\right).
$$
\begin{proof}
Observe that $\ZZ_2$ is a reductive group so the GIT quotient of $\XD{\Xf{n}}$ by $\ZZ_2$ exists and is a pseudo-quotient. Hence, by Proposition \ref{prop:core} and the results of Section \ref{section:stratification}, we have that $\DelHod{\Xred{\Xf{n}} \sslash \SL{2}} = \DelHod{\XD{\Xf{n}} \sslash \ZZ_2}$. For computing this last polynomial, observe that $\XD{\Xf{n}} = (\CC^*)^n$ and use the calculations of \cite{MM} (see also Example \ref{ex:equiv-product}).
\end{proof}
\end{cor}

After this analysis, we have computed the Deligne-Hodge polynomials of both strata so, adding up the contributions we obtain that
\begin{empheq}{align*}
	\DelHod{\Xf{n} \sslash \SL{2}} = \frac{1}{2} \, {\left(q + 1\right)}^{n - 1} q + \frac{1}{2} \, {\left(q
- 1\right)}^{n - 1} q - {\left(q - 1\right)}^{n - 1} q^{n - 1} +
{\left(q^{3} - q\right)}^{n - 1}.
\end{empheq}

\begin{rmk}
This result agrees with the computations of \cite{Lawton-Munoz:2016} and \cite{Cavazos-Lawton:2014}. Actually, the argument is, somehow, parallel to the one in the former paper.
\end{rmk}

\begin{rmk}
From the analysis of this section, it is also possible to study the case in which $\Gamma = \ZZ^n$ is the free abelian group with $n$ generators and $G = \SL{m}$. Observe that, in this case, $\Rep{\SL{m}}(\ZZ^n)$ is the set of tuples of $n$ pairwise commutating matrices of $\SL{m}$. Since commutating matrices share a common eigenvector, all the representations of $\Rep{\SL{m}}(\ZZ^n)$ are reducible so $\Rep{\SL{m}}(\ZZ^n) = \Xred{\Rep{\SL{m}}}(\ZZ^n)$. Hence, analogously to Corollary \ref{cor:core-diagonal}, $(\XD{\Rep{\SL{m}}}(\ZZ^n), S_m)$ is a core for $\Xred{\Rep{\SL{m}}}(\ZZ^n)$, where $S_m$ acts on $\XD{\Rep{\SL{m}}}(\ZZ^n) = (\CC^*)^{n(m-1)}$ by permutation of the eigenvalues. Hence, we obtain that
$$
	\Rep{\SL{m}}(\ZZ^n) \sslash \SL{m} = (\CC^*)^{n(m-1)} \sslash S_m.
$$
Analogously, for $G = \GL{m}(\CC)$, we obtain that $\Rep{\GL{m}(\CC)}(\ZZ^n) \sslash \GL{m}(\CC) = (\CC^*)^{nm} \sslash S_m$.
This reproves Theorem 5.1 of \cite{Florentino-Silva:2017}.
In the case $m=2$, the Deligne-Hodge polynomial of these character varietes can be computed by means of Example \ref{ex:equiv-product}, obtaining that
$$
	\DelHod{\Rep{\SL{2}}(\ZZ^n) \sslash \SL{2}} = \DelHod{(\CC^*)^n \sslash \ZZ_2} = \frac{1}{2} \left((q-1)^{n} + (q+1)^{n}\right),
$$
$$
	\DelHod{\Rep{\GL{2}}(\ZZ^n) \sslash \GL{2}} = \DelHod{(\CC^*)^{2n} \sslash \ZZ_2} = \frac{1}{2} \left((q-1)^{2n} + (q+1)^{2n}\right).
$$
In the higher rank case, we need to use stronger results about equivariant cohomology in order to compute the corresponding quotient by $S_m$. This is, precisely, the strategy accomplished in \cite{Florentino-Silva:2017}.
\end{rmk}

\subsection{Surface groups}
\label{sec:surface-groups-nonpar}

In this section, we will consider the case $\Gamma = \pi_1(\Sigma_g)$, where $\Sigma_g$ is the genus $g \geq 1$ compact surface. Recall that, for short, we will denote the associated $\SL{2}$-representation variety by $\Xs{g} = \mathfrak{X}_{\SL{2}}(\pi_1(\Sigma_g))$ and its GIT quotient $\cR_g = \Xs{g} \sslash \SL{2}$. We have that $\Xs{g} \subseteq \Xf{2g}$ as a closed subvariety.

In order to identify the elements of $\Xs{g}$, let us denote the set of upper triangular matrices of $\Xf{n}$ by $\XU{\Xf{n}} \cong (\CC^*)^n \times \CC^n$ and $\XU{\Xs{g}} = \XU{\Xf{2g}} \cap \Xs{g}$. Given $A \in \XU{\Xf{2g}}$, say
$$
	A= \left(
	\begin{pmatrix}
		\lambda_1 & \alpha_1 \\
		0 & \lambda_1^{-1}
	\end{pmatrix},
	\begin{pmatrix}
		\mu_1^{-1} & \beta_1 \\
		0 & \mu_1
	\end{pmatrix}, \ldots,
	\begin{pmatrix}
		\lambda_{g} & \alpha_{g} \\
		0 & \lambda_{g}
	\end{pmatrix},
	\begin{pmatrix}
		\mu_{g} & \beta_{g} \\
		0 & \mu_{g}
	\end{pmatrix}
	\right)
$$
with $\lambda_i, \mu_i \in \CC^*$ and $\alpha_i, \beta_i \in \CC$, a straighforward computation shows that $A \in \XU{\Xs{g}}$ if and only if
\begin{equation*}\label{eq:cond-upper}
	\sum_{i=1}^g \lambda_i\mu_i \big(\left(\mu_i-\mu_i^{-1}\right)\beta_i - \left(\lambda_i-\lambda_i^{-1}\right)\alpha_i\big)= 0.
\end{equation*}
In particular, this implies that, for the strata given by Proposition \ref{prop:isotropy-of-matrices}, we have $\XP{\Xs{g}} = \XP{\Xf{2g}}$, $\XD{\Xs{g}} = \XD{\Xf{2g}}$ and $\XI{\Xs{g}} = \XI{\Xf{2g}}$. Let $\pi \subseteq \CC^{2g}$ be the $(\alpha_i, \beta_i)$-plane defined by the previous equation for fixed $(\lambda_i, \mu_i)$.

\begin{itemize}
	\item For $\Xred{\Xs{g}}$, observe that, using the equality $\XD{\Xs{g}} = \XD{\Xf{2g}}$ and Corollary \ref{cor:core-diagonal}, we obtain that $(\XD{\Xs{g}}, \ZZ_2)$ is a core for the action. Therefore, since $\XD{\Xf{2g}} = (\CC^*)^n$, we have
$$
	\DelHod{\Xred{\Xs{g}} \sslash \SL{2}} = \DelHod{(\CC^*)^{2g} \sslash \ZZ_2} = \frac{1}{2} \left((q-1)^{2g} + (q+1)^{2g}\right).
$$
	\item For $\Xirred{\Xs{g}}$, observe that, since $\Xirred{\Xf{2g}}$ is an open set of $\Xf{2g}$ in which the action of $\PGL{2}$ is closed and free, and $\Xs{g} \subseteq \Xf{2g}$ is closed, then $\Xirred{\Xs{g}} = \Xirred{\Xf{2g}} \cap \Xs{g}$ is an open subset of $\Xs{g}$ with a closed and free action. Therefore, Corollary \ref{cor:luna-thm-epol} gives us that
	$$
		\DelHod{\Xirred{\Xs{g}} \sslash \SL{2}} = \frac{\DelHod{\Xirred{\Xs{g}}}}{\DelHod{\PGL{2}}} = \frac{\DelHod{\Xirred{\Xs{g}}}}{q^3-q}.
	$$
In order to complete the calculation, it is enough to compute $\DelHod{\Xirred{\Xs{g}}}$. For this purpose, we will use that $\Xirred{\Xs{g}} = \Xs{g} - \XPh{\Xs{g}} - \XDh{\Xs{g}} - \XI{\Xs{g}} - \XTilde{\Xs{g}}$ and count.

\begin{itemize}
	\item $\XPh{\Xs{g}}$. In this case, since $\XP{\Xs{g}} = \XP{\Xf{2g}}$ and $\Xs{g}$ is $\SL{2}$-invariant, we have that $\XPh{\Xs{g}} = \XPh{\Xf{2g}}$. Therefore,
$$
    \DelHod{\XPh{\Xs{g}}} = \DelHod{\XPh{\Xf{2g}}} = 2^{2g}(q^2-1)\frac{q^{2g} -1}{q-1}.
$$
	\item $\XDh{\Xs{g}}$. Again, $\XD{\Xs{g}} = \XD{\Xf{2g}}$ and, thus, $\XDh{\Xs{g}} = \XDh{\Xf{2g}}$. Therefore,
$$
	\DelHod{\XDh{\Xs{g}}} = \DelHod{\XDh{\Xf{2g}}} = \frac{q^3-q}{2} \left((q-1)^{2g-1} + (q+1)^{2g-1}\right) - 2^{2g}q^2.
$$

	\item $\XI{\Xs{g}}$. Again, $\XI{\Xs{g}} = \XI{\Xf{2g}}$, which are $2^{2g}$ matrices so $\DelHod{\XI{\Xs{g}}} = 2^{2g}$.
	\item $\XTilde{\Xs{g}}$. In this case, any element is conjugated to one of the form
$$
	\left(
	\begin{pmatrix}
		\lambda_1 & \alpha_1 \\
		0 & \lambda_1^{-1}
	\end{pmatrix},
	\begin{pmatrix}
		\mu_1^{-1} & \beta_1 \\
		0 & \mu_1
	\end{pmatrix}, \ldots,
	\begin{pmatrix}
		\lambda_{g} & \alpha_{g} \\
		0 & \lambda_{g}
	\end{pmatrix},
	\begin{pmatrix}
		\mu_{g} & \beta_{g} \\
		0 & \mu_{g}
	\end{pmatrix}
	\right),
$$
with $(\lambda_1, \mu_1, \ldots, \lambda_g, \mu_g) \in (\CC^*)^{2g} - \left\{(\pm 1, \ldots, \pm 1)\right\}$ and $(\alpha_1, \beta_1, \ldots, \alpha_{g}, \beta_g) \in \pi - l$ where $l$ is the line spanned by $(\lambda_1 - \lambda_1^{-1}, \mu_1 - \mu_1^{-1}, \ldots, \lambda_g - \lambda_g^{-1}, \mu_g - \mu_g^{-1})$. Thus, we have a fibration
$$
	U \longrightarrow \PGL{2} \times \Omega \longrightarrow \XTilde{\Xs{g}},
$$
where $U = \CC \times \CC^*$ and $\Omega$ is a fibration with trivial monodromy $\pi - l \to \Omega \to (\CC^*)^{2g} - \left\{(\pm 1, \ldots, \pm 1)\right\}$. Using that $\DelHod{\pi} = q^{2g-1}$ and $\DelHod{l} = q$, the Deligne-Hodge polynomial is
$$
	\DelHod{\XTilde{\Xs{g}}} = \frac{q^3-q}{(q-1)q} \left((q-1)^{2g} - 2^{2g}\right)\left(q^{2g-1}-q\right).
$$
\end{itemize}
Therefore, putting all together we have that
$$
	\DelHod{\Xred{\Xs{g}}} = (q+1) (q-1)^{2g}\left(q^{2g-1}-q\right) + \frac{q^3-q}{2} \left((q-1)^{2g-1} + (q+1)^{2g-1}\right) - 2^{2g}(q^2-1).
$$
\end{itemize}
From \cite{MM}, Proposition 11, we know that the Deligne-Hodge polynomial of the total space is
\begin{align*}
	\DelHod{\Xs{g}} =&\, 2^{2g - 1} {\left(q - 1\right)}^{2g - 1} {\left(q + 1\right)}
q^{2g - 1} + 2^{2g - 1} {\left(q + 1\right)}^{2g - 1}
{\left(q - 1\right)} q^{2g - 1} \\
&+ \frac{1}{2} \, {\left(q +
1\right)}^{2g - 1} {\left(q - 1\right)}^{2} q^{2g - 1} +
\frac{1}{2} \, {\left(q - 1\right)}^{2g - 1} {\left(q + 1\right)}
{\left(q - 3\right)} q^{2g - 1} \\&+ {{\left(q + q^{2 \,
g-1}\right)} {\left(q^2 - 1\right)}^{2 \,
g - 1}}.
\end{align*}
So we finally find that
\begin{empheq}{align*}
	\DelHod{\cR_g} =&\, \frac{1}{2} \, {\left({\left(2^{2g} + 2 \, {\left(q - 1\right)}^{2g - 2} + q - 1\right)} q^{2g - 2} + q^{2} + 2 \, {\left(q -1\right)}^{2g - 2} + q\right)} {\left(q + 1\right)}^{2g - 2} \\
	&+\frac{1}{2} \, {\left({\left(2^{2g} - 1\right)} {\left(q -
1\right)}^{2g - 2} - {\left(q - 1\right)}^{2g - 2} q - 2^{2g
+ 1}\right)} q^{2g - 2} + \frac{1}{2} \, {\left(q - 1\right)}^{2 \,
g - 1} q.
\end{empheq}

\begin{rmk}
This result agrees with the one given in \cite{MM}, Theorem 14, and in \cite{Baraglia-Hekmati:2016}, Theorem 1.3.
\end{rmk}

%%%%%%%%%%%%%%%%%%%%%%%%%%%%%%%%%%%%%%%%%%%%%%%%
%              PARABOLIC CLASS J_+
%%%%%%%%%%%%%%%%%%%%%%%%%%%%%%%%%%%%%%%%%%%%%%%%

\section{Parabolic $\SL{2}(\CC)$-character varieties}
\label{sec:parabolic-rep}

In this section, we will discuss the case of parabolic representation varieties. Recall from Section \ref{sec:parabolic-repr-var} that, given a finitely generated group $\Gamma$ and an algebraic group $G$, a \emph{parabolic structure} $Q$ is a finite set of pairs $(\gamma, \lambda)$ where $\gamma \in \Gamma$ and $\lambda \subseteq G$ is a locally closed subset which is closed under conjugation. The corresponding representation variety is
$$
	\Rep{G}(\Gamma, Q) = \left\{\rho \in \Rep{G}(\Gamma)\,|\, \rho(\gamma) \in \lambda\,\,\textrm{for all } (\gamma, \lambda) \in Q\right\}.
$$
As a subvariety of $\Rep{G}(\Gamma)$, all the previous results about stability translate to parabolic representation varieties. In this way, all the points of $\Rep{G}(\Gamma, Q)$ are semi-stable for this action, the stable points are $\Xirred{\Rep{G}}(\Gamma, Q) = \Xirred{\Rep{G}}(\Gamma) \cap\Rep{G}(\Gamma, Q)$, in which the action of $G/G^0$ is closed and free.

In this Section, we are going to focus in the following cases (cf.\ Example \ref{ex:par-struc-repr}):
\begin{itemize}
	\item Fix some elements $h_1, \ldots, h_s \in G$ and let $[h_i]$ be their conjugacy classes. Then, we take $\Gamma = F_{n+s}$ and $Q = \left\{(\gamma_{1}, [h_1]), \ldots, (\gamma_s, [h_s])\right\}$, where $\gamma_1, \ldots, \gamma_s \in F_{n+s}$ is an independent set. 
	\item Let $\Sigma = \Sigma_g - \left\{p_1, \ldots, p_s\right\}$ with $p_i \in \Sigma_g$ distinct points. As parabolic structure, we will take $Q = \left\{(\gamma_1, [h_1]), \ldots, (\gamma_s, [h_s])\right\}$ for some fixed elements $h_i \in G$, where $\gamma_i$ are the positive oriented loops around $p_i$. Observe that the epimorphism $F_{2g +s} \to \pi_1(\Sigma)$ gives an inclusion $\Rep{G}(\pi_1(\Sigma), Q) \subseteq \Rep{G}(F_{2g+s}, Q)$.
\end{itemize}

As in the previous section, we will focus on the case $k = \CC$ and $G = \SL{2}$. In this paper, we are going to study the parabolic character varieties whose parabolic data lie in the conjugacy classes of the matrices
$$
	J_+ = \begin{pmatrix}
	1 & 1\\
	0 & 1\\
\end{pmatrix} \hspace{1cm} J_- = \begin{pmatrix}
	-1 & 1\\
	0 & -1\\
\end{pmatrix}
\hspace{1cm} - \Id = \begin{pmatrix}
	-1 & 0\\
	0 & -1\\
\end{pmatrix}.
$$
As we will see, the most important case to be considered is to take the parabolic structure $Q_s^+ = \left\{(\gamma_1, [J_+]), \ldots, (\gamma_s, [J_+])\right\}$. The rest of cases can be easily obtained from this one, as shown in Section \ref{subsec:parabolic-data-gen}.

\subsection{Free groups}
For short, let us denote $\Xfp{n}{s} = \mathfrak{X}_{\SL{2}}(F_{n+s}, Q_s^+)$. As before, we have a stratification $\Xfp{n}{s} = \Xirred{\Xfp{n}{s}} \sqcup \Xred{\Xfp{n}{s}}$. Hence, by the results of Section \ref{section:stratification}, we have that
$$
	\DelHod{\Xfp{n}{s} \sslash \SL{2}} = \DelHod{\Xirred{\Xfp{n}{s}} \sslash \SL{2}} + \DelHod{\Xred{\Xfp{n}{s}} \sslash \SL{2}}.
$$
Let us analyze each stratum separately.

\begin{itemize}
	\item $\Xirred{\Xfp{n}{s}}$ is again an open subvariety in which the the action of $\PGL{2}$ is closed and free, so Corollary \ref{cor:luna-thm-epol} gives us that
	$$
		\DelHod{\Xirred{\Xfp{n}{s}} \sslash \SL{2}} = \frac{\DelHod{\Xirred{\Xfp{n}{s}}}}{\DelHod{\PGL{2}}} = \frac{\DelHod{\Xirred{\Xfp{n}{s}}}}{q^3-q}.
	$$
In order to compute $\DelHod{\Xirred{\Xfp{n}{s}}}$ we count:

\begin{itemize}
	\item $\XPh{\Xfp{n}{s}}$. In this case, $\XP{\Xfp{n}{s}} = \left\{\pm 1\right\}^n \times \left(\CC^{n} \times (\CC^*)^s\right)$ and the action restricts to projectivizing the second factor. Therefore,
$$
    \DelHod{\XPh{\Xfp{n}{s}}} = 2^{n}(q^2-1)\frac{q^{n}(q-1)^s}{q-1}.
$$
	
	\item $\XDh{\Xfp{n}{s}}$. In this case, since $J_+$ is not diagonalizable, $\XD{\Xfp{n}{s}} = \emptyset$ so $\XDh{\Xfp{n}{s}} = \emptyset$ and it makes no contribution.

	\item $\XI{\Xfp{n}{s}}$. Again, it is empty so it makes no contribution.
	\item $\XTilde{\Xfp{n}{s}}$. In this case, any element is conjugated to one of the form
$$
	\left(
	\begin{pmatrix}
		\lambda_1 & \alpha_1 \\
		0 & \lambda_1^{-1}
	\end{pmatrix}, \ldots,
	\begin{pmatrix}
		\lambda_n & \alpha_n \\
		0 & \lambda_n^{-1}
	\end{pmatrix},
	\begin{pmatrix}
		1 & c_{1} \\
		0 & 1
	\end{pmatrix}, \ldots,
	\begin{pmatrix}
		1 & c_{s} \\
		0 & 1
	\end{pmatrix}
	\right)
$$
with $(\lambda_1, \ldots, \lambda_n) \in (\CC^*)^n - \left\{(\pm 1, \ldots, \pm 1)\right\}$ and $(\alpha_1, \ldots, \alpha_{n}, c_1, \ldots, c_s) \in \CC^{n} \times (\CC^*)^s$. Thus, we have a fibration
$$
	U \longrightarrow \PGL{2} \times \Omega \longrightarrow \XTilde{\Xfp{n}{s}},
$$
where $\Omega = \big((\CC^*)^n - \left\{(\pm 1, \ldots, \pm 1)\right\}\big) \times \big(\CC^{n} \times (\CC^*)^s\big)$. Observe that we do not need to remove any antidiagonal value since the intersection of the line spanned by $(\lambda_1-\lambda_1^{-1}, \ldots, \lambda_n - \lambda_n^{-1}, 0, \ldots, 0)$ with $\CC^{n} \times (\CC^*)^s$ is empty.
Hence, the Deligne-Hodge polynomial is
$$
	\DelHod{\XTilde{\Xfp{n}{s}}} = \frac{q^3-q}{(q-1)q} \left((q-1)^n - 2^n\right)q^{n}(q-1)^s.
$$
\end{itemize}
Therefore, putting all together we have that
$$
	\DelHod{\Xred{\Xfp{n}{s}}} = {\left(q - 1\right)}^{n} {\left(q - 1\right)}^{s} {\left(q + 1\right)}
q^{n}.
$$
Using that $\Xfp{n}{s} = \SL{2}^n \times [J_+]^s$ and $\DelHod{[J_+]} = (q^2-1)$, we obtain that $\DelHod{\Xfp{n}{s}} = (q^3-q)^n(q^2-1)^s$, so we finally find
$$
	\DelHod{\Xirred{\Xfp{n}{s}}} =  (q-1)^n(q-1)^s\left((q^2+q)^n(q+1)^s -(q+1)q^n\right).
$$
	\item For $\Xred{\Xfp{n}{s}}$ the situation becomes completely different from the previous ones. As we have shown, $\XD{\Xfp{n}{s}} = \emptyset$ so it can be no longer a core for the action. The key point now is that, precisely for this reason, the action of $\SL{2}$ on $\Xred{\Xfp{n}{s}}$ is closed. To be precise, let $A \in \Xred{\Xfp{n}{s}} \subseteq \Xred{\Xf{n + s}}$ and let $\overline{[A]}$ be the closure of its orbit in $\Xf{n + s}$. The difference $\overline{[A]} - [A]$ lies in $\XD{\Xf{n + s}}$ so, since $\XD{\Xf{n + s}} \cap \Xfp{n}{s} = \XD{\Xfp{n}{s}} = \emptyset$, the orbit of $A$ is closed in $\Xred{\Xfp{n}{s}}$. However, the action of $\PGL{2}$ on $\Xred{\Xfp{n}{s}}$ is not free everywhere so we have to distinguish between two strata:
	\begin{itemize}
		\item $\XTilde{\Xfp{n}{s}}$. Here, the action of $\PGL{2}$ is free so, by Corollary \ref{cor:luna-thm-epol}, we have
		$$
			\DelHod{\XTilde{\Xfp{n}{s}} \sslash \SL{2}} = \frac{\DelHod{\XTilde{\Xfp{n}{s}}}}{q^3-q} = \left((q-1)^n - 2^n\right)q^{n-1}(q-1)^{s-1}.
		$$
		\item $\XPh{\Xfp{n}{s}}$. Here, the action of $\PGL{2}$ is not free, but it has stabilizer isomorphic to $\Stab\,J_+ \cong \CC$. The fact that $[J_+] \cong \SL{2}/\Stab\,J_+$ implies that the GIT quotient
		$
			\XPh{\Xfp{n}{s}} \to \XPh{\Xfp{n}{s}} \sslash \SL{2}
		$
is a locally trivial fibration with fiber $\SL{2}/\Stab\,J_+$ and trivial monodromy. Hence, we have that
		$$
			\DelHod{\XPh{\Xfp{n}{s}} \sslash \SL{2}} = \frac{\DelHod{\XPh{\Xfp{n}{s}}}}{q^2-1} = 2^{n}q^{n}(q-1)^{s-1}.
		$$
	Observe that we have used that $\DelHod{\SL{2}/\Stab\,J_+} = \DelHod{\SL{2}}/\DelHod{\Stab\,J_+} = q^2-1$.
	\end{itemize}
We have a stratification $\Xred{\Xfp{n}{s}} = \XPh{\Xfp{n}{s}} \sqcup \XTilde{\Xfp{n}{s}}$, where $\XTilde{\Xfp{n}{s}}$ is open orbitwise-closed. Thus, we have
$$
    \DelHod{\Xred{\Xfp{n}{s}} \sslash \SL{2}} = \DelHod{\XPh{\Xfp{n}{s}} \sslash \SL{2}} + \DelHod{\XTilde{\Xfp{n}{s}} \sslash \SL{2}} =  
{\left(2^{n} + {\left(q - 1\right)}^{n - 1}\right)} {\left(q -
1\right)}^{s} q^{n - 1}.
$$
\end{itemize}
Summarizing, the analysis above shows that
\begin{empheq}{align*}
	\DelHod{\Xfp{n}{s} \sslash \SL{2}} = 2^{n} {\left(q - 1\right)}^{s} q^{n - 1} + {\left(q^{3} - q\right)}^{n -
1} {\left(q^{2} - 1\right)}^{s}.
\end{empheq}

\subsection{Surface groups}
\label{sec:surface-groups-par}

Let $\Sigma$ be the genus $g$ compact surface with $s$ punctures and let us denote $\Xsp{g}{s} = \mathfrak{X}_{\SL{2}}(\pi_1(\Sigma), Q_s^+)$. As we mentioned in Section \ref{section:stratification}, the decomposition into reducible representations and irreducible ones gives an equality
$$
	\DelHod{\Xsp{g}{s} \sslash \SL{2}} = \DelHod{\Xred{\Xsp{g}{s}} \sslash \SL{2}} + \DelHod{\Xirred{\Xsp{g}{s}} \sslash \SL{2}}.
$$
In order to understand this stratification, observe that, as closed subvarieties of the one for the free case, we have $\XD{\Xsp{g}{s}} = \XDh{\Xsp{g}{s}} = \XI{\Xsp{g}{s}} = \emptyset$. Let us define $\XU{\Xfp{n}{s}} = \XU{\Xf{n+s}} \cap \Xfp{n}{s}$ and $\XU{\Xsp{g}{s}} = \XU{\Xf{2g+s}} \cap \Xsp{g}{s}$ and take $A \in \XU{\Xfp{2g}{s}}$, say
$$
	A = \left(
	\begin{pmatrix}
		\lambda_1 & \alpha_1 \\
		0 & \lambda_1^{-1}
	\end{pmatrix},
	\begin{pmatrix}
		\mu_1^{-1} & \beta_1 \\
		0 & \mu_1
	\end{pmatrix}, \ldots,
	\begin{pmatrix}
		\lambda_{g} & \alpha_{g} \\
		0 & \lambda_{g}
	\end{pmatrix},
	\begin{pmatrix}
		\mu_{g} & \beta_{g} \\
		0 & \mu_{g}
	\end{pmatrix},
	\begin{pmatrix}
		1 & c_1 \\
		0 & 1
	\end{pmatrix}, \ldots,
	\begin{pmatrix}
		1 & c_s \\
		0 & 1
	\end{pmatrix}
	\right),
$$
with $\lambda_i, \mu_i \in \CC^*$, $\alpha_i, \beta_i \in \CC$ and $c_i \in \CC^*$. Then, we have that $A \in \XU{\Xsp{g}{s}}$ if and only if
\begin{equation*}\label{eq:cond-upper:parabolic-jplus}
	\sum_{i=1}^g \lambda_i\mu_i \big(\left(\mu_i-\mu_i^{-1}\right)\beta_i - \left(\lambda_i-\lambda_i^{-1}\right)\alpha_i\big) + \sum_{i=1}^s c_i= 0.
\end{equation*}
In particular, $\XD{\Xsp{g}{s}} = \XD{\Xfp{2g}{s}}$. Let us analyze each stratum separately.

\begin{itemize}
	\item $\Xirred{\Xsp{g}{s}}$ is again an open subvariety in which the the action of $\PGL{2}$ is closed and free so Corollary \ref{cor:luna-thm-epol} gives us that
	$$
		\DelHod{\Xirred{\Xsp{g}{s}} \sslash \SL{2}} = \frac{\DelHod{\Xirred{\Xsp{g}{s}}}}{\DelHod{\PGL{2}}} = \frac{\DelHod{\Xsp{g}{s}} - \DelHod{\Xred{\Xsp{g}{s}}}}{q^3-q}.
	$$
The calculation of each of the strata of $\Xred{\Xsp{g}{s}}$ mimics the corresponding one for $\Xred{\Xfp{2g}{s}}$ but taking care of the previous equation.

\begin{itemize}
	\item $\XPh{\Xsp{g}{s}}$. In this case $\XP{\Xsp{g}{s}} = \left\{\pm 1\right\}^{2g} \times \big(\CC^{2g} \times \pi_s\big)$, where
$$
    \pi_s = \left\{\left.\sum_{j = 1}^s c_j = 0\;\;\right|\, c_j \neq 0\right\}.
$$
In order to compute the Deligne-Hodge polynomial of this space, just observe that $\pi_s = (\CC^*)^{s-1} - \pi_{s-1}$. Therefore, using the base case $\pi_1 = \emptyset$, we have
\begin{align*}
    \DelHod{\pi_s} &= (q-1)^{s-1} - \DelHod{\pi_{s-1}} = \sum_{k=1}^{s-1} (-1)^{k+1}(q-1)^{s-k} \\&= (-1)^s \left(\frac{(1-q)^s - 1}{q} + 1\right).
\end{align*}
Observe that, as in the free case, the action of $\SL{2}$ on $\XP{\Xsp{g}{s}}$ is just the projectivization on $\CC^{n} \times \pi_s$ so we obtain a fibration with trivial monodromy
$$
    \CC^* \longrightarrow \frac{\SL{2}}{\Stab\,J_+} \times \left\{\pm 1\right\}^{2g} \times \big(\CC^{2g} \times \pi_s\big) \longrightarrow \XPh{\Xsp{g}{s}}.
$$
Therefore, we have
$$
    \DelHod{\XPh{\Xsp{g}{s}}} = 2^{2g}(q^2-1)\frac{q^{2g}}{q-1}\left((-1)^s \left(\frac{(1-q)^s - 1}{q} + 1\right)\right).
$$
	
	\item $\XDh{\Xsp{g}{s}}$. This case makes no contribution since $\XD{\Xfp{n}{s}} = \emptyset$.

	\item $\XI{\Xsp{g}{s}}$. Again, it is empty so it makes no contribution.
	\item $\XTilde{\Xsp{g}{s}}$. In this case, any element is conjugated to one of the form
$$
	\left(
	\begin{pmatrix}
		\lambda_1 & \alpha_1 \\
		0 & \lambda_1^{-1}
	\end{pmatrix}, \ldots,
	\begin{pmatrix}
		\mu_g & \beta_n \\
		0 & \mu_g^{-1}
	\end{pmatrix},
	\begin{pmatrix}
		1 & c_{1} \\
		0 & 1
	\end{pmatrix}, \ldots,
	\begin{pmatrix}
		1 & c_{s} \\
		0 & 1
	\end{pmatrix}
	\right)
$$
with $(\lambda_1, \ldots, \mu_g) \in (\CC^*)^{2g} - \left\{(\pm 1, \ldots, \pm 1)\right\}$ and $(\alpha_1, \ldots, \beta_{g}, c_1, \ldots, c_s) \in  \Pi_s$, where we denote
$$
    \Pi_s = \left\{	\sum_{i=1}^g \lambda_i\mu_i \big(\left(\mu_i-\mu_i^{-1}\right)\beta_i - \left(\lambda_i-\lambda_i^{-1}\right)\alpha_i\big) + \sum_{i=1}^s c_i= 0
\right\},
$$
for fixed $(\lambda_i, \mu_i)$.
In order to compute the Deligne-Hodge polynomial of $\Pi_s$ observe that $\Pi_s = \CC^{2g} \times (\CC^*)^{s-1} - \Pi_{s-1}$. Using as base case that $\Pi_1$ is $\CC^{2g}$ minus a hyperplane, we have
\begin{align*}
    \DelHod{\Pi_s} &= q^{2g}(q-1)^{s-1} - \DelHod{\Pi_{s-1}} \\
    &= q^{2g}\sum_{k=1}^{s} (-1)^{k+1}(q-1)^{s-k} + (-1)^s q^{2g-1} = q^{2g-1}(q-1)^s.
\end{align*}
Thus, we have a fibration
$$
	U \longrightarrow \PGL{2} \times \Omega \longrightarrow \XTilde{\Xsp{g}{s}},
$$
where $\Pi_s \to \Omega \to (\CC^*)^{2g} - \left\{(\pm 1, \ldots, \pm 1)\right\}$ is a fibration with trivial monodromy. Therefore, the Deligne-Hodge polynomial is
$$
	\DelHod{\XTilde{\Xsp{g}{s}}} = \frac{q^3-q}{(q-1)q} \left((q-1)^{2g} - 2^{2g}\right)q^{2g-1}(q-1)^s.
$$
\end{itemize}
Therefore, putting all together, we obtain that
\begin{align*}
	\DelHod{\Xred{\Xsp{g}{s}}} &=          
2^{2g} \left(-1\right)^{s} {\left(q + 1\right)} q^{2g}
{\left(\frac{{\left(-q + 1\right)}^{s} - 1}{q} + 1\right)} \\&\;\;\;\;- {\left(2^{2
\, g} - {\left(q - 1\right)}^{2g}\right)} {\left(q - 1\right)}^{s}
{\left(q + 1\right)} q^{2g - 1}.
\end{align*}
In \cite{GP-2018}, it is proven that the Deligne-Hodge polynomial of the whole representation variety is
\begin{align*}
\DelHod{\Xsp{g}{s}} =& \,{\left(q^2 - 1\right)}^{2g + s - 1} q^{2g - 1} +
\frac{1}{2} \, {\left(q -
1\right)}^{2g + s - 1}q^{2g -
1}(q+1){\left({2^{2g} + q - 3}\right)} 
\\ &+ \frac{\left(-1\right)^{s}}{2} \,
{\left(q + 1\right)}^{2g + s - 1} q^{2g - 1} (q-1){\left({2^{2g} +q -1}\right)}.
\end{align*}
Therefore, substracting the contribution of the previous strata, we obtain
\begin{align*}
	\DelHod{\Xirred{\Xsp{g}{s}}} =&\, 2^{2g - 1} \left(-1\right)^{s} {\left(q + 1\right)}^{2g + s - 1}
{\left(q - 1\right)} q^{2g - 1} \\&- 2^{2g} \left(-1\right)^{s} {\left(q + 1\right)} q^{2g}
{\left(\frac{{\left(1-q\right)}^{s} - 1}{q} + 1\right)}\\&+ \frac{1}{2} \, \left(-1\right)^{s}
{\left(q + 1\right)}^{2g + s - 1} {\left(q - 1\right)}^{2} q^{2g
- 1}  + {\left(2^{2
\, g} - {\left(q - 1\right)}^{2g}\right)} {\left(q - 1\right)}^{s}
{\left(q + 1\right)} q^{2g - 1} \\
&+ \frac{1}{2} \, {\left(q -
1\right)}^{2g + s - 1} {\left(q + 1\right)} {\left(q - 3\right)}
q^{2g - 1}\\
	& + \frac{{\left(2^{2g} q^{2} + 2^{2g + 1} q + 2^{2
\, g} + 2 \, {\left(q + 1\right)}^{2g + s}\right)} {\left(q -
1\right)}^{2g + s - 1} q^{2g - 1}}{2 \, {\left(q + 1\right)}}.
\end{align*}

	\item For $\Xred{\Xsp{g}{s}}$, the situation absolutely analoguous for the one of $\Xred{\Xfp{2g}{s}}$. Since the action there was closed, it is also on $\Xred{\Xsp{g}{s}} = \Xred{\Xfp{2g}{s}} \cap \Xsp{g}{s}$. Therefore, stratifying in terms of the stabilizers for the action, we have:
	\begin{itemize}
		\item $\XTilde{\Xsp{g}{s}}$. Here, the action of $\PGL{2}$ is free so, by Corollary \ref{cor:luna-thm-epol}, we have
		$$
			\DelHod{\XTilde{\Xsp{g}{s}} \sslash \SL{2}} = \frac{\DelHod{\XTilde{\Xsp{g}{s}}}}{q^3-q} = (-1)^s\left((q-1)^{2g} - 2^{2g}\right)q^{2g-2}(1-q)^{s-1}.
		$$
		\item $\XPh{\Xsp{g}{s}}$. Here, the action of $\PGL{2}$ is not free, but it has stabilizer isomorphic to $\Stab\,J_+ \cong \CC$. As in the free case, the GIT quotient is a locally trivial fibration with trivial monodromy so
		$$
			\DelHod{\XPh{\Xsp{g}{s}} \sslash \SL{2}} = \frac{\DelHod{\XPh{\Xsp{g}{s}}}}{q^2-1} = (-1)^s2^{2g}\frac{q^{2g}}{q-1}\left(\frac{(1-q)^s - 1}{q} + 1\right).
		$$
	\end{itemize}
Therefore, using the stratification $\Xred{\Xsp{g}{s}} = \XPh{\Xsp{g}{s}} \sqcup \XTilde{\Xsp{g}{s}}$ with $\XTilde{\Xsp{g}{s}}$ an open orbitwise-closed set, Theorem \ref{prop:decomposition-quotient} and Section \ref{section:stratification} give us
\begin{align*}
    \DelHod{\Xred{\Xsp{g}{s}} \sslash \SL{2}} &= \DelHod{\XPh{\Xsp{g}{s}} \sslash \SL{2}} + \DelHod{\XTilde{\Xsp{g}{s}} \sslash \SL{2}} \\&= 
\frac{2^{2g}q^{2g} \left(-1\right)^{s} {}
{\left(\frac{{\left(1-q\right)}^{s} - 1}{q} + 1\right)}}{q - 1} -
{{\left(2^{2g} - {\left(q - 1\right)}^{2g}\right)} {\left(q
- 1\right)}^{s-1} q^{2g - 2}}{}.
\end{align*}

\end{itemize}
Summarizing, the analysis above shows that
\begin{empheq}{align*}
	\DelHod{\Xsp{g}{s} \sslash \SL{2}} =&\,  {\left(q^2 - 1\right)}^{2g + s -
2} q^{2g - 2} +\left(-1\right)^{s} 2^{2g}  {\left(q - 1\right)} q^{2g - 2}
{\left({1-\left(1-q\right)}^{s - 1}\right)}\\
&+ \frac{1}{2}{\left(q - 1\right)}^{2g +s - 2} q^{2g - 2} \, {\left(2^{2g} + q - 3\right)}  \\
&+ \frac{1}{2} {\left(q + 1\right)}^{2g +
s - 2} q^{2g - 2}\,
\left(2^{2g} + q - 1\right).
\end{empheq}

\subsection{General parabolic data of Jordan type}
\label{subsec:parabolic-data-gen}

On the genus $g$ closed orientable surface, $\Sigma_g$, sonsider the parabolic structure $Q = \left\{(\gamma_1, [C_1]), \ldots, (\gamma_s, [C_s])\right\}$, where $C_i = J_+, J_-$ or $- \Id$. Let $r_+$ be the number of $J_+$, $r_-$ the number of $J_-$ and $t$ the number of $- \Id$ (so that $r_+ + r_- + t = s$), and let $\sigma = (-1)^{r_- + t}$. Observe that $J_+ \in [-J_-]$ and $[- \Id] = \left\{ -  \Id\right\}$ so, depending on $\sigma$, we have:
\begin{itemize}
	\item If $\sigma = 1$, then we have $\mathfrak{X}_{\SL{2}}(\Sigma_g, Q) = \mathfrak{X}_{\SL{2}}(\Sigma_g, Q_{r}^+)$ where $r = r_+ + r_-$. Hence, we have $\DelHod{\mathfrak{X}_{\SL{2}}(\Sigma_g, Q) \sslash \SL{2}} = \DelHod{\Xsp{g}{r} \sslash \SL{2}}$ and the polynomial follows from the computations above.
	\item If $\sigma = -1$, then we have $\mathfrak{X}_{\SL{2}}(\Sigma_{g}, Q) = \mathfrak{X}_{\SL{2}}(\Sigma_{g}, Q_{r}^-)$ where $r = r_+ + r_-$ and consider $Q_r^- = \left\{(\gamma_1, [J_+]), \ldots, (\gamma_r, [J_+]), (\gamma_{r+1}, \left\{- \Id\right\}) \right\}$. This is the so-called twisted representation variety.
This variety does not contain reducible representations. To check that,
let $A = (A_1, B_1, \ldots, C_1, \ldots, C_r, - \Id) \in \mathfrak{X}_{\SL{2}}(\Sigma_{g}, Q_{r}^-)$ so that
$$
	\prod_{i=1}^{2g} [A_i, B_i] \prod_{j=1}^{s} C_i = - \Id.
$$
If $v \in \CC^2 - \left\{0\right\}$ is a common eigenvector of $A$, then, since all the eigenvalues of the commutators $[A_i, B_i]$ and the $C_i$ are equal to $1$, the left hand side of the previous equation fixes $v$ but the right hand side does not. This proves that such a $v$ cannot exists. Therefore, the action of $\PGL{2}$ on $\mathfrak{X}_{\SL{2}}(\Sigma_{g}, Q_{r}^-)$ is closed and free.

As proven in \cite{GP-2018}, the representation variety has Deligne-Hodge polynomial
\begin{align*}
	\DelHod{\mathfrak{X}_{\SL{2}}(\Sigma_{g}, Q_{r}^-)} = &\, {\left(q - 1\right)}^{2g + r - 1} (q+1)q^{2g - 1}{{\left( {\left(q + 1\right)}^{2 \,
g + r-2}+2^{2g-1}-1\right)} } \\
&+ \left(-1\right)^{r + 1}2^{2g - 1}  {\left(q + 1\right)}^{2g + r
- 1} {\left(q - 1\right)} q^{2g - 1}.
\end{align*}
Therefore, Corollary \ref{cor:luna-thm-epol} gives us
\begin{align*}
	\DelHod{\mathfrak{X}_{\SL{2}}(\Sigma_{g}, Q) \sslash \SL{2}}  =&\, \left(-1\right)^{r-1}2^{2g - 1}  {\left(q + 1\right)}^{2g + r - 2} q^{2g - 2} \\ &+ {\left(q - 1\right)}^{2g + r - 2} q^{2g - 2}\left( {\left(q + 1\right)}^{2g + r - 2} + 2^{2g - 1} - 1\right).
\end{align*} 
\end{itemize}

Summarizing, we have proven the following result.

\begin{thm} With the notations above, the Deligne-Hodge polynomials of parabolic $\SL{2}(\CC)$-character varieties are:
\begin{itemize}
	\item If $\sigma = 1$, then
	\begin{align*}
	\DelHod{\Char{\SL{2}(\CC)}(\Sigma_g, Q)} =& \,{\left(q^2 - 1\right)}^{2g + r -
2} q^{2g - 2} +\left(-1\right)^{r} 2^{2g}  {\left(q - 1\right)} q^{2g - 2}
{\left({1-\left(1-q\right)}^{r - 1}\right)}\\
&+ \frac{1}{2}{\left(q - 1\right)}^{2g +r - 2} q^{2g - 2} \, {\left(2^{2g} + q - 3\right)}  \\
&+ \frac{1}{2} {\left(q + 1\right)}^{2g +
r - 2} q^{2g - 2}\,
\left(2^{2g} + q - 1\right).
	\end{align*}
	\item If $\sigma = -1$, then
	\begin{align*}
	\hspace{-1.75cm}\DelHod{\Char{\SL{2}(\CC)}(\Sigma_g, Q)} =& \left(-1\right)^{r-1}2^{2g - 1}  {\left(q + 1\right)}^{2g + r - 2} q^{2g - 2} \\ &+ {\left(q - 1\right)}^{2g + r - 2} q^{2g - 2}\left( {\left(q + 1\right)}^{2g + r - 2} + 2^{2g - 1} - 1\right).
\end{align*} 
\end{itemize}
\end{thm}

Recall that the union of the conjugacy classes of $J_+, J_-$ and $- \Id \in \SL{2}$, together with $ \Id$, is exactly the set of $C \in \SL{2}$ such that $\tr C = \pm 2$. Therefore, with this result, we have finally computed the Deligne-Hodge polynomial of the parabolic character variety of a punctured closed orientable surface with any parabolic structure of the form $Q = \left\{(\gamma_i, [C_i])\right\}$ where $C_i \in \SL{2} - D$.

%% file: Chapters/FutureWork.tex
\label{future-work} % For referencing the chapter elsewhere, use \ref{Chapter1} 

\lhead{Future work} % This is for the header on each page - perhaps a shortened title

%----------------------------------------------------------------------------------------

\subsection*{Parabolic $\SL{2}(\CC)$-character varietes with semi-simple punctures}

As we showed in Sections \ref{sec:hodge-theory-sl2}, \ref{sec:tube-J+} and \ref{sec:genus-tube}, in the case $G=\SL{2}(\CC)$ the core submodule $\cW = \langle T_{\pm 1}, T_\pm, T_{\Bt}, S_{\pm 2}, S_2 \otimes S_{-2}\rangle \subseteq \KM{[\SL{2}(\CC)/\SL{2}(\CC)]}$ is an invariant submodule for $\Zg{\SL{2}}(L)$ and $\Zg{\SL{2}}(L_\lambda)$ for $\lambda = [J_\pm] \subseteq \SL{2}(\CC)$. However, some preliminary computations show that, if we want to extend the previous result to semi-simple conjugacy classes, the submodule $\cW$ is no longer invariant.

To be precise, pick $t \in \CC-\left\{\pm 2\right\}$ and consider the class $\lambda_t$ of the elements of $\SL{2}(\CC)$ with trace $t$ i.e.\ $\lambda_t = [D_\mu]$ with $\mu + \mu^{-1}=t$. In that case, we have that $\Zg{\SL{2}}(L_{\lambda_t})(\cW) \subseteq \langle \cW, S_{t}, S_{-t}\rangle$, where $S_{t}$ (resp.\ $S_{-t}$) denotes the monodromy representation of order two on $\CC-\left\{\pm 2, \pm t\right\}$ that is non-trivial on the loop around $t$ (resp.\ around $-t$).

Nevertheless, this homomorphism can be computed explicitly using the techniques of Section \ref{sec:tube-J+} without further considerations. Moreover, if we now consider $\lambda_{t'}$ with $t' \neq \pm t$, then the same calculation can be performed and the image submodule can be controlled. However, if we want to add a new tube with $t' = \pm t$, a non-trivial interference between the two puctures occurs. This effect is new with respect to the Jordan type case and reflects an arithmetic dependence previously observed for parabolic Higgs bundles.

For this reason, a prospective future work is to study parabolic $\SL{2}(\CC)$-representation varieties with semi-simple punctures, both in the generic and non-generic setting. In addition, it may be important to analyze the stability of the geometric structures of this parabolic representation varieties with respect to the punctures in a (kind of) chambers-and-walls framework as in \cite{GP-Gothen-Munoz}. These cases are especially important since the parabolic $\SL{2}(\CC)$-character variety with two semi-simple punctures over a torus is diffeomorphic to the moduli space of doubly periodic instantons by means of the Nahm transform.  

\subsection*{Higher rank character varietes}

The constructed TQFT, $\Zs{G}: \Bordp{n} \to \Modt{\K{\MHSq}}$, is valid for any complex group $G$ and any dimension $n$. For this reason, it would be interesting to considered different groups that $\SL{2}(\CC)$, for example, the higher rank groups $\SL{r}(\CC)$ with $r \geq 2$. As a first step, it would be useful to study the rank $3$ case and to compute the Hodge structures of $\SL{3}(\CC)$-parabolic representation varieties by means of the TQFT $\Zs{\SL{3}(\CC)}$ and its geometric counterpart $\Zg{\SL{3}(\CC)}$.

Recall that, as shown in Section \ref{sec:sl2-repr-var}, the key space in the calculation of $\SL{2}(\CC)$-representation varieties was the piecewise quotient $[\SL{2}(\CC) / \SL{2}(\CC)]$ by inner automorphisms. In this case, we have that the GIT quotient $\SL{2}(\CC) \sslash \SL{2}(\CC) = \CC$, with quotient map given by the trace. Moreover, this GIT quotient collapses the orbits of the two Jordan type matrices to the completely reducible ones. In this way, $[\SL{2}(\CC) / \SL{2}(\CC)]$ can be understood as $\CC$ with the points $\pm 2$ doubled. These doubled points are precisely the origin of the monodromy interaction captured by the Hodge monodromy representation.

Analogously, for the rank $3$ case, the computations of $\Zg{\SL{3}(\CC)}$ are based on the piecewise quotient $[\SL{3}(\CC) / \SL{3}(\CC)]$. In general, the GIT quotient is given by the coefficients of the characteristic polynomial so, now, $\SL{3}(\CC) \sslash \SL{3}(\CC) = \CC^2$ and the quotient map controls both the trace and the minors of the matrix. Hence, $[\SL{3}(\CC) / \SL{3}(\CC)]$ is just $\CC^2$ with some doubled curves corresponding to the collapsing Jordan type orbits. Thus, we expect that the computations of the geometric TQFT will be similar to the rank $2$ case but with Hodge monodromy representations supported on curves.

Another exciting prospective work would be to consider higher dimensional manifolds for the TQFT. For example, it would be interesting to analyze the case of $3$-manifolds and to link it with Thurston's geometrisation program.

\subsection*{Relation with other constructions of TQFTs}

The building method of TQFTs given in Section \ref{sec:lax-monoidal-tqft} is, formally, similar to other existing constructions in the literature so it would be of primary interest to link them. An inspiring starting point is the paper \cite{Dijkgraaf-Witten:1990}, in which a TQFT is constructed for character varieties with finite gauge group. The quantisation method there is, somehow, similar to our approach since, in both cases, there is a kind of weighted `counting' of the character varieties (number of points in \cite{Dijkgraaf-Witten:1990} and by means of $K$-theory in this thesis). In this direction, in \cite{Ben-Zvi-Gunningham-Nadler} (see also \cite{Ben-Zvi-Francis-Nadler:2010} and \cite{Ben-Zvi-Nadler}) it was proposed an extended TQFT for the homology of character varieties and related with geometric Langlands duality. There, the TQFT is constructed indirectly via the Cobordism Hypothesis (see \cite{Lurie:2009}) by proving that the monoidal category of Harish-Chandra bimodules is a 2-dualizable Calabi-Yau category. It would be fascinating to be able to restate all these proposals in a common framework and to compare them.

Further than the particular scope of application, there are also strong parallelisms in the construction methods of the literature that are based on the same philosophical principle: TQFTs can be constructed as a composition of a field theory and a quantisation, in the spirit of a mathematical formulation of Lagrangian field theories. In this way, it would be very interesting to unravel the relation between the sheaf theoretic formalism of \cite{Gaitsgory-Rozenblyum}, based in correspondences, and the $\cC$-algebras developed in Section \ref{sec:C-algebras}. As explained in Remark \ref{rmk:C-algebra-as-correspondence}, both are models of Grothendieck's six functors formalism and can be used for constructing quantisations. Another alternative point of view is the work \cite{Rovi-Schoenbauer:2018}, in which cut-and-paste invariants are related with TQFTs. In this direction, a prospective future work would be to recast the Seifert-van Kampen property (Section \ref{sec:field-theory}) with these cut-and-paste invariants.

\subsection*{TQFTs for the confluence scheme of Painlev\'e equations}

Another interesting scope of application of the ideas developed in this thesis is to formulate the results of \cite{Chekhov-Mazzocco-Rubtsov:2017} and \cite{Chekhov-Mazzocco-Rubtsov:2017b} about monodromy manifolds and decorated character varieties in the language of TQFTs. The point is that the behaviour of these manifolds with respect to the `chewing-gum' operations strongly suggests that the confluence scheme of Painlev\'e equations and their relation with decorated character varieties should underlie a TQFT that governs the degeneration process.

As mentioned in Section \ref{sec:TQFT-over-sheaf}, we may consider an extra $2$-category structure on $\Bord{n}$ to allow non-reversible deformations on bordisms. In this way, after extending the category $\Bord{2}$ to comprise also bordered cusped surfaces, the `chewing-gum' operations of \cite{Chekhov-Mazzocco-Rubtsov:2017} may be understood as non-invertible $2$-morphisms between bordisms. Under this interpretation, the confluence scheme would correspond to a commutative diagram in $\Bord{2}$. Moreover, seeing decorated representation varieties as complexifications of bordered cusped Teichm\"uller spaces, these $2$-cells would also have a combinatorial interpretation in terms of cusped fat graphs (see \cite{Fock-Goncharov}). 

For that, it is necessary to analyze what kind of field theory should be used. It seems plausible that a modification of the generalized representation varieties considered in Section \ref{sec:TQFT-for-repr} may work. An heuristic argument for so is that use of shear coordinates in \cite{Chekhov-Mazzocco-Rubtsov:2017} and \cite{Chekhov-Mazzocco-Rubtsov:2017b} is very similar to the description of the TQFT given in Section \ref{almost-TQFT-strategy}. Moreover, the natural conjugacy action on the representation varieties does not preserve the field theory considered in this thesis, but the action by unipotent Borel subgroups, as in \cite{Chekhov-Mazzocco-Rubtsov:2017}, might do. That suggests this should be right framework. Moreover, it is interesting to look for a suitable quantisation giving rise to a complete TQFT. In this direction, a first try would be to consider the Cherednik and Zhedanov algebras that were used in \cite{Mazzocco:2016} to relate Painlev\'e equations with the degenerations of the $q$-Askey scheme.

\subsection*{Extension of pseudo-quotients}

In Section \ref{subsec:reducible-rep}, we have intensively used pseudo-quotients for computing the Deligne-Hodge polynomial of $\SL{2}(\CC)$-character varieties from the corresponding representation varieties. A panoramic view of the techniques used in that chapter shows that the stratifications of the representation variety considered were determined by the stabilizers of the action of $\SL{2}(\CC)$. This brings us back to the so-called Luna's stratification \cite{Luna:1973} that, roughly speaking, stratifies an algebraic variety in terms of the stabilizers. Therefore, a deep study of the relation between these two techniques is required. It is expectable that Luna's stratification would be the key for considering quotients of higher rank representation varieties.

Furthermore, as proposed  to us by L. \'Alvarez-Consul, the construction of the quotient of a stack by a group action is, formally, similar to the one of pseudo-quotients as the quotient stack takes into account the action of the stabilizers by the action. For this reason, it would be interesting to analyze the relation of pseudo-quotients with quotients of stacks. This exploration can also be very useful to address mirror symmetry conjectures (see last objective) since, as explained in \cite{Martinez:2017}, in the parabolic setting a correction term must be expected coming from the stabilizers for the action.

\subsection*{Equivariant Hodge monodromy representation}

In Section \ref{sec:equivariant-hodge-mono}, we sketched how the Hodge monodromy representation behaves when we have an underlying action of $\ZZ_2$. This is just the first step in the more ambitious goal of studying equivariant Hodge monodromy representations.

As a first step, it would be interesting to extend the results of Section \ref{sec:equivariant-hodge-mono} to the case of an action of a finite group $F$ on an algebraic variety $X$. This study is necessary, for example, for addressing the problem of $\SL{r}(\CC)$-varieties. The idea is that, for $\SL{2}(\CC)$, the action of $\ZZ_2$ appeared as the action of interchanging the two eigenvalues of a matrix. For this reason, in the $\SL{r}(\CC)$ case, it can be expected to have an underlying action of $S_r$ by interchanging the $r$ eigenvalues (see \cite{Lawton-Munoz:2016} for the $r=3$ case).

As shown in \cite{Florentino-Silva:2017}, Section 4.2, the mixed Hodge structure on $H_c^\bullet(X/F;\QQ)$ is precisely the invariant part of $H_c^\bullet(X;\QQ)$ under the action of $F$ on cohomology. Obviously, it is in agreement with the general principle in GIT that the rings associated to $X \sslash F$ are precisely the $F$-invariant parts of the ones on $X$. For this reason, we could expect that similar results must hold for Hodge monodromy representations. Moreover, it would be very interesting to link it with the theory of equivariant mixed Hodge modules \cite{Tanisaki:1987}.

\subsection*{Construction of strict TQFTs for character varieties}

The TQFT for $G$-representation varieties is not a strict monoidal functor, but only a lax monoidal functor. Recall that, from the Zorro's property (Lemma \ref{lem:zorro}), it is absolutely necessary to relax the monoidality condition if we want to use Grothendieck rings of mixed Hodge modules as quantisations, since they are very far from being finitely generated. Hence, this finiteness condition is a major obstruction to $\Zs{G}$ to be strict monoidal.

Nonetheless, as we showed in \ref{sec:sl2-repr-var}, for the $2$-TQFT we do not need the whole module $\KM{[G/G]}$ and we can restrict our attention to the submodule $\cV = \cV_{(S^1, \star)} \subseteq \KM{[G/G]}$. Moreover as byproduct of the computations of this thesis, \cite{MM} and \cite{Martinez:2017}, it was shown that, actually, $\cV$ is finitely generated, at least for $G=\SL{2}(\CC)$ and $G=\PGL{2}(\CC)$.
For this reason, we expect that, restricting $\Zs{G}$ to a submodule of the target, we could promote it to a monoidal functor. Hence, the objective in this direction would be to find out under which conditions the submodule $\cV$ is finitely generated and to construct strict monoidal TQFTs from it. 

For this purpose, a good idea would be to base on the calculations of $\SL{3}(\CC)$-character varieties previously mentioned. The point is that, for $G=\SL{2}(\CC)$, this finiteness condition might have appeared as a consequence of the fact that it has finitely many non-polystable orbits and the monodromy always arises around these orbits. However, this no longer holds for $G=\SL{3}(\CC)$ so it may be expected that the computations of this case will be helpful for clarifying the situation. Even in the case that the result was not true, it would give the first example of an almost-TQFT that cannot be extended to a strict TQFT.

This kind of behaviour usually means that the whole $\Zs{G}$ actually is a monoidal functor, but with a different monoidal structure on the source and target categories (see \cite{Wasserman:2017} and \cite{Wasserman-phd} for an example of this phenomenon). Hence, another objective in this direction would be to design an appropiate tensor product that makes the TQFT $\Zs{G}$ monoidal.

Finally, a step further would be to improve the construction method developed in this thesis and to obtain an extended TQFT (see \cite{Freed:1994} and \cite{Lawrence:1993}). For that purpose, we will need to bind the gluing properties of character varieties and Hodge modules with the presentations of the category of extended $2$-bordisms (see \cite{Schommer-Pries:2009}) and once-extended $3$-dimensional TQFTs (see \cite{Bartlett-Douglas-Schommer-Pries-Vicary:TBAa}, \cite{Bartlett-Douglas-Schommer-Pries-Vicary:2014} and \cite{Bartlett-Douglas-Schommer-Pries-Vicary:TBAb}). The later case is particulary exiciting since these TQFTs are known to be equivalent to modular tensor categories (see \cite{Bakalov-Kirillov} and \cite{Bartlett-Douglas-Schommer-Pries-Vicary:2015}).

\subsection*{TQFTs accross the non-abelian Hodge correspondence}

So far, we have focused our study on the character varieties part of the non-abelian Hodge theory. However, as explained in \ref{sec:lax-monoidal-tqft}, the construction of the TQFT $\Zs{G}$ requires the use of two ingredients: the field theory, that captures which space we are studying; and the quantisation, that captures which algebraic invariant we are modelling. Hence, if we still use mixed Hodge modules as quantisation but we vary the field theory, we can study Hodge structures on different spaces. In particular, a prospective work would be to construct field theories adapted to moduli spaces of Higgs bundles and of flat connections accross the non-abelian Hodge theory (see Section \ref{sec:non-abelian-hodge}). However, these moduli spaces depend not only on the topological structure of the underlying surface (as for character varieties) but on its complex structure. For this reason, it will be necessary to endow the surface with an extra structure that encodes the complex structure as, for example, a conformal structure, and to addapt the TQFT to this context.

This understanding is relevant for the following reason. The correspondences of non-abelian Hodge theory are far from being holomorphic and, actually, each of the moduli spaces on this correspondence is equipped with its own complex structure. In other words, there is a single object from the topological point of view but it admits three different (but compatible) K\"ahler structures, in a structure called hyperk\"ahler. For this reason, the construction of this triple of TQFTs will allow us to fully understand this hyperk\"ahler structure.

\subsection*{Mirror symmetry}

Finally, another framework in which character varieties are central is the geometric Langlands program (see \cite{Beilinson-Drinfeld}). In \cite{Hausel-Thaddeus}, Hausel and Thaddeus proved that the moduli space of $\SL{2}(\CC)$-Higgs bundles is the first non-trivial example of the Strominger-Yau-Zaslow conditions for mirror symmetry for Calabi-Yau manifolds (see \cite{Strominger-Yau-Zaslow}). As shown in \cite{Hausel:2005} and \cite{Hausel-Rodriguez-Villegas:2008}, this conjecturally would imply that a collection of conjectures concerning Hodge structures of character varieties must hold.

In particular, it can be expected some very astonishing symmetries between Hodge structures of $G$-character varieties and ${^L}G$-character varieties, where ${^L}G$ is the Langlands dual group of $G$. The validity of these conjectures has been discussed in some cases as in \cite{LMN} and \cite{Martinez:2017} for $G=\SL{2}(\CC)$ and ${^L}G = \textrm{PGL}_2(\CC)$. Despite of that, the general case remains unsolved.

At a long term, we hope that the ideas introduced in this work by means of the constructed TQFT for representation varieties might be useful to shed light into these questions. In particular, it would be very interesting to study whether the TQFTs of $G$ and ${^L}G$ are somehow related and what kind of implications it has for mirror symmetry.